# EQUIVARIANT LOCALIZATION IN FACTORIZATION HOMOLOGY AND APPLICATIONS IN MATHEMATICAL PHYSICS II: GAUGE THEORY APPLICATIONS

DYLAN BUTSON

ABSTRACT. We give an account of the theory of factorization spaces, categories, functors, and algebras, following the approach of [Ras15a]. We apply these results to give geometric constructions of factorization $\mathbb{E}_n$ algebras describing mixed holomorphic-topological twists of supersymmetric gauge theories in low dimensions. We formulate and prove several recent predictions from the physics literature in this language:

We recall the Coulomb branch construction of [BFN18] from this perspective. We prove a conjecture from [CG18] that the Coulomb branch factorization $\mathbb{E}_1$ algebra $\mathcal{A}(G, N)$ acts on the factorization algebra of chiral differential operators $\mathcal{D}^{\mathrm{ch}}(Y)$ on the quotient stack $Y = N/G$. We identify the latter with the semi-infinite cohomology of $\mathcal{D}^{\mathrm{ch}}(N)$ with respect to $\hat{\mathfrak{g}}$, following the results of [Ras20b]. Both these results require the hypothesis that $Y$ admits a Tate structure, or equivalently that $\mathcal{D}^{\mathrm{ch}}(N)$ admits an action of $\hat{\mathfrak{g}}$ at level $\kappa = -\mathrm{Tate}$.

We construct an analogous factorization $\mathbb{E}_2$ algebra $\mathcal{F}(Y)$ describing the local observables of the mixed holomorphic-B twist of four dimensional $\mathcal{N} = 2$ gauge theory. We apply the theory of equivariant factorization algebras of the prequel [But20a] in this example: we identify $S^1$ equivariant structures on $\mathcal{F}(Y)$ with Tate structures on $Y = N/G$, and prove that the corresponding filtered quantization of $\iota^! \mathcal{F}(Y)$ is given by the two-periodic Rees algebra of chiral differential operators on $Y$. This gives a mathematical account of the results of [BLL+15]. Finally, we apply the equivariant cigar reduction principle of [But20a] to explain the relationship between these results and our account of the results of [CG18] described above.

## CONTENTS

# 1. INTRODUCTION

We begin with some overarching remarks about the background and motivation for the present work, as well as its companion paper and formal prequel [But20a], which we call Parts II and I, respectively. In Section 1.3 we give a brief review of the foundational results of Part I. In Section 1.4 we give an overview of the results of the present work.

## 1.1. **Background: Factorization algebras, representation theory, and quantum field theory.** Factorization algebras were introduced by Beilinson and Drinfeld in [BD04] as a model for algebras of observables in two dimensional chiral conformal quantum field theories, defined in the language of algebraic geometry. Factorization algebras in this setting generalize vertex algebras to global objects defined over algebraic curves, vaguely analogous to sheaves on them. A generalization of the theory of factorization algebras to higher dimensional varieties was also given in [FG11], analogously modeling holomorphic quantum field theories in higher dimensions, which by definition generalize the holomorphic behaviour of observables in chiral conformal field theories in two real dimensions. From the beginning, the development of this theory was motivated by the essential connection between chiral conformal field theory and representation theory of affine Lie algebras.

An analogue of factorization algebras defined over smooth manifolds in the language of algebraic topology was proposed by Lurie in [Lur08], as an example of a class of extended topological field theories in the mathematical sense defined therein, and pursued by Ayala, Francis, Lurie, and collaborators in [AF15, AFR15, AFT16, Lur09a, Lur12]. Factorization algebras in the topological setting analogously generalize algebras over the little $n$-discs operad, and again describe the algebras of observables in topological quantum field theories of dimension $n$. In the case $n = 1$ these are equivalent to usual (homotopy) associative algebras, a central topic of study in classical representation theory.

Thus, there is a natural dictionary between predictions of quantum field theory or string theory, which have led to groundbreaking ideas in a variety of areas of mathematics, and statements in representation theory phrased in terms of factorization algebras. This dictionary is both a primary motivation and the main source of new ideas for this series of papers.

The work of Costello [Cos11] and Costello-Gwilliam [CG16] established such a dictionary in a more analytic context, constructing a variant of factorization algebras defined over smooth manifolds in terms of the differential geometric input data of a Lagrangian classical field theory satisfying certain ellipticity requirements together with a choice of renormalization scheme. These ideas were very influential for the present series of papers, and have led to many other developments following this paradigm [Cos13, CS15, BY16, GW18, ES19, SW19, ESW20].

The present series of papers also closely follows the program of Ben-Zvi, Nadler, and collaborators, which gives approaches to many facets of geometric representation theory in terms of extended topological field theory and derived algebraic geometry [BZN09, BZFN10, BZN13, BZG17, BZN18]. In particular, the use of sheaf theory in constructing extended topological field theories from geometry is a central theme of the present series of papers, which is borrowed from *loc. cit.*. Further, the derived stacks and sheaf theories defined on them which are relevant for our constructions can often be predicted from statements about the shifted symplectic geometry of the spaces of solutions to the Euler-Lagrange equations in the relevant classical field theories. This relies on a family of ideas about functoriality of shifted geometric quantization, closely related to those in *loc. cit.*, which I learned from Pavel Safronov [Saf20].

Finally, the circle of ideas and mathematical technology around the local geometric Langlands correspondence [ABC+18], derived geometric Satake correspondence [BF08], and Coulomb branch construction [BFM05, BFN18, BFN19b], provided a collection of mathematically well-understood examples and established techniques which were crucial for the technical underpinning for the present series of papers. In particular, we follow the sheaf theory foundations given in [GR14a, GR14b, Gai15, GR17a, GR17b, Ras15a, Ras15b, Ras20b] and references therein. These ideas can naturally be interpreted in certain holomorphic-topological twists of supersymmetric quantum field theories, as we explain below. These interpretations have also been studied in a more mathematical context, for example in [EY18, BZN18, EY19, EY20, RY19], and I have benefited greatly from ongoing discussions with Justin Hilburn and Philsang Yoo about these ideas. In particular, the forthcoming papers [HY], [GY], and [HR] will also contain some of their ideas that we follow in the present work.

In terms of the various perspectives we have just discussed, we can summarize an underlying goal of this series of papers as follows:

We develop a dictionary between factorization algebras and quantum field theory in the mixed holomorphic-topological setting, using a synthesis of the chiral and topological variants of factorization algebras; examples of interest are given by factorization compatible sheaf theory constructions, motivated by shifted geometric quantization of spaces of solutions to equations of motion in supersymmetric gauge theories, and using tools from geometric representation theory and derived algebraic geometry.

## 1.2. Motivation: Holomorphic-topological twists of supersymmetric quantum field theories and $\Omega$-backgrounds.

The more broad goal of this series of papers is to use this dictionary to formulate and prove results from a particular family of interconnected predictions of string theory, at the intersections of affine representation theory [KW07, GW09, Gai18, BPRR15], enumerative geometry [AGT10, Nek16, NP17, GR19], low-dimensional topology [GGP16, DGP18, Wit12], and integrable systems [Nek03, NS10, NW10]. These ideas are centred around the six dimensional $\mathcal{N} = (2, 0)$ superconformal field theory, sometimes called "theory X", which is an elusive, non-Lagrangian quantum field theory that morally describes fluctuations of M5 branes in M theory (which we remind the reader are geometric objects supported on six dimensional spaces). This theory is considered on a spacetime of the form $C \times M$, for $C$ a smooth algebraic curve and $M$ a smooth four manifold, and this gives rise to natural predictions relating chiral factorization algebras over the curve $C$ with the differential topology of the four manifold $M$, or the enumerative geometry of sheaves in the case $M = S$ is a smooth algebraic surface over $\mathbb{C}$.

As an intermediate step, we establish analogous predictions from three and four dimensional gauge theories following [BDG17, CG18, BLL$^+$15], which correspondingly relate to the representation theory of classical Lie algebras, and of quantizations of symplectic singularities more generally [BDGH16], as well as to more classical aspects of enumerative geometry [BDG$^+$18] and integrable systems [NS09, CWY18].

Similar ideas have been studied extensively in mathematics already in both of the above contexts, often explicitly motivated by the same physics considerations; for example [FG06, Ara18, Bra04, SV13, MO19, BFN14, Neg17, RSYZ19, FG20, BD99, BFN18, BZG17, Cos13] are a few which have been influential in our understanding of this family of ideas, ordered roughly corresponding to the physics references above.

The preceding predictions are nominally phrased in terms of string theory and supersymmetric quantum field theory, which are notoriously difficult to understand and often not yet defined mathematically, but an important common feature of these results from the physical perspective is that they often factor through mixed holomorphic-topological twists of the relevant quantum field theories. As a result, these theories are expected to be amenable to descriptions in terms of algebraic geometry and topology, and in particular the algebras of observables of these theories are expected to correspond to objects in the synthesis of chiral and topological factorization algebras mentioned above that we study in the present work. This is the fundamental reason for the effectiveness of the mathematical tools considered in the present work in the relevant physics context.

However, there is another salient feature of many of the physical constructions and corresponding mathematical interpretations mentioned above, which has not been codified mathematically in our explanation so far: in the seminal paper [Nek03], Nekrasov introduced a construction in quantum field theory called an $\Omega$-*background*, an additional structure on a partially topological quantum field theory which (when it exists) deforms the given theory in a way that enforces rotational equivariance with respect to a fixed $S^1$ action on the underlying spacetime. A primary consequence is that cohomological calculations in $\Omega$-deformed topological field theories are given by the analogous calculations in equivariant cohomology.

Moreover, motivated by the localization theorem in equivariant cohomology, we expect these calculations should in some sense localize to the fixed points of the underlying $S^1$ action, after passing to an appropriate localization $\mathbb{K}[\varepsilon][f^{-1}]$ of the base ring $\mathbb{K}[\varepsilon] := H^{\bullet}_{S^1}(\mathrm{pt}; \mathbb{K})$. In fact, calculations in the algebras of observables of the $\Omega$-deformed theories localize to calculations in (families over $\mathbb{K}[\varepsilon][f^{-1}]$ of) algebras of observables over the fixed point locus. Furthermore, such families of algebras of observables have been observed in [NS09, NW10] to define filtered quantizations of the algebra specialized at the central fibre over $\mathbb{K}[\varepsilon]$.

In the companion paper [But20a], formally Part I of the series, we establish the foundations of the theory of equivariant factorization algebras in the mixed chiral-topological setting. Moreover, in this language we give an account of the equivariant localization and quantization phenomena associated with the $\Omega$-background construction in holomorphic-topological quantum field theory described above. We briefly review these results presently in Section 1.3.

In the present work, formally Part II of the series, we develop methods for constructing examples of equivariant factorization algebras corresponding to holomorphic-topological twists of supersymmetric gauge theories, and apply the results of Part I in these examples. In Section 1.4.1 we give an overview and summary of the results of Chapter 1 of the present work, and in Section 1.4.2 we give the same for Chapter 2.

### 1.3. **Review of Part I.**

In this subsection, we give review the results of the companion paper and formal prequel [But20a]. To begin, we develop a theory of equivariant factorization algebras $A \in \mathrm{Alg}^{\mathrm{fact}}(X)^G$ on algebraic varieties $X$ with the action of a connected algebraic group $G$ and define an equivariant analogue of factorization homology

$$\int_X^G : \mathrm{Alg}^{\mathrm{fact}}(X)^G \to H^{\bullet}_G(\mathrm{pt})\text{-Mod} .$$

In the case $G = (\mathbb{C}^{\times})^n$, we prove an equivariant localization theorem in this context:

*Theorem* 1.3.1. [But20a] Let $A \in \mathrm{Alg}^{\mathrm{fact}}(X)^G$ be an equivariant factorization algebra. The natural map

$$\int_{X^G}^G \iota^! A \overset{\cong}{\Longrightarrow} \int_X^G A$$

induces an equivalence over the localization $H^{\bullet}_G(\mathrm{pt})[f_k^{-1}]$.

Once correctly formulated, the proof of this statement follows straightforwardly from the results of [GKM97]. Nonetheless, it provides an important link between higher dimensional factorization algebras on $X$, which are often subtle to understand algebraically due to their homotopical nature, and lower dimensional factorization algebras on $X^G$, which can be identified with more familiar objects such as associative algebras or vertex algebras.

Next, we study the algebraic structure of equivariant factorization algebras in the simplest examples, explain relations to algebras over variants of the framed little $n$-disks operad, and give an account in this language of the relationship to deformation quantization predicted in the physics literature, as described above. The latter is given as follows:

*Theorem* 1.3.2. [But20a] There are equivalences of categories

$$\mathrm{Alg}^{\mathrm{fact}}(X \times \mathbb{A}^1)^{\mathbb{G}_a \rtimes \mathbb{G}_m} \cong \mathrm{Alg}^{\mathrm{fact}}_{\mathbb{E}^{S^1}_2}(X) \cong \mathrm{Alg}^{\mathrm{fact}}_{\mathbb{BD}^u_0}(X) \ .$$

The latter category is equivalent to that of (two-periodic) filtered quantizations of (shifted) Coisson algebras to chiral factorization algebras on $X$, by the chiral Poisson additivity theorem of Rozenblyum. Thus, factorization $\mathbb{E}^{S^1}_2$ algebras $A \in \mathrm{Alg}^{\mathrm{fact}}_{\mathbb{E}^{S^1}_2}(X)$ induce quantizations $\iota^! A \in \mathrm{Alg}^{\mathrm{fact}}(X)_{/\mathbb{K}[\varepsilon]}$.

Finally, we explain the manifestation in this language of the physical principle of equivariant cigar reduction, which plays a central role in our applications of interest in Part II [But20b], as we explain below. Consider an $S^1$ equivariant factorization $\mathbb{E}_2$ algebra

$$A \in \mathrm{Alg}^{\mathrm{fact}}_{\mathbb{E}^{S^1}_2}(X) \qquad \text{and define} \qquad A_0 := \mathrm{oblv}^{\mathbb{E}_0}_{\mathbb{E}^{S^1}_2} A \ \in \mathrm{Alg}^{\mathrm{fact}}(X) \ .$$

Note that $A_0$ is canonically a module over $A$ in the $\mathbb{E}_2$ sense, so that there is a module structure

$$(1.3.1) \qquad A_0 \in \mathrm{CH}_\bullet(A)\text{-Mod}(\mathrm{Alg}^{\mathrm{fact}}(X)) \qquad \text{or equivalently a map} \qquad \mathrm{CH}_\bullet(A) \to \mathrm{CH}^\bullet(A_0)$$

in the category $\mathrm{Alg}^{\mathrm{fact}}_{\mathbb{E}_1}(X)$ of factorization $\mathbb{E}_1$ algebras. In these terms, we have the following additional structure relating the factorization $\mathbb{E}^{S^1}_2$ algebra and the corresponding factorization $\mathbb{BD}^u_0$ algebra, which has a geometric interpretation in physics as the *equivariant cigar reduction principle*, pictured in Figure 1:

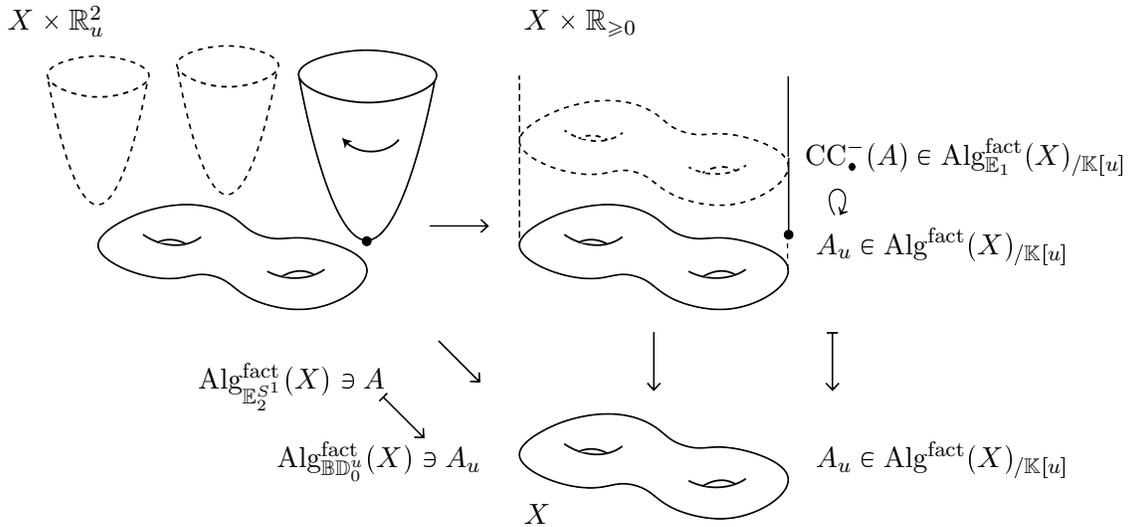

FIGURE 1. The equivariant cigar reduction principle

*Proposition* 1.3.3. [But20a] The family of factorization algebras $A_u \in \mathrm{Alg}^{\mathrm{fact}}(X)_{/\mathbb{K}[u]}$ underlying the factorization $\mathbb{B}\mathbb{D}_0^u$ algebra corresponding to $A \in \mathrm{Alg}^{\mathrm{fact}}_{\mathbb{E}_2^{S^1}}(X)$ under Theorem 1.3.2, admits a canonical module structure

$$A_u \in \mathrm{CC}_{\bullet}^-(A)\text{-Mod}(\mathrm{Alg}^{\mathrm{fact}}(X)_{/\mathbb{K}[u]}) \qquad \text{such that} \qquad A_u|_{\{u=0\}} = A_0 \in \mathrm{CH}_{\bullet}(A)\text{-Mod}(\mathrm{Alg}^{\mathrm{fact}}(X)) \, ,$$

its restriction to the central fibre agrees with the module structure of Equation 1.3.1, where $\mathrm{CC}_{\bullet}^-(A) \in \mathrm{Alg}^{\mathrm{fact}}_{\mathbb{E}_1}(X)_{/\mathbb{K}[u]}$ denotes the negative cyclic chains on $A$, considered as a family of factorization $\mathbb{E}_1$ algebras over $\mathbb{K}[u] = \mathrm{H}^{\bullet}_{S^1}(\mathrm{pt})$ with central fibre $\mathrm{CH}_{\bullet}(A) \in \mathrm{Alg}^{\mathrm{fact}}_{\mathbb{E}_1}(X)$ .

In the present work, we use this to explain the relationship between the construction of chiral algebras corresponding to four dimensional $\mathcal{N} = 2$ superconformal field theories in [BLL+15], which we give a mathematical account of in terms of equivariant factorization algebras, and the construction of boundary chiral algebras for (holomorphic-)topological twists of three dimensional $\mathcal{N} = 4$ theories following [CG18], which is the other central topic of this paper.

1.4. **Overview of Part II.** We now give an overview of the present work. There are two overarching goals of this paper, which are accomplished in Chapters 1 and 2, respectively:

(1) We develop tools to give geometric constructions of factorization $\mathbb{E}_n$ algebras.
(2) We apply the theory of equivariant factorization algebras, as developed in the prequel [But20a], in examples that describe holomorphic-topological twists of supersymmetric gauge theories, constructed as in the preceding chapter. In this language, we formulate and prove several variants of predictions from the physics literature about these quantum field theories.

We now give a more detailed overview of these chapters:

1.4.1. *Overview of Chapter 1.* Chapter 1 is of a primarily technical nature: we review the notions of factorization categories and functors, as well as factorization spaces and maps of such, following the main construction of [Ras15a]. Moreover, we explain various linearization constructions which use sheaf theory to produce the former from the latter. In particular, in Section 8, we give a construction of factorization $\mathbb{E}_n$ algebras in terms of spaces that are equipped with a compatible factorization and convolution structure. This construction is the main source of the examples in Chapter 2.

1.4.2. *Overview of Chapter 2.* In Chapter 2, we study various examples of equivariant factorization $\mathbb{E}_n$ algebras corresponding to holomorphic-topological twists of supersymmetric gauge theories equipped with an $\Omega$ background. We begin with a breif review of the constructions of the $\mathbb{E}_2$ and $\mathbb{E}_3$ algebras corresponding to the two and three dimensional topological $B$ model, as a warm up for the more involved examples which are the primary focus of the present work.

In sections 12 and 13, we recall the construction of the factorization algebra corresponding to the three dimensional $A$ model, following [BFN18]: Let $C$ be a smooth algebraic curve, $\mathbb{D}$ the formal disk, $G$ an affine, reductive algebraic group with Lie algebra $\mathfrak{g}$, and $N$ a finite dimensional $G$ representation. The Coulomb branch construction introduced in *loc. cit.* gives a factorization algebra on $C$ with compatible associative (or $\mathbb{E}_1$) algebra structure, defined by

$$(1.4.1) \qquad \mathcal{A}(G, N) := \mathrm{C}^{\mathrm{BM}}_{\bullet}(Y_{\mathcal{O}} \times_{Y_{\mathcal{K}}} Y_{\mathcal{O}}) \in \mathrm{Alg}^{\mathrm{fact}}_{\mathbb{E}_1}(C) \cong \mathrm{Alg}_{\mathbb{E}_1}(\mathrm{Alg}^{\mathrm{fact}}(C)) \, ,$$

where $Y = N/G$ as a stack, $Y_{\mathcal{O}} = \mathrm{Maps}(\mathbb{D}, Y)$, $Y_{\mathcal{K}} = \mathrm{Maps}(\mathbb{D}^{\circ}, Y)$, and $\mathrm{C}^{\mathrm{BM}}_{\bullet}$ are the Borel-Moore chains. The factorization $\mathbb{E}_1$ algebra $\mathcal{A}(G, N) \in \mathrm{Alg}^{\mathrm{fact}}_{\mathbb{E}_1}(C)$ describes the local observables of the three dimensional $A$ model gauge theory on $C \times \mathbb{R}$ with gauge group $G$ and matter representation $T^\vee N$. Equivalently, we can think of this theory as a sigma model with target space the cotangent stack $T^\vee Y$.

In [BFN18], it is explained that this construction also gives a filtered quantization of a graded Poisson algebra, which in good cases describes a quantization of the symplectic singularity which is dual to $T^\vee Y$ in the sense of symplectic duality [BLPW14], or three dimensional mirror symmetry [IS96]; this was the original motivation for the construction. The relationship between these results is an example of the equivalence of Theorem 1.3.1:

*Theorem* 1.4.1. [BFN18] For $C = \mathbb{A}^1$ the factorization $\mathbb{E}_1$ algebra defined above

$$\mathcal{A}(G, N) \in \mathrm{Alg}^{\mathrm{fact}}_{\mathbb{E}_1}(\mathbb{A}^1)^{\mathbb{G}_a \rtimes \mathbb{G}_m} \cong \mathrm{Alg}_{\mathbb{E}_3^{S^1}}(\mathrm{Vect}_{\mathbb{K}}) \cong \mathrm{Alg}_{\mathbb{B}\mathbb{D}_1^u}(\mathrm{D}(\mathbb{K}[u]))$$

admits a canonical $\mathbb{G}_a \rtimes \mathbb{G}_m$ equivariant structure and thus, under the equivalence Theorem 1.3.1, defines a filtered quantization of a (2-shifted) Poisson algebra to an associative (or $\mathbb{E}_1$) algebra.

Concretely, passing to $\mathbb{G}_m$ equivariant (with respect to loop rotation) Borel-Moore homology in the expression from Equation 1.4.1 gives a quantization of the homology $\mathbb{P}_3$ algebra, viewed as a graded commutative algebra with Poisson bracket of degree $-2$, as in [BFN18]. We also explain an analogous construction of the three dimensional $B$ model in [But20b], which gives a filtered quantization of $T^\vee Y$ itself by this mechanism.

In sections 14, 15, 16 and 17, we study the factorization algebra on $C$ of chiral differential operators $\mathcal{D}^{\mathrm{ch}}(Y) \in \mathrm{Alg}^{\mathrm{fact}}(C)$ on $Y = N/G$, and their relationship to the three dimensional $A$ model above, culminating in a proof in this language of a prediction of Costello-Gaiotto from [CG18].

On a smooth curve $C$, the factorization algebra of chiral differential operators on a scheme $Y$ is defined, following [GMS00, BD04, KV06], by

$$\mathcal{D}^{\mathrm{ch}}(Y) := \mathcal{H}\mathrm{om}_{D(Y_{\mathbb{K}})}(\mathcal{D}_{Y_{\mathbb{K}}}, \iota_* \omega_{Y_\mathbb{O}}) \cong \Gamma(Y_{\mathbb{K}}, \mathrm{oblv}\, \iota_* \omega_{Y_\mathbb{O}}) \in \mathrm{Alg}^{\mathrm{fact}}(C)$$

where oblv $: D(Y_{\mathbb{K}}) \to \mathrm{IndCoh}(Y_{\mathbb{K}})$ denotes the forgetful functor from $D$ modules to coherent sheaves. In the case that $Y = N/G$ is a quotient stack, we show that the chiral differential operators are given by

$$\mathcal{D}^{\mathrm{ch}}(Y) \cong \mathrm{C}^{\frac{\infty}{2}}(\hat{\mathfrak{g}}, \mathfrak{g}_\mathbb{O}, G_\mathbb{O}; \mathcal{D}^{\mathrm{ch}}(N))$$

the $G_\mathbb{O}$-equivariant semi-infinite cohomology, or BRST reduction, of $\mathcal{D}^{\mathrm{ch}}(N) \in \mathrm{Alg}^{\mathrm{fact}}(C)$ with respect to the Lie algebra $\hat{\mathfrak{g}}$. The approach to sheaf theory on the relevant infinite dimensional varieties and stacks, as well as the identification with semi-infinite cohomology, follows [Ras15b, Ras20b] and references therein.

*Remark* 1.4.2. For $Y$ a scheme this construction requires a trivialization of the determinant gerbe [KV06], which for $Y = N/G$ we identify with a lift of the $G_\mathbb{O}$ action on $\mathcal{D}^{\mathrm{ch}}(N)$ to an action of $\hat{\mathfrak{g}}$ at level $-$Tate. Physically, this corresponds to the requirement that the corresponding four dimensional $\mathcal{N} = 2$ theory is superconformal.

Under the hypothesis of the preceding remark, we prove the prediction of Costello-Gaiotto [CG18]:

*Theorem* 1.4.3. The chiral differential operators on $Y = N/G$ admits a canonical module structure

$$\mathcal{D}^{\mathrm{ch}}(Y) \in \mathcal{A}(G, N)\text{-Mod}(\mathrm{Alg}^{\mathrm{fact}}(C))$$

over the factorization $\mathbb{E}_1$ algebra $\mathcal{A}(G, N) \in \mathrm{Alg}^{\mathrm{fact}}_{\mathbb{E}_1}(C)$ on $C$ constructed in [BFN18].

This result corresponds to the statement that the three dimensional $A$ model to $Y$ admits a chiral boundary condition, so that the algebra of local observables $\mathcal{A}(G, N)$ of the three dimensional theory on $C \times \mathbb{R}_{\geqslant 0}$ acts on the chiral algebra $\mathcal{D}^{\mathrm{ch}}(Y)$ of boundary observables on $C$, as pictured on the right.

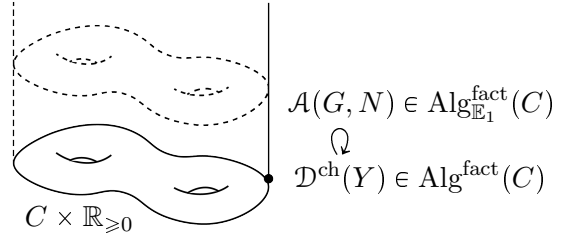

$$\mathcal{A}(G, N) \in \mathrm{Alg}_{\mathbb{E}_1}^{\mathrm{fact}}(C)$$
$$\circlearrowright$$
$$\mathcal{D}^{\mathrm{ch}}(Y) \in \mathrm{Alg}^{\mathrm{fact}}(C)$$

The preceding theorem is equivalent to the existence of a map $\mathcal{A}(G, N) \to \mathrm{HC}^{\bullet}(\mathcal{D}^{\mathrm{ch}}(Y))$ of factorization $\mathbb{E}_1$ algebras, where the latter denotes the Hochschild cochains (or derived centre) of the factorization algebra $\mathcal{D}^{\mathrm{ch}}(Y)$. It was suggested in [CG18] that this is often an isomorphism, and since $\mathcal{D}^{\mathrm{ch}}(G, N)$ is manifestly independent of the Lagrangian splitting of $V = T^{\vee} N$, we obtain:

*Corollary* 1.4.4. In such situations, $\mathcal{A}(G, N)$ is canonically independent of the choice of splitting $V = T^{\vee} N$.

A related construction which also implies this corollary was suggested independently by Raskin.

In Section 18, we also construct a family of factorization $\mathbb{E}_1$ algebras $\mathcal{C}(Y)^{\hbar} \in \mathrm{Alg}_{\mathbb{E}_1}^{\mathrm{fact}}(X)_{/\mathbb{K}[\hbar]}$ over $\mathbb{K}[\hbar]$, with generic fibre $\mathcal{C}(Y)^{\hbar}|_{\{\hbar=1\}} = \mathcal{A}(Y)$ the factorization $\mathbb{E}_1$ algebra of the three dimensional $A$ model, and central fibre that of the holomorphic-B twist. Moreover, we show that the module structure on $\mathcal{D}^{\mathrm{ch}}(Y)$ over $\mathcal{A}(Y)$ given above extends to a module structure over $\mathcal{C}(Y)^{\hbar}$,

(1.4.2) $$\mathcal{D}^{\mathrm{ch}}(Y)_{\hbar} \in \mathcal{C}(Y)^{\hbar}\text{-Mod}(\mathrm{Alg}^{\mathrm{fact}}(C)_{/\mathbb{K}[\hbar]}) \ ,$$

where $\mathcal{D}^{\mathrm{ch}}(Y)_{\hbar} \in \mathrm{Alg}^{\mathrm{fact}}(C)_{/\mathbb{K}[\hbar]}$ is the filtered quantization of chiral differential operators to $Y$.

Finally, in Section 19, we use the theory of equivariant factorization algebras developed in Part I to give a mathematical account of the construction of chiral algebras corresponding to 4d $\mathcal{N} = 2$ gauge theories introduced in [BLL$^+$15]. Further, we explain the relation of this construction with our formulation of the predictions of [CG18] explained above.

Let $\mathcal{Z}(Y) = Y_{\mathcal{O}} \times_{Y_{\mathcal{K}}} Y_{\mathcal{O}}$ be the stack underlying the construction of $\mathcal{A}(G, N)$ from [BFN18] recalled above. The category $\mathrm{IndCoh}(\mathcal{Z}(Y)) \in \mathrm{Cat}_{\mathbb{E}_1,\mathrm{un}}^{\mathrm{fact}}(C)$ of (ind)coherent sheaves on $\mathcal{Z}(Y)$, studied in [CW19] for $Y = BG$, defines a factorization $\mathbb{E}_1$ category, with unit factorization algebra $\mathrm{u}_{\bullet} \mathcal{O}_{Y_{\mathcal{O}}} \in \mathrm{Alg}_{\mathbb{E}_1}^{\mathrm{fact}}(\mathrm{IndCoh}(\mathcal{Z}(Y)))$, so that

$$\mathcal{F}(Y) = \mathcal{H}om_{\mathrm{IndCoh}(\mathcal{Z}(Y))}(\mathrm{u}_{\bullet} \mathcal{O}_{Y_{\mathcal{O}}}, \mathrm{u}_{\bullet} \mathcal{O}_{Y_{\mathcal{O}}}) \ \in \ \mathrm{Alg}_{\mathbb{E}_2}^{\mathrm{fact}}(C)$$

defines a factorization $\mathbb{E}_2$ algebra on $C$, which describes the local operators of the four dimensional mixed holomorphic-B model on $C \times \mathbb{R}^2$, which occurs as the corresponding mixed twist of four dimensional $\mathcal{N} = 2$ gauge theory. For $Y$ satisfying the hypotheses of Remark 1.4.2, we have:

*Theorem* 1.4.5. There is a canonical $S^1$ equivariant structure on $\mathcal{F}(Y)$ such that

$$\mathcal{F}(Y) \mapsto \mathcal{D}^{\mathrm{ch}}(Y)_u \qquad \text{under the equivalence} \qquad \mathrm{Alg}_{\mathbb{E}_2^{S^1}}^{\mathrm{fact}}(C) \cong \mathrm{Alg}_{\mathbb{BD}_0}^{\mathrm{fact}}(C)$$

of Theorem 1.3.2, where $\mathcal{D}^{\mathrm{ch}}(Y)_u \in \mathrm{Alg}_{\mathbb{BD}_0^u}^{\mathrm{fact}}(X)$ is the (two-periodic) filtered quantization of the factorization algebra of chiral differential operators to $Y$.

Let $\mathcal{C}(Y)^u \in \mathrm{Alg}_{\mathbb{E}_1}^{\mathrm{fact}}(C)_{/\mathbb{K}[u]}$ denote (the two-periodic variant of) the family of factorization $\mathbb{E}_1$ algebras claimed above, which describes the deformation from the holomorphic to the A twist of 3d $\mathcal{N} = 4$ gauge theory. The following result identifies this with the equivariant reduction on $S^1$ of holomorphic-B twisted 4d $\mathcal{N} = 2$ gauge theory, and thus relates the results of [CG18] and [BLL$^+$15] via their fomulations in the present work:

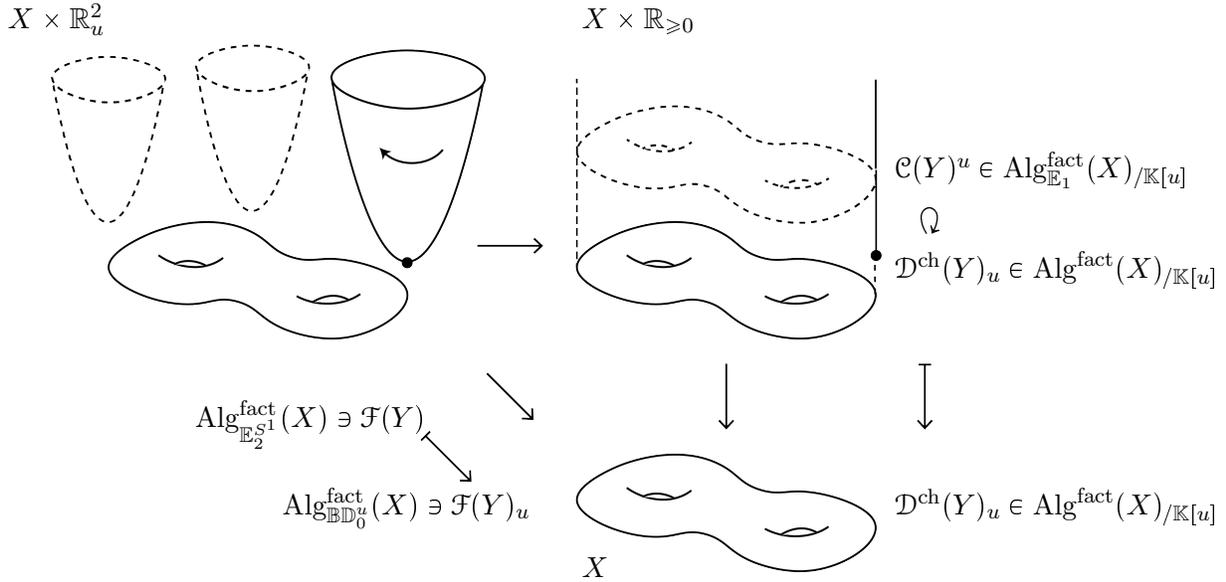

FIGURE 2. The equivariant cigar reduction principle relates the four dimensional holomorphic-B model to the three dimensional A model

*Theorem* 1.4.6. There is an equivalence of families of factorization $\mathbb{E}_1$ algebras on $X$

$$\mathrm{CC}_\bullet^-(\mathcal{F}(Y)) \xrightarrow{\cong} \mathcal{C}(Y)^u \in \mathrm{Alg}_{\mathbb{E}_1}^{\mathrm{fact}}(C)_{/\mathbb{K}[u]} \ ,$$

such that under the equivalence of Theorem 1.4.5, the module structure of the preceding proposition

$$\mathcal{F}(Y)_u \in \mathrm{CC}_\bullet^-(\mathcal{F}(Y))\text{-}\mathrm{Mod}(\mathrm{Alg}^{\mathrm{fact}}(C)_{/\mathbb{K}[u]}) \qquad \text{identifies with} \qquad \mathcal{D}^{\mathrm{ch}}(Y)_u \in \mathcal{C}(Y)^u\text{-}\mathrm{Mod}(\mathrm{Alg}^{\mathrm{fact}}(C)_{/\mathbb{K}[u]}) \ ,$$

equipped with the module structure recalled in Equation 1.4.2 above.

In fact, the preceding theorem does more than just identify the equivariant $S^1$ reduction $\mathrm{CC}_\bullet^-(\mathcal{F}(Y))$ of the holomorphic-B twist of four dimensional $\mathcal{N}=2$ gauge theory with the deformation from the holomorphic-B twist to the A twist of three dimensional $\mathcal{N}=4$ gauge theory. It also identifies the canonical family of modules over the former corresponding to the family of boundary condition induced by the 'cigar tip' as in Proposition 1.3.3, with the family of boundary conditions for the latter whose observables form the filtered quantization of the factorization algebra of chiral differential operators explained above. The reader should compare Figure 2 above with the general case illustrated in Figure 1 of the preceding section.

1.5. **Acknowledgements.** I am especially thankful to Kevin Costello, for first introducing me to so many of the ideas underlying this series of papers, as well as to Sam Raskin, for patient explanations about various more technical aspects. I would also like to thank David Ben-Zvi, Alexander Braverman, Davide Gaiotto, Justin Hilburn, Surya Raghavendran, Pavel Safronov, Brian Williams, and Philsang Yoo for useful discussions.

1.6. **References to the companion paper.** As we have mentioned, this paper is the second among at least two papers in this series, and we will systematically refer to propositions, definitions, etc. from each of the companion papers in the other. In this paper, references to the present text are given by red hyperlinks, such as 4.2.3 which refers to the remark of that enumeration. References to Part I are given by turquoise hyperlinks (which link to the companion pdf if both files are in the same folder) and their enumeration is prefaced by a I, such as I-4.2.3.

# Chapter 1
# Factorization Objects

## 2. Factorization objects

**2.1. Overview.** In this section we outline a construction that defines uniformly a large class of factorization *objects*, including factorization algebras, spaces, and categories. The construction was implicitly introduced in [Ras15a] in order to give a homotopy coherent definition of the latter, and we thank Sam Raskin for many helpful discussions about this material; of course, any inaccuracies or errors in the present work are a result of the Author's misunderstanding of the material.

*Warning* 2.1.1. There are certain 2-categorical details required in the main construction of this section, especially in its application to factorization categories and their variants, which are treated in [Ras15a]. We will ignore these details throughout the general discussion of the present section, but in each of the applications in the following sections we give the definitions and results which follow from the careful treatment of *loc. cit.*.

The minimal input data for the construction is a lax symmetric monoidal functor $\mathrm{F}: \mathrm{Sch}_{\mathrm{aff}}^{\mathrm{op}} \to \mathcal{C}$, where $\mathcal{C} = \mathrm{Cat}^{\times}$ or some variant thereof and $\mathrm{Sch}_{\mathrm{aff}}$ is given the cartesian monoidal structure. For an affine scheme $S \in \mathrm{Sch}_{\mathrm{aff}}$, the object $\mathrm{F}(S) \in \mathcal{C}$ is interpreted as the category of 'F-objects over $S$'.

For example, if $\mathrm{F} = \mathrm{QCoh}: \mathrm{Sch}_{\mathrm{aff}}^{\mathrm{op}} \to \mathrm{DGCat}_{\mathrm{cont}}$ then $\mathrm{F}(S) = \mathrm{QCoh}(S) \in \mathrm{DGCat}$ is the DG category of quasicoherent sheaves on $S$, or if $F = \mathrm{ShvCat}: \mathrm{Sch}_{\mathrm{aff}}^{\mathrm{op}} \to 2\mathrm{Cat}$ then $\mathrm{F}(S) = \mathrm{ShvCat}(S)$ is the (2-)category of quasicoherent sheaves of categories on $S$. The condition that $F$ is lax symmetric monoidal is interpreted as the existence of a notion of exterior product of F-objects

$$\boxtimes: \mathrm{F}(S) \otimes \mathrm{F}(T) \to \mathrm{F}(S \times T)$$

given by the usual exterior product operation in the previous examples. The requirement of functoriality from $\mathrm{Sch}^{\mathrm{op}}$ is interpreted as the existence of a 'pullback functor on F-objects'

$$f^* := \mathrm{F}(f): \mathrm{F}(S) \to \mathrm{F}(T) \qquad \text{for each} \qquad f: T \to S .$$

We take the Kan extension of F to define $\mathrm{F}: \mathrm{PreStk}^{\mathrm{op}} \to \mathcal{C}$ which is again symmetric monoidal with respect to the cartesian monoidal structure, so that in particular we obtain a category $\mathrm{F}(\mathrm{Ran}_X) \in \mathcal{C}$ of F-objects on $\mathrm{Ran}_X$.

Recall from the introduction to Subsection I-5.3 that the operation of disjoint union of finite subsets defines a correspondence of prestacks

$$\mathrm{Ran}_X^{\times 2} \xleftarrow{j_{\mathrm{disj}}} (\mathrm{Ran}_X^{\times 2})_{\mathrm{disj}} \xrightarrow{\sqcup} \mathrm{Ran}_X ,$$

which underlies the definition of the chiral tensor structure on $D(\mathrm{Ran}_X)$, and in turn the definition of factorization algebra. More generally, this will be the primary structure on $\mathrm{Ran}_X$ used in the definition of factorization F-objects.

We can now give the heuristic definition of a non-unital factorization F object:

*Tentative Definition* 2.1.2. A factorization F-object $\Psi$ on $X$ is:

- an object $\Psi \in \mathrm{F}(\mathrm{Ran}_X)$, and
- an isomorphism $\alpha: j_{\mathrm{disj}}^* \Psi^{\boxtimes 2} \xrightarrow{\cong} \sqcup^* \Psi$ of objects in $\mathrm{F}((\mathrm{Ran}_X^{\times 2})_{\mathrm{disj}})$.

The isomorphism $\alpha$ encodes the factorization structure maps for the object $\Psi$. In particular, in the case $F := D^! : \mathrm{Sch}^{\mathrm{op}}_{\mathrm{aff}} \to \mathrm{DGCat}_{\mathrm{cont}}$ as in Definition B.6.3, the preceding definition agrees with the definition I-6.0.2 of non-unital factorization algebra.

2.1.1. *Unital structures and lax prestacks.* In order to define *unital* factorization objects, an additional technical generalization is required: the structure maps

$$\Delta(\pi)^! A_I \to A_J \qquad \text{for each} \qquad \pi : I \twoheadrightarrow J$$

required in the definition I-4.3.1 of a unital $D$ module $A \in D(\mathrm{Ran}_{X,\mathrm{un}})$ on $\mathrm{Ran}_X$ are only required to be isomorphisms in the case $\pi : I \twoheadrightarrow J$ is a surjection, and otherwise are typically not equivalences unless $A = \omega_{\mathrm{Ran}_{X,\mathrm{un}}}$ is the unit factorization algebra. As such, these maps can not be encoded as the gluing data for a geometric object (in this case, a $D$ module) over the moduli space $\mathrm{Ran}_X \sqcup \{\varnothing\}$ of (possibly empty) subsets of $X$, viewed as an ordinary prestack.

To address this issue, [Ras15a] employs the notion of *lax* prestack, which is essentially just a functor $\mathcal{Y} : \mathrm{Sch}^{\mathrm{op}}_{\mathrm{aff}} \to \mathrm{Cat}$ valued in $\mathrm{Cat}$ rather than $\mathrm{Grpd}$; it is explained in *loc. cit.* that a functor $F : \mathrm{Sch}^{\mathrm{op}}_{\mathrm{aff}} \to \mathcal{C}$ as above can also be appropriately Kan extended to lax prestacks, so that we can define a category $F(\mathcal{Y}) \in \mathcal{C}$ of $F$-objects on a lax prestack $\mathcal{Y}$, such that the gluing data for an object in $F(\mathcal{Y})$ can require structure morphisms which are not necessarily invertible.

Moreover, *loc. cit.* constructs a lax prestack $\mathrm{Ran}_{X,\mathrm{un}} \in \mathrm{PreStk}^{\mathrm{lax}}$ which encodes the moduli space of (possibly empty) finite subsets of $X$ together with the additional data of inclusion relations between subsets, in such a way that for $F = D^!$ as above, the category $D(\mathrm{Ran}_{X,\mathrm{un}})$ of $D$ modules on $\mathrm{Ran}_{X,\mathrm{un}}$ agrees with the same-denoted category of Definition I-4.3.1.

The corresponding heuristic definition of a *unital* factorization $F$ object is as follows:

*Tentative Definition* 2.1.3. A unital factorization $F$-object $\Psi$ on $X$ is:

- an object $\Psi \in F(\mathrm{Ran}_{X,\mathrm{un}})$,
- an isomorphism $\alpha : j^*_{\mathrm{disj}} \Psi^{\boxtimes 2} \xrightarrow{\cong} \sqcup^* \Psi$ of objects in $F((\mathrm{Ran}^{\times 2}_{X,\mathrm{un}})_{\mathrm{disj}})$, and
- an isomorphism $\beta : u_{F(\mathrm{pt})} \xrightarrow{\cong} \iota^*_{\varnothing} \Psi$ of objects in $F(\mathrm{pt})$,

where $\iota_{\varnothing} : \mathrm{pt} \hookrightarrow \mathrm{Ran}_{X,\mathrm{un}}$ is the map of lax prestacks corresponding to the inclusion of the empty set, and $u_{F(\mathrm{pt})} \in F(\mathrm{pt})$ is the unit for the canonical symmetric monoidal structure on $F(\mathrm{pt})$.

## 2.2. $F$-**Objects on lax prestacks.**

Let $\mathcal{C} = \mathrm{Cat}^{\otimes}_{\mathrm{pres}}$, $\mathrm{Cat}^{\times}$, or $\mathrm{DGCat}^{\otimes}_{\mathrm{cont}}$ equipped with the indicated symmetric monoidal structure; see the conventions in Section I-2. Fix a lax symmetric monoidal functor $F : \mathrm{Sch}^{\mathrm{op}}_{\mathrm{aff}} \to \mathcal{C}$, interpreted as assigning to $Y \in \mathrm{Sch}_{\mathrm{aff}}$ the category $F(Y) \in \mathcal{C}$ of $F$-objects on $\mathcal{Y}$.

*Example* 2.2.1. The functor $F = \mathrm{QCoh}^{\bullet} : \mathrm{Sch}^{\mathrm{op}}_{\mathrm{aff}} \to \mathrm{DGCat}^{\otimes}_{\mathrm{cont}}$ of quasicoherent sheaves under usual pullback is lax symmetric monoidal via the usual exterior product of quasicoherent sheaves.

*Example* 2.2.2. The functor $F = D^! : \mathrm{Sch}^{\mathrm{op}}_{\mathrm{aff}} \to \mathrm{DGCat}^{\otimes}_{\mathrm{cont}}$ of $D$ modules under usual pullback is lax symmetric monoidal with respect to the usual exterior product of $D$ modules; see Definition B.6.2.

*Example* 2.2.3. The functor $F = \mathrm{PreStk}_{/(\cdot)} : \mathrm{Sch}^{\mathrm{op}}_{\mathrm{aff}} \to \mathrm{Cat}^{\times}$ which assigns to $S \in \mathrm{Sch}_{\mathrm{aff}}$ the category $\mathrm{PreStk}_{/S}$ of prestacks over $S$ is lax symmetric monoidal with respect to the product of prestacks

$$\times : \mathrm{PreStk}_{/S} \times \mathrm{PreStk}_{/T} \to \mathrm{PreStk}_{/T \times S} \ .$$

See Appendix B.1 and references therein for background on prestacks.

*Example* 2.2.4. The functor $F = \mathrm{ShvCat} : \mathrm{Sch}_{\mathrm{aff}}^{\mathrm{op}} \to \mathrm{Cat}^{\times}$ which assigns to $S \in \mathrm{Sch}_{\mathrm{aff}}$ the category

$$\mathrm{ShvCat}(S) := \mathrm{QCoh}(S)\text{-}\mathrm{Mod}(\mathrm{DGCat}_{\mathrm{cont}}^{\otimes})$$

of (quasicoherent) sheaves of categories on $S$. This functor is lax symmetric monoidal with respect to the exterior product of sheaves of categories, which is a direct generalization of that of quasicoherent sheaves; see Appendix A in [Ras15a] or [Gai15] for background on sheaves of categories.

*Remark* 2.2.5. Although we require only the stated hypotheses on $\mathrm{F} : \mathrm{Sch}_{\mathrm{aff}}^{\mathrm{op}} \to \mathcal{C}$ in order to define the notion of factorization F object below, several further properties of the functor F will be necessary in order to implement non-trivial constructions with factorization F objects. For example, several of the above functors are part of a (partially defined) six functors formalism in the sense of Subappendix A.3, and this structure will play a crucial role in most applications. This structure can be homotopy coherently summarized by the requirement that the functors naturally extend to variants of the correspondence category $\mathrm{Sch}_{\mathrm{corr}}^{\mathrm{aff}}$ in a way that encodes the six functors compatibilities; this perspective is explained in detail in [GR17a].

*Warning* 2.2.6. In the following sections will carry out several constructions using hypotheses as in the preceeding remark. Although we do indicate which properties are required in each case, we do not give a systematic, homotopy coherent treatment in terms of the formalism of correspondence categories developed [GR17a].

Recall there is a canonical functor $\mathrm{Sch}_{\mathrm{aff}} \hookrightarrow \mathrm{PreStk}^{\mathrm{lax}}$; see Appendix B.1 and references therein for a review of prestacks.

*Definition* 2.2.7. The functor $\mathrm{F} : (\mathrm{PreStk}^{\mathrm{lax}})^{\mathrm{op}} \to \mathcal{C}$ assigning to a lax prestack $\mathcal{Y}$ the category $\mathrm{F}(\mathcal{Y})$ of F-objects on $\mathcal{Y}$ is defined by right Kan extension of $\mathrm{F} : \mathrm{Sch}_{\mathrm{aff}}^{\mathrm{op}} \to \mathcal{C}$ along $\mathrm{Sch}_{\mathrm{aff}}^{\mathrm{op}} \hookrightarrow (\mathrm{PreStk}^{\mathrm{lax}})^{\mathrm{op}}$.

*Remark* 2.2.8. Concretely, by pointwise evaluation of the opposite left Kan extension as a weighted colimit, an object of $M \in \mathrm{F}(\mathcal{Y})$ is given by an assignment

$$S \mapsto \big[(\cdot)^*(M): \ \mathrm{Maps}(S,\mathcal{Y}) \longrightarrow \mathrm{F}(S) \ \big] \qquad [\varphi : S \to T] \mapsto \begin{array}{c} \mathrm{Maps}(T,\mathcal{Y}) \xrightarrow{(\cdot)^*(M)} \mathrm{F}(T) \\ {\scriptstyle (\cdot)\circ\varphi} \big\downarrow \qquad\quad \varphi^* \big\downarrow \\ \mathrm{Maps}(S,\mathcal{Y}) \xrightarrow{(\cdot)^*(M)} \mathrm{F}(S) \end{array} ,$$
$$\text{``}\, [f : S \to \mathcal{Y}] \longmapsto f^*(M) \,\text{``}$$

defined for each $S \in \mathrm{Sch}_{\mathrm{aff}}$ and each morphism $\varphi : S \to T$ of affine schemes. Note that $\mathrm{Maps}(S,\mathcal{Y}) \in \mathrm{Cat}$ can contain non-invertible morphisms; this is the primary additional generality provided by the use of lax prestacks here, as explained above.

## 2.3. Multiplicative F-objects on commutative monoids in lax prestacks.
Let $\mathrm{PreStk}_{\mathrm{corr}}$ be the category of lax prestacks under correspondences. The cartesian monoidal structure on $\mathrm{PreStk}$ extends canonically to $\mathrm{PreStk}_{\mathrm{corr}}$ to define a symmetric monoidal structure. The disjoint union operation on $\mathrm{Ran}_X$ underlying the definition of factorization, as discussed in the overview of this section, as well as the geometric structure on $\mathrm{Ran}_{X,\mathrm{un}}$ underlying the unital variant, are encoded in terms of a (unital) monoid structure on $\mathrm{Ran}_X$ ($\mathrm{Ran}_{X,\mathrm{un}}$) in the correspondence category. The factorization structure on an F-object on $\mathrm{Ran}_X$ is then encoded as 'multiplicativity' data with respect to this monoid structure.

Let $\mathrm{Comm}(\mathrm{PreStk}_{\mathrm{corr}}^{\mathrm{lax}})$ denote the category of commutative monoid objects in lax prestacks under correspondences.

*Remark* 2.3.1. Concretely, a commutative monoid object $\mathcal{Y} \in \mathrm{Comm}(\mathrm{PreStk}^{\mathrm{lax}}_{\mathrm{corr}})$ in lax prestacks is given by an underlying prestack $\mathcal{Y}$, together with correspondences of such

$$
\begin{array}{ccc}
& \mathrm{mult}_{\mathcal{Y}} & \\
{}^{m_1}\swarrow & & \searrow^{m_2} \\
\mathcal{Y}^{\times 2} & & \mathcal{Y}
\end{array}
\qquad
\begin{array}{ccc}
& \mathrm{unit}_{\mathcal{Y}} & \\
{}^{e_1}\swarrow & & \searrow^{e_2} \\
\mathrm{pt} & & \mathcal{Y}
\end{array}
\quad,
$$

and compatible higher arity analogues of the multiplication map, satisfying the relations of a unital commutative monoid.

*Example* 2.3.2. The unital Ran space $\mathrm{Ran}_{X,\mathrm{un}} \in \mathrm{Comm}(\mathrm{PreStk}^{\mathrm{lax}}_{\mathrm{corr}})$ is a commutative monoid in lax prestacks, with structure maps given by

$$
\begin{array}{ccc}
& (\mathrm{Ran}^{\times 2}_{X,\mathrm{un}})_{\mathrm{disj}} & \\
{}^{j_{\mathrm{disj}}}\swarrow & {\scriptstyle \sqcup}\downarrow & \\
\mathrm{Ran}^{\times 2}_{X,\mathrm{un}} & & \mathrm{Ran}_{X,\mathrm{un}}
\end{array}
\qquad
\begin{array}{ccc}
& \mathrm{pt} & \\
\swarrow & & \searrow^{\iota_\varnothing} \\
\mathrm{pt} & & \mathrm{Ran}_{X,\mathrm{un}}
\end{array}
\quad,
$$

where $\iota_\varnothing : \mathrm{pt} \hookrightarrow \mathrm{Ran}_{X,\mathrm{un}}$ is the map corresponding to the inclusion of the empty set.

Toward giving the formal definition of a multiplicative object on a commutative lax prestack, we give a heuristic recollection of the main construction from Section 5 of [Ras15a].

*Definition* 2.3.3. Let $\mathcal{I} \in \mathrm{Cat}$ be an abstract index category and $\mathrm{F} : \mathcal{I}^{\mathrm{op}} \to \mathcal{C}$. Define the correspondence Grothendieck construction $\mathrm{Groth}_{\mathrm{corr}}(\mathrm{F}) \in \mathrm{Cat}$ of $\mathrm{F}$ as the category with:

* objects $\mathrm{ob}(\mathrm{Groth}_{\mathrm{corr}}(\mathrm{F}) \in \mathrm{Cat})$ given by pairs $i \in \mathcal{I}$ and $\Psi_i \in F(i)$, and
* maps $\mathrm{Hom}_{\mathrm{Groth}_{\mathrm{corr}}(\mathrm{F}) \in \mathrm{Cat}}((i, \Psi_i), (j, \Psi_j))$ from $((i, \Psi_i)$ to $(j, \Psi_j))$ given by correspondences

$$
\begin{array}{ccc}
& h & \\
{}^{\alpha}\swarrow & & \searrow^{\beta} \\
i & & j
\end{array}
\qquad \text{together with a map} \qquad \alpha^*(\Psi_i) \to \beta^*(\Psi_j) \qquad \text{in } F(j).
$$

*Remark* 2.3.4. There is a natural forgetful functor $\mathrm{Groth}_{\mathrm{corr}}(\mathrm{F}) \to \mathcal{I}$ defined by $(i, \Psi_i) \mapsto i$.

*Proposition* 2.3.5. Suppose $\mathcal{I}$ is symmetric monoidal and that $F : \mathcal{I} \to \mathcal{C}$ is lax monoidal. Then there is a natural symmetric monoidal structure

$$\mathrm{Groth}_{\mathrm{corr}}(\mathrm{F}) \times \mathrm{Groth}_{\mathrm{corr}}(\mathrm{F}) \to \mathrm{Groth}_{\mathrm{corr}}(\mathrm{F}) \qquad \text{defined on objects by} \qquad ((i, \Psi_i), (j, \Psi_j)) \mapsto (i \otimes j, \Psi_i \boxtimes \Psi_j)$$

where $i \otimes j$ denotes the monoidal structure on $\mathcal{I}$ and $\Psi_i \boxtimes \Psi_j$ denotes the lax monoidal structure map.

*Remark* 2.3.6. Under the hypotheses of the preceeding proposition, the forgetful functor $\mathrm{Groth}_{\mathrm{corr}}(\mathrm{F}) \to \mathcal{I}$ of Remark 2.3.4 is symmetric monoidal.

Let $\mathcal{I} = \mathrm{PreStk}^{\mathrm{lax}}$ be the category of lax prestacks, $F : \mathrm{PreStk}^{\mathrm{op}}_{\mathrm{lax}} \to \mathcal{C}$ be lax monoidal, and $\mathrm{PreStk}^{\mathrm{lax}, F}_{\mathrm{corr}} = \mathrm{Groth}_{\mathrm{corr}}(F)$ denote the symmetric monoidal category given by the correspondence Grothendieck construction above.

*Definition* 2.3.7. Let $\mathcal{Y} \in \mathrm{Comm}(\mathrm{PreStk}_{\mathrm{corr}}^{\mathrm{lax}})$ be a commutative monoid in lax prestacks under correspondences. A weakly multiplicative F-object $\Psi$ on $\mathcal{Y}$ is a commuative monoid object $(\mathcal{Y}, \Psi) \in \mathrm{Comm}(\mathrm{PreStk}_{\mathrm{corr}}^{\mathrm{lax},F})$ with image under the forgetful functor of Remark 2.3.4 given by $\mathcal{Y} \in \mathrm{Comm}(\mathrm{PreStk}_{\mathrm{corr}}^{\mathrm{lax}})$.

Let $\mathrm{F}^{\mathrm{mult},w}(\mathcal{Y})$ denote the category of weakly multiplicative F-objects on $\mathcal{Y}$. Throughout, we will abuse notation and refer to the underlying object $\Psi \in \mathrm{F}(\mathcal{Y})$ as the weakly multiplicative object over $\mathcal{Y}$.

*Remark* 2.3.8. Concretely, for $\mathcal{Y} \in \mathrm{Comm}(\mathrm{PreStk}_{\mathrm{corr}}^{\mathrm{lax}})$ with structure maps denoted as in Remark 2.3.1, a weakly multiplicative F-object $\Psi$ on $\mathcal{Y}$ is given by

- An object $\Psi \in \mathrm{F}(\mathcal{Y})$,
- a multiplication map $\eta_m : m_1^* \Psi^{\boxtimes 2} \to m_2^* \Psi$, and
- a unit map $\eta_e : e_1^* u_{\mathrm{F}(\mathrm{pt})} \to e_2^* \Psi$,

together with compatible higher arity analogues of the multiplication map, satisfying the relations of a unital commutative monoid, where $u_{\mathrm{F}(\mathrm{pt})} \in \mathrm{F}(\mathrm{pt})$ is the unit for the canonical symmetric monoidal structure on $\mathrm{F}(\mathrm{pt})$.

*Definition* 2.3.9. Let $\mathcal{Y} \in \mathrm{Comm}(\mathrm{PreStk}_{\mathrm{corr}}^{\mathrm{lax}})$. A multiplicative F-object on $\mathcal{Y}$ is a weakly multiplicative F-object $\Psi \in \mathrm{F}^{\mathrm{mult},w}(\mathcal{Y})$ such that the unit map $\eta_e$, multiplication map $\eta_m$, and its higher arity analogues, are isomorphisms.

Let $\mathrm{F}^{\mathrm{mult}}(\mathcal{Y})$ denote the category of multiplicative F-objects on $\mathcal{Y}$.

## 2.4. **Factorization F-objects.** We can now state the main definition of this section:

*Definition* 2.4.1. Let $F : \mathrm{Sch}_{\mathrm{aff}}^{\mathrm{op}} \to \mathcal{C}$ be lax monoidal. A (unital) factorization F-object $\Psi$ on $X$ is a multiplicative F-object $\Psi \in \mathrm{F}^{\mathrm{mult}}(\mathrm{Ran}_{X,\mathrm{un}})$ on the unital Ran space $\mathrm{Ran}_{X,\mathrm{un}}$ with respect to the commutative monoid structure of Example 2.3.2.

Similarly, a non-unital factorization F-object is a (non-unital) multiplicative F-object on $\mathrm{Ran}_X$.

Although we will not discuss it further until Section 4, we also define the following variant: Let $(\cdot)^{\mathrm{op}} : \mathrm{Cat}^\times \to \mathrm{Cat}^\times$ be the endofunctor of taking the opposite category, so that for $\mathrm{F} : \mathrm{Sch}_{\mathrm{aff}}^{\mathrm{op}} \to \mathrm{Cat}^\times$ we obtain another functor $\mathrm{F}_{\mathrm{op}} := (\cdot)^{\mathrm{op}} \circ \mathrm{F} : \mathrm{Sch}_{\mathrm{aff}}^{\mathrm{op}} \to \mathrm{Cat}^\times$.

*Definition* 2.4.2. A co-unital factorization F-object on $X$ is a unital factorization $\mathrm{F}_{\mathrm{op}}$-object on $X$.

Let $\mathrm{F}_{\mathrm{un}}^{\mathrm{fact}}(X) = \mathrm{F}^{\mathrm{mult}}(\mathrm{Ran}_{X,\mathrm{un}})$ denote the category of (unital) factorization F-objects on $X$, $\mathrm{F}^{\mathrm{fact}}(X) = \mathrm{F}^{\mathrm{mult}}(\mathrm{Ran}_X)$ the category of non-unital factorization F objects, and $\mathrm{F}_{\mathrm{co\text{-}un}}^{\mathrm{fact}}(X)$ the category of co-unital factorization F objects.

*Remark* 2.4.3. Concretely, following Remark 2.3.8, a factorization F-object $\Psi \in \mathrm{F}_{\mathrm{un}}^{\mathrm{fact}}(X)$ on $X$ is

- an object $\Psi \in \mathrm{F}(\mathrm{Ran}_{X,\mathrm{un}})$,
- an isomorphism $\alpha : j_{\mathrm{disj}}^* \Psi^{\boxtimes 2} \xrightarrow{\cong} \sqcup^* \Psi$ of objects in $\mathrm{F}((\mathrm{Ran}_{X,\mathrm{un}}^{\times 2})_{\mathrm{disj}})$, and
- an isomorphism $\beta : u_{\mathrm{F}(\mathrm{pt})} \xrightarrow{\cong} \iota_\emptyset^* \Psi$ of objects in $\mathrm{F}(\mathrm{pt})$,

together with compatible higher arity analogues of the multiplication map $\alpha$ satisfying the relations of a unital commutative monoid. Here $\iota_\emptyset : \mathrm{pt} \hookrightarrow \mathrm{Ran}_{X,\mathrm{un}}$ is the map of lax prestacks corresponding to the inclusion of the empty set, and $u_{\mathrm{F}(\mathrm{pt})} \in \mathrm{F}(\mathrm{pt})$ is the unit for the canonical symmetric monoidal structure on $\mathrm{F}(\mathrm{pt})$.

*Example* 2.4.4. Let $F = D^! : \mathrm{Sch}_{\mathrm{aff}}^{\mathrm{op}} \to \mathrm{DGCat}_{\mathrm{cont}}^{\otimes}$ be the functor of $D$ modules under usual pullback, as in Example 2.2.2. Then the category of factorization F objects $\mathcal{D}_{\mathrm{un}}^{\mathrm{fact}}(X) = \mathrm{Alg}_{\mathrm{un}}^{\mathrm{fact}}(X)$ is the category of unital factorization algebras from Definition I-6.0.2.

In particular, the above construction provides the formal definition of the compatibility between the unit and factorization structure maps, deferred from Remark I-6.0.3.

*Example* 2.4.5. Let $F = \mathrm{QCoh}^{\bullet} : \mathrm{Sch}_{\mathrm{aff}}^{\mathrm{op}} \to \mathrm{DGCat}_{\mathrm{cont}}^{\otimes}$ be the functor of quasicoherent sheaves under usual pullback, as in Example 2.2.1. Then the category of factorization F objects $\mathrm{QCoh}_{\mathrm{un}}^{\mathrm{fact}}(X)$ is the category of unital factorization quasicoherent sheaves; this is a quasicoherent variant of the usual notion of factorization algebra.

*Remark* 2.4.6. For $X = \tilde{X}_{\mathrm{dR}}$, the category $\mathrm{QCoh}_{\mathrm{un}}^{\mathrm{fact}}(X) \cong \mathrm{Alg}_{\mathrm{un}}^{\mathrm{fact}}(\tilde{X})$ of factorization quasicoherent sheaves on $X_{\mathrm{dR}}$ is equivalent to the usual category of factorization algebras on $X$.

*Example* 2.4.7. Let $F = \mathrm{PreStk}_{/(\cdot)} : \mathrm{Sch}_{\mathrm{aff}}^{\mathrm{op}} \to \mathrm{Cat}^{\times}$ be the functor of prestacks over a space, as in Example 2.2.3. Then the category of factorization F objects $\mathrm{PreStk}_{\mathrm{un}}^{\mathrm{fact}}(X)$ is the category of unital factorization prestacks, also called simply unital factorization spaces. This notion is discussed in detail in Section 4.

*Example* 2.4.8. Let $F = \mathrm{ShvCat}(\cdot) : \mathrm{Sch}_{\mathrm{aff}}^{\mathrm{op}} \to \mathrm{Cat}^{\times}$ be the functor of sheaves of categories on a space, as in Example 2.2.4. Then the category of factorization F objects $\mathrm{ShvCat}_{\mathrm{un}}^{\mathrm{fact}}(X)$ is the category of unital factorization categories, also called chiral categories. This notion is discussed in detail in Section 3.

There are some additional technical subtleties required to correctly treat this Example; see Warning 2.1.1.

## 2.5. Universal linearization.

In this subsection, we explain a mechanism for producing factorization G-objects from factorization F-objects on $X$, for distinct lax symmetric monoidal functors $F, G : \mathrm{Sch}_{\mathrm{aff}}^{\mathrm{op}} \to \mathcal{C}$, in a way which can be applied universally over all $X$. The primary application of this procedure is in Subsection 5, and we recommend the reader skip this subsection, and return to it only as necessary while reading *loc. cit.*.

Throughout, let $F, G : \mathrm{Sch}_{\mathrm{aff}}^{\mathrm{op}} \to \mathcal{C}$ be lax symmetric monoidal functors, as above.

*Definition* 2.5.1. A monoidal compatible natural transformation $\eta : F \to G$ is a morphism in the category of lax monoidal functors $\mathrm{Sch}_{\mathrm{aff}}^{\mathrm{op}}$ to $\mathcal{C}$.

*Remark* 2.5.2. Concretely, a monoidal compatible natural transformation is a usual natural transformation $\eta : F \to G$ together with commutativity data for the diagram

$$
\begin{array}{ccc}
F(S) \times F(T) & \xrightarrow{\eta_S \otimes \eta_T} & G(S) \otimes G(T) \\
\Big\downarrow{\boxtimes} & \overset{\sim}{\nearrow} & \Big\downarrow{\boxtimes} \\
F(S \times T) & \xrightarrow{\eta_{S \times T}} & G(S \times T)
\end{array}
$$

natural in $S, T \in \mathrm{Sch}_{\mathrm{aff}}$, and its higher arity analogues, where $\otimes$ denotes the symmetric monoidal structure on $\mathcal{C}$.

*Proposition* 2.5.3. A monoidal compatible natural transformation $\eta : F \to G$ induces a symmetric monoidal functor $\mathrm{PreStk}_{\mathrm{corr}}^{\mathrm{lax}, F} \to \mathrm{PreStk}_{\mathrm{corr}}^{\mathrm{lax}, G}$, intertwining the forgetful functors to $\mathrm{PreStk}_{\mathrm{corr}}^{\mathrm{lax}}$.

*Corollary* 2.5.4. A monoidal compatible natural transformation $\eta : F \to G$ induces a natural functor $\eta_X : F_{un}^{fact}(X) \to G_{un}^{fact}(X)$ from the category of factorization F-objects on $X$ to factorization G-objects.

Recall that a co-unital factorization F object is a unital factorization $F_{op}$-object where $F_{op} : Sch_{aff}^{op} \to \mathcal{C}$, denotes the functor F composed with taking the opposite category of the resulting category.

*Example* 2.5.5. Let $\mathcal{C} = Cat^\times$, $F = PreStk_{/(\cdot)}$, and $G = QCoh(\cdot)$. Define $\eta : F \to G_{op}$ by

$$\eta(S) : PreStk_{/S} \to QCoh(S)^{op} \qquad (\mathcal{Y} \xrightarrow{p} S) \mapsto p_* \mathcal{O}_\mathcal{Y} .$$

Then the above provides a functor from (unital) factorization spaces on $X$ to (counital) factorization quasicoherent sheaves on $X$; see Section 5 for more on this construction.

*Example* 2.5.6. Taking $X = \tilde{X}_{dR}$ in the previous example, we obtain a functor from factorization $D$ spaces on $\tilde{X}$ to usual factorization algebras on $\tilde{X}$.

*Example* 2.5.7. Let $\mathcal{C} = Cat^\times$, $F = PreStk_{/(\cdot)}$, and $G = ShvCat(\cdot)$. Define $\eta : F \to G_{op}$ by

$$\eta(S) : PreStk_{/S} \to ShvCat(S)^{op} \qquad (\mathcal{Y} \xrightarrow{p} S) \mapsto p_* QCoh_{\mathcal{Y}} .$$

Then the above provides a functor from (unital) factorization spaces over $X$ to (counital) factorization categories over $X$; see Section 5 for more on this construction.

## 3. Factorization categories and functors

In this section, we review the definition of factorization categories in terms of the general structure of Section 2, following [Ras15a] throughout, and outline a few related constructions from *loc. cit.*.

3.1. **Factorization categories.** Let $ShvCat : Sch_{aff}^{op} \to Cat^\times$ be the functor of sheaves of categories on a space, as in example 2.2.4. Recall from Example 2.4.8 we have the following definition:

*Definition* 3.1.1. A (unital) factorization category is a (unital) factorization sheaf of categories $C \in ShvCat_{un}^{fact}(X)$.

Similarly, a non-unital factorization category is a non-unital factorization sheaf of categories $C \in ShvCat^{fact}(X)$.

Let $Cat^{fact}(X) = ShvCat^{fact}(X)$ and $Cat_{un}^{fact}(X) = ShvCat_{un}^{fact}(X)$ denote the (2-)categories of non-unital and unital factorization categories.

*Remark* 3.1.2. There is a natural forgetful functor $Cat_{un}^{fact}(X) \to Cat^{fact}(X)$.

*Remark* 3.1.3. Concretely, following remark 2.4.3, a factorization category $C \in Cat_{un}^{fact}(X)$ on $X$ is

- a sheaf of categories $C \in ShvCat(Ran_{X,un})$ on the unital Ran space of $X$,
- an isomorphism $\alpha : j_{disj}^* C^{\boxtimes 2} \xrightarrow{\cong} \sqcup^* C$ of sheaves of categories over $(Ran_{X,un}^{\times 2})_{disj}$, and
- an isomorphism $\beta : Vect \xrightarrow{\cong} \iota_\emptyset^* C$ of objects in $DGCat = ShvCat(pt)$,

together with compatible higher arity analogues of the multiplication map $\alpha$ satisfying the relations of a unital commutative monoid, where $\iota_\emptyset : pt \hookrightarrow Ran_{X,un}$ is the map corresponding to the inclusion of the empty set.

*Remark* 3.1.4. The data of a factorization category $\mathrm{C} \in \mathrm{Cat}^{\mathrm{fact}}_{\mathrm{un}}(X)$ on $X$ can be further unpacked in terms of its restrictions to the strata of the Ran space, as in the case of a usual factorization algebra:

- The object $\mathrm{C} \in \mathrm{ShvCat}(\mathrm{Ran}_{X,\mathrm{un}})$ is given by an assignment

$$(3.1.1) \qquad I \mapsto \mathrm{C}_I \in \mathrm{ShvCat}(X^I) \qquad [\pi : I \to J] \mapsto [\Phi_\pi : \Delta(\pi)^* \, \mathrm{C}_I \to \mathrm{C}_J]$$

  for each $I \in \mathrm{fSet}$ and $\pi : I \to J$, such that the maps corresponding to $\pi$ surjective are isomorphisms;
- the factorization data $\alpha$ defines equivalences

$$(3.1.2) \qquad j(\pi)^* \, \mathrm{C}_I \xrightarrow{\cong} j(\pi)^* \boxtimes_{j \in J} \mathrm{C}_{I_j}$$

  of sheaves of categories on $U(\pi) \subset X^I$ for each $I$ and $\pi : I \twoheadrightarrow J$, in analogy with Definition I-6.0.2; and
- the unit data $\beta$, together with the gluing data for $\mathrm{C} \in \mathrm{ShvCat}(\mathrm{Ran}_{X,\mathrm{un}})$ in Equation 3.1.1 above, defines the analogues of the factorization unit maps of Remark I-4.3.2. In particular, for $\pi : I \hookrightarrow J$ injective, or further in particular for $\pi : \emptyset \to J$, this gives maps

$$(3.1.3) \qquad \mathrm{unit}^\pi_{\mathrm{C}} : \mathrm{C}_I \boxtimes \mathrm{QCoh}_{X^{I_\pi}} \to \mathrm{C}_J \qquad \text{and} \qquad \mathrm{unit}^J_{\mathrm{C}} : \mathrm{QCoh}_{X^J} \to \mathrm{C}_J$$

  of sheaves of categories on $X^J$, in analogy with Remark I-6.0.3; recall the notation $I_\pi$ and related conventions around partitions are given in Subsection I-2.2.

*Remark* 3.1.5. Recall that

$$\mathrm{Hom}_{\mathrm{ShvCat}(X)}(\mathrm{QCoh}_X, \mathrm{C}) = \Gamma(X, \mathrm{C})$$

so that the data of the unit maps in the preceeding Remark can be interpretted as giving objects $\mathrm{unit}_I \in \Gamma(X^I, \mathrm{C}^I)$. The compatibility with the gluing data for $\mathrm{C} \in \mathrm{ShvCat}(\mathrm{Ran}_{X,\mathrm{un}})$ implies these define an object $\mathrm{unit}_{\mathrm{C}} \in \Gamma(\mathrm{Ran}_{X,\mathrm{un}}, \mathrm{C})$, called the factorization unit; we discuss this further in Example 3.2.7 and Definition 3.2.8.

*Example* 3.1.6. There is a canonical (unital) factorization category $\mathrm{C} = \mathrm{QCoh}_{\mathrm{Ran}_{X,\mathrm{un}}} \in \mathrm{Cat}^{\mathrm{fact}}_{\mathrm{un}}(X)$ called the unit factorization category on $X$, defined by taking $\mathrm{C}_I = \mathrm{QCoh}_{X^I}$ for each $I \in \mathrm{fSet}$ together with the canonical factorization and unit structure maps; this is the categorical analogue of Example I-6.0.4, but in the quasicoherent setting.

*Proposition* 3.1.7. The symmetric monoidal structure $\otimes^*$ on $\mathrm{ShvCat}(\mathrm{Ran}_{X,\mathrm{un}})$ induces a canonical symmetric monoidal structure on $\mathrm{Cat}^{\mathrm{fact}}_{\mathrm{un}}(X)$, such that $\mathrm{QCoh}_{\mathrm{Ran}_{X,\mathrm{un}}}$ is the tensor unit.

Let $\mathrm{Cat}^{\mathrm{fact}}_{\mathrm{un}}(X)^*$ denote the symmetric monoidal (2-)category of unital factorization categories in the $\otimes^*$ monoidal structure.

## 3.2. Factorization functors and factorization algebras internal to a factorization category.

In this section, we define the notion of a factorization functor $\varphi : \mathrm{C} \to \mathrm{D}$ between (unital) factorization categories $\mathrm{C}, \mathrm{D} \in \mathrm{Cat}^{\mathrm{fact}}_{\mathrm{un}}(X)$ and define a factorization algebra internal to a factorization category $\mathrm{C}$ as a factorization functor from the unit factorization category $\mathrm{QCoh}_{\mathrm{Ran}_{X,\mathrm{un}}} \to \mathrm{C}$.

*Definition* 3.2.1. For (unital) factorization categories $\mathrm{C}, \mathrm{D} \in \mathrm{Cat}^{\mathrm{fact}}_{\mathrm{un}}(X)$, a (unital) factorization functor $\varphi : \mathrm{C} \to \mathrm{D}$ is a morphism in the category $\mathrm{Cat}^{\mathrm{fact}}_{\mathrm{un}}(X)$. Similarly, for non-unital factorization categories, a factorization functor is a morphism in $\mathrm{Cat}^{\mathrm{fact}}(X)$.

*Remark* 3.2.2. Concretely, in terms of the description in Remark 3.1.3 a factorization functor $\varphi : C \to D$ is given by

- a map of sheaves of categories $\varphi : C \to D$ over $\mathrm{Ran}_{X,\mathrm{un}}$,
- commutativity data for the diagram witnessing the compatibility of $\varphi$ with $\alpha_C$ and $\alpha_D$, and
- commutativity data for the diagram witnessing the compatibility of $\varphi$ with $\beta_C$ and $\beta_D$.

*Remark* 3.2.3. Further, in terms of the description in Remark 3.1.4, this yields:

- an assignment

$$(3.2.1) \qquad I \mapsto [\varphi_I : C_I \to D_I] \qquad [\pi : I \to J] \mapsto \begin{array}{ccc} \Delta(\pi)^* \, C_I & \longrightarrow & C_J \\ {\scriptstyle f_J} \downarrow \ {\scriptstyle \eta_\pi} \nearrow & & \downarrow {\scriptstyle \Delta(\pi)^*(f_I)} \\ \Delta(\pi)^* D_I & \longrightarrow & D_J \end{array}$$

  such that natural transformation $\eta_\pi$ is required to be an equivalence for $\pi : I \twoheadrightarrow J$ surjective;
- commutativity data for the diagram

$$(3.2.2) \qquad \begin{array}{ccc} j(\pi)^* \, C_i & \longrightarrow & j(\pi)^* \boxtimes_{j \in J} C_{I_j} \\ {\scriptstyle j(\pi)^* \varphi_I} \downarrow \ {\scriptstyle \sim} \nearrow & & \downarrow {\scriptstyle j(\pi)^* \boxtimes_{j \in J} \varphi_{I_j}} \\ j(\pi)^* D_I & \longrightarrow & j(\pi)^* \boxtimes_{j \in J} D_{I_j} \end{array}$$

  in $\mathrm{ShvCat}(X^I)$ for each $I$ and $\pi : I \twoheadrightarrow J$; and
- commutativity data for the diagram

$$(3.2.3) \qquad \begin{array}{ccc} \mathrm{Vect} & \longrightarrow & C_\emptyset \\ \| \ \ {\scriptstyle \sim} \nearrow & & \downarrow \\ \mathrm{Vect} & \longrightarrow & D_\emptyset \end{array} \ .$$

*Remark* 3.2.4. Note that the unital structure on a factorization functor allows for the flexibility of lax commutativity data in the diagram in Equation 3.2.1 above. Alternatively, a co-unital factorization functor would be defined analogously, but would instead allow for oplax commutativity data in the analogous diagram.

We now give the definition of a factorization algebra internal to a factorization category, which plays a central role in our desired applications of the results of the present section and those that follow.

*Definition* 3.2.5. A (unital) factorization algebra $A$ internal to a (unital) factorization category $C \in \mathrm{Cat}^{\mathrm{fact}}_{\mathrm{un}}(X)$ is a (unital) factorization functor $\mathrm{QCoh}_{\mathrm{Ran}_{X,\mathrm{un}}} \to C$.

Similarly, a non-unital factorization algebra internal to a non-unital factorization category $C \in \mathrm{Cat}^{\mathrm{fact}}(X)$ is given by a non-unital factorization functor from $\mathrm{QCoh}_{\mathrm{Ran}_X} \to C$.

Let $\mathrm{Alg}^{\mathrm{fact}}_{\mathrm{un}}(C)$ denote the category of factorization algebras internal to $C$ and similarly $\mathrm{Alg}^{\mathrm{fact}}(C)$ the category of non-unital factorization algebras internal to $C$.

*Remark* 3.2.6. Concretely, following Remark 3.2.3, a factorization algebra $A \in \mathrm{Alg}^{\mathrm{fact}}_{\mathrm{un}}(C)$ internal to a factorization category $C \in \mathrm{Cat}^{\mathrm{fact}}_{\mathrm{un}}(X)$ is given by:

- an object $A \in \Gamma(\mathrm{Ran}_{X,\mathrm{un}}, \mathrm{C})$, given by an assignment

$$(3.2.4) \qquad I \mapsto A_I \in \Gamma(X^I, \mathrm{C}_I) \qquad\qquad [\pi : I \to J] \mapsto [\eta_\pi : \Phi_\pi(\Delta(\pi)^\bullet A_I) \to A_J] \ ,$$

such that for $\pi$ surjective the corresponding maps in $\Gamma(X^J, \mathrm{C}_J)$ are isomorphisms, where $\Delta(\pi)^\bullet$ is the functor on sections induced by the unit $\Delta(\pi)^\bullet : \mathrm{C}_I \to \Delta(\pi)_* \Delta(\pi)^* \mathrm{C}$ of the $(\Delta(\pi)^*, \Delta(\pi)_*)$ adjunction, and $\Phi_\pi$ is the functor on sections induced by structure map of Equation 3.1.1;

- an equivalence

$$j(\pi)^\bullet A_I \xrightarrow{\cong} j(\pi)^\bullet (\boxtimes_j A_{I_j}) \qquad \text{in the category} \qquad \Gamma(U(\pi), j(\pi)^* \mathrm{C}_I) \cong \Gamma(U(\pi), j(\pi)^* \boxtimes_{j \in J} \mathrm{C}_{I_j})$$

for each $I$ and $\pi : I \twoheadrightarrow J$, where the equivalence of categories is that induced from the equivalence of sheaves of categories stipulated in Equation 3.1.2; and

- an equivalence $A_\varnothing \cong \mathbb{K}$ in $\mathrm{C}_\varnothing \cong \mathrm{Vect}$, which together with the structure maps of Equation 3.2.4 above, determine maps

$$(3.2.5) \qquad \mathrm{unit}_\mathrm{C}^\pi(A_I \boxtimes \mathcal{O}_{X^{I\pi}}) \to A_J \qquad \text{and} \qquad \mathrm{unit}_\mathrm{C}^J(\mathcal{O}_{X^J}) \to A_J$$

for $\pi : I \hookrightarrow J$ injective, and in particular for $\pi : \varnothing \hookrightarrow J$, respectively, where the functors $\mathrm{unit}_\mathrm{C}$ are those induced by the maps of sheaves of categories in Equation 3.1.3.

*Example* 3.2.7. The object $\mathrm{unit}_\mathrm{C} \in \Gamma(\mathrm{Ran}_{X,\mathrm{un}}, \mathrm{C})$ of Example 3.1.5 is canonically a factorization algebra object internal to $\mathrm{C}$.

*Definition* 3.2.8. Let $\mathrm{C}$ be a unital factorization category. The factorization unit object $\mathrm{unit}_\mathrm{C} \in \mathrm{Alg}_{\mathrm{un}}^{\mathrm{fact}}(\mathrm{C})$ is the factorization algebra internal to $\mathrm{C}$ defined in the preceeding example.

*Example* 3.2.9. Let $\mathrm{QCoh}_{\mathrm{Ran}_{X,\mathrm{un}}} \in \mathrm{Cat}_{\mathrm{un}}^{\mathrm{fact}}(X)$ be the factorization category of Example 3.1.6. The category $\mathrm{Alg}_{\mathrm{un}}^{\mathrm{fact}}(\mathrm{QCoh}_{\mathrm{Ran}_{X,\mathrm{un}}})$ of factorization algebras internal to $\mathrm{QCoh}_{\mathrm{Ran}_{X,\mathrm{un}}}$ is equivalent to the category $\mathrm{QCoh}_{\mathrm{un}}^{\mathrm{fact}}(X)$ of unital factorization quasicoherent sheaves, defined in Example 2.4.5. This follows tautologically from the identifications $\Gamma(X^I, \mathrm{QCoh}_{X^I}) \cong \mathrm{QCoh}(X^I)$.

*Example* 3.2.10. Taking $X = \tilde{X}_{\mathrm{dR}}$ in the preceding example, as in Remark 2.4.6, gives the factorization category

$$D_{\mathrm{Ran}_{\tilde{X},\mathrm{un}}} := \mathrm{QCoh}_{\mathrm{Ran}_{X,\mathrm{un}}} \in \mathrm{Cat}_{\mathrm{un}}^{\mathrm{fact}}(X) \qquad \text{such that} \qquad \mathrm{Alg}_{\mathrm{un}}^{\mathrm{fact}}(D_{\mathrm{Ran}_{\tilde{X},\mathrm{un}}}) \cong \mathrm{Alg}_{\mathrm{un}}^{\mathrm{fact}}(\tilde{X})$$

is equivalent to the usual category of unital factorization algebras.

*Remark* 3.2.11. Throughout this section, all the definitions and results were stated in the quasicoherent format, and only recovered the more familiar factorization algebras by taking $X = \tilde{X}_{\mathrm{dR}}$ as in the preceding example. This pattern will also be used in the remaining sections in Part 1.

## 4. Factorization spaces

4.1. **Factorization spaces.** Let $\mathrm{PreStk}_{/(\cdot)} : \mathrm{Sch}_{\mathrm{aff}}^{\mathrm{op}} \to \mathrm{Cat}^\times$ be the functor of prestacks over a scheme, as in Example 2.2.3. Recall from Example 2.4.7 we have the following definition:

*Definition* 4.1.1. A non-unital factorization space is a non-unital factorization prestack $\mathcal{Y} \in \mathrm{PreStk}^{\mathrm{fact}}(X)$.

Similarly, a (unital) factorization space is a (unital) factorization prestack $\mathcal{Y} \in \mathrm{PreStk}_{\mathrm{un}}^{\mathrm{fact}}(X)$.

*Remark* 4.1.2. There is a natural forgetful functor $\mathrm{PreStk}_{\mathrm{un}}^{\mathrm{fact}}(X) \to \mathrm{PreStk}^{\mathrm{fact}}(X)$.

*Remark* 4.1.3. We could take $X = \tilde{X}_{\mathrm{dR}}$ throughout to render all of the following definitions in the flat, as opposed to quasicoherent, format; see also Example 3.2.10 and Remark 3.2.11.

*Remark* 4.1.4. Concretely, following remark 2.4.3, a factorization space $\mathcal{Y} \in \mathrm{PreStk}_{\mathrm{un}}^{\mathrm{fact}}(X)$ on $X$ is

- a relative prestack $\mathcal{Y} \in \mathrm{PreStk}_{\mathrm{Ran}_{X,\mathrm{un}}}$ over the unital Ran space $\mathrm{Ran}_{X,\mathrm{un}} \in \mathrm{PreStk}^{\mathrm{lax}}$ of $X$,
- an equivalence

$$\alpha : (\mathrm{Ran}_{X,\mathrm{un}}^{\times 2})_{\mathrm{disj}} \times_{\mathrm{Ran}_{X,\mathrm{un}}^{\times 2}} \mathcal{Y}^{\times 2} \xrightarrow{\cong} (\mathrm{Ran}_{X,\mathrm{un}}^{\times 2})_{\mathrm{disj}} \times_{\mathrm{Ran}_{X,\mathrm{un}}} \mathcal{Y}$$

  of relative prestacks over $(\mathrm{Ran}_{X,\mathrm{un}}^{\times 2})_{\mathrm{disj}}$, and
- an isomorphism $\beta : \mathrm{pt} \xrightarrow{\cong} \{\varnothing\} \times_{\mathrm{Ran}_{X,\mathrm{un}}} \mathcal{Y}$ of prestacks,

together with compatible higher arity analogues of the multiplication map $\alpha$, satisfying the relations of a unital commutative monoid.

*Remark* 4.1.5. The data of a factorization space $\mathcal{Y} \in \mathrm{PreStk}_{\mathrm{un}}^{\mathrm{fact}}(X)$ on $X$ can be further unpacked in terms of its restrictions to the strata of the Ran space, as in the case of a usual factorization algebra:

- The object $\mathcal{Y} \in \mathrm{PreStk}_{/\mathrm{Ran}_{X,\mathrm{un}}}$ is given by an assignment

$$(4.1.1) \qquad I \mapsto \mathcal{Y}_I \in \mathrm{PreStk}_{/X^I} \qquad [\pi : I \to J] \mapsto [\phi_\pi : X^J \times_{X^I} \mathcal{Y}_I \to \mathcal{Y}_J]$$

  for each $I \in \mathrm{fSet}_\varnothing$ and $\pi : I \to J$, such that the maps corresponding to $\pi$ surjective are isomorphisms;
- The factorization data $\alpha$ defines equivalences

$$(4.1.2) \qquad U(\pi) \times_{X^I} \mathcal{Y}_I \xrightarrow{\cong} U(\pi) \times_{X^I} (\times_{j \in J} \mathcal{Y}_{I_j})$$

  of prestacks over $X^I$ for each $I$ and $\pi : I \twoheadrightarrow J$, in analogy with Definition I-6.0.2; and
- The unit data $\beta$, together with the gluing data for $\mathcal{Y} \in \mathrm{PreStk}_{/\mathrm{Ran}_{X,\mathrm{un}}}$ above, defines the analogues of the factorization unit maps of Remark I-6.0.3: for $\pi : I \hookrightarrow J$ injective, or further in particular for $\pi : \varnothing \to J$, this gives maps

$$(4.1.3) \qquad \mathrm{unit}_\mathcal{Y}^\pi : X^{I_\pi} \times \mathcal{Y}_I \to \mathcal{Y}_J \qquad \text{and} \qquad \mathrm{unit}_\mathcal{Y}^J : X^J \to \mathcal{Y}_J$$

  of prestacks over $X^J$, in analogy with Remark I-6.0.3; recall the notation $I_\pi$ and related conventions around partitions are given in Subsection I-2.2.

*Example* 4.1.6. The space $\mathrm{Ran}_{X,\mathrm{un}} \in \mathrm{PreStk}_{\mathrm{lax}}$ defines a trivial unital factorization prestack $\mathrm{Ran}_{X,\mathrm{un}} \in \mathrm{PreStk}_{\mathrm{un}}^{\mathrm{fact}}(X$ which is analogous to the trivial factorization category $\mathrm{QCoh}_{\mathrm{Ran}_{X,\mathrm{un}}} \in \mathrm{Cat}_{\mathrm{un}}^{\mathrm{fact}}(X)$.

*Remark* 4.1.7. In analogy with Remark 3.1.5, the unit maps in Equation 4.1.3 glue together to define a section $\mathrm{unit}_\mathcal{Y} \in \Gamma(\mathrm{Ran}_{X,\mathrm{un}}, \mathcal{Y})$ of the relative prestack $\mathcal{Y} \twoheadrightarrow \mathrm{Ran}_{X,\mathrm{un}}$.

Towards defining more interesting factorization spaces, we introduce the following notation:

*Remark* 4.1.8. Let $S \in \mathrm{Sch}_{\mathrm{aff}}$ and $x = (x_i) : S \to X^I$ denote an $S$ point of $X^I$. Let $\Gamma_{x_i} \subset S \times X$ denote the graph of $x_i : S \to X$ and define the spaces

$$\Gamma_x = \bigcup_{i \in I} \Gamma_{x_i} \qquad \mathbb{D}_x = (S \times X)_{\Gamma_x/S}^\wedge \qquad \mathbb{D}_x^\circ = \mathbb{D}_x \backslash \Gamma_x$$

over $S \times X$, which we call the graph of $x$, the formal disk around $x$, and the punctured formal disc around $x$, respectively. In the case where $S = \mathrm{pt}$, we will also write $\{x\}$ for $\Gamma_x$ to denote the union of the points $x_i \in X(\mathbb{K})$; note that we take the union in the set-theoretic sense, so that in particular $\{x\} \subset X$ does not keep track of the information of multiplicities of repeated points $x_{i_0} = x_{i_1}$ for $i_0 \neq i_1 \in I$.

*Example* 4.1.9. The Beilinson-Drinfeld Grassmannian is the unital factorization space $\mathrm{Gr}_{G,\mathrm{Ran}_{X,\mathrm{un}}} \in \mathrm{PreStk}^{\mathrm{fact}}_{\mathrm{un}}(X)$ defined by

$$\mathrm{Gr}_{G,I}(S) = \{x : S \to X^I, \ P \in \mathrm{Bun}_G(S \times X), \ \varphi : P|_{S \times X \backslash \Gamma_x} \xrightarrow{\cong} P^{\mathrm{triv}}\}$$

together with the natural gluing data over $\mathrm{Ran}_X$ and factorization structure maps. In particular, the fibre of $\mathrm{Gr}_{G,\{1\}}$ over $x \in X$ is given by the usual affine Grassmannian $\mathrm{Gr}_{G,x} = G(\mathcal{K}_x)/G(\mathcal{O}_x)$.

The unital structure on $\mathrm{Gr}_{G,\mathrm{Ran}_{X,\mathrm{un}}}$ is defined, for example for $\pi : I \hookrightarrow J$ injective, by the map

$$X^{I_\pi} \times \mathrm{Gr}_{G,I} \to \mathrm{Gr}_{G,J} \qquad (x = (x_1, x_2), \ P, \ \varphi : P|_{X \times S \backslash \Gamma_{x_2}} \to P^{\mathrm{triv}}) \mapsto (x, P, \varphi : P|_{X \times S \backslash \Gamma_x} \to P^{\mathrm{triv}})$$

restricting the trivialization $\varphi$ along $X \times S \backslash \Gamma_x \hookrightarrow X \times S \backslash \Gamma_{x_2}$, where $x = (x_1, x_2) \in X^J = X^{I_\pi} \times X^I$, and $\Gamma_x$ is as defined in Remark 4.1.8. In particular, the map on $\mathrm{Gr}_{G,\{1\}}$ over $x \in X$ is given by the usual inclusion $\{[t^0]\} \hookrightarrow \mathrm{Gr}_G$.

*Remark* 4.1.10. We may also consider non-unital or co-unital factorization spaces. The former is described concretely as in remarks 4.1.4 and 4.1.5 but only requiring the data assigned to non-empty subsets and $\pi : I \twoheadrightarrow J$ surjective. The latter is also described concretely as in these two remarks but requiring maps

$$I \mapsto \mathcal{Y}_I \in \mathrm{PreStk}_{/X^I} \qquad [\pi : I \to J] \mapsto \big[\psi_\pi : \mathcal{Y}_J \to X^J \times_{X^I} \mathcal{Y}_I\big]$$

in place of those in Equation 4.1.1; note that this data only differs in the case $\pi : I \to J$ is not injective, and in particular induces maps

$$\mathrm{counit}^\pi_{\mathcal{Y}} : \mathcal{Y}_J \to X^{I_\pi} \times \mathcal{Y}_I \qquad \text{and} \qquad \mathrm{counit}^J_{\mathcal{Y}} : \mathcal{Y}_J \to X^J$$

in place of those in Equation 4.1.3.

*Example* 4.1.11. Let $Y \in \mathrm{PreStk}_{/X}$ be a prestack over $X$. The jet space $\mathcal{J}(Y)_{\mathrm{Ran}_{X,\mathrm{un}}} \in \mathrm{PreStk}^{\mathrm{fact}}_{\mathrm{co\text{-}un}}(X)$ is the counital factorization space defined by

$$\mathcal{J}(Y)_I(S) = \{x : S \to X^I, a : \mathbb{D}_x \to Y \text{ a map over } X\}$$

together with the natural gluing data over $\mathrm{Ran}_X$ and factorization structure maps, where $\mathbb{D}_x$ is as defined in Remark 4.1.8. In particular, the fibre of $\mathcal{J}(Y)_{\{1\}}$ over $x \in X$ is given by the usual arc space $\mathcal{J}(Y)_x = Y_{\mathcal{O}_x} = \mathrm{Maps}(\mathbb{D}_x, Y)$ of $Y$.

The counital structure on $\mathcal{J}(Y)_{\mathrm{Ran}_{X,\mathrm{un}}}$ is defined, for example for $\pi : I \hookrightarrow J$ injective, by the map

$$\mathcal{J}(Y)_J \to X^{I_\pi} \times \mathcal{J}(Y)_I \qquad (x = (x_1, x_2), a : \mathbb{D}_x \to Y) \mapsto (x_1, (x_2, a|_{\mathbb{D}_{x_2}} : \mathbb{D}_{x_2} \to Y))$$

given by restricting the map $a$ along $\mathbb{D}_{x_2} \hookrightarrow \mathbb{D}_x$, where $x = (x_1, x_2) \in X^J = X^{I_\pi} \times X^I$ and $\mathbb{D}_x$ is as defined in Remark 4.1.8. In particular, the map on $\mathcal{J}(Y)_x$ over $x \in X$ is given by the projection to the point $Y_{\mathcal{O}_x} \to \{x\}$.

4.2. **Correspondences and unital structures.** More generally, consider as a starting point the functor $\mathrm{PreStk}_{\mathrm{corr}/(\cdot)} : \mathrm{Sch}^{\mathrm{op}}_{\mathrm{aff}} \to \mathrm{Cat}^\times$ of prestacks under (arbitrary) correspondences over a scheme, which also admits a natural lax symmetric monoidal structure. We use the same notation $\mathrm{PreStk}^{\mathrm{fact}}_{\mathrm{un}}(X)$ and term 'unital factorization space' to refer to the resulting more general variant discussed below, unless otherwise stated. The previous unital and counital variants occur as special cases; see Definition 4.2.2 and Remark 4.2.3.

*Remark* 4.2.1. In the category $\mathrm{PreStk}^{\mathrm{fact}}_{\mathrm{un}}(X) := \mathrm{PreStk}^{\mathrm{fact}}_{\mathrm{corr,un}}(X)$ of unital factorization spaces over $X$, the equivalence $\alpha$, or equivalently the maps of Equation 4.1.2, would still required to be invertible and thus genuine morphisms. However, the maps of Equation 4.1.1 for $\pi : I \to J$ not surjective would be given more generally by correspondences

$$(4.2.1) \qquad [\pi : I \to J] \mapsto \qquad \begin{array}{c} \mathrm{unit}^\pi \\ {}^{\psi_\pi}\swarrow \qquad \searrow{}^{\phi_\pi} \\ X^J \times_{X^I} \mathcal{Y}^I \qquad\qquad \mathcal{Y}_J \end{array} \quad .$$

In particular, the maps of Equation 4.1.3, would be given by

$$(4.2.2) \qquad \begin{array}{c} \mathrm{unit}^\pi_{\mathcal{Y}} \\ {}^{\psi_\pi}\swarrow \qquad \searrow{}^{\phi_\pi} \\ X^{I_\pi} \times \mathcal{Y}_I \qquad\qquad \mathcal{Y}_J \end{array} \qquad \text{and} \qquad \begin{array}{c} \mathrm{unit}^J_{\mathcal{Y}} \\ {}^{\psi_J}\swarrow \qquad \searrow{}^{\phi_J} \\ X^J \qquad\qquad \mathcal{Y}_J \end{array}$$

for $\pi : I \hookrightarrow J$ injective, and in particular for $\pi : \emptyset \to J$, respectively.

Note that in this more general format, the notions of unital and counital factorization prestack coincide: the category $\mathrm{PreStk}_{\mathrm{corr}/S}$ is equivalent to its opposite category, naturally in $S \in \mathrm{Sch}_{\mathrm{aff}}$, so that there is a canonical equivalence of functors $\mathrm{PreStk}_{\mathrm{corr}/(\cdot),\mathrm{op}} \cong \mathrm{PreStk}_{\mathrm{corr}/(\cdot)}$, and thus a natural equivalence between the space of unital and counital factorization prestacks.

*Definition* 4.2.2. A unital factorization space $\mathcal{Y}$ such that the maps $\psi_\pi$ (or $\phi_\pi$) of Equation 4.2.1 are isomorphisms for each $\pi : I \to J$ is called strictly unital (or counital).

*Remark* 4.2.3. The preceeding definitions of strictly unital and counital factorization spaces agree with the previously defined notions of Definition 4.1.1 and Remark 4.1.10.

*Example* 4.2.4. Let $Y \in \mathrm{PreStk}_{/X}$ be a prestack over $X$. The meromorphic jet space $\mathcal{J}^{\mathrm{mer}}(Y) \in \mathrm{PreStk}^{\mathrm{fact}}_{\mathrm{un}}(X)$ of $Y$ is the unital factorization space (in the more general sense of the preceeding remark) defined by

$$\mathcal{J}^{\mathrm{mer}}(Y)_I(S) = \{x : S \to X^I, a : \mathbb{D}^\circ_x \to Y \text{ a map over } X\}$$

together with the natural gluing data over $\mathrm{Ran}_X$ and factorization structure maps. In particular, the fibre of $\mathcal{J}^{\mathrm{mer}}(Y)_{\{1\}}$ over $x \in X$ is given by the usual algebraic loop space $\mathcal{J}^{\mathrm{mer}}(Y)_x = Y_{\mathcal{K}_x} = \mathrm{Maps}(\mathbb{D}^\circ_x, Y)$ of $Y$.

The unital structure on $\mathcal{J}^{\mathrm{mer}}_{\mathrm{Ran}_{X,\mathrm{un}}}$ is defined, for example for $\pi : I \hookrightarrow J$ injective, by the correspondence

$$\mathrm{unit}^\pi_{\mathcal{J}^{\mathrm{mer}}(Y)}(S) = \{x = (x_1, x_2) : S \to X^J \; , \; a : \mathbb{D}_x \backslash \Gamma_{x_2} \to Y \text{ a map over } X\}$$

together with the maps

$$\begin{array}{c} \mathrm{unit}^\pi_{\mathcal{J}^{\mathrm{mer}}(Y)} \\ \swarrow \qquad\qquad \searrow \\ X^{I_\pi} \times \mathcal{J}^{\mathrm{mer}}(Y)_I \qquad\qquad \mathcal{J}^{\mathrm{mer}}(Y)_J \end{array} \quad \text{defined using} \quad \begin{array}{c} \mathbb{D}_x \backslash \Gamma_{x_2} \\ \nearrow \qquad\qquad \nwarrow \\ \mathbb{D}_{x_2} \backslash \Gamma_{x_2} \qquad\qquad \mathbb{D}_x \backslash \Gamma_x \end{array}$$

to restrict $a : \mathbb{D}_x \backslash \Gamma_{x_2} \to Y$, giving the desired maps, where $x = (x_1, x_2) \in X^J = X^{I_\pi} \times X^I$, and $\mathbb{D}_x$, $\Gamma_x$ and $\mathbb{D}^\circ_x = \mathbb{D}_x \backslash \Gamma_x$ are as defined in Remark 4.1.8. In particular, the correspondence corresponding to $\varnothing \hookrightarrow J$, and its restriction to each $x \in X$, are given by

$$(4.2.3) \qquad \begin{array}{ccc} & \mathfrak{J}(Y)_J & \\ & \swarrow \qquad \searrow & \\ X^J & & \mathfrak{J}^{\mathrm{mer}}(Y)_J \end{array} \qquad \text{and} \qquad \begin{array}{ccc} & Y_{\mathcal{O}_x} & \\ & \swarrow \qquad \searrow & \\ \mathrm{pt} & & Y_{\mathcal{K}_x} \end{array} \;\; .$$

*Remark* 4.2.5. Canonical descent to de Rham stack of $X$.

4.3. **Maps of factorization spaces and factorizable sections.** A map of factorization spaces $f : \mathcal{X} \to \mathcal{Y}$ is by definition a map in the category $\mathrm{PreStk}^{\mathrm{fact}}_{\mathrm{un}}(X)$ of unital factorization prestacks (under correspondences, unless stated otherwise). We breifly unpack the definition, in close analogy with remarks 3.2.2 and 3.2.3 describing factorization functors.

*Remark* 4.3.1. Concretely, in terms of the descriptions in remarks 4.1.4 and 4.1.5, a map of unital factorization spaces (under correspondences) from $\mathcal{X}$ to $\mathcal{Y}$ is given by

- a map $f : \mathcal{X} \to \mathcal{Y}$ of relative prestacks over $\mathrm{Ran}_{X,\mathrm{un}}$, specified by

$$(4.3.1) \qquad I \mapsto [f_I : \mathcal{X}_I \to \mathcal{Y}_I] \qquad\qquad [\pi : I \to J] \mapsto \begin{array}{ccc} & & \mathrm{unit}^{\pi}_{\mathcal{X}} \\ & \swarrow \quad \searrow & \\ X^J \times_{X^I} \mathcal{X}_I \quad & \mathrm{unit}^{\pi}_{\mathcal{Y}} & \quad \mathcal{X}_J \\ \downarrow \quad\quad \swarrow & & \downarrow \\ X^J \times_{X^I} \mathcal{Y}_I & & \mathcal{Y}_J \end{array}$$

together with commutativity data for the latter diagram analogous to that required in Equation 3.2.1 in the setting of factorization categories, for each $I \in \mathrm{fSet}$ and $\pi : I \to J$; and

- compatibility data for the above correspondence with the factorization and unital structures on $\mathcal{X}$ and $\mathcal{Y}$, exactly analogous to that required in the second and third items of Remark 3.2.2 and given concretely in equations 3.2.2 and 3.2.3 in the setting of factorization categories, but defined in terms of the factorization space data described in the second and third items of Remark 4.1.4, and given concretely in equations 4.1.2 and 4.1.3, respectively.

*Example* 4.3.2. Every (unital) factorization space has a canonical map to the trivial (unital) factorization space $\mathrm{Ran}_{X,\mathrm{un}}$.

*Remark* 4.3.3. Concretely, in terms of the descriptions in remarks 4.1.4 and 4.1.5, a map of unital factorization spaces (under correspondences) from $\mathcal{X}$ to $\mathcal{Y}$ is given by

- a correspondence of relative prestacks over $\mathrm{Ran}_{X,\mathrm{un}}$

(4.3.2)

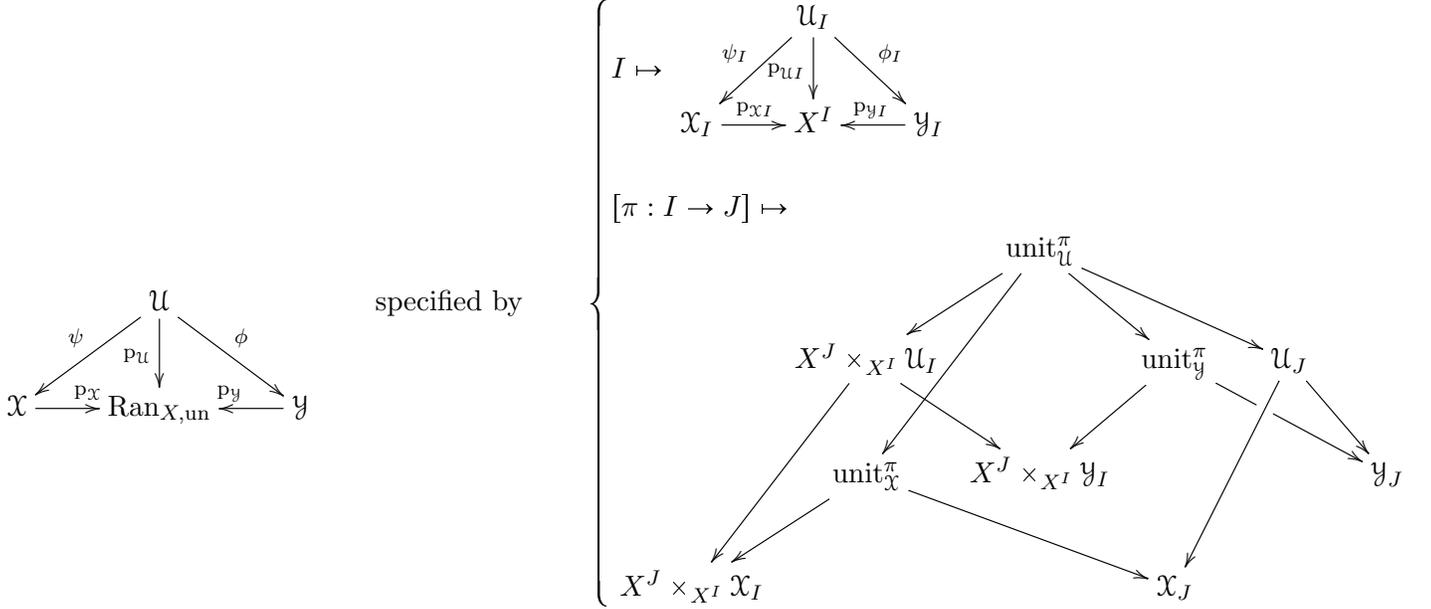

together with commutativity data, for each $I \in \mathrm{fSet}$ and $\pi : I \to J$; and
- compatibility data for the above correspondence with the factorization and unital structures on $\mathcal{X}$ and $\mathcal{Y}$, exactly analogous to that required in the second and third items of Remark 3.2.2 and given concretely in equations 3.2.2 and 3.2.3 in the setting of factorization categories, but defined in terms of the factorization space data described in the second and third items of Remark 4.1.4 (though in the correspondence setting), and given concretely in equations 4.1.2 and 4.2.2, respectively.

*Remark* 4.3.4. The preceeding data of a correspondence of (unital) factorization spaces defines in particular a unital factorization structure on the correspondence space $\mathcal{U} \in \mathrm{PreStk}_{\mathrm{un}}^{\mathrm{fact}}(X)$.

4.3.1. *Factorizable sections.* We now introduce the analogue of a factorization algebra object internal to a factorization category in the setting of factorization spaces. This is simply another factorization space with a strict map to the given one, as we now explain:

*Definition* 4.3.5. The space of (unital) factorizable sections $\Gamma_{\mathrm{un}}^{\mathrm{fact}}(X, \mathcal{Y})$ of a (unital) factorization space $\mathcal{Y}$ is the space of maps $\mathrm{Ran}_{X,\mathrm{un}} \to \mathcal{Y}$ of (unital) factorization spaces on $X$.

*Remark* 4.3.6. Concretely, a factorizable section $\mathcal{S} \in \Gamma^{\mathrm{fact}}(X, \mathcal{Y})$ is given by

- a correspondence of relative prestacks over $\mathrm{Ran}_{X,\mathrm{un}}$

(4.3.3)

$$
\begin{cases}
I \mapsto & \begin{array}{c}
\mathcal{S}_I \\
{}^{\psi_I}\swarrow \ {}^{\mathrm{p}_{\mathcal{S}I}}\downarrow \ {}^{\phi_I}\searrow \\
X^I =\!\!= X^I \xleftarrow{\mathrm{p}_{\mathcal{Y}I}} \mathcal{Y}_I
\end{array} \\
\\
[\pi : I \to J] \mapsto & \begin{array}{c}
\mathrm{unit}_{\mathcal{S}}^{\pi} \\
\swarrow \qquad\qquad \searrow \\
X^J \times_{X^I} \mathcal{S}_I \qquad\qquad \mathrm{unit}_{\mathcal{Y}}^{\pi} \qquad \mathcal{S}_J \\
\searrow \qquad \swarrow \qquad\qquad \searrow \\
X^J \times_{X^I} \mathcal{Y}_I \qquad\qquad \mathcal{Y}_J
\end{array}
\end{cases}
$$

$$
\begin{array}{c}
\mathcal{S} \\
{}^{\psi}\swarrow \ {}^{\mathrm{p}_{\mathcal{S}}}\downarrow \ {}^{\phi}\searrow \\
\mathrm{Ran}_{X,\mathrm{un}} =\!\!= \mathrm{Ran}_{X,\mathrm{un}} \xleftarrow{\mathrm{p}_{\mathcal{Y}}} \mathcal{Y}
\end{array}
\qquad \text{specified by}
$$

for each $I \in \mathrm{fSet}$ and $\pi : I \to J$, together with commutativity data, as in the first item of Remark 4.3.3 above; and

- compatibility data for the above correspondences with the factorization and unital structure on $\mathcal{Y}$, as in the second item of *loc. cit.*.

*Remark* 4.3.7. A factorizable section $\mathcal{S} \in \Gamma^{\mathrm{fact}}(X, \mathcal{Y})$ is canonically itself a unital factorization space $\mathcal{S} \in \mathrm{PreStk}_{\mathrm{un}}^{\mathrm{fact}}(X)$ on $X$ by Remark 4.3.4, the map $\mathcal{S} \to \mathrm{Ran}_{X,\mathrm{un}}$ identified with the canonical such map of Example 4.3.2, and $\mathcal{S} \to \mathcal{Y}$ defines a strict map of (unital) factorization spaces. Indeed, we have:

*Proposition* 4.3.8. The category $\Gamma^{\mathrm{fact}}(X, \mathcal{Y})$ of unital is identified with that of (unital) factorization spaces over $\mathcal{Y}$.

*Example* 4.3.9. The factorization unit correspondence $\mathrm{unit}_{\mathcal{Y}}$ is canonically a factorizable section. Thus, it is itself a (strictly counital) factorization space with a strict map of unital factorization spaces $\phi_{\mathcal{Y}} : \mathrm{unit}_{\mathcal{Y}} \to \mathcal{Y}$.

## 5. Linearization of factorization spaces

5.1. **Factorization categories associated to factorization spaces.** Throughout, we use the general format for linearization discussed in subsection 2.5, though we elaborate on the details of these examples of the construction.

*Example* 5.1.1. Following Example 2.5.7, we let $\mathcal{C} = \mathrm{Cat}^{\times}$, $F = \mathrm{PreStk}_{/(\cdot)} : \mathrm{Sch}_{\mathrm{aff}}^{\mathrm{op}} \to \mathcal{C}$, and $G = \mathrm{ShvCat}(\cdot) : \mathrm{Sch}_{\mathrm{aff}}^{\mathrm{op}} \to \mathcal{C}$, and explained that there was a monoidal compatible natural transformation $\eta : F \to G_{\mathrm{op}}$ defined by

$$
\eta(S) : \mathrm{PreStk}_{/S} \to \mathrm{ShvCat}(S)^{\mathrm{op}} \quad
\begin{cases}
(\mathcal{Y} \xrightarrow{\mathrm{p}} S) & \mapsto \mathrm{p}_* \mathrm{QCoh}_{\mathcal{Y}} \\
\\
\begin{array}{c}
\mathcal{X} \xrightarrow{\phi} \mathcal{Y} \\
{}_{\mathrm{p}_{\mathcal{X}}}\searrow \ \downarrow {}^{\mathrm{p}_{\mathcal{Y}}} \\
S
\end{array} & \mapsto [\mathrm{p}_{\mathcal{Y}*}(\phi^{\bullet}) : \mathrm{p}_{\mathcal{Y}*} \mathrm{QCoh}_{\mathcal{Y}} \to \mathrm{p}_{\mathcal{Y}*} \phi_* \phi^* \mathrm{QCoh}_{\mathcal{Y}} \cong \mathrm{p}_{\mathcal{X}*} \mathrm{QCoh}_{\mathcal{X}}]
\end{cases} ,
$$

where the map $\phi^\bullet : \mathrm{QCoh}_{\mathcal{Y}} \to \phi_* \phi^* \mathrm{QCoh}_{\mathcal{Y}}$ is given by the unit of the $(\phi^*, \phi_*)$ adjunction, and the naturality of $\eta$ is equivalent to the requirement of base change for ShvCat

$$f^* \mathrm{p}_* \mathrm{QCoh}_{\mathcal{Y}} \xrightarrow{\cong} \tilde{\mathrm{p}}_* \tilde{f}^* \mathrm{QCoh}_{\mathcal{Y}} \cong \tilde{\mathrm{p}}_* \mathrm{QCoh}_{T \times_S \mathcal{Y}} \qquad \text{along} \qquad \begin{array}{ccc} T \times_S \mathcal{Y} & \xrightarrow{\tilde{f}} & \mathcal{Y} \\ {\scriptstyle \tilde{\mathrm{p}}} \downarrow & & \downarrow {\scriptstyle \mathrm{p}} \\ T & \xrightarrow{f} & S \end{array} .$$

From the preceeding example, we obtain a functor from the category of unital factorization spaces over $X$ to that of counital factorization categories over $X$, by corollary 2.5.4.

*Definition* 5.1.2. Let $\mathcal{Y} \in \mathrm{PreStk}_{\mathrm{un}}^{\mathrm{fact}}(X)$ be a unital factorization space. The counital factorization category $\mathrm{QCoh}_{\mathcal{Y}} \in \mathrm{Cat}_{\mathrm{co\text{-}un}}^{\mathrm{fact}}(X)$ of quasicoherent sheaves on $\mathcal{Y}$ is defined as the factorization category resulting from applying Corollary 2.5.4 to the monoidal compatible natural transformation from Example 2.5.7.

*Remark* 5.1.3. Concretely, in terms of the decsriptions in remarks 3.2.3 and 4.1.5, the factorization category $\mathrm{QCoh}_{\mathcal{Y}} \in \mathrm{Cat}_{\mathrm{co\text{-}un}}^{\mathrm{fact}}(X)$ is given by the following data:

- The object $\mathrm{p}_* \mathrm{QCoh}_{\mathcal{Y}} \in \mathrm{ShvCat}(\mathrm{Ran}_{X,\mathrm{un}})$, determined by the assignment

$$I \mapsto \mathrm{p}_{I,*} \mathrm{QCoh}_{\mathcal{Y}_I} \in \mathrm{ShvCat}(X^I)$$

(5.1.1) $\qquad [\pi : I \to J] \mapsto \big[ \mathrm{p}_{J,*}(\phi_\pi^\bullet) : \mathrm{p}_{J,*} \mathrm{QCoh}_{\mathcal{Y}_J} \to \mathrm{p}_{J,*} \phi_{\pi *} \phi_\pi^* \mathrm{QCoh}_{\mathcal{Y}_J} \cong \Delta(\pi)^* \mathrm{p}_{I,*} \mathrm{QCoh}_{\mathcal{Y}_I} \big]$

for each $I \in \mathrm{fSet}$ and $\pi : I \to J$; the latter isomorphism is defined by base change along $X^J \times_{X^I} \mathcal{Y}_I$, and the map is given by the unit $\phi^\bullet$ of the $(\phi_\pi^*, \phi_{\pi *})$ adjunction, where $g_\pi$ is as in Equation 4.1.1;

- the factorization data $\alpha$ is given by the equivalences

$$j(\pi)^* \mathrm{p}_{I,*} \mathrm{QCoh}_{\mathcal{Y}_I} \xrightarrow{\cong} j(\pi)^* \boxtimes_{j \in J} \mathrm{p}_{I_j,*} \mathrm{QCoh}_{\mathcal{Y}_{I_j}}$$

induced by the equivalences in Equation 4.1.2; and

- the factorization counit maps are given by the maps

$$\mathrm{unit}_{\mathcal{Y}}^{\pi,\bullet} : \mathrm{p}_{J,*} \mathrm{QCoh}_{\mathcal{Y}_J} \to \mathrm{p}_{I,*} \mathrm{QCoh}_{\mathcal{Y}_I} \boxtimes \mathrm{QCoh}_{X^{I_\pi}} \qquad \text{and} \qquad \mathrm{unit}_{\mathcal{Y}}^{J,\bullet} : \mathrm{QCoh}_{\mathcal{Y}_J} \to \mathrm{p}_{J,*} \mathrm{QCoh}_{X^J}$$

defined, as above, as (the image under $\mathrm{p}_{J*}$ of) the unit $\mathrm{unit}_{\mathcal{Y}}^\bullet$ of the $(\mathrm{unit}_{\mathcal{Y}}^*, \mathrm{unit}_{\mathcal{Y}*})$ adjunctions for the unit maps in Equation 4.1.3.

*Remark* 5.1.4. There is a non-unital variant of the preceeding definitions and remarks, given as usual by not requiring the data associated to $I = \emptyset$ nor $\pi : I \to J$ non-injective.

*Remark* 5.1.5. There is also a variant which produces a unital factorization category from a counital factorization space; the variant explained in the following subsection in terms of correspondences simeltaneously generalizes this and the above.

*Example* 5.1.6.

## 5.2. Linearization of correspondences and (co)unital structures.
Following Remark 4.2.1, there is an analogue of the preceeding definition of $\mathrm{QCoh}_{\mathcal{Y}}$ in the case where the unit maps for the factorization space $\mathcal{Y} \in \mathrm{PreStk}_{\mathrm{un}}^{\mathrm{fact}}(X)$ are given by correspondences. However, it will be necessary to restrict the class of factorization spaces under consideration in order to carry out this construction, as well as many other constructions throughout the remainder of the paper. A simple format for this is the following:

*Remark* 5.2.1. Let $\mathrm{F}_P \subset \mathrm{PreStk}_{\mathrm{corr}/(\cdot)} : \mathrm{Sch}_{\mathrm{aff}}^{\mathrm{op}} \to \mathrm{Cat}^\times$ be a lax symmetric monoidal subfunctor of $\mathrm{PreStk}_{\mathrm{corr}/(\cdot)}$; this data is interpreted as assigning to $S \in \mathrm{Sch}_{\mathrm{aff}}$ the subcategory $\mathrm{F}_P(S) \to \mathrm{PreStk}_{/S} \in \mathrm{Cat}^\times$ of the category of prestacks over $S$ 'satisfying property $P$', though note that we also allow to restrict the class of morphisms in the category.

We introduce the following formal language for the resulting notion:

*Definition* 5.2.2. A factorization space $\mathcal{Y} \in \mathrm{PreStk}_{P,\mathrm{un}}^{\mathrm{fact}}(X)$ on $X$ with property $P$ is a factorization $\mathrm{F}_P$-object for $\mathrm{F}_P$ a lax symmetric monoidal subfunctor of $\mathrm{PreStk}_{/(\cdot)}^{\mathrm{corr}}$, as in the preceeding remark.

*Remark* 5.2.3. The inclusion of the subfunctor $\eta : \mathrm{F}_P \to \mathrm{PreStk}_{/(\cdot)}$ is a canonically a lax compatible natural transformation, and thus induces a functor $\mathrm{PreStk}_{P,\mathrm{un}}^{\mathrm{fact}}(X) \to \mathrm{PreStk}_{\mathrm{un}}^{\mathrm{fact}}(X)$, interpreted as the inclusion of the subcategory of unital factorization spaces on $X$ satisfying property $P$.

*Example* 5.2.4. One general class of subfunctors are given by $\mathrm{PreStk}_{(\mathrm{corr};\alpha,\beta)/(\cdot)} : \mathrm{Sch}_{\mathrm{aff}}^{\mathrm{op}} \to \mathrm{Cat}^\times$ which assigns to an affine scheme $S \in \mathrm{Sch}^{\mathrm{aff}}$ the category with objects the same as those in $\mathrm{PreStk}_{/S}$ but with morphisms

$$\mathcal{X} \rightsquigarrow \mathcal{Y} \qquad \text{given by correspondences} \qquad \begin{array}{ccc} & \mathcal{U} & \\ {\scriptstyle\psi}\swarrow & & \searrow{\scriptstyle\phi} \\ \mathcal{X} & & \mathcal{Y} \end{array}$$

such that $\psi : \mathcal{U} \to \mathcal{X}$ is a map 'of class $\alpha$' and $\phi : \mathcal{U} \to \mathcal{Y}$ a map 'of class $\beta$'. We call a choice of a pair of such classes $(\alpha, \beta)$ compatible if this construction defines a (lax symmetric monoidal) subfunctor of $\mathrm{PreStk}_{\mathrm{corr}/(\cdot)}$.

The resulting class of factorization spaces $\mathrm{PreStk}_{(\alpha,\beta),\mathrm{un}}^{\mathrm{fact}}(X)$ are those $\mathcal{Y} \in \mathrm{PreStk}_{\mathrm{un}}^{\mathrm{fact}}(X)$ which admit a description as in Remark 4.2.1 such that the correspondences in equations 4.2.1 and 4.2.2 required to be composed of maps $f$ of class $\alpha$ and $g$ of class $\beta$ as above.

Note that for $\alpha \neq \beta$, the notions of unital and counital factorization spaces of type $(\alpha, \beta)$ are no longer equivalent.

*Example* 5.2.5. Let $\alpha = $ all be the class of all morphisms, $\beta = $ sch-qcqs to be the class of schematic, quasi-compact and quasi-separated morphisms, $\mathrm{F} = \mathrm{PreStk}_{(\mathrm{corr};\alpha,\beta)/(\cdot)} : \mathrm{Sch}_{\mathrm{aff}}^{\mathrm{op}} \to \mathrm{Cat}^\times$ and $\mathrm{G} = \mathrm{ShvCat}(\cdot) : \mathrm{Sch}_{\mathrm{aff}}^{\mathrm{op}} \to \mathrm{Cat}^\times$. Then there exists a monoidal compatible natural transformation $\eta : \mathrm{F} \to \mathrm{G}$ defined by

(5.2.1)

$$\eta(S) : \mathrm{PreStk}_{(\mathrm{corr};\alpha,\beta)/S} \to \mathrm{ShvCat}(S)^{\mathrm{op}} \begin{cases} (\mathcal{Y} \xrightarrow{\mathrm{p}} S) & \mapsto \mathrm{p}_* \, \mathrm{QCoh}_{\mathcal{Y}} \\ \begin{array}{ccc} & \mathcal{U} & \\ {\scriptstyle\psi}\swarrow \,{\scriptstyle\mathrm{p}_{\mathcal{U}}}\downarrow & & \searrow{\scriptstyle\phi} \\ \mathcal{X} \xrightarrow[\mathrm{p}_{\mathcal{X}}]{} & S & \xleftarrow[\mathrm{p}_{\mathcal{Y}}]{} \mathcal{Y} \end{array} & \mapsto \big[ \mathrm{p}_{\mathcal{Y}*}(\phi_\bullet) \circ \mathrm{p}_{\mathcal{X}*}(\psi^\bullet) : \mathrm{p}_{\mathcal{X}*} \, \mathrm{QCoh}_{\mathcal{X}} \to \mathrm{p}_{\mathcal{Y}*} \, \mathrm{QCoh}_{\mathcal{Y}} \big] \end{cases} ,$$

where the map $\psi^\bullet : \mathrm{QCoh}_{\mathcal{X}} \to \psi_* \psi^* \, \mathrm{QCoh}_{\mathcal{X}}$ is the unit of the $(\psi_*, \psi^*)$ adjunction and

$$\phi_\bullet : \phi_* \phi^* \, \mathrm{QCoh}_{\mathcal{Y}} \to \mathrm{QCoh}_{\mathcal{Y}}$$

is the right adjoint to $\phi^\bullet$, which is always defined and satisfies base change for $\phi$ a sch-qcqs map, so that the desired map is given by the composition

$$\mathrm{p}_{\mathcal{X}*} \, \mathrm{QCoh}_{\mathcal{X}} \xrightarrow{\mathrm{p}_{\mathcal{X}*}(\psi^\bullet)} \mathrm{p}_{\mathcal{X}*} \psi_* \psi^* \, \mathrm{QCoh}_{\mathcal{X}} \cong \mathrm{p}_{\mathcal{U}*} \, \mathrm{QCoh}_{\mathcal{U}} \cong \mathrm{p}_{\mathcal{Y}*} \phi_* \phi^* \, \mathrm{QCoh}_{\mathcal{Y}} \xrightarrow{\mathrm{p}_{\mathcal{Y}*}(\phi_\bullet)} \mathrm{p}_{\mathcal{Y}*} \, \mathrm{QCoh}_{\mathcal{Y}} \ .$$

*Example* 5.2.6. Similarly, reversing the roles of $\alpha$ and $\beta$ in the preceeding example, there is a monoidal compatible natural transformation $\eta : \mathrm{PreStk}_{(\mathrm{corr};\beta,\alpha)} \to \mathrm{ShvCat}(\cdot)_{\mathrm{op}}$, generalizing that of Example 5.1.1. It is defined as in Equation 5.2.1 of the preceeding example, except that it assigns to a correspondence of *loc. cit.* the map

$$\mathrm{p}_{\mathfrak{X}*}(\psi_\bullet) \circ \mathrm{p}_{\mathfrak{Y}*}(\phi^\bullet) : \mathrm{p}_{\mathfrak{Y}*} \mathrm{QCoh}_\mathfrak{Y} \to \mathrm{p}_{\mathfrak{X},*} \mathrm{QCoh}_\mathfrak{X} \ .$$

From the preceeding examples, we obtain a functors from categories of unital factorization spaces to those of (co)unital factorization categories, by Corollary 2.5.4:

*Definition* 5.2.7. Let $\mathcal{Y} \in \mathrm{PreStk}^{\mathrm{fact}}_{(\mathrm{all},\mathrm{sch\text{-}qcqs}),\mathrm{un}}(X)$ a unital factorization spaces over $X$ of type $(\mathrm{all},\mathrm{sch\text{-}qcqs})$. Then $\mathrm{QCoh}_\mathcal{Y} \in \mathrm{Cat}^{\mathrm{fact}}_{\mathrm{un}}(X)$ is the unital factorization category resulting from applying Corollary 2.5.4 to the monoidal compatible natural transformation of Example 5.2.5.

Similarly, for $\mathcal{Y} \in \mathrm{PreStk}^{\mathrm{fact}}_{(\mathrm{sch\text{-}qcqs},\mathrm{all}),\mathrm{un}}(X)$, $\mathrm{QCoh}_\mathcal{Y} \in \mathrm{Cat}^{\mathrm{fact}}_{\mathrm{co\text{-}un}}(X)$ is the counital factorization category resulting the monoidal compatible natural transformation of Example 5.2.6.

*Remark* 5.2.8. Concretely, in analogy with Remark 5.1.3, in terms of the decsriptions in remarks 3.2.3 and 4.1.5, the factorization category $\mathrm{QCoh}_\mathcal{Y} \in \mathrm{Cat}^{\mathrm{fact}}(X)$ is given by the following data:

- The object $\mathrm{p}_* \mathrm{QCoh}_\mathcal{Y} \in \mathrm{ShvCat}(\mathrm{Ran}_{X,\mathrm{un}})$ is determined by the assignment

$$I \mapsto \mathrm{p}_{I,*} \mathrm{QCoh}_{\mathcal{Y}_I} \in \mathrm{ShvCat}(X^I)$$

(5.2.2)
$$[\pi : I \to J] \mapsto \Big[ \mathrm{p}_{J,*}(\phi_{\pi\bullet}) \circ \tilde{\mathrm{p}}_{I*}(\psi^\bullet_\pi) : \Delta(\pi)^* \mathrm{p}_{I,*} \mathrm{QCoh}_{\mathcal{Y}_I} \cong \tilde{\mathrm{p}}_{I,*} \mathrm{QCoh}_{X^J \times_{X^I} \mathcal{Y}_I} \to \mathrm{p}_{J,*} \mathrm{QCoh}_{\mathcal{Y}_J} \Big]$$

for each $I \in \mathrm{fSet}$ and $\pi : I \to J$, where the map in the latter datum is given by base change along $X^J \times_{X^I} \mathcal{Y}_I$, together with the map defined by the construction in Equation 5.2.1 for the correspondence structure data of Equation 4.2.1, as follows:

For $\psi^\bullet_\pi$ the unit of the $(\psi_{\pi*}, \psi^*_\pi)$ adjunction and $\phi_{\pi,\bullet}$ the right adjoint to $\phi^\bullet_\pi$, the above composition gives a map

$$\tilde{\mathrm{p}}_{I,*} \mathrm{QCoh}_{X^J \times_{X^I} \mathcal{Y}_I} \xrightarrow{\mathrm{p}_{I*}(\psi^\bullet_\pi)} \tilde{\mathrm{p}}_{I,*} \psi_{\pi*} \psi^*_\pi \mathrm{QCoh}_{X^J \times_{X^I} \mathcal{Y}_I} \cong \mathrm{p}_{J*}\phi_{\pi*}\phi^*_\pi \mathrm{QCoh}_{\mathcal{Y}_J} \xrightarrow{\mathrm{p}_{J*}(\phi_{\pi\bullet})} \mathrm{p}_{J,*} \mathrm{QCoh}_{\mathcal{Y}_J} \ ,$$

as desired.

- the factorization data $\alpha$ is given by the equivalences

(5.2.3)
$$j(\pi)^* \mathrm{p}_{I,*} \mathrm{QCoh}_{\mathcal{Y}_I} \xrightarrow{\cong} j(\pi)^* \boxtimes_{j \in J} \mathrm{p}_{I_j,*} \mathrm{QCoh}_{\mathcal{Y}_{I_j}}$$

induced by the equivalences in Equation 4.1.2; and

- the factorization unit maps are given by the maps

(5.2.4)
$$\phi_{\pi,\bullet}\psi^\bullet_\pi : \mathrm{p}_{I,*} \mathrm{QCoh}_{\mathcal{Y}_I} \boxtimes \mathrm{QCoh}_{X^{I_\pi}} \to \mathrm{p}_{J,*} \mathrm{QCoh}_{\mathcal{Y}_J} \qquad \text{and} \qquad \phi_{J,\bullet}\psi^\bullet_J : \mathrm{QCoh}_{X^J} \to \mathrm{p}_{J,*} \mathrm{QCoh}_{\mathcal{Y}_J}$$

defined, as above, as (the images under $\mathrm{p}_*$ of) the unit $\psi^\bullet$ of the $(\psi^*, \psi_*)$ adjunction, and the right adjoint $\phi_\bullet$ of the analogous adjunction unit $\phi^\bullet$, for the maps in the unit correspondences of Equation 4.2.2.

*Example* 5.2.9. quasicoh of loop space

*Remark* 5.2.10. As in Remark 3.2.11, all of the preceeding definitions and remarks are given in the quasicoherent format, but there are several natural flat variants by considering the cases $X = \tilde{X}_{\mathrm{dR}}$ and/or $\mathcal{Y} = \tilde{Y}_{\mathrm{dR}}$.

### 5.3. Linearization of maps and factorization functors.

The statement of Proposition 2.5.3 is that a symmetric monoidal compatible natural transformation $\eta : \mathrm{F} \to \mathrm{G}$ of lax monoidal functors $\mathrm{F}, \mathrm{G} : \mathrm{Sch}_{\mathrm{aff}}^{\mathrm{op}} \to \mathcal{C}$ gives rise to a *functor* $\mathrm{F}_{\mathrm{un}}^{\mathrm{fact}}(X) \to \mathrm{G}_{\mathrm{un}}^{\mathrm{fact}}(X)$. Thus, in addition to the procedure which takes a factorization space $\mathcal{Y}$ to the corresponding factorization category $\mathrm{QCoh}_{\mathcal{Y}}$ desribed in the preceeding subsection, a map (or correspondence) of factorization spaces $f : \mathcal{X} \to \mathcal{Y}$ naturally defines a factorization functor $f_\bullet : \mathrm{QCoh}_{\mathcal{X}} \to \mathrm{QCoh}_{\mathcal{Y}}$. We describe the latter construction explicitly below, noting that these results do not quite follow formally from Proposition 2.5.3, due to the subtleties of Warning 2.1.1:

For concreteness, we restrict to the format of Example 5.2.5, taking $\mathrm{F} = \mathrm{PreStk}_{(\mathrm{corr};\mathrm{all},\mathrm{sch\text{-}qcqs})/(\cdot)} :$ $\mathrm{Sch}_{\mathrm{aff}}^{\mathrm{op}} \to \mathrm{Cat}^\times$, $G = \mathrm{ShvCat} : \mathrm{Sch}_{\mathrm{aff}}^{\mathrm{op}} \to \mathrm{Cat}^\times$, and $\eta : \mathrm{F} \to \mathrm{G}$ as defined in Equation 5.2.1 of *loc. cit.*. We begin by unpacking the resulting construction of the pushforward and pullback functors for a strict map of such factorization spaces, and then explain the functor induced by a general correspondence of factorization spaces as the composition of the former:

*Example* 5.3.1. Let $f : \mathcal{X} \to \mathcal{Y}$ be a strict map in the category of (type (all, sch-qcqs), unital) factorization spaces, and recall from Remark 4.3.1 that this is in particular described by an assignment

$$(5.3.1) \quad I \mapsto [f_I : \mathcal{X}_I \to \mathcal{Y}_I] \qquad \text{and} \qquad [\pi : I \to J] \mapsto$$

together with commutativity of the latter diagram, for each $\pi : I \to J$.

The induced (unital) factorization functor $f_\bullet : \mathrm{QCoh}_{\mathcal{X}} \to \mathrm{QCoh}_{\mathcal{Y}}$ is given, in terms of the descriptions in Remark 5.2.8, by the assignment
$$(5.3.2)$$

$$I \mapsto [f_{I,\bullet} : \mathrm{p}_{\mathcal{X}_{I*}} \mathrm{QCoh}_{\mathcal{X}_I} \to \mathrm{p}_{\mathcal{Y}_{I*}} \mathrm{QCoh}_{\mathcal{Y}_I}] \qquad [\pi : I \to J] \mapsto$$

for each $I \in \mathrm{fSet}_\varnothing$ and $\pi : I \to J$, where the natural transformation $\eta_\pi$ is defined by

$$\phi_{\mathcal{Y}\bullet}\psi_{\mathcal{Y}}^\bullet \tilde{f}_{I,\bullet} \cong \phi_{\mathcal{Y}\bullet}(\tilde{f}_I)_\bullet \tilde{\psi}_{\mathcal{Y}}^\bullet \to \phi_{\mathcal{Y}\bullet}(\tilde{f}_I)_\bullet \chi_\bullet \chi^\bullet \tilde{\psi}_{\mathcal{Y}}^\bullet \cong f_{J\bullet}\phi_{\mathcal{Y}\bullet}\psi_{\mathcal{X}}^\bullet \qquad \text{where}$$

summarizes the maps used: $\widetilde{\mathrm{unit}}^\pi_{\mathcal{Y}} = \mathrm{unit}^\pi_{\mathcal{Y}} \underset{X^J \times_{X^I} \mathcal{Y}_I}{\times} X^J \times_{X^I} \mathcal{X}_I$ is the pullback, $\chi : \mathrm{unit}^\pi_{\mathcal{X}} \to \widetilde{\mathrm{unit}}^\pi_{\mathcal{Y}}$ is the canonical such map, and we apply base change, the unit of the $(\chi^\bullet, \chi_\bullet)$ adjunction, and commutativity of the diagram in Equation 5.3.1.

*Example* 5.3.2. Let $f : \mathcal{X} \to \mathcal{Y}$ be a strict map of factorization spaces, as in the preceeding example. The induced (co)unital factorization functor $f^\bullet : \mathrm{QCoh}_{\mathcal{Y}} \to \mathrm{QCoh}_{\mathcal{X}}$ is given, in terms of the descriptions in Remark 5.2.8, by the assignment

(5.3.3)

$$I \mapsto [f_{I,\bullet} : \mathrm{p}_{\mathcal{Y}_I *} \mathrm{QCoh}_{\mathcal{Y}_I} \to \mathrm{p}_{\mathcal{X}_I *} \mathrm{QCoh}_{\mathcal{X}_I}] \qquad [\pi : I \to J] \mapsto$$

$$\begin{array}{ccc}
\tilde{\mathrm{p}}_{I*} \mathrm{QCoh}_{X^J \times_{X^I} \mathcal{Y}_I} & \xrightarrow{\phi_{\mathcal{Y}_\bullet} \psi_{\mathcal{Y}}^\bullet} & \mathrm{p}_{\mathcal{X}_{J*}} \mathrm{QCoh}_{\mathcal{Y}_J} \\
{\scriptstyle \tilde{f}_I^\bullet} \downarrow & {\scriptstyle \eta_\pi} \Downarrow & \downarrow {\scriptstyle f_J^\bullet} \\
\tilde{\mathrm{p}}_{I*} \mathrm{QCoh}_{X^J \times_{X^I} \mathcal{X}_I} & \xrightarrow{\phi_{\mathcal{X}_\bullet} \psi_{\mathcal{X}}^\bullet} & \mathrm{p}_{\mathcal{Y}_{J*}} \mathrm{QCoh}_{\mathcal{X}_J}
\end{array}$$

for each $I \in \mathrm{fSet}_{\varnothing}$ and $\pi : I \to J$, where the natural transformation $\eta_\pi$ is defined by

$$f_J^\bullet \phi_{\mathcal{Y}_\bullet} \psi_{\mathcal{Y}}^\bullet \cong \tilde{\phi}_{\mathcal{Y}_\bullet} \tilde{f}_J^\bullet \psi_{\mathcal{Y}}^\bullet \to \tilde{\phi}_{\mathcal{Y}_\bullet} \xi_\bullet \xi^\bullet \tilde{f}_J^\bullet \psi_{\mathcal{Y}}^\bullet \cong \phi_{\mathcal{X}_\bullet} \psi_{\mathcal{X}}^\bullet \tilde{f}_I^\bullet \qquad \text{where}$$

$$\begin{array}{ccc}
& \mathrm{unit}^\pi_{\mathcal{X}} & \\
{\scriptstyle \xi} \swarrow \quad {\scriptstyle \phi_{\mathcal{X}}} & & \\
{\scriptstyle f_\pi} \downarrow \quad \mathrm{unit}^\pi_{\mathcal{X}} \times_{\mathcal{X}_J} \mathcal{Y}_{J-} & \xrightarrow{\quad} & \mathcal{X}_J \\
& {\scriptstyle \tilde{f}_J} \downarrow \qquad {\scriptstyle \phi_{\mathcal{Y}}} & \downarrow {\scriptstyle f_J} \\
& \mathrm{unit}^\pi_{\mathcal{Y}} & \xrightarrow{\phi_{\mathcal{Y}}} \mathcal{Y}_J
\end{array}$$

summarizes the maps used: $\xi : \mathrm{unit}^\pi_{\mathcal{X}} \to \mathrm{unit}^\pi \times_{\mathcal{X}_J} \mathcal{Y}_J$ is the canonical map such map, and we use base change together with the unit of the $(\xi^\bullet, \xi_\bullet)$ adjunction, and commutativity of the diagram in Equation 5.3.1.

*Remark* 5.3.3. In the case that $\mathcal{X}$ and $\mathcal{Y}$ are strictly co-unital, the map $\xi$ is an isomorphism, and thus the map $\eta_\pi$ of the preceeding example is also an equivalence, so that it can also be interpreted as data defining a unital factorization functor.

*Example* 5.3.4. Now, consider a general morphism of unital factorization spaces over $X$ under correspondences, as in Equation 4.3.2.

5.4. **Factorization compatible sheaves.** Throughout, let $\mathcal{Y} \in \mathrm{PreStk}_{\mathrm{un}}^{\mathrm{fact}}(X)$ be a (unital) factorization space admitting a (unital) factorization category of quasicoherent sheaves $\mathrm{QCoh}_{\mathcal{Y}} \in \mathrm{Cat}_{\mathrm{un}}^{\mathrm{fact}}(X)$, as for example in Definition 5.2.7.

*Definition* 5.4.1. A (unital) factorization compatible quasicoherent sheaf $\mathcal{F} \in \mathrm{Alg}_{\mathrm{un}}^{\mathrm{fact}}(\mathrm{QCoh}_{\mathcal{Y}})$ on $\mathcal{Y}$ is a (unital) factorization algebra internal to the factorization category $\mathrm{QCoh}_{\mathcal{Y}} \in \mathrm{Cat}_{\mathrm{un}}^{\mathrm{fact}}(X)$.

*Remark* 5.4.2. Concretely, following the descriptions in remarks 3.2.6 and 5.2.8, a factorization compatible quasicoherent sheaf $\mathcal{F} \in \mathrm{Alg}_{\mathrm{un}}^{\mathrm{fact}}(\mathrm{QCoh}_{\mathcal{Y}})$ on $\mathcal{Y}$ is given by

- an object $\mathcal{F} \in \mathrm{QCoh}(\mathcal{Y})$, given by an assignment

(5.4.1) $$I \mapsto \mathcal{F}_I \in \mathrm{QCoh}(\mathcal{Y}_I) \qquad [\pi : I \to J] \mapsto \left[ \eta_\pi : \phi_{\pi\bullet} \psi_\pi^\bullet \tilde{\Delta}(\pi)^\bullet \mathcal{F}_I \to \mathcal{F}_J \right],$$

where the latter are maps in $\mathrm{QCoh}(\mathcal{Y}_J)$, such that for $\pi$ surjective the corresponding maps are isomorphisms; here $\tilde{\Delta}(\pi) : X^J \times_{X^I} \mathcal{Y}_I \to \mathcal{Y}_I$ is the canonical map covering $\Delta(\pi) : X^J \to X^I$, and the maps $\phi_{\pi\bullet}$ and $\psi_\pi^\bullet$ as in Equation 5.2.2;

- an equivalence

$$\tilde{j}(\pi)^\bullet \mathcal{F}_I \xrightarrow{\cong} \tilde{j}'(\pi)^\bullet (\boxtimes_j \mathcal{F}_{I_j}) \qquad \text{in the category} \qquad \mathrm{QCoh}(U(\pi) \times_{X^I} \mathcal{Y}_I) \cong \mathrm{QCoh}(U(\pi) \times_{X^I} (\times_{j \in J} \mathcal{Y}_{I_j}))$$

for each $I$ and $\pi : I \twoheadrightarrow J$, where the equivalence of categories is that induced from the equivalence of Equation 5.2.4, and $\tilde{j}(\pi) : U(\pi) \times_{X^I} \mathcal{Y}_I \hookrightarrow \mathcal{Y}_I$ is the map lifting $j(\pi) : U(\pi) \hookrightarrow X^I$, and similarly for $\tilde{j}'(\pi)$; and

- an equivalence $\mathcal{F}_\varnothing \cong \mathbb{K}$ in $\mathrm{QCoh}(\mathcal{Y}_\varnothing) \cong \mathrm{Vect}$, which together with the structure maps of Equation 5.4.1 above, determines maps

$$\phi_{\pi \bullet} \psi^\bullet_\pi (\mathcal{F}_I \boxtimes \mathcal{O}_{X^{I_\pi}}) \to \mathcal{F}_J \qquad \text{and} \qquad \phi_{J \bullet} \psi^\bullet_J (\mathcal{O}_{X^J}) \to \mathcal{F}_J$$

for $\pi : I \hookrightarrow J$ injective, and in particular for $\pi : \varnothing \hookrightarrow J$, respectively, where $\phi_\bullet$ and $\psi^\bullet$ are the functors on objects induced by the maps of sheaves of categories in Equation 5.2.4.

*Example* 5.4.3. The factorization unit object $\mathrm{unit}_{\mathrm{QCoh}_\mathcal{Y}} \in \mathrm{Alg}^{\mathrm{fact}}_{\mathrm{un}}(\mathrm{QCoh}_\mathcal{Y})$ defines a factorization compatible sheaf: By Example 4.3.9, there are maps of unital factorization spaces

$$\mathrm{p}_{\mathrm{unit}_\mathcal{Y}} : \mathrm{unit}_\mathcal{Y} \to \mathrm{Ran}_{X,\mathrm{un}} \qquad \text{and} \qquad \phi_\mathcal{Y} : \mathrm{unit}_\mathcal{Y} \to \mathcal{Y} \ .$$

If $\mathcal{Y}$ is as in Definition 5.2.7, these induce unital factorization functors

$$\mathrm{p}^\bullet_{\mathrm{unit}_\mathcal{Y}} : \mathrm{QCoh}_{\mathrm{Ran}_{X,\mathrm{un}}} \to \mathrm{QCoh}_{\mathrm{unit}_\mathcal{Y}} \qquad \text{and} \qquad \phi_{\mathcal{Y}\bullet} : \mathrm{QCoh}_{\mathrm{unit}_\mathcal{Y}} \to \mathrm{QCoh}_\mathcal{Y} \ ,$$

by Examples 5.3.2 and 5.3.1, respectfully. In terms of these functors, the unit factorization compatible sheaf

$$\mathrm{unit}_{\mathrm{QCoh}_\mathcal{Y}} \in \mathrm{Alg}^{\mathrm{fact}}_{\mathrm{un}}(\mathrm{QCoh}_\mathcal{Y}) \qquad \text{is given by} \qquad \mathrm{unit}_{\mathrm{QCoh}_\mathcal{Y}} = \phi_{\mathcal{Y}\bullet} \mathrm{p}^\bullet_{\mathrm{unit}_\mathcal{Y}} \mathcal{O}_{\mathrm{Ran}_{X,\mathrm{un}}} = \phi_{\mathcal{Y}\bullet} \mathcal{O}_{\mathrm{unit}_\mathcal{Y}} \ .$$

Concretely, the factorization unit is given by the assignment

$$I \mapsto \phi_{I \bullet} \psi^\bullet_I \mathcal{O}_{X^I} \cong \phi_{I \bullet} \mathcal{O}_{\mathrm{unit}^I_\mathcal{Y}}$$

together with the canonical compatibility data determined on the factorization unit.

*Example* 5.4.4. example: the unit for the BD grassmannian

*Example* 5.4.5. Let $\mathcal{Y} \in \mathrm{PreStk}^{\mathrm{fact}}(X)$ be a (unital) factorization space such that $\mathrm{p}_\mathcal{Y} : \mathcal{Y} \to \mathrm{Ran}_{X,\mathrm{un}}$ induces a unital factorization functor $\mathrm{p}^\bullet_\mathcal{Y} : \mathrm{QCoh}_{\mathrm{Ran}_{X,\mathrm{un}}} \to \mathrm{QCoh}_\mathcal{Y}$, as in Remark 5.3.3. Then the structure sheaf $\mathcal{O}_\mathcal{Y} = \mathrm{p}^\bullet_\mathcal{Y} \mathcal{O}_{\mathrm{Ran}_{X,\mathrm{un}}} \in \mathrm{Alg}^{\mathrm{fact}}(\mathrm{QCoh}_\mathcal{Y})$ is canonically a non-unital factorization compatible sheaf. Concretely, it is given by the assignment $I \mapsto \mathcal{O}_{\mathcal{Y}_I}$ together with the natural structure maps. In the case that $\mathcal{Y}$ is strictly counital, the resulting object is just the factorization unit.

*Remark* 5.4.6. For $f : \mathcal{X} \to \mathcal{Y}$ a map of factorization spaces inducing a unital factorization functor $f_\bullet : \mathrm{QCoh}_\mathcal{Y} \to \mathrm{QCoh}_\mathcal{Y}$, there is an induced functor

$$\mathrm{Alg}^{\mathrm{fact}}_{\mathrm{un}}(\mathcal{X}) \to \mathrm{Alg}^{\mathrm{fact}}_{\mathrm{un}}(\mathcal{Y}) \ ,$$

of pushforward of factorization compatible quasicoherent sheaves.

*Example* 5.4.7. More generally, let $\mathcal{X} \in \mathrm{PreStk}^{\mathrm{fact}}_{\mathrm{un}}(X)$ be a factorization space such that $\mathrm{p}_\mathcal{X} : \mathcal{X} \to \mathrm{Ran}_{X,\mathrm{un}}$ induces a factorization functor $\mathrm{p}^\bullet_\mathcal{X} : \mathrm{QCoh}_{\mathrm{Ran}_{X,\mathrm{un}}} \to \mathrm{QCoh}_\mathcal{X}$, as in Example 5.4.5, and let $f : \mathcal{X} \to \mathcal{Y}$ be a map of (unital) factorization spaces inducing a unital factorization functor $f_\bullet : \mathrm{QCoh}_\mathcal{Y} \to \mathrm{QCoh}_\mathcal{Y}$. Then there is a factorization compatible quasicoherent sheaf

$$f_\bullet \mathrm{p}^\bullet_\mathcal{X} \mathcal{O}_{\mathrm{Ran}_{X,\mathrm{un}}} \in \mathrm{Alg}^{\mathrm{fact}}_{\mathrm{un}}(\mathrm{QCoh}_\mathcal{Y}) \qquad \text{given concretely by} \qquad I \mapsto f_{I \bullet} \mathrm{p}^\bullet_{\mathcal{X}_I} \mathcal{O}_{X^I} \cong f_{I \bullet} \mathcal{O}_{\mathcal{X}_I}$$

together with the induced compatibility data.

In the special case when $\mathfrak{X} = \mathrm{unit}_\mathcal{Y}$, this induces the factorization unit object $\mathrm{unit}_\mathcal{Y} \in \mathrm{Alg}_\mathrm{un}^\mathrm{fact}(\mathrm{QCoh}_\mathcal{Y})$ by Example 5.4.3 above.

### 5.5. Factorization linearization functors.

Throughout, let $\mathcal{Y} \in \mathrm{PreStk}_\mathrm{un}^\mathrm{fact}(X)$ be a (unital) factorization space admitting a (unital) factorization category of quasicoherent sheaves $\mathrm{QCoh}_\mathcal{Y} \in \mathrm{Cat}_\mathrm{un}^\mathrm{fact}(X)$, as for example in Definition 5.2.7.

*Definition* 5.5.1. A factorizable linearization functor on $\mathcal{Y}$ is a factorization functor $\varphi : \mathrm{QCoh}_\mathcal{Y} \to \mathrm{QCoh}_{\mathrm{Ran}_{X,\mathrm{un}}}$.

*Remark* 5.5.2. Given a factorizable linearization functor $\varphi : \mathrm{QCoh}_\mathcal{Y} \to \mathrm{QCoh}_{\mathrm{Ran}_{X,\mathrm{un}}}$, we obtain a functor

$$\varphi \circ (\cdot) : \mathrm{Alg}_\mathrm{un}^\mathrm{fact}(\mathrm{QCoh}_\mathcal{Y}) \to \mathrm{QCoh}_\mathrm{un}^\mathrm{fact}(X)$$

from the category of factorization compatible quasicoherent sheaves on $\mathcal{Y}$ to that of factorization quasicoherent sheaves on $X$, by Example 3.2.9.

*Remark* 5.5.3. Similarly, following Example 3.2.10, taking $X = \tilde{X}_\mathrm{dR}$ in the above, a factorizable linearization functor $\varphi : \mathrm{p}_* \mathrm{QCoh}_\mathcal{Y} \to \mathrm{QCoh}_{\mathrm{Ran}_{X,\mathrm{un}}}$ determines a functor

$$\varphi \circ (\cdot) : \mathrm{Alg}_\mathrm{un}^\mathrm{fact}(\mathrm{QCoh}_\mathcal{Y}) \to \mathrm{Alg}_\mathrm{un}^\mathrm{fact}(\tilde{X})$$

from the category of factorization compatible quasicoherent sheaves on $\mathcal{Y}$ to that of factorization algebras on $\tilde{X}$.

Note that in this variant, the factorization space $\mathcal{Y}$ is equipped with a flat structure along $\mathrm{Ran}_{\tilde{X},\mathrm{un}}$ so that the quasicoherent direct image of an object in $\mathrm{QCoh}(\mathcal{Y})$ determines a $D$ module on $\mathrm{Ran}_{\tilde{X},\mathrm{un}}$, but we have not used the category $D(\mathcal{Y})$ of $D$ modules on $\mathcal{Y}$.

*Example* 5.5.4. Suppose $\mathcal{Y} \in \mathrm{PreStk}^\mathrm{fact}(X)$ is such that the factorization category $\mathrm{QCoh}_\mathcal{Y} \in \mathrm{Cat}_{\mathrm{co}\text{-}\mathrm{un}}^\mathrm{fact}(X)$ is counital. Then the counital structure defines a factorization linearization functor, in analogy with Example 5.4.3. Concretely, if $\mathcal{Y}$ is as in Definition 5.2.7, the factorization linearization functor is given by the assignment

$$I \mapsto [\psi_{I,\bullet}\phi_I^\bullet : \mathrm{p}_{I*}\,\mathrm{QCoh}_{\mathcal{Y}_I} \to \mathrm{QCoh}_{X^I}]$$

together with the induced compatibility data.

## 6. Factorization $\mathcal{P}$ algebras

Let $\mathcal{P} \in \mathrm{HOp}^\mathrm{co}(\mathrm{Vect}_\mathbb{K})$ be a cocommutative Hopf operad over $\mathbb{K}$. Recall from Appendix I-C.2 that the cocommutative Hopf structure on $\mathcal{P}$ implies that for any $\mathbb{K}$-linear symmetric monoidal category $\mathcal{C}$, the category of $\mathcal{P}$ algebras $\mathrm{Alg}_\mathcal{P}(\mathcal{C})$ internal to $\mathcal{C}$ is symmetric monoidal with respect to the tensor product $\otimes_\mathcal{C}$ of the underlying objects in $\mathcal{C}$.

Let $\mathrm{F}_\mathcal{P} := \mathrm{Alg}_\mathcal{P}(D^!(\cdot)) : \mathrm{Sch}_\mathrm{aff}^\mathrm{op} \to \mathrm{Cat}^\times$ denote the functor

$$X \mapsto \mathrm{Alg}_\mathcal{P}(D(X)) \qquad [f : X \to Y] \mapsto f^! : \mathrm{Alg}_\mathcal{P}(D(Y)) \to \mathrm{Alg}_\mathcal{P}(D(X))$$

of $\mathcal{P}$ algebra objects internal to $D(X)$, together with the pullback functor on such objects induced by the symmetric monoidal pullback functor $f^! : D(Y) \to D(X)$.

*Proposition* 6.0.1. There is a canonical lax symmetric monoidal structure on $\mathrm{F}_\mathcal{P} : \mathrm{Sch}_\mathrm{aff}^\mathrm{op} \to \mathrm{Cat}^\times$ given by

$$\boxtimes : \mathrm{Alg}_\mathcal{P}(D(X)) \times \mathrm{Alg}_\mathcal{P}(D(Y)) \to \mathrm{Alg}_\mathcal{P}(D(X \times Y)) \qquad (A, B) \mapsto A \boxtimes B$$

where $\boxtimes : D(X) \times D(Y) \to D(X \times Y)$ is the usual exterior product, which admits a lift to algebra objects via the cocommutative Hopf structure on $\mathcal{P}$.

*Definition* 6.0.2. A (unital) factorization $\mathcal{P}$ algebra on $X$ is a factorization $F_{\mathcal{P}}$-object $A \in F_{\mathcal{P},\mathrm{un}}^{\mathrm{fact}}(X)$ on $X$, for $F_{\mathcal{P}}$ as defined above together with the lax symmetric monoidal structure of the previous proposition.

Similarly, a non-unital factorization $\mathcal{P}$ algebra $A \in F_{\mathcal{P}}^{\mathrm{fact}}(X)$ is a non-unital factorization $F_{\mathcal{P}}$-object.

Let $\mathrm{Alg}_{\mathcal{P},\mathrm{un}}^{\mathrm{fact}}(X) = F_{\mathcal{P},\mathrm{un}}^{\mathrm{fact}}(X)$ denote the category of factorization $\mathcal{P}$ algebras, and $\mathrm{Alg}_{\mathcal{P}}^{\mathrm{fact}}(X) = F_{\mathcal{P}}^{\mathrm{fact}}(X)$ the category of non-unital factorization $\mathcal{P}$ algebras.

*Remark* 6.0.3. Concretely, following remark 2.4.3, a factorization $\mathcal{P}$ algebra $A \in \mathrm{Alg}_{\mathcal{P},\mathrm{un}}^{\mathrm{fact}}(X)$ on $X$ is

- an $\mathcal{P}$ algebra object $A \in \mathrm{Alg}_{\mathcal{P}}(D(\mathrm{Ran}_{X,\mathrm{un}}))$ internal to $D$ modules on the unital Ran space of $X$,
- an isomorphism $\alpha : j_{\mathrm{disj}}^! A^{\boxtimes 2} \xrightarrow{\cong} \sqcup^! A$ of $\mathcal{P}$ algebras internal to $D((\mathrm{Ran}_{X,\mathrm{un}}^{\times 2})_{\mathrm{disj}})$, and
- an isomorphism $\beta : \mathbb{K} \xrightarrow{\cong} A_{\varnothing}$ of $\mathcal{P}$ algebras in $\mathrm{Vect}_{\mathbb{K}}$,

together with compatible higher arity analogues of the multiplication map $\alpha$ satisfying the relations of a unital commutative monoid, where $\iota_{\varnothing} : \mathrm{pt} \hookrightarrow \mathrm{Ran}_{X,\mathrm{un}}$ is the map corresponding to the inclusion of the empty set.

*Remark* 6.0.4. The data of a factorization $\mathcal{P}$ algebra $A \in \mathrm{Alg}_{\mathcal{P},\mathrm{un}}^{\mathrm{fact}}(X)$ on $X$ can be further unpacked in terms of its restrictions to the strata of the Ran space, as in the case of a usual factorization algebra:

- The object $A \in \mathrm{Alg}_{\mathcal{P}}(D(\mathrm{Ran}_{X,\mathrm{un}}))$ is given by an assignment

$$I \mapsto A_I \in \mathrm{Alg}_{\mathcal{P}}(D(X^I)) \qquad [\pi : I \to J] \mapsto \Delta(\pi)^! A_I \to A_J$$

  for each $I \in \mathrm{fSet}$ and $\pi : I \to J$, such that the maps corresponding to $\pi$ surjective are isomorphisms;
- The factorization data $\alpha$ defines equivalences

$$j(\pi)^! A_I \xrightarrow{\cong} j(\pi)^! \boxtimes_{j \in J} A_{I_j}$$

  of $\mathcal{P}$ algebras internal to $D(X^I)$ for each $I$ and $\pi : I \twoheadrightarrow J$, in analogy with Definition I-6.0.2;
- The unit data $\beta$, together with the gluing data for $A \in \mathrm{Alg}_{\mathcal{P}}(D(\mathrm{Ran}_{X,\mathrm{un}}))$ above, defines the analogues of the factorization unit maps of Remark I-4.3.2. In particular, for $\pi : I \hookrightarrow J$ injective, or further in particular for $\pi : \varnothing \to J$, this gives maps

$$A_I \boxtimes \omega_{X^{I_\pi}} \to A_J \qquad \text{and} \qquad \omega_{X^J} \to \mathrm{C}_J$$

  of $\mathcal{P}$ algebras in $D(X^J)$, in analogy with Remark I-6.0.3; recall the notation $I_\pi$ and related conventions around partitions are given in Subsection I-2.2.

*Remark* 6.0.5. The preceeding proposition, definition, and remarks apply equal well in the quasicoherent format, and the $D$ module setting given above could be recovered as usual by taking $X = \tilde{X}_{\mathrm{dR}}$. However, since we do not heavily use the quasicoherent variant, we give the definitions in the $D$ module setting throughout this section.

*Remark* 6.0.6. Recall from Example I-C.1.15 that a (unital, cocommutative, Hopf) operad $\mathcal{P} \in \mathrm{HOp}^{\mathrm{co}}(\mathrm{Top})$ in spaces gives rise to such an object internal to $\mathrm{Vect}_{\mathbb{K}}$, which we denote by $\mathcal{P} = C_{\bullet}(\mathcal{P}; \mathbb{K}) \in \mathrm{HOp}^{\mathrm{co}}(\mathrm{Vect}_{\mathbb{K}})$ again, by abuse of notation. For such operads, we can define the category of algebras internal to a (not-necessarily stable or $\mathbb{K}$-linear) category. In this case, we have an alternative description of factorization $\mathcal{P}$ algebras, given by the following proposition.

Recall that $\mathrm{Alg}_{\mathrm{un}}^{\mathrm{fact}}(X)^{\otimes^!}$ denotes the category of unital factorization algebras equipped with the $\otimes^!$ tensor structure of Proposition I-6.0.5.

*Proposition* 6.0.7. Let $\mathcal{P} \in \mathrm{HOp}_{\mathrm{un}}^{\mathrm{co}}(\mathrm{Top})$ be a (unital) cocommutative Hopf algebra in spaces. There is an equivalence of categories

$$\mathrm{Alg}_{\mathcal{P},\mathrm{un}}^{\mathrm{fact}}(X) \xrightarrow{\cong} \mathrm{Alg}_{\mathcal{P}}(\mathrm{Alg}_{\mathrm{un}}^{\mathrm{fact}}(X)^{\otimes^!}) ,$$

which intertwines the natural forgetful functors to $\mathrm{Alg}_{\mathrm{un}}^{\mathrm{fact}}(X)$ and $\mathrm{Alg}_{\mathcal{P}}(D(\mathrm{Ran}_{X,\mathrm{un}}))$.

*Proof.* First, we recall the proof of Proposition I-6.0.5, which was previously postponed: given $A, B \in \mathrm{Alg}_{\mathrm{un}}^{\mathrm{fact}}(X)$, the tensor product of their underlying $D$ modules $A \otimes^! B \in D(\mathrm{Ran}_{X,\mathrm{un}})$ has canonical factorization structure defined by

$$j_{\mathrm{disj}}^!((A \otimes^! B)^{\boxtimes 2}) \cong j_{\mathrm{disj}}^!(A^{\boxtimes 2} \otimes^! B^{\boxtimes 2}) \cong j_{\mathrm{disj}}^!(A^{\boxtimes 2}) \otimes^! j_{\mathrm{disj}}^!(B^{\boxtimes 2}) \xrightarrow{\alpha_A \otimes^! \alpha_B} \sqcup^!(A) \otimes^! \sqcup^!(B) \cong \sqcup^!(A \otimes^! B) .$$

All of the above equivalences are natural in $A, B$, such that the above construction, together with the compatibility data of the underlying symmetric monoidal structure on $D(\mathrm{Ran}_{X,\mathrm{un}})$, defines a symmetric monoidal structure on $\mathrm{Alg}_{\mathrm{un}}^{\mathrm{fact}}(X)$.

We now give the proof of the proposition: the equivalence is fixed by the requirement that it intertwines the forgetful functors as stated, so it suffices to check that the factorization and unit data stipulated by the definitions of each of the categories are canonically equivalent.

The factorization structure on a factorization $\mathcal{P}$ algebra is given in terms of the underlying $\mathcal{P}$ algebra $A \in \mathrm{Alg}_{\mathcal{P}}(D(\mathrm{Ran}_X))$ by an isomorphism $\alpha : j_{\mathrm{disj}}^! A^{\boxtimes 2} \xrightarrow{\cong} \sqcup^!$ of $\mathcal{P}$ algebras internal to $D((\mathrm{Ran}_{X,\mathrm{un}}^{\times 2})_{\mathrm{disj}})$. The data witnessing that these equivalences are maps of $\mathcal{P}$ algebras is given by compatible commutativity data for

$$
\begin{array}{ccc}
j_{\mathrm{disj}}^!(A^{\boxtimes 2})^{\otimes^! n} & \xrightarrow{\cong} & (\sqcup^! A)^{\otimes^! n} \\
\downarrow & & \downarrow \\
j_{\mathrm{disj}}^!(A^{\boxtimes 2}) & \xrightarrow{\cong} & \sqcup^! A
\end{array} \quad ,
$$

for each $n \in \mathbb{N}$. Alternatively, the multiplication structure maps for a $\mathcal{P}$ algebra object in $\mathrm{Alg}_{\mathrm{un}}^{\mathrm{fact}}(X)^{\otimes^!}$ are maps of factorization algebras $A^{\otimes^! n} \to A$. The data witnessing that these structure maps respect the factorization structure is compatible commutativity data for

$$
\begin{array}{ccc}
j_{\mathrm{disj}}^!(A^{\otimes^! n})^{\boxtimes 2} & \xrightarrow{\cong} & \sqcup^!(A^{\otimes^! n}) \\
\downarrow & & \downarrow \\
j_{\mathrm{disj}}^!(A^{\boxtimes 2}) & \xrightarrow{\cong} & \sqcup^! A
\end{array} \quad .
$$

Now, by the preceeding construction of the factorization data on $\otimes^!$ tensor products of factorization algebras, there exists canonical commutativity data for

$$
\begin{array}{ccc}
j^!_{\mathrm{disj}}(A^{\otimes^! n})^{\boxtimes 2} & \xrightarrow{\;\cong\;} & \sqcup^!(A^{\otimes^! n}) \\
\downarrow & & \downarrow \\
j^!_{\mathrm{disj}}(A^{\boxtimes 2})^{\otimes^! n} & \xrightarrow{\;\cong\;} & (\sqcup^! A)^{\otimes^! n}
\end{array}
,
$$

so that the preceeding two diagrams are canonically equivalent. Similarly, the unit data in each case are tautologically equivalent. □

## 7. Factorization $\mathbb{E}_n$ algebras and factorization $\mathbb{E}_n$ categories

### 7.1. **Factorization $\mathbb{E}_n$-algebras.**

Let $\mathbb{E}_n = C_\bullet(\mathbb{E}_n; \mathbb{K}) \in \mathrm{HOp}^{\mathrm{co}}(\mathrm{Vect}_{\mathbb{K}})$ be the $\mathbb{E}_n$ operad, the cocommutative Hopf operad in Vect as in Definition I-C.4, and consider the category $\mathrm{Alg}^{\mathrm{fact}}_{\mathbb{E}_n,\mathrm{un}}(X)$ of unital factorization $\mathbb{E}_n$ algebras on $X$, following the general framework of Section 6 in the case $\mathcal{P} = \mathbb{E}_n$.

*Remark* 7.1.1. A primary motivation for the present work is the application of the category of unital factorization $\mathbb{E}_n$ algebras on a variety $X$ as an approximation to the notion of holomorphic-topological quantum field theory on the 'space-time manifold' $M = X \times \mathbb{R}^n$; the results in Part 2 of the present work should be understood as examples of this paradigm.

Due to the central role of this definition, we formally repeat the exposition of the previous subsection in this case for the reader's convenience:

*Definition* 7.1.2. A (unital) factorization $\mathbb{E}_n$ algebra $A \in \mathrm{Alg}^{\mathrm{fact}}_{\mathbb{E}_n,\mathrm{un}}(X)$ is a unital factorization $\mathrm{F}_{\mathbb{E}_n}$-object for $\mathrm{F}_{\mathbb{E}_n} := \mathrm{Alg}_{\mathbb{E}_n}(D^!(\cdot)) : \mathrm{Sch}^{\mathrm{op}}_{\mathrm{aff}} \to \mathrm{Cat}^\times$ together with the lax symmetric monoidal structure of Proposition 6.0.1.

*Remark* 7.1.3. Concretely, following remark 2.4.3, a factorization $\mathbb{E}_n$ algebra $A \in \mathrm{Alg}^{\mathrm{fact}}_{\mathbb{E}_n,\mathrm{un}}(X)$ on $X$ is

- an $\mathbb{E}_n$ algebra object $A \in \mathrm{Alg}_{\mathbb{E}_n}(D(\mathrm{Ran}_{X,\mathrm{un}}))$ internal to $D$ modules on the unital Ran space of $X$,
- an isomorphism $\alpha : j^!_{\mathrm{disj}} A^{\boxtimes 2} \xrightarrow{\cong} \sqcup^! A$ of $\mathbb{E}_n$ algebras internal to $D((\mathrm{Ran}^{\times 2}_{X,\mathrm{un}})_{\mathrm{disj}})$, and
- an isomorphism $\beta : \mathbb{K} \xrightarrow{\cong} A_\emptyset$ of $\mathbb{E}_n$ algebras in $\mathrm{Vect}_{\mathbb{K}}$,

together with compatible higher arity analogues of the multiplication map $\alpha$ satisfying the relations of a unital commutative monoid, where $\iota_\emptyset : \mathrm{pt} \hookrightarrow \mathrm{Ran}_{X,\mathrm{un}}$ is the map corresponding to the inclusion of the empty set.

*Remark* 7.1.4. The data of a factorization $\mathbb{E}_n$ algebra $A \in \mathrm{Alg}^{\mathrm{fact}}_{\mathbb{E}_n,\mathrm{un}}(X)$ on $X$ can be further unpacked in terms of its restrictions to the strata of the Ran space, as in the case of a usual factorization algebra:

- The object $A \in \mathrm{Alg}_{\mathbb{E}_n}(D(\mathrm{Ran}_{X,\mathrm{un}}))$ is given by an assignment

$$
I \mapsto A_I \in \mathrm{Alg}_{\mathbb{E}_n}(D(X^I)) \qquad [\pi : I \to J] \mapsto \Delta(\pi)^! A_I \to A_J
$$

for each $I \in \mathrm{fSet}$ and $\pi : I \to J$, such that the maps corresponding to $\pi$ surjective are isomorphisms;

- The factorization data $\alpha$ defines equivalences

$$j(\pi)^! A_I \overset{\cong}{\Rightarrow} j(\pi)^! \boxtimes_{j \in J} A_{I_j}$$

of $\mathbb{E}_n$ algebras internal to $D(X^I)$ for each $I$ and $\pi : I \twoheadrightarrow J$, in analogy with Definition I-6.0.2;

- The unit data $\beta$, together with the gluing data for $A \in \mathrm{Alg}_{\mathbb{E}_n}(D(\mathrm{Ran}_{X,\mathrm{un}}))$ above, defines the analogues of the factorization unit maps of Remark I-4.3.2. In particular, for $\pi : I \hookrightarrow J$ injective, or further in particular for $\pi : \emptyset \rightarrow J$, this gives maps

$$A_I \boxtimes \omega_{X^{I_\pi}} \rightarrow A_J \qquad \text{and} \qquad \omega_{X^J} \rightarrow C_J$$

of $\mathbb{E}_n$ algebras in $D(X^J)$, in analogy with Remark I-6.0.3; recall the notation $I_\pi$ and related conventions around partitions are given in Subsection I-2.2.

*Definition* 7.1.5. Let $C \in \mathrm{Cat}^{\mathrm{fact}}_{\mathrm{un}}(X)$ be a factorization category. A (unital) factorization algebra $\eta \in \mathrm{Alg}^{\mathrm{fact}}_{\mathrm{un}}(C)$ internal to C admits internal Hom objects if there exists a (unital) factorization functor

(7.1.1)

$$\mathcal{H}om_C(\eta, \cdot) : C \rightarrow \mathrm{QCoh}_{\mathrm{Ran}_{X,\mathrm{un}}} \qquad \text{with commutativity data for}$$

$$\begin{array}{ccc} C & \overset{\mathcal{H}om(\eta,\cdot)}{\longrightarrow} & \mathrm{QCoh}_{\mathrm{Ran}_{X,\mathrm{un}}} \\ & \searrow_{\mathrm{Hom}(\eta,\cdot)} & \downarrow_{\pi_\bullet} \\ & & \mathrm{Vect}_{\mathbb{K}} \end{array} \quad .$$

*Remark* 7.1.6. In particular, for $\eta \in \mathrm{Alg}^{\mathrm{fact}}_{\mathrm{un}}(C)$ admitting internal Hom objects over $X = \tilde{X}_{\mathrm{dR}}$ and $\xi \in \mathrm{Alg}^{\mathrm{fact}}_{\mathrm{un}}(C)$ an arbitrary factorization algebra internal to C, the internal Hom space

$$\mathcal{H}om_C(\eta, \xi) \in \mathrm{Alg}^{\mathrm{fact}}_{\mathrm{un}}(\tilde{X})$$

defines a unital factorization algebra.

*Example* 7.1.7. Let $C \in \mathrm{Cat}^{\mathrm{fact}}_{\mathrm{un}}(X)$ be a factorization category and $\eta \in \mathrm{Alg}^{\mathrm{fact}}_{\mathrm{un}}(C)$ a factorization algebra internal to C that admits internal Hom objects over $X = \tilde{X}_{\mathrm{dR}}$, as in the preceeding Definition. Then the internal Hom object

$$A = \mathcal{H}om_C(\eta, \eta) \in \mathrm{Alg}^{\mathrm{fact}}_{\mathbb{E}_1,\mathrm{un}}(\tilde{X})$$

is equipped with a natural factorization $\mathbb{E}_1$ structure, given by composition of endomorphisms.

As a special case of Proposition 6.0.7, we have:

*Corollary* 7.1.8. There is an equivalence of categories

$$\mathrm{Alg}^{\mathrm{fact}}_{\mathbb{E}_n,\mathrm{un}}(X) \overset{\cong}{\Rightarrow} \mathrm{Alg}_{\mathbb{E}_n}\big(\mathrm{Alg}^{\mathrm{fact}}_{\mathrm{un}}(X)^{\otimes^!}\big) \,,$$

which intertwines the natural forgetful functors to $\mathrm{Alg}^{\mathrm{fact}}_{\mathrm{un}}(X)$ and $\mathrm{Alg}_{\mathbb{E}_n}(D(\mathrm{Ran}_{X,\mathrm{un}}))$.

In terms of the preceeding description of factorization $\mathbb{E}_n$ algebras, we make the following definition:

*Definition* 7.1.9. Let $A \in \mathrm{Alg}^{\mathrm{fact}}_{\mathbb{E}_1,\mathrm{un}}(X) \cong \mathrm{Alg}_{\mathbb{E}_1}\big(\mathrm{Alg}^{\mathrm{fact}}_{\mathrm{un}}(X)^{\otimes^!}\big)$ be a (unital) factorization $\mathbb{E}_1$ algebra. A (unital) factorization $\mathbb{E}_1$ module $M \in A\text{-}\mathrm{Mod}(\mathrm{Alg}^{\mathrm{fact}}_{\mathrm{un}}(X)^{\otimes^!})$ over $A$ is a module object for $A$ internal to $\mathrm{Alg}^{\mathrm{fact}}_{\mathrm{un}}(X)^{\otimes^!}$.

*Example* 7.1.10. Let $C \in \mathrm{Cat}^{\mathrm{fact}}_{\mathrm{un}}(X)$ be a factorization category and $\eta \in \mathrm{Alg}^{\mathrm{fact}}_{\mathrm{un}}(C)$ a factorization algebra internal to C that admits internal Hom objects over $X = \tilde{X}_{\mathrm{dR}}$, as in Example 7.1.7, so that $A = \mathcal{H}om_C(\eta, \eta) \in \mathrm{Alg}^{\mathrm{fact}}_{\mathbb{E}_1, \mathrm{un}}(\tilde{X})$. Then for any factorization functor $\varphi : C \to \mathrm{QCoh}_{\mathrm{Ran}_{X,\mathrm{un}}}$ we have

$$\varphi(\eta) \in A\text{-}\mathrm{Mod}(\mathrm{Alg}^{\mathrm{fact}}_{\mathrm{un}}(\tilde{X})) \ .$$

In particular, for another factorization algebra $\xi \in \mathrm{Alg}^{\mathrm{fact}}_{\mathrm{un}}(C)$ admiting internal Hom objects, there is a natural $A$ module structure

$$\mathcal{H}om_C(\xi, \eta) \in A\text{-}\mathrm{Mod}(\mathrm{Alg}^{\mathrm{fact}}_{\mathrm{un}}(\tilde{X}))$$

on the internal Hom object $\mathcal{H}om_C(\xi, \eta) \in \mathrm{Alg}^{\mathrm{fact}}_{\mathrm{un}}(X)$, given by composition of maps.

## 7.2. Factorization $\mathbb{E}_n$-categories and factorization $\mathbb{E}_n$-algebra objects.

*Example* 7.2.1. Let $\mathrm{ShvCat}_{\mathbb{E}_n}(\cdot) = \mathrm{Alg}_{\mathbb{E}_n}(\mathrm{ShvCat}(\cdot)) : \mathrm{Sch}^{\mathrm{op}}_{\mathrm{aff}} \to \mathrm{Cat}^\times$ be the functor of $\mathbb{E}_n$-monoidal sheaves of categories. This functor admits a natural lax symmetric monoidal structure directly generalizing that of Proposition 6.0.1 by replacing the role of the lax monoidal structure of Example 2.2.2 with that of Example 2.2.4.

*Definition* 7.2.2. A (unital) factorization $\mathbb{E}_n$ category $C^\otimes$ on $X$ is a (unital) factorization $\mathrm{ShvCat}_{\mathbb{E}_n}$-object for $\mathrm{ShvCat}_{\mathbb{E}_n} : \mathrm{Sch}^{\mathrm{op}}_{\mathrm{aff}} \to \mathrm{Cat}^\times$ the lax monoidal functor defined in the preceeding Remark.

Let $\mathrm{Cat}^{\mathrm{fact}}_{\mathbb{E}_n}(X) = \mathrm{ShvCat}^{\mathrm{fact}}_{\mathbb{E}_n}(X)$ and $\mathrm{Cat}^{\mathrm{fact}}_{\mathbb{E}_n, \mathrm{un}}(X) = \mathrm{ShvCat}^{\mathrm{fact}}_{\mathbb{E}_n, \mathrm{un}}(X)$ denote the (2-)categories of non-unital and unital factorization $\mathbb{E}_n$-categories.

*Remark* 7.2.3. The preceeding definition is the categorical analogue of the quasicoherent variant of factorization $\mathbb{E}_n$ algebra; see also remark 6.0.5. A factorization $\mathbb{E}_n$ category on $X = \tilde{X}_{\mathrm{dR}}$ is the categorical analogue of an object in $\mathrm{Alg}^{\mathrm{fact}}_{\mathbb{E}_n, \mathrm{un}}(X)$.

*Remark* 7.2.4. Concretely, following Remark 3.1.4, a (unital) factorization $\mathbb{E}_n$ category $C \in \mathrm{Cat}^{\mathrm{fact}}_{\mathbb{E}_n, \mathrm{un}}(X)$ on $X$ can be unpacked in terms of its restrictions to the strata of the Ran space, as follows:

- An object $C \in \mathrm{Alg}_{\mathbb{E}_n}(\mathrm{ShvCat}(\mathrm{Ran}_{X, \mathrm{un}})^{\otimes *})$ is given by an assignment

$$(7.2.1) \qquad I \mapsto C_I \in \mathrm{Alg}_{\mathbb{E}_n}(\mathrm{ShvCat}(X^I)^{\otimes *}) \qquad [\pi : I \to J] \mapsto [\Phi_\pi : \Delta(\pi)^* \, C_I \to C_J]$$

for each $I \in \mathrm{fSet}$ and $\pi : I \to J$, such that the maps corresponding to $\pi$ surjective are isomorphisms;
- the factorization data $\alpha$ defines equivalences

$$(7.2.2) \qquad j(\pi)^* \, C_I \xrightarrow{\cong} j(\pi)^* \boxtimes_{j \in J} C_{I_j}$$

of sheaves of $\mathbb{E}_n$ categories on $U(\pi) \subset X^I$ for each $I$ and $\pi : I \twoheadrightarrow J$; and
- the unit data $\beta$, together with the gluing data for $C \in \mathrm{Alg}_{\mathbb{E}_n}(\mathrm{ShvCat}(\mathrm{Ran}_{X, \mathrm{un}}))$ above, defines the analogues of the factorization unit maps of Remark I-4.3.2. In particular, for $\pi : I \hookrightarrow J$ injective, or further in particular for $\pi : \emptyset \to J$, this gives maps

$$(7.2.3) \qquad \mathrm{unit}^\pi_C : C_I \boxtimes \mathrm{QCoh}_{X^{I\pi}} \to C_J \qquad \text{and} \qquad \mathrm{unit}^J_C : \mathrm{QCoh}_{X^J} \to C_J$$

of sheaves of $\mathbb{E}_n$ categories on $X^J$.

*Example* 7.2.5. Let $\mathcal{Y} \in \mathrm{PreStk}^{\mathrm{fact}}_{\mathrm{un}}(X)$ be a factorization prestack on $X$ such that the category $\mathrm{QCoh}_{\mathcal{Y}} \in \mathrm{Cat}^{\mathrm{fact}}_{\mathrm{un}}(X)$ admits the structure of a unital factorization category, as for example in Definition 5.2.7. Then $\mathrm{QCoh}_{\mathcal{Y}} \in \mathrm{Cat}^{\mathrm{fact}}_{\mathbb{E}_\infty, \mathrm{un}}(X)$ is canonically a factorization symmetric monoidal category, under the usual tensor product of quasicoherent sheaves.

*Example* 7.2.6. The unit factorization category $\mathrm{QCoh}_{\mathrm{Ran}_{X,\mathrm{un}}} \in \mathrm{Cat}^{\mathrm{fact}}_{\mathbb{E}_\infty,\mathrm{un}}(X)$ is canonically a factorization symmetric monoidal category, under the usual tensor product of quasicoherent sheaves, by the preceeding example. In particular, applying the symmetric monoidal forgetful functor $\mathrm{Alg}_{\mathbb{E}_\infty} \to \mathrm{Alg}_{\mathbb{E}_n}$ as in Remark I-C.4.12, $\mathrm{QCoh}_{\mathrm{Ran}_{X,\mathrm{un}}} \in \mathrm{Cat}^{\mathrm{fact}}_{\mathbb{E}_n,\mathrm{un}}(X)$ defines a factorization $\mathbb{E}_n$-category for each $n \in \mathbb{N}$.

*Definition* 7.2.7. Let $\mathrm{C} \in \mathrm{Cat}^{\mathrm{fact}}_{\mathbb{E}_n,\mathrm{un}}(X)$ be a (unital) factorization $\mathbb{E}_n$ category. A (unital) factorization $\mathbb{E}_n$-algebra $A \in \mathrm{Alg}^{\mathrm{fact}}_{\mathbb{E}_n,\mathrm{un}}(\mathrm{C})$ internal to $\mathrm{C}$ is a map of (unital) factorization $\mathbb{E}_n$ categories $\mathrm{QCoh}_{\mathrm{Ran}_{X,\mathrm{un}}} \to \mathrm{C}$.

*Remark* 7.2.8. The data of such an object, or more generally a functor between factorization $\mathbb{E}_n$ categories, can be unpacked explicitly in terms of the description of Remark 7.2.4, following Remarks 3.2.6, and 3.2.2 and 3.2.3, respectfully.

In particular, a factorization $\mathbb{E}_n$-algebra $A \in \mathrm{Alg}^{\mathrm{fact}}_{\mathbb{E}_n,\mathrm{un}}(\mathrm{C})$ internal to $\mathrm{C}$ is given by:

- an $\mathbb{E}_n$ algebra object $A \in \mathrm{Alg}_{\mathbb{E}_n}(\Gamma(\mathrm{Ran}_{X,\mathrm{un}}, \mathrm{C}))$, given by an assignment

$$(7.2.4) \qquad I \mapsto A_I \in \mathrm{Alg}_{\mathbb{E}_n}(\Gamma(X^I, \mathrm{C}_I)) \qquad [\pi : I \to J] \mapsto [\eta_\pi : \Phi_\pi(\Delta(\pi)^\bullet A_I) \to A_J] \ ,$$

such that for $\pi$ surjective the corresponding maps in $\mathrm{Alg}_{\mathbb{E}_n}(\Gamma(X^J, \mathrm{C}_J))$ are isomorphisms;
- an equivalence

$j(\pi)^\bullet A_I \xrightarrow{\cong} j(\pi)^\bullet(\boxtimes_j A_{I_j}) \qquad$ in the category $\qquad \mathrm{Alg}_{\mathbb{E}_n}(\Gamma(U(\pi), j(\pi)^* \mathrm{C}_I)) \cong \mathrm{Alg}_{\mathbb{E}_n}(\Gamma(U(\pi), j(\pi)^* \boxtimes_{j \in J} \mathrm{C}_{I_j}))$

for each $I$ and $\pi : I \twoheadrightarrow J$; and
- an equivalence $A_\emptyset \cong \mathbb{K}$ in $\mathrm{Alg}_{\mathbb{E}_n}(\mathrm{C}_\emptyset) \cong \mathrm{Alg}_{\mathbb{E}_n}(\mathrm{Vect})$, which together with the structure maps of Equation 7.2.4 above, determine maps

$$\mathrm{unit}^\pi_{\mathrm{C}}(A_I \boxtimes \mathbb{O}_{X^{I\pi}}) \to A_J \qquad \text{and} \qquad \mathrm{unit}^J_{\mathrm{C}}(\mathbb{O}_{X^J}) \to A_J$$

for $\pi : I \hookrightarrow J$ injective, and in particular for $\pi : \emptyset \hookrightarrow J$, respectively.

*Example* 7.2.9. The factorization unit $\mathrm{unit}_{\mathrm{C}} \in \mathrm{Alg}^{\mathrm{fact}}_{\mathbb{E}_n,\mathrm{un}}(\mathrm{C})$ of a unital factorization $\mathbb{E}_n$-category is canonically a factorization $\mathbb{E}_n$-algebra internal to $\mathrm{C}$.

*Example* 7.2.10. Following Remark 6.0.5, the analogue of Example 3.2.9 is as follows: the category of unital factorization $\mathbb{E}_n$ algebras $\mathrm{Alg}^{\mathrm{fact}}_{\mathbb{E}_n,\mathrm{un}}(\mathrm{QCoh}_{\mathrm{Ran}_{X,\mathrm{un}}})$ internal to $\mathrm{QCoh}_{\mathrm{Ran}_{X,\mathrm{un}}}$ is equivalent to the usual category of factorization $\mathbb{E}_n$ quasicoherent sheaves. In particular, taking $X = \tilde{X}_{\mathrm{dR}}$ as in Example 3.2.10 and following Remark 3.2.11, there is an equivalence

$$\mathrm{Alg}^{\mathrm{fact}}_{\mathbb{E}_n,\mathrm{un}}(\mathrm{QCoh}_{\mathrm{Ran}_{X,\mathrm{un}}}) \cong \mathrm{Alg}^{\mathrm{fact}}_{\mathbb{E}_n,\mathrm{un}}(\tilde{X}) \ .$$

In analogy with Proposition 7.1.8, we have

*Proposition* 7.2.11. There is an equivalence of categories

$$\mathrm{Cat}^{\mathrm{fact}}_{\mathbb{E}_n,\mathrm{un}}(X) \xrightarrow{\cong} \mathrm{Alg}_{\mathbb{E}_n}(\mathrm{Cat}^{\mathrm{fact}}_{\mathrm{un}}(X)^*) \ ,$$

which intertwines the natural forgetful functors to $\mathrm{Cat}^{\mathrm{fact}}_{\mathrm{un}}(X)$ and $\mathrm{Alg}_{\mathbb{E}_n}(\mathrm{ShvCat}(\mathrm{Ran}_{X,\mathrm{un}})^{\otimes *})$.

In analogy with Definition 7.1.9, we make the following definition:

*Definition* 7.2.12. Let $\mathrm{C} \in \mathrm{Cat}^{\mathrm{fact}}_{\mathbb{E}_1,\mathrm{un}}(X) \cong \mathrm{Alg}_{\mathbb{E}_1}(\mathrm{Cat}^{\mathrm{fact}}_{\mathrm{un}}(X)^*)$ be a (unital) factorization $\mathbb{E}_1$ category. A (unital) factorization $\mathbb{E}_1$ module category $\mathrm{D} \in \mathrm{C}\text{-}\mathrm{Mod}(\mathrm{Cat}^{\mathrm{fact}}_{\mathrm{un}}(X)^*)$ over $\mathrm{C}$ is a module object for $\mathrm{C}$ internal to $\mathrm{Cat}^{\mathrm{fact}}_{\mathrm{un}}(X)^*$.

*Example* 7.2.13. A factorization $\mathbb{E}_1$ category $C \in \mathrm{Cat}^{\mathrm{fact}}_{\mathbb{E}_1,\mathrm{un}}(X)$ is canonically an $\mathbb{E}_1$ module category over itself $C \in C\text{-Mod}(\mathrm{Cat}^{\mathrm{fact}}_{\mathrm{un}}(X)^\star)$ via the monoidal structure multiplication map.

In this context, we extend the definition of factorization $\mathbb{E}_1$ module to the internal setting:

*Definition* 7.2.14. Let $C \in \mathrm{Cat}^{\mathrm{fact}}_{\mathbb{E}_n,\mathrm{un}}(X)$ be a (unital) factorization $\mathbb{E}_n$ category, $D \in C\text{-Mod}(\mathrm{Cat}^{\mathrm{fact}}_{\mathrm{un}}(X)^\star)$ a (unital) factorization $\mathbb{E}_1$ module category over $C$, and $A \in \mathrm{Alg}^{\mathrm{fact}}_{\mathbb{E}_n,\mathrm{un}}(C)$ a (unital) factorization $\mathbb{E}_n$-algebra internal to $C$.

A (unital) factorization $\mathbb{E}_1$ module object $M \in A\text{-Mod}(\mathrm{Alg}^{\mathrm{fact}}_{\mathrm{un}}(D))$ over $A$ internal to $D$ is a module object for $A$ internal to $\mathrm{Alg}^{\mathrm{fact}}_{\mathrm{sun}}(D)$.

*Remark* 7.2.15. Concretely, a factorization $\mathbb{E}_1$ module object $M \in A\text{-Mod}(\mathrm{Alg}^{\mathrm{fact}}_{\mathrm{un}}(D))$ is given by a factorization algebra $M \in \mathrm{Alg}^{\mathrm{fact}}_{\mathrm{un}}(D)$ internal to $D$, together with a module structure map

$$\rho : A \star M \to M \qquad \text{where} \qquad (\cdot) \star (\cdot) : \mathrm{Alg}^{\mathrm{fact}}_{\mathrm{sun}}(C) \times \mathrm{Alg}^{\mathrm{fact}}_{\mathrm{sun}}(D) \to \mathrm{Alg}^{\mathrm{fact}}_{\mathrm{sun}}(D)$$

denotes the map on factorization algebras induced by the $C$ module structure on $D$, together with its higher arity analogues and compatibility data with the $\mathbb{E}_1$ algebra structure data for $A$.

*Remark* 7.2.16. Further, in terms of the description of $A \in \mathrm{Alg}^{\mathrm{fact}}_{\mathbb{E}_1,\mathrm{un}}(C)$ in Remark 7.2.8, a factorization $\mathbb{E}_1$ module object $M \in A\text{-Mod}(\mathrm{Alg}^{\mathrm{fact}}_{\mathrm{un}}(D))$ over $A$ internal to $D$ is given by:

- an $A$ module object $M \in A\text{-Mod}(\Gamma(\mathrm{Ran}_{X,\mathrm{un}},D))$, given by an assignment

$$(7.2.5) \qquad I \mapsto M_I \in A_I\text{-Mod}(\Gamma(X^I,D_I)) \qquad\qquad [\pi : I \to J] \mapsto [\eta_\pi : \Phi_\pi(\Delta(\pi)^\bullet M_I) \to M_J] \ ,$$

  such that for $\pi$ surjective the corresponding maps in $A_J\text{-Mod}(\Gamma(X^J, D_J))$ are isomorphisms;
- an equivalence

$$j(\pi)^\bullet M_I \xrightarrow{\cong} j(\pi)^\bullet(\boxtimes_j M_{I_j}) \qquad \text{in} \qquad j(\pi)^\bullet A_I\text{-Mod}(\Gamma(U(\pi), j(\pi)^* D_I))$$

$$\cong j(\pi)^\bullet(\boxtimes_j A_{I_j})\text{-Mod}(\Gamma(U(\pi), j(\pi)^* \boxtimes_{j \in J} D_{I_j}))$$

  for each $I$ and $\pi : I \twoheadrightarrow J$; and
- an equivalence $M_\emptyset \cong \mathbb{K}$ in $A_\emptyset\text{-Mod}(D_\emptyset) \cong \mathrm{Vect}$, which together with the structure maps of Equation 7.2.5 above, determine maps

$$\mathrm{unit}^\pi_C(M_I \boxtimes \mathcal{O}_{X^{I_\pi}}) \to M_J \qquad \text{and} \qquad \mathrm{unit}^J_C(\mathcal{O}_{X^J}) \to M_J$$

  for $\pi : I \hookrightarrow J$ injective, and in particular for $\pi : \emptyset \hookrightarrow J$, respectively.

## 8. Convolution constructions for factorization $\mathbb{E}_n$ algebras and categories

### 8.1. Convolution on factorization spaces and constructions of factorization $\mathbb{E}_n$ algebras and categories.
In this subsection, we apply the format of Subappendix A.6 to outline a general construction of factorization $\mathbb{E}_1$ algebras, and factorization $\mathbb{E}_1$ categories, via convolution.

*Proposition* 8.1.1. Let $\mathcal{X}, \mathcal{Y} \in \mathrm{PreStk}^{\mathrm{fact}}_{\mathrm{un}}(X)$ be (unital) factorization spaces and $f : \mathcal{X} \to \mathcal{Y}$ a strict map of (unital) factorization spaces. Then $\mathcal{Z} = \mathcal{X} \times_{\mathcal{Y}} \mathcal{X} \in \mathrm{PreStk}^{\mathrm{fact}}_{\mathrm{un}}(X)$ is naturally a (unital) factorization space, and the canonical projection maps $\pi_{\mathcal{X},i} : \mathcal{Z} \to \mathcal{X}$ are maps of (unital) factorization spaces for $i = 1, 2$.

*Proof.* Recall that the data defining the strict map $f : \mathcal{X} \to \mathcal{Y}$ of underlying relative prestacks over $\mathrm{Ran}_{X,\mathrm{un}}$ is given in Equation 4.3.1. We define a relative prestack $\mathcal{Z} = \mathcal{X} \times_{\mathcal{Y}} \mathcal{X}$ over $\mathrm{Ran}_{X,\mathrm{un}}$ by the assignments

$$
\begin{cases}
I & \mapsto \mathcal{Z}_I := \mathcal{X}_I \times_{\mathcal{Y}_I} \mathcal{X}_I \\[2em]
[\pi : I \twoheadrightarrow J] & \mapsto
\end{cases}
$$

$$
\xymatrix{ & \mathrm{unit}_{\mathcal{X}}^{\pi} \underset{\mathrm{unit}_{\mathcal{Y}}^{\pi}}{\times} \mathrm{unit}_{\mathcal{X}}^{\pi} & \\
{}^{\psi_{\mathcal{Z}}^{\pi}} \swarrow & & \searrow {}^{\phi_{\mathcal{Z}}^{\pi}} \\
X^J \times_{X^I} (\mathcal{X}_I \times_{\mathcal{Y}_I} \mathcal{X}_I) & & \mathcal{X}_J \times_{\mathcal{Y}_J} \mathcal{X}_J }
$$

where the maps are those induced on the given fibre products from the corresponding maps in Equation 5.3.1. The factorization data is defined for each $\pi : I \twoheadrightarrow J$ by the natural equivalences

$$
U(\pi) \times_{X^I} \mathcal{Z}_I = U(\pi) \times_{X^I} (\mathcal{X}_I \times_{\mathcal{Y}_I} \mathcal{X}_I) \xrightarrow{\cong} U(\pi) \times_{X^I} \left( \prod_{j \in J} \mathcal{X}_{I_j} \underset{\left( \prod_{j \in J} \mathcal{Y}_{I_j} \right)}{\times} \prod_{j \in J} \mathcal{X}_{I_j} \right) = U(\pi) \times_{X^I} \left( \prod_{j \in J} \mathcal{Z}_{I_j} \right)
$$

given by the fibre product of those defined for $\mathcal{X}$ and $\mathcal{Y}$, and similarly for the unit data. The canonical projections $\pi_{\mathcal{X}} : \mathcal{Z} \to \mathcal{X}$ evidently admit natural factorization compatibility data in the sense of Equation 4.3.1 and thus define maps of factorization spaces.    $\square$

*Remark* 8.1.2. The preceeding proposition and its proof evidently apply to general fibre products of distict factorization spaces $\mathcal{X} \times_{\mathcal{Y}} \mathcal{W}$, and is the definition of fibre products in the category of (unital) factorization spaces.

*Remark* 8.1.3. More generally, the iterated fibre products

$$
\mathcal{Z}_{(n)} = \mathcal{X} \times_{\mathcal{Y}} \mathcal{X} \times_{\mathcal{Y}} \ldots \times_{\mathcal{Y}} \mathcal{X} \quad \in \mathrm{PreStk}_{\mathrm{un}}^{\mathrm{fact}}(X) \qquad \text{and projections} \qquad \pi_{ij} : \mathcal{Z}_{(n)} \to \mathcal{Z}
$$

define factorization spaces and maps of such for each $n \in \mathbb{N}$ and $i, j \in \{1, ..., n\}$.

Towards stating the analogue of Proposition A.6.9 in the (quasicoherent) factorization setting, we restrict to the following format:

*Remark* 8.1.4. Throughout the remainder of this subsection, let $\mathcal{X}, \mathcal{Y} \in \mathrm{PreStk}_{\mathrm{un}}^{\mathrm{fact}}(X)$, $f : \mathcal{X} \to \mathcal{Y}$ and $\mathcal{Z} \in \mathrm{PreStk}_{\mathrm{un}}^{\mathrm{fact}}(X)$ be as in the preceeding proposition, and such that $\mathrm{QCoh}_{\mathcal{Z}_{(n)}} \in \mathrm{Cat}_{\mathrm{un}}^{\mathrm{fact}}(X)$ define (unital) factorization categories for each $n \in \mathbb{N}$; for example, consider $\mathcal{X}, \mathcal{Y} \in \mathrm{PreStk}_{(\mathrm{all,sch\text{-}qcqs}),\mathrm{un}}^{\mathrm{fact}}$ as in Example 5.2.5 and $f : \mathcal{X} \to \mathcal{Y}$ a map in this category.

Further, suppose that for each $i, j$ the functors $\pi_{ij}^{\bullet} : \mathrm{QCoh}_{\mathcal{Z}} \to \mathrm{QCoh}_{\mathcal{Z}_{(n)}}$ of pullback along the maps $\pi_{ij} : \mathcal{Z}_{(n)} \to \mathcal{Z}$ induce (unital) factorization functors; for example, as in Example 5.3.2 and Remark 5.3.3.

*Proposition* 8.1.5. The factorization category $\mathrm{QCoh}_{\mathcal{Z}} \in \mathrm{Cat}_{\mathbb{E}_1,\mathrm{un}}^{\mathrm{fact}}(X)$ of quasicoherent sheaves on $\mathcal{Z}$ is naturally a (unital) $\mathbb{E}_1$-factorization category with respect to the convolution monoidal structure $(\cdot) \star (\cdot) : \mathrm{QCoh}_{\mathcal{Z}} \otimes^{\star} \mathrm{QCoh}_{\mathcal{Z}} \to \mathrm{QCoh}_{\mathcal{Z}}$ defined by the composition

$$
\mathrm{QCoh}_{\mathcal{Z}} \otimes^{\star} \mathrm{QCoh}_{\mathcal{Z}} \xrightarrow{\pi_{12}^{\bullet} \boxtimes \pi_{23}^{\bullet}} \mathrm{QCoh}_{\mathcal{Z}_{(3)}^{\times 2}} \xrightarrow{\Delta^{\bullet}} \mathrm{QCoh}_{\mathcal{Z}_{(3)}} \xrightarrow{\pi_{13,\bullet}} \mathrm{QCoh}_{\mathcal{Z}} \ .
$$

Further, the pushforward functor $\mathrm{p}_{\mathcal{Z}\bullet} : \mathrm{QCoh}_{\mathcal{Z}}^{\star} \to \mathrm{QCoh}_{\mathrm{Ran}_{X,\mathrm{un}}}^{\otimes !}$ is a unital, lax $\mathbb{E}_1$-monoidal, factorization functor, with respect to the convolution monoidal structure on $\mathrm{QCoh}_{\mathcal{Z}}$ defined above, and the symmetric monoidal structure on $\mathrm{QCoh}_{\mathrm{Ran}_{X,\mathrm{un}}}$ of Example 7.2.6. In particular, there exist natural maps of factorization quasicoherent sheaves

$$\mathrm{p}_{\mathcal{Z}\bullet}\mathcal{M} \otimes \mathrm{p}_{\mathcal{Z}\bullet}\mathcal{N} \to \mathrm{p}_{\mathcal{Z}\bullet}(\mathcal{M} \star \mathcal{N})$$

for each $\mathcal{M}, \mathcal{N} \in \mathrm{Alg}^{\mathrm{fact}}_{\mathrm{un}}(\mathrm{QCoh}_{\mathcal{Z}})$.

*Proof.* By Proposition 7.2.11 above, it suffices to prove the analogous statements for $\mathbb{E}_1$ algebras internal to factorization categories. Using the strategy outlined in Remark A.6.11, it suffices to apply the proof of Proposition A.6.9 in this setting. The hypotheses of the preceding remark guarantee that the required factorization functors exist, and that the underlying functors of quasicoherent pushforward and pullback satisfy the required adjointness properties of the $*$ pushforward and pullback functors of Proposition A.6.9. □

Following propositions A.6.2 and A.6.10, we also have:

*Corollary* 8.1.6. Let $\mathcal{A} \in \mathrm{Alg}^{\mathrm{fact}}_{\mathbb{E}_1,\mathrm{un}}(\mathrm{QCoh}_{\mathcal{Z}}^{\star})$ be a factorization $\mathbb{E}_1$-algebra internal to $\mathrm{QCoh}_{\mathcal{Z}}$ over $X = \tilde{X}_{\mathrm{dR}}$. Then $\mathrm{p}_{\mathcal{Z}\bullet}\mathcal{A} \in \mathrm{Alg}^{\mathrm{fact}}_{\mathbb{E}_1,\mathrm{un}}(\tilde{X})$ is naturally an factorization $\mathbb{E}_1$-algebra on $\tilde{X}$, so that we obtain a functor

$$\mathrm{p}_{\mathcal{Z}\bullet} : \mathrm{Alg}^{\mathrm{fact}}_{\mathbb{E}_1,\mathrm{un}}(\mathrm{QCoh}_{\mathcal{Z}}^{\star}) \to \mathrm{Alg}^{\mathrm{fact}}_{\mathbb{E}_1,\mathrm{un}}(\tilde{X}) \ .$$

*Remark* 8.1.7. Concretely, following Remark 7.2.8, a factorization $\mathbb{E}_1$-algebra $\mathcal{A} \in \mathrm{Alg}^{\mathrm{fact}}_{\mathbb{E}_1,\mathrm{un}}(\mathrm{QCoh}_{\mathcal{Z}})$ internal to $\mathrm{QCoh}_{\mathcal{Z}}$ is given by:

- an $\mathbb{E}_1$ algebra object $\mathcal{A} \in \mathrm{Alg}_{\mathbb{E}_1}(\mathrm{QCoh}(\mathcal{Z}))$, given by an assignment

$$(8.1.1) \qquad I \mapsto A_I \in \mathrm{Alg}_{\mathbb{E}_1}(\mathrm{QCoh}(\mathcal{Z}_I)^{\star}) \qquad [\pi : I \to J] \mapsto \left[ \eta_\pi : \phi^{\pi}_{\mathcal{Z}\bullet}\psi^{\pi,\bullet}_{\mathcal{Z}}\tilde{\Delta}(\pi)^{\bullet}\mathcal{A}_I \to \mathcal{A}_J \right] \ ,$$

  such that for $\pi$ surjective the corresponding maps in $\mathrm{Alg}_{\mathbb{E}_1}(\mathrm{QCoh}(\mathcal{Z}_J))$ are isomorphisms, where the former are given by an underlying object $\mathcal{A}_I \in \mathrm{QCoh}(\mathcal{Z}_I)$ together with a map

$$(8.1.2) \qquad \pi_{13,I\bullet}(\pi^{\bullet}_{12,I}A_I \otimes \pi^{\bullet}_{23,I}A_I) \to A_I \quad \text{in} \quad \mathrm{QCoh}(\mathcal{Z}_I)$$

  and its higher arity analogues, as in Proposition A.6.2, where $\pi_{ij,I} : \mathcal{Z}_{(3)I} \to \mathcal{Z}_I$ are the components of the maps of factorization spaces in Remark 8.1.3 over $X^I$;

- an equivalence

$$\tilde{j}(\pi)^{\bullet}A_I \xrightarrow{\cong} \tilde{j}(\pi)^{\bullet}(\boxtimes_j A_{I_j}) \qquad \text{in the category} \qquad \mathrm{Alg}_{\mathbb{E}_1}(\mathrm{QCoh}(U(\pi)\times_{X^I}\mathcal{Z}_I)) \cong \mathrm{Alg}_{\mathbb{E}_1}(\mathrm{QCoh}(U(\pi)\times_{X^I}(\times_{j\in J}\mathcal{Z}_{I_j})))$$

  for each $I$ and $\pi : I \to J$; and

- an equivalence $A_{\varnothing} \cong \mathbb{K}$ in $\mathrm{Alg}_{\mathbb{E}_1}(\mathrm{QCoh}(\mathcal{Z}_{\varnothing})^{\star}) \cong \mathrm{Alg}_{\mathbb{E}_1}(\mathrm{Vect})$, which together with the structure maps of Equation 8.1.1 above, determine maps

$$(8.1.3) \qquad \phi_{\mathcal{Z},\pi\bullet}\psi^{\bullet}_{\mathcal{Z},\pi}(A_I\boxtimes\mathcal{O}_{X^{I\pi}}) \to A_J \qquad \text{and} \qquad \phi_{\mathcal{Z},J\bullet}\psi^{\bullet}_{\mathcal{Z},J}(\mathcal{O}_{X^J}) \to A_J$$

  for $\pi : I \hookrightarrow J$ injective, and in particular for $\pi : \varnothing \hookrightarrow J$, respectively.

*Remark* 8.1.8. In particular, following the preceding remark, the preceding corollary implies the factorization analogue of Proposition A.6.2, analogously to the fact that *loc. cit.* follows from Proposition A.6.10.

*Definition* 8.1.9. Let $\mathcal{N} \in \mathrm{Alg}_{\mathrm{un}}^{\mathrm{fact}}(\mathrm{QCoh}_{\mathcal{Y}})$ be a (unital) factorization compatible quasicoherent sheaf on $\mathcal{Y}$ that admits internal Hom objects in the sense of Definition 8.1.9. Then $\mathcal{N}$ admits internal Hom objects over $\mathcal{Y}$ if there is a lift of the functor of Equation 7.1.1 to a unital factorization functor

$$\tilde{\mathcal{Hom}}_{\mathrm{QCoh}(\mathcal{Y})}(\mathcal{N}, \cdot) : \mathrm{QCoh}_{\mathcal{Y}} \to \mathrm{QCoh}_{\mathcal{Y}} \qquad \text{with commutativity data for}$$

*Remark* 8.1.10. Concretely, for $\mathcal{N} \in \mathrm{Alg}_{\mathrm{un}}^{\mathrm{fact}}(\mathrm{QCoh}_{\mathcal{Y}})$ admitting internal Hom objects over $X$, we have

$$\mathcal{Hom}_{\mathrm{QCoh}(\mathcal{Y})}(\mathcal{N}, \mathcal{M}) \in \mathrm{QCoh}_{\mathrm{un}}^{\mathrm{fact}}(X) \qquad \text{defined by} \qquad I \mapsto \mathrm{p}_{\mathcal{Y}_I \bullet} \underline{\mathrm{Hom}}_{\mathrm{QCoh}(\mathcal{Y}_I)}(N_I, M_I) \in \mathrm{QCoh}(X^I) \,,$$

together with the natural induced compatibility data, for each $\mathcal{M} \in \mathrm{Alg}_{\mathrm{un}}^{\mathrm{fact}}(\mathcal{Y})$.

Similarly, for $\mathcal{N} \in \mathrm{Alg}_{\mathrm{un}}^{\mathrm{fact}}(\mathrm{QCoh}_{\mathcal{Y}})$ admitting internal Hom objects over $\mathcal{Y}$, we have

$$\tilde{\mathcal{Hom}}_{\mathrm{QCoh}(\mathcal{Y})}(\mathcal{N}, \mathcal{M}) \in \mathrm{Alg}_{\mathrm{un}}^{\mathrm{fact}}(\mathrm{QCoh}_{\mathcal{Y}}) \qquad \text{defined by} \qquad I \mapsto \underline{\mathrm{Hom}}_{\mathrm{QCoh}(\mathcal{Y}_I)}(N_I, M_I) \in \mathrm{QCoh}(\mathcal{Y}_I) \,.$$

We now state the analogues of Remark A.6.5 and examples A.6.6 and A.6.7 in the quasicoherent, factorization setting. To begin, we recall the quasicoherent analogue of Remark A.6.5:

*Remark* 8.1.11. Let $X, Y$ and $W \in \mathrm{Sch}_{\mathrm{ft}}$ be finite type schemes, $M \in \mathrm{QCoh}(X)$, $N \in \mathrm{Coh}(W)$, $f : X \to Y$, $g : W \to Y$ proper, and $Z = X \times_Y W$. Then the internal Hom object

$$\underline{\mathrm{Hom}}_{\mathrm{QCoh}(Y)}(g_\bullet N, f_\bullet M) \in \mathrm{QCoh}(Y)$$

can be calculated as

$$\begin{aligned}
\underline{\mathrm{Hom}}_{\mathrm{QCoh}(Y)}(g_\bullet N, f_\bullet M) &\cong g_\bullet \underline{\mathrm{Hom}}_{\mathrm{QCoh}(W)}(N, g^! f_\bullet M) \\
&\cong g_\bullet \underline{\mathrm{Hom}}_{\mathrm{QCoh}(W)}(N, \pi_{W \bullet} \pi_X^\bullet M \otimes \omega_{W/Y}) \\
&\cong g_\bullet \pi_{W \bullet} \underline{\mathrm{Hom}}_{\mathrm{QCoh}(X \times_Y W)}(\pi_W^\bullet N, \pi_X^\bullet M \otimes \omega_{Z/X})
\end{aligned}$$

so that we obtain a lift of the internal Hom object to $\mathrm{QCoh}(X \times_Y W)$, which we denote

$$\tilde{\underline{\mathrm{Hom}}}_{\mathrm{QCoh}(Y)}(g_\bullet N, f_\bullet M) = M \boxtimes_Y (N^\vee \otimes \omega_{W/Y}) := (\pi_X^\bullet M) \otimes_{\mathcal{O}_Z} \pi_W^\bullet (N^\vee \otimes \omega_{W/Y}) \in \mathrm{QCoh}(X \times_Y W) \,.$$

In keeping with the simplifying assumptions of the preceeding Remark, we introduce the category $\mathrm{Sch}_{\mathrm{ft},\mathrm{un}}^{\mathrm{fact}}(X)$ of finite type, schematic (unital) factorization spaces on $X$ following the general Definition 5.2.2:

*Example* 8.1.12. The functor $\mathrm{Sch}_{\mathrm{ft},\mathrm{corr}/(\cdot)} : \mathrm{Sch}_{\mathrm{aff}}^{\mathrm{op}} \to \mathrm{Cat}^\times$, which assigns to $S \in \mathrm{Sch}_{\mathrm{aff}}$ the category $\mathrm{Sch}_{\mathrm{ft},\mathrm{corr}/S}$ of schemes of finite type over $S$ under correspondences, is naturally a lax symmetric monoidal subfunctor of $\mathrm{PreStk}_{\mathrm{corr}/(\cdot)}$.

Let $\mathrm{Sch}_{\mathrm{ft},\mathrm{un}}^{\mathrm{fact}}(X)$ denote the category of finite type factorization schemes over $X$, defined by the subfunctor of the preceeding example following Definition 5.2.2.

*Example* 8.1.13. Let $\mathcal{X}, \mathcal{Y}$ and $\mathcal{W} \in \mathrm{Sch}_{\mathrm{ft},\mathrm{un}}^{\mathrm{fact}}(X)$, $M \in \mathrm{Alg}_{\mathrm{un}}^{\mathrm{fact}}(\mathrm{QCoh}_{\mathcal{X}})$, $N \in \mathrm{Alg}_{\mathrm{un}}^{\mathrm{fact}}(\mathrm{QCoh}_{\mathcal{W}})$ coherent, $f : \mathcal{X} \to \mathcal{Y}$, $g : \mathcal{W} \to \mathcal{Y}$ proper, and $\mathcal{Z} = \mathcal{X} \times_{\mathcal{Y}} \mathcal{W} \in \mathrm{Sch}_{\mathrm{ft},\mathrm{un}}^{\mathrm{fact}}(X)$ over $X = \tilde{X}_{\mathrm{dR}}$. Then following Remark 8.1.10, the internal Hom object over $\mathcal{Y}$ is given by

$$\tilde{\mathcal{Hom}}_{\mathrm{QCoh}(\mathcal{Y})}(g_\bullet \mathcal{N}, f_\bullet \mathcal{M}) \in \mathrm{Alg}_{\mathrm{un}}^{\mathrm{fact}}(\mathrm{QCoh}_{\mathcal{Y}}) \qquad \text{defined by} \qquad I \mapsto \underline{\mathrm{Hom}}_{\mathrm{QCoh}(\mathcal{Y}_I)}(g_{I \bullet} N_I, f_{I \bullet} M_I) \in \mathrm{QCoh}(\mathcal{Y}_I) \,.$$

Moreover, following Remark 8.1.11, there is a canonical lifted internal Hom object

$$\mathcal{H}\tilde{\mathrm{om}}_{\mathrm{QCoh}(\mathcal{Y})}(g_\bullet \mathcal{N}, f_\bullet \mathcal{M}) \in \mathrm{Alg}^{\mathrm{fact}}_{\mathrm{un}}(\mathrm{QCoh}_{\mathcal{Z}}) \qquad \text{defined by} \qquad I \mapsto \underline{\mathrm{Hom}}_{\mathrm{QCoh}(\mathcal{Y}_I)}(g_{I\bullet} N_I, f_{I\bullet} M_I) \in \mathrm{QCoh}(\mathcal{Z}_I) \,,$$

in the sense that we have a canonical equivalences

$$g_\bullet \pi_{\mathcal{W}_\bullet} \mathcal{H}\tilde{\mathrm{om}}_{\mathrm{QCoh}(\mathcal{Y})}(g_\bullet \mathcal{N}, f_\bullet \mathcal{M}) \cong \mathcal{H}\tilde{\mathrm{om}}_{\mathrm{QCoh}(\mathcal{Y})}(g_\bullet \mathcal{N}, f_\bullet \mathcal{M}) \in \mathrm{Alg}^{\mathrm{fact}}_{\mathrm{un}}(\mathrm{QCoh}_\mathcal{Y}) \qquad \text{, and}$$

$$\mathrm{p}_{\mathcal{Z}\bullet} \mathcal{H}\tilde{\mathrm{om}}_{\mathrm{QCoh}(\mathcal{Y})}(g_\bullet \mathcal{N}, f_\bullet \mathcal{M}) \cong \mathcal{H}\tilde{\mathrm{om}}_{\mathrm{QCoh}(\mathcal{Y})}(g_\bullet \mathcal{N}, f_\bullet \mathcal{M}) \in \mathrm{Alg}^{\mathrm{fact}}_{\mathrm{un}}(\tilde{X}) \qquad \text{.}$$

The quasicoherent, factorization analogue of the construction in Example A.6.6 is given by:

*Example* 8.1.14. Consider the special case of the preceeding example where $\mathcal{X} = \mathcal{W}$ and $M = N$. Then the lifted internal Hom object

$$\tilde{\mathcal{A}} = \mathcal{H}\tilde{\mathrm{om}}_{\mathrm{QCoh}(\mathcal{Y})}(f_\bullet \mathcal{M}, f_\bullet \mathcal{M}) = \mathcal{M} \boxtimes_\mathcal{Y} (\mathcal{M}^\vee \otimes \omega_{\mathcal{X}/\mathcal{Y}}) \in \mathrm{Alg}^{\mathrm{fact}}_{\mathbb{E}_1, \mathrm{un}}(\mathrm{QCoh}^\star_\mathcal{Z})$$

is naturally a factorization $\mathbb{E}_1$-algebra internal to $\mathrm{QCoh}^\star_\mathcal{Z} \in \mathrm{Cat}^{\mathrm{fact}}_{\mathbb{E}_1, \mathrm{un}}(X)$ as defined in Proposition 8.1.5.

Concretely, in terms of the description of Remark 8.1.7, the $\mathbb{E}_1$-algebra $A_I \in \mathrm{Alg}_{\mathbb{E}_1}(\mathrm{QCoh}(\mathcal{Z}_I)^\star)$ assigned to each $I \in \mathrm{fSet}_\varnothing$ is given by

$$A_I = \underline{\mathrm{Hom}}_{\mathrm{QCoh}(\mathcal{Y}_I)}(f_{I\bullet} M_I, f_{I\bullet} M_I) = M_I \boxtimes_{\mathcal{Y}_I} (M_I^\vee \otimes \omega_{\mathcal{X}_I/\mathcal{Y}_I}) \in \mathrm{QCoh}(\mathcal{Z}_I)$$

with algebra structure maps as in Equation 8.1.2 are given by

$$\pi_{13, I\bullet}(\pi^\bullet_{12, I} A_I \otimes \pi^\bullet_{23, I} A_I) \to A_I$$

corresponding under the $(\pi_{13\bullet}, \pi^!_{13})$ adjunction to

$$\begin{aligned}
\pi^\bullet_{12, I} A_I \otimes \pi^\bullet_{23, I} A_I &= (M_I \boxtimes_{\mathcal{Y}_I} (M_I^\vee \otimes \omega_{\mathcal{X}_I/\mathcal{Y}_I}) \boxtimes_{\mathcal{Y}_I} \mathcal{O}_{X_I}) \otimes (\mathcal{O}_{X_I} \boxtimes_{\mathcal{Y}_I} M_I \boxtimes_{\mathcal{Y}_I} (M_I^\vee \otimes \omega_{\mathcal{X}_I/\mathcal{Y}_I})) \\
&\cong M_I \boxtimes_{\mathcal{Y}_I} (M_I^\vee \otimes \omega_{\mathcal{X}_I/\mathcal{Y}_I} \otimes M_I) \boxtimes_{\mathcal{Y}_I} (M_I^\vee \otimes \omega_{\mathcal{X}_I/\mathcal{Y}_I}) \\
&\to M_I \boxtimes_{\mathcal{Y}_I} \omega_{\mathcal{X}_I/\mathcal{Y}_I} \boxtimes_{\mathcal{Y}_I} (M_I^\vee \otimes \omega_{\mathcal{X}_I/\mathcal{Y}_I}) \\
&\cong \pi^!_{13, I} A_I
\end{aligned}$$

where the map is given by the duality pairing on $\mathcal{M}_I \in \mathrm{Coh}(\mathcal{X}_I)$.

Further, the analogue of Example A.6.7 is given by:

*Example* 8.1.15. Consider the special case of the preceding example where the factorization algebras are given by $\mathcal{M} = \mathcal{N} = \mathcal{O}_\mathcal{X} := \mathrm{p}^\bullet_\mathcal{X} \mathcal{O}_{\mathrm{Ran}_{X, \mathrm{un}}} \in \mathrm{Alg}^{\mathrm{fact}}_{\mathrm{un}}(\mathrm{QCoh}_\mathcal{X})$ over $X = \tilde{X}_{\mathrm{dR}}$, so that

$$\mathcal{A} = \mathcal{H}\tilde{\mathrm{om}}(f_\bullet \mathcal{O}_\mathcal{X}, f_\bullet \mathcal{O}_\mathcal{X}) \cong \mathrm{p}_{\mathcal{X}\bullet} f^! f_\bullet \mathrm{p}^\bullet_\mathcal{X} \mathcal{O}_{\mathrm{Ran}_{X, \mathrm{un}}} \cong \mathrm{p}_{\mathcal{X}\bullet} \pi_{\mathcal{X}\bullet} \pi^!_\mathcal{X} \mathrm{p}^\bullet_\mathcal{X} \mathcal{O}_{\mathrm{Ran}_{X, \mathrm{un}}} = \mathrm{p}_{\mathcal{Z}\bullet} \omega_{\mathcal{Z}/\mathcal{X}} \in \mathrm{Alg}^{\mathrm{fact}}_{\mathrm{un}}(\tilde{X}).$$

Correspondingly, the lifted internal Hom object is given by

$$\tilde{\mathcal{A}} = \mathcal{H}\tilde{\mathrm{om}}(f_\bullet \mathcal{O}_\mathcal{X}, f_\bullet \mathcal{O}_\mathcal{X}) = \mathcal{O}_\mathcal{X} \boxtimes_\mathcal{Y} \omega_{\mathcal{X}/\mathcal{Y}} \cong \omega_{\mathcal{Z}/\mathcal{X}} \in \mathrm{Alg}^{\mathrm{fact}}_{\mathbb{E}_1, \mathrm{un}}(\mathrm{QCoh}^\star_\mathcal{Z}) \,.$$

Concretely, for each $I \in \mathrm{fSet}_\varnothing$, we have

$$\tilde{A}_I = \underline{\mathrm{Hom}}_{\mathrm{QCoh}(\mathcal{Y}_I)}(f_{I\bullet} \mathcal{O}_{\mathcal{X}_I}, f_{I\bullet} \mathcal{O}_{\mathcal{X}_I}) = \mathcal{O}_{\mathcal{X}_I} \boxtimes_{\mathcal{Y}_I} \omega_{\mathcal{X}_I/\mathcal{Y}_I} \cong \omega_{\mathcal{Z}_I/\mathcal{X}_I} \in \mathrm{Alg}_{\mathbb{E}_1}(\mathrm{QCoh}(\mathcal{Z}_I)^\star) \,.$$

This construction can be summarized by the diagram in factorization spaces, and induced diagram in factorization categories, given by

(8.1.4)

$$
\begin{array}{ccc}
& \mathcal{X} \times_{\mathcal{Y}} \mathcal{X} & \\
{\scriptstyle \pi_{\mathcal{X}}} \swarrow & & \searrow {\scriptstyle \pi_{\mathcal{X}}} \\
\mathcal{X} & & \mathcal{X} \\
{\scriptstyle p_{\mathcal{X}}} \swarrow \quad {\scriptstyle f} \searrow & \swarrow {\scriptstyle f} \quad \searrow {\scriptstyle p_{\mathcal{X}}} & \\
\mathrm{Ran}_{X,\mathrm{un}} \quad & \mathcal{Y} & \quad \mathrm{Ran}_{X,\mathrm{un}}
\end{array}
\qquad \text{and} \qquad
\begin{array}{ccc}
& \mathrm{QCoh}_{\mathcal{X} \times_{\mathcal{Y}} \mathcal{X}} & \\
{\scriptstyle \pi_{\mathcal{X}}^!} \nearrow & & \nwarrow {\scriptstyle \pi_{\mathcal{X}\bullet}} \\
\mathrm{QCoh}_{\mathcal{X}} & & \mathrm{QCoh}_{\mathcal{X}} \\
{\scriptstyle p_{\mathcal{X}}^\bullet} \nearrow \quad {\scriptstyle f_\bullet} \searrow & \nearrow {\scriptstyle f^!} \quad \nwarrow {\scriptstyle p_{\mathcal{X}\bullet}} & \\
\mathrm{QCoh}_{\mathrm{Ran}_{X,\mathrm{un}}} & \mathrm{QCoh}_{\mathcal{Y}} & \mathrm{QCoh}_{\mathrm{Ran}_{X,\mathrm{un}}}
\end{array}
.
$$

*Remark* 8.1.16. As in Remark A.6.8, defined much more generally. IndCoh analogue, de Rham analogue.

## 8.2. Convolution constructions of modules over factorization $\mathbb{E}_n$-algebras and categories.
In this subsection, we apply the format of Subappendix A.6.1 to outline a general construction of factorization module categories and factorization algebra module objects over the factorization $\mathbb{E}_1$ categories and factorization $\mathbb{E}_1$ algebra objects of the preceeding subsection.

*Example* 8.2.1. Let $\mathcal{X}, \mathcal{W}, \mathcal{Y} \in \mathrm{PreStk}_{\mathrm{un}}^{\mathrm{fact}}(X)$ be (unital) factorization spaces, and $f : \mathcal{X} \to \mathcal{Y}$ and $g : \mathcal{W} \to \mathcal{Y}$ maps of such. Following Remark 8.1.2, we have that $\mathcal{Z}^{\mathcal{W}} = \mathcal{X} \times_{\mathcal{Y}} \mathcal{W} \in \mathrm{PreStk}_{\mathrm{un}}^{\mathrm{fact}}(X)$ is a (unital) factorization space, specified as a relative prestack over $\mathrm{Ran}_{X,\mathrm{un}}$ by the assignments

$$
\begin{cases}
I & \mapsto \mathcal{Z}_I^{\mathcal{W}} := \mathcal{X}_I \times_{\mathcal{Y}_I} \mathcal{W}_I \\
& \qquad\qquad\qquad\qquad \mathrm{unit}_{\mathcal{X}}^{\pi} \underset{\mathrm{unit}_{\mathcal{Y}}^{\pi}}{\times} \mathrm{unit}_{\mathcal{W}}^{\pi} \\
[\pi : I \to J] & \mapsto \qquad\qquad {\scriptstyle \psi_{\mathcal{Z}^{\mathcal{W}}}^{\pi}} \swarrow \qquad\qquad\qquad \searrow {\scriptstyle \phi_{\mathcal{Z}^{\mathcal{W}}}^{\pi}} \\
& \quad X^J \times_{X^I} (\mathcal{X}_I \times_{\mathcal{Y}_I} \mathcal{W}_I) \qquad\qquad\qquad \mathcal{X}_J \times_{\mathcal{Y}_J} \mathcal{W}_J
\end{cases}
.
$$

Further, the canonical projection maps $\pi_{\mathcal{X}} : \mathcal{Z}^{\mathcal{W}} \to \mathcal{X}$ and $\pi_{\mathcal{W}} : \mathcal{Z}^{\mathcal{W}} \to \mathcal{W}$ are maps of (unital) factorization spaces.

*Remark* 8.2.2. More generally, the iterated fibre products

$$\mathcal{Z}_{(n)}^{\mathcal{W}} = \mathcal{X} \times_{\mathcal{Y}} \mathcal{X} \times_{\mathcal{Y}} \ldots \times_{\mathcal{Y}} \mathcal{W} \quad \in \mathrm{PreStk}_{\mathrm{un}}^{\mathrm{fact}}(X) \qquad \text{and projections} \qquad \pi_{ij} : \mathcal{Z}_{(n)}^{\mathcal{W}} \to \mathcal{Z}^{(\mathcal{W})}$$

define factorization spaces and maps of such for each $n \in \mathbb{N}$ and $i, j \in \{1, ..., n\}$.

*Remark* 8.2.3. Throughout the remainder of this subsection, let $\mathcal{X}, \mathcal{W}, \mathcal{Y} \in \mathrm{PreStk}_{\mathrm{un}}^{\mathrm{fact}}(X)$, $f : \mathcal{X} \to \mathcal{Y}$ and $g : \mathcal{W} \to \mathcal{Y}$, and $\mathcal{Z}^{\mathcal{W}} \in \mathrm{PreStk}_{\mathrm{un}}^{\mathrm{fact}}(X)$ be as in the preceeding example. Further, suppose that $\mathcal{X}$, $\mathcal{Y}$ and $f : \mathcal{X} \to \mathcal{Y}$ satisfy the hypotheses of Remark 8.1.4, and analogously that $\mathcal{W}$ and $g : \mathcal{W} \to \mathcal{Y}$ are such that $\mathrm{QCoh}_{\mathcal{Z}_{(n)}^{\mathcal{W}}} \in \mathrm{Cat}_{\mathrm{un}}^{\mathrm{fact}}(X)$ define (unital) factorization categories for each $n \in \mathbb{N}$, and the maps $\pi_{ij} : \mathcal{Z}_{(n)}^{\mathcal{W}} \to \mathcal{Z}^{(\mathcal{W})}$ induce (unital) factorization functors $\pi_{ij}^{\bullet} : \mathrm{QCoh}_{\mathcal{Z}^{(\mathcal{W})}} \to \mathrm{QCoh}_{\mathcal{Z}_{(n)}^{\mathcal{W}}}$ for each $i, j$. For example, we can take factorization spaces and maps of such as in examples 5.2.5 and 5.3.2, and Remark 5.3.3.

The following is the analogue of Proposition A.6.13 compatible with Proposition 8.1.5, recalling the latter is the analogue of Proposition A.6.9:

*Proposition* 8.2.4. The factorization category of quasicoherent sheaves on $\mathcal{Z}^{\mathcal{W}}$

$$\mathrm{QCoh}_{\mathcal{Z}^{\mathcal{W}}} \in \mathrm{QCoh}_{\mathcal{Z}}^{\star}\text{-Mod}(\mathrm{Cat}_{\mathrm{un}}^{\mathrm{fact}}(X))$$

is naturally a factorization module category over $\mathrm{QCoh}_{\mathcal{Z}} \in \mathrm{Cat}_{\mathbb{E}_1,\mathrm{un}}^{\mathrm{fact}}(X)$, with respect to the convolution module structure $(\cdot) \star (\cdot) : \mathrm{QCoh}_{\mathcal{Z}} \otimes^{\ast} \mathrm{QCoh}_{\mathcal{Z}^{\mathcal{W}}} \to \mathrm{QCoh}_{\mathcal{Z}^{\mathcal{W}}}$ defined by the composition

$$\mathrm{QCoh}_{\mathcal{Z}} \otimes^{\ast} \mathrm{QCoh}_{\mathcal{Z}^{\mathcal{W}}} \xrightarrow{\pi_{12}^{\ast}\boxtimes\pi_{23}^{\bullet}} \mathrm{QCoh}_{(\mathcal{Z}_{(3)}^{\mathcal{W}})^{\times 2}} \xrightarrow{\Delta^{\bullet}} \mathrm{QCoh}_{\mathcal{Z}_{(3)}^{\mathcal{W}}} \xrightarrow{\pi_{13,\bullet}} \mathrm{QCoh}_{\mathcal{Z}^{\mathcal{W}}} \ .$$

Further, the pushforward functors $\mathrm{p}_{\mathcal{Z}\bullet} : \mathrm{QCoh}_{\mathcal{Z}} \to \mathrm{QCoh}_{\mathrm{Ran}_{X,\mathrm{un}}}^{\otimes!}$ and $\mathrm{p}_{\mathcal{Z}^{\mathcal{W}}\bullet} : \mathrm{QCoh}_{\mathcal{Z}^{\mathcal{W}}} \to \mathrm{QCoh}_{\mathrm{Ran}_{X,\mathrm{un}}}$ define unital, lax compatible factorization functors, with respect to the above module structure. In particular, there exist natural maps of factorization quasicoherent sheaves

$$\mathrm{p}_{\mathcal{Z}\bullet}\mathcal{A} \otimes \mathrm{p}_{\mathcal{Z}\bullet}\mathcal{R} \to \mathrm{p}_{\mathcal{Z}\bullet}(\mathcal{A} \star \mathcal{R})$$

for each $\mathcal{A} \in \mathrm{Alg}_{\mathrm{un}}^{\mathrm{fact}}(\mathrm{QCoh}_{\mathcal{Z}})$ and $\mathcal{R} \in \mathrm{Alg}_{\mathrm{un}}^{\mathrm{fact}}(\mathrm{QCoh}_{\mathcal{Z}^{\mathcal{W}}})$.

Similarly, following propositions A.6.14 and A.6.12, compatible with Corollary 8.1.6 we have:

*Corollary* 8.2.5. Let $\mathcal{A} \in \mathrm{Alg}_{\mathbb{E}_1,\mathrm{un}}^{\mathrm{fact}}(\mathrm{QCoh}_{\mathcal{Z}}^{\star})$ and $\mathcal{R} \in A\text{-Mod}(\mathrm{Alg}_{\mathrm{un}}^{\mathrm{fact}}(\mathrm{QCoh}_{\mathcal{Z}^{\mathcal{W}}})$ a factorization $\mathbb{E}_1$ module object internal to $\mathrm{QCoh}_{\mathcal{Z}^{\mathcal{W}}}$ over $X = \tilde{X}_{\mathrm{dR}}$. Then $\mathrm{p}_{\mathcal{Z}^{\mathcal{W}}\bullet}\mathcal{R} \in (\mathrm{p}_{\mathcal{Z}\bullet}A)\text{-Mod}(\mathrm{Alg}_{\mathrm{un}}^{\mathrm{fact}}(\tilde{X}))$ is naturally a factorization $\mathbb{E}_1$ module over $\mathrm{p}_{\mathcal{Z}\bullet}\mathcal{A} \in \mathrm{Alg}_{\mathbb{E}_1,\mathrm{un}}^{\mathrm{fact}}(\tilde{X})$, so that we obtain a functor

$$\mathrm{p}_{\mathcal{Z}^{\mathcal{W}}\bullet} : \mathcal{A}\text{-Mod}(\mathrm{Alg}_{\mathrm{un}}^{\mathrm{fact}}(\mathrm{QCoh}_{\mathcal{Z}^{\mathcal{W}}})) \to (\mathrm{p}_{\mathcal{Z}\bullet}\mathcal{A})\text{-Mod}(\mathrm{Alg}_{\mathrm{un}}^{\mathrm{fact}}(\tilde{X})) \ .$$

*Remark* 8.2.6. Concretely, following Remark 7.2.16 and in terms of the description of Remark 8.1.7, a factorization $\mathbb{E}_1$ module object $\mathcal{R} \in A\text{-Mod}(\mathrm{Alg}_{\mathrm{un}}^{\mathrm{fact}}(\mathrm{QCoh}_{\mathcal{Z}^{\mathcal{W}}})$ over $A$ internal to $\mathrm{QCoh}_{\mathcal{Z}^{\mathcal{W}}}$ is given by:

- an $\mathcal{A}$ module object $\mathcal{R} \in \mathcal{A}\text{-Mod}(\mathrm{QCoh}(\mathcal{Z}^{\mathcal{W}}))$, given by an assignment

$$(8.2.1) \quad I \mapsto R_I \in A_I\text{-Mod}(\mathrm{QCoh}(\mathcal{Z}_I^{\mathcal{W}})) \qquad [\pi : I \to J] \mapsto \left[\eta_\pi : \phi_{\mathcal{Z}^{\mathcal{W}}\bullet}^{\pi}\psi_{\mathcal{Z}^{\mathcal{W}}}^{\pi,\bullet}\tilde{\Delta}(\pi)^{\bullet}R_I \to R_J\right] \ ,$$

  such that for $\pi$ surjective the corresponding maps in $A_J\text{-Mod}(\mathrm{QCoh}(\mathcal{Z}_J^{\mathcal{W}})$ are isomorphisms, where the former are given by an underlying object $R_I \in \mathrm{QCoh}(\mathcal{Z}_I^{\mathcal{W}})$ together with a map

$$(8.2.2) \quad \pi_{13,I\bullet}(\pi_{12,I}^{\bullet}A_I \otimes \pi_{23,I}^{\bullet}R_I) \to R_I \quad \text{in} \quad \mathrm{QCoh}(\mathcal{Z}_I^{\mathcal{W}})$$

  and its higher arity analogues, as in Proposition A.6.12, where $\pi_{ij,I} : \mathcal{Z}_{(3)I}^{\mathcal{W}} \to \mathcal{Z}_I^{(\mathcal{W})}$ are the components of the maps of factorization spaces in Remark 8.2.2 over $X^I$;
- an equivalence

$$j(\pi)^{\bullet}R_I \xrightarrow{\cong} j(\pi)^{\bullet}(\boxtimes_j R_{I_j}) \qquad \text{in} \qquad (j(\pi)^{\bullet}A_I)\text{-Mod}(\mathrm{QCoh}(U(\pi) \times_{X^I} \mathcal{Z}_I^{\mathcal{W}})))$$

$$\cong (j(\pi)^{\bullet}(\boxtimes_j A_{I_j}))\text{-Mod}(\mathrm{QCoh}(U(\pi) \times_{X^I} (\times_{j \in J}\mathcal{Z}_{I_j}^{\mathcal{W}})))$$

  for each $I$ and $\pi : I \to J$; and
- an equivalence $R_{\emptyset} \cong \mathbb{K}$ in $A_{\emptyset}\text{-Mod}(\mathrm{QCoh}_{\mathcal{Z}_{\emptyset}^{\mathcal{W}}}) \cong \mathrm{Vect}$, which together with the structure maps of Equation 8.2.1 above, determine maps

$$\mathrm{unit}_C^{\pi}(R_I\boxtimes\mathcal{O}_{X^{I_\pi}}) \to R_J \qquad \text{and} \qquad \mathrm{unit}_C^{J}(\mathcal{O}_{X^J}) \to R_J$$

  for $\pi : I \hookrightarrow J$ injective, and in particular for $\pi : \emptyset \hookrightarrow J$, respectively.

*Remark* 8.2.7. In particular, following Remark 8.1.8 and the preceding remark, the preceeding corollary implies the factorization analogue of Proposition A.6.12, analogously to the fact that *loc. cit.* follows from Proposition A.6.14.

Following Example 8.1.14, the quasicoherent, factorization analogue of Example A.6.16 is given by:

*Example* 8.2.8. Let $\mathcal{X}, \mathcal{Y}$ and $\mathcal{W} \in \mathrm{Sch}^{\mathrm{fact}}_{\mathrm{ft,un}}(X)$, $M \in \mathrm{Alg}^{\mathrm{fact}}_{\mathrm{un}}(\mathrm{QCoh}_{\mathcal{X}})$, $N \in \mathrm{Alg}^{\mathrm{fact}}_{\mathrm{un}}(\mathrm{QCoh}_{\mathcal{W}})$ coherent, $f : \mathcal{X} \to \mathcal{Y}$, $g : \mathcal{W} \to \mathcal{Y}$ proper, and $\mathcal{Z} = \mathcal{X} \times_{\mathcal{Y}} \mathcal{X}$ and $\mathcal{Z}^{\mathcal{W}} = \mathcal{X} \times_{\mathcal{Y}} \mathcal{W} \in \mathrm{Sch}^{\mathrm{fact}}_{\mathrm{ft,un}}(X)$ over $X = \tilde{X}_{\mathrm{dR}}$. Further, as in examples 8.1.14 and 8.1.13, consider

$$\tilde{\mathcal{A}} = \tilde{\mathcal{H}\mathrm{om}}_{\mathrm{QCoh}(\mathcal{Y})}(f_{\bullet}\mathcal{M}, f_{\bullet}\mathcal{M}) = \mathcal{M} \boxtimes_{\mathcal{Y}} (\mathcal{M}^{\vee} \otimes \omega_{\mathcal{X}/\mathcal{Y}}) \qquad \in \mathrm{Alg}^{\mathrm{fact}}_{\mathbb{E}_1,\mathrm{un}}(\mathrm{QCoh}^{\star}_{\mathcal{Z}}) \text{ , and}$$

$$\tilde{\mathcal{R}} = \tilde{\mathcal{H}\mathrm{om}}_{\mathrm{QCoh}(\mathcal{Y})}(g_{\bullet}\mathcal{N}, f_{\bullet}\mathcal{M}) = \mathcal{M} \boxtimes_{\mathcal{Y}} (\mathcal{N}^{\vee} \otimes \omega_{\mathcal{W}/\mathcal{Y}}) \qquad \in \mathrm{Alg}^{\mathrm{fact}}_{\mathrm{un}}(\mathrm{QCoh}_{\mathcal{Z}^{\mathcal{W}}}).$$

Then $\tilde{\mathcal{R}}$ is naturally a factorization $\mathbb{E}_1$ module object

$$\tilde{\mathcal{R}} \in \tilde{\mathcal{A}}\text{-}\mathrm{Mod}(\mathrm{Alg}^{\mathrm{fact}}_{\mathrm{un}}(\mathrm{QCoh}_{\mathcal{Z}^{\mathcal{W}}}))$$

over $\tilde{\mathcal{A}} \in \mathrm{Alg}^{\mathrm{fact}}_{\mathbb{E}_1,\mathrm{un}}(\mathrm{QCoh}^{\star}_{\mathcal{Z}})$ internal to $\mathrm{QCoh}_{\mathcal{Z}^{\mathcal{W}}}\text{-}\mathrm{Mod}(\mathrm{Cat}^{\mathrm{fact}}_{\mathrm{un}}(X))$.

Concretely, for each $I \in \mathrm{fSet}_{\varnothing}$ we have

$$\tilde{R}_I = \tilde{\mathcal{H}\mathrm{om}}_{\mathrm{QCoh}(\mathcal{Y}_I)}(g_{I,\bullet}N_I, f_{I,\bullet}M_I) = M_I \boxtimes_{\mathcal{Y}_I} (N_I^{\vee} \otimes \omega_{\mathcal{W}_I/\mathcal{Y}_I}) \ \in A_I\text{-}\mathrm{Mod}(\mathrm{QCoh}(\mathcal{Z}^{\mathcal{W}}_I))$$

with module structure maps as in Equation 8.2.2 given by

$$\pi_{13,I\bullet}(\pi^{\bullet}_{12,I}A_I \otimes \pi^{\bullet}_{23,I}R_I) \to R_I$$

corresponding under the $(\pi_{13\bullet}, \pi^!_{13})$ adjunction to

$$\pi^{\bullet}_{12,I}A_I \otimes \pi^{\bullet}_{23,I}R_I = (M_I \boxtimes_{\mathcal{Y}_I} (M_I^{\vee} \otimes \omega_{\mathcal{X}_I/\mathcal{Y}_I}) \boxtimes_{\mathcal{Y}_I} \mathcal{O}_{X_I}) \otimes (\mathcal{O}_{X_I} \boxtimes_{\mathcal{Y}_I} M_I \boxtimes_{\mathcal{Y}_I} (N_I^{\vee} \otimes \omega_{\mathcal{X}_I/\mathcal{Y}_I}))$$
$$\cong M_I \boxtimes_{\mathcal{Y}_I} (M_I^{\vee} \otimes \omega_{\mathcal{X}_I/\mathcal{Y}_I} \otimes M_I) \boxtimes_{\mathcal{Y}_I} (N_I^{\vee} \otimes \omega_{\mathcal{X}_I/\mathcal{Y}_I})$$
$$\to M_I \boxtimes_{\mathcal{Y}_I} \omega_{\mathcal{X}_I/\mathcal{Y}_I} \boxtimes_{\mathcal{Y}_I} (N_I^{\vee} \otimes \omega_{\mathcal{X}_I/\mathcal{Y}_I})$$
$$\cong \pi^!_{13,I}R_I$$

where the map is given by the duality pairing on $M_I \in \mathrm{Coh}(\mathcal{X}_I)$.

Following 8.1.15, the quasicoherent, factorization analogue of Example A.6.17 is given by:

*Example* 8.2.9. Consider the special case of the preceding example where the factorization algebras are given by $\mathcal{M} = \mathcal{O}_{\mathcal{X}} := \mathrm{p}^{\bullet}_{\mathcal{X}}\mathcal{O}_{\mathrm{Ran}_{X,\mathrm{un}}} \in \mathrm{Alg}^{\mathrm{fact}}_{\mathrm{un}}(\mathrm{QCoh}_{\mathcal{X}})$ and $\mathcal{N} = \mathcal{O}_{\mathcal{W}} := \mathrm{p}^{\bullet}_{\mathcal{W}}\mathcal{O}_{\mathrm{Ran}_{X,\mathrm{un}}} \in \mathrm{Alg}^{\mathrm{fact}}_{\mathrm{un}}(\mathrm{QCoh}_{\mathcal{W}})$ over $X = \tilde{X}_{\mathrm{dR}}$, so that

$$\mathcal{R} = \mathcal{H}\mathrm{om}(g_{\bullet}\mathcal{O}_{\mathcal{W}}, f_{\bullet}\mathcal{O}_{\mathcal{X}}) \cong \mathrm{p}_{\mathcal{W}\bullet}g^! f_{\bullet}\mathrm{p}^{\bullet}_{\mathcal{X}}\mathcal{O}_{\mathrm{Ran}_{X,\mathrm{un}}} \cong \mathrm{p}_{\mathcal{W}\bullet}\pi_{\mathcal{W}\bullet}\pi^!_{\mathcal{X}}\mathrm{p}^{\bullet}_{\mathcal{X}}\mathcal{O}_{\mathrm{Ran}_{X,\mathrm{un}}} = \mathrm{p}_{\mathcal{Z}^{\mathcal{W}}\bullet}\omega_{\mathcal{Z}^{\mathcal{W}}/\mathcal{W}} \in \mathcal{A}\text{-}\mathrm{Mod}(\mathrm{Alg}^{\mathrm{fact}}_{\mathrm{un}}(\tilde{X})).$$

Correspondingly, the lifted internal Hom object is given by

$$\tilde{\mathcal{R}} = \tilde{\mathcal{H}\mathrm{om}}(g_{\bullet}\mathcal{O}_{\mathcal{W}}, f_{\bullet}\mathcal{O}_{\mathcal{X}}) = \mathcal{O}_{\mathcal{X}} \boxtimes_{\mathcal{Y}} \omega_{\mathcal{W}/\mathcal{Y}} \cong \omega_{\mathcal{Z}^{\mathcal{W}}/\mathcal{X}} \in \tilde{\mathcal{A}}\text{-}\mathrm{Mod}(\mathrm{Alg}^{\mathrm{fact}}_{\mathrm{un}}(\mathrm{QCoh}_{\mathcal{Z}^{\mathcal{W}}})) \text{ .}$$

Concretely, for each $I \in \mathrm{fSet}_{\varnothing}$, we have

$$\tilde{R}_I = \underline{\mathrm{Hom}}_{\mathrm{QCoh}(\mathcal{Y}_I)}(g_{I\bullet}\mathcal{O}_{\mathcal{W}_I}, f_{I\bullet}\mathcal{O}_{\mathcal{X}_I}) = \mathcal{O}_{\mathcal{X}_I} \boxtimes_{\mathcal{Y}_I} \omega_{\mathcal{W}_I/\mathcal{Y}_I} \cong \omega_{\mathcal{Z}^{\mathcal{W}}_I/\mathcal{X}_I} \ \in A_I\text{-}\mathrm{Mod}(\mathrm{QCoh}(\mathcal{Z}_I)^{\star}) \text{ .}$$

This construction can be summarized by the diagram in factorization schemes, and induced diagram in factorization categories, given by
(8.2.3)

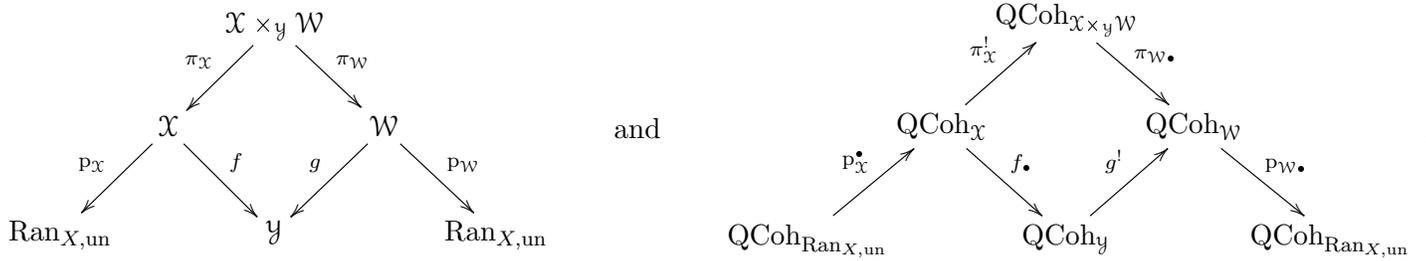

*Remark* 8.2.10. As in Remark 8.1.16, defined much more generally. IndCoh analogue, de Rham analogue.

# Chapter 2

# Holomorphic-topological field theory in the $\Omega$-background and equivariant factorization $\mathbb{E}_n$ algebras

## 9. Overview of Chapter 2

### 9.1. General overview.
In this Chapter, we apply the tools developed in Chapter 1 to construct examples of factorization $\mathbb{E}_n$ algebras $A \in \mathrm{Alg}^{\mathrm{fact}}_{\mathbb{E}_n, \mathrm{un}}(X)$ describing mixed chiral-topological field theories on $X \times \mathbb{R}^n$ coming from holomorphic-topological twists of supersymmetric gauge theories. Moreover, we apply the theory of equivariant factorization algebras from Chapter I-2 (which we recall refers to Chapter 2 of the prequel [But20a]) in these examples to formulate and prove statements in representation theory which codify predictions from physics about the corresponding field theories, as outlined in the introduction. Recall that these results were introduced and summarized already in Section 1.4.2.

The results of Section 8 provide the main technical mechanism for constructing examples. Recall the diagrams of Equation 8.1.4, which summarize the geometric construction of factorization $\mathbb{E}_1$ algebras given in Example 8.1.15. In each of the following sections, we construct factorization $\mathbb{E}_n$ algebras of local observables in chiral-topological field theories on $X \times \mathbb{R}^n$ using variants of this format, and deduce properties of these objects and relationships between them in terms of such geometric constructions.

In the sections dedicated to introducing individual examples, we proceed by defining the following algebraic data, which is typically constructed geometrically as listed in each case:

- a unital factorization $\mathbb{E}_{n-1}$ category $C \in \mathrm{Cat}^{\mathrm{fact}}_{\mathbb{E}_{n-1}, \mathrm{un}}(X)$ on $X$, which corresponds physically to the category of line operators of the field theory, typically constructed geometrically by some variant of $\mathrm{QCoh}_{\mathcal{Y}} \in \mathrm{Cat}^{\mathrm{fact}}_{\mathbb{E}_{n-1}, \mathrm{un}}(X)$ for an appropriate factorization space $\mathcal{Y}$ with convolution structure maps inducing a compatible $\mathbb{E}_{n-1}$ monoidal structure;
- a factorization $\mathbb{E}_n$ algebra on $X$, corresponding to the algebra of local observables, given by

$$\mathcal{A} = \mathcal{H}om_{\mathrm{QCoh}(\mathcal{Y})}(\mathrm{unit}_{\mathrm{QCoh}_{\mathcal{Y}}}, \mathrm{unit}_{\mathrm{QCoh}_{\mathcal{Y}}}) \in \mathrm{Alg}^{\mathrm{fact}}_{\mathbb{E}_n, \mathrm{un}}(X) \ ,$$

  the endomorphism algebra of the unit factorization $\mathbb{E}_{n-1}$ algebra $\mathrm{unit}_{\mathrm{QCoh}_{\mathcal{Y}}} \in \mathrm{Alg}^{\mathrm{fact}}_{\mathbb{E}_{n-1}, \mathrm{un}}(\mathrm{QCoh}_{\mathcal{Y}})$, typically constructed geometrically by the analogous variant of $\mathrm{unit}_{\mathrm{QCoh}_{\mathcal{Y}}} = f_\bullet \mathrm{p}_X^\bullet \mathcal{O}_{\mathrm{Ran}_{X, \mathrm{un}}}$;
- a unital factorization $\mathbb{E}_n$ category $\tilde{C} \in \mathrm{Cat}^{\mathrm{fact}}_{\mathbb{E}_n, \mathrm{un}}(X)$ on $X$ together with a factorization functor $\varphi : \tilde{C} \to \mathrm{QCoh}_{\mathrm{Ran}_{X, \mathrm{un}}}$, typically given geometrically by the analogous variant of $\mathrm{p}_\bullet : \mathrm{QCoh}_{\mathcal{X} \times_{\mathcal{Y}} \mathcal{X}} \to \mathrm{QCoh}_{\mathrm{Ran}_{X, \mathrm{un}}}$ the coherent pushforward functor; and
- a lift $\tilde{\mathcal{A}} \in \mathrm{Alg}^{\mathrm{fact}}_{\mathbb{E}_n, \mathrm{un}}(\tilde{C})$ of $\mathcal{A}$ to a factorization $\mathbb{E}_n$ algebra internal to $\tilde{C}$ such that $\varphi(\tilde{\mathcal{A}}) \cong \mathcal{A}$, defining the internal to $\tilde{C}$ endomorphism object of $\mathrm{unit}_{\mathrm{QCoh}(\mathcal{Y})}$, which exists in the typical geometric setting by the analogous variant of coherent base change.

In the remaining sections, we prove a number of results relating these factorization $\mathbb{E}_n$ algebras and categories in terms of the geometric constructions explained above. These can mostly be understood as examples of the following general types of results:

- concrete computations of the factorization $\mathbb{E}_n$ algebras constructed geometrically as above, in terms of sheaf theory on the relevant varieties or stacks;
- analogous geometric constructions of module objects over these factorization $\mathbb{E}_n$ algebras, following the results of Section A.6.1;
- analogous concrete computations of these module objects in terms of sheaf theory;
- conditions for the existence of equivariant structures on factorization $\mathbb{E}_n$ algebras;
- analogous computations of the families of factorization algebras on fixed point subvarieties induced by equivariant factorization algebras, following the results of Section I-18; and
- applications of the equivariant cigar reduction principle to deduce relationships between various factorization $\mathbb{E}_n$ algebras, following the results of sections I-25.2 and I-25.4.

*Warning* 9.1.1. Many of the constructions in this chapter require the input data of a space $Y$, the cotangent stack of which defines the 'sigma-model target space' of the relevant field theory. For each example, we begin with an overview describing the general features expected from the construction, without formulating precise enough hypotheses on $Y \in \mathrm{PreStk}$ to give complete proofs of the claims. We still write the expected statements as 'Propositions', 'Examples' and 'Theorems', including similar warnings to this one throughout; we hope that this will not be a source of confusion for the reader.

In the examples relevant for our main results, we will restrict to the case that $Y$ is a specific class of space, typically either a smooth, finite type, affine scheme or a global quotient stack $Y = N/G$ of a finite-type, linear $G$ representation $N$ by a reductive, affine algebraic group $G$. At this level of generality, we give careful accounts of the relevant constructions, and complete the proofs of the main results outlined in the overview sections.

## 9.2. Summary. We now summarize the results of each section in the remainder of this chapter.

9.2.1. *The two dimensional B model.* In Section 10, as a warm-up we recall the well-known construction of the two dimensional B model to $Y$ as an $\mathbb{E}_2$ algebra $\mathcal{V}(Y) \in \mathrm{Alg}_{\mathbb{E}_2}(\mathrm{Vect})$. We emphasize the role of the $\mathbb{E}_1$ monoidal category of line operators $\mathrm{QCoh}(Y^{\times 2}) \in \mathrm{Alg}_{\mathbb{E}_1}(\mathrm{DGCat})$ in this construction, which plays the role of the factorization category C described above in this example. Further, we recall the relationship between $S^1$ equivariant structures on the two dimensional B model $\mathbb{E}_2$ algebra and Calabi-Yau structures on $Y$.

9.2.2. *The three dimensional B model.* In Section 11, we recall an analogous construction of the three dimensional B model to $Y$ as an $\mathbb{E}_3$ algebra $\mathcal{B}(Y) \in \mathrm{Alg}_{\mathbb{E}_3}(\mathrm{Vect})$, following [BZ14] for example. In this case, the category of line operators is given by the $\mathbb{E}_2$ monoidal category $\mathrm{QCoh}(\mathcal{L}Y) \in \mathrm{Alg}_{\mathbb{E}_2}(\mathrm{DGCat})$ of coherent sheaves on the topological loop space $\mathcal{L}Y = Y \times_{Y \times Y} Y$ of $Y$, which plays the role of the factorization category in this example. Further, we explain that there is a canonical $S^1$ equivariant structure on $\mathcal{B}(Y)$ such that the $\mathbb{BD}_1^u$ algebra corresponding to it under Proposition I-25.1.1 is given by the (two-periodic) Rees algebra of differential operators on $Y$, following [BZ14, BZN12]. Finally, we explain the result of the equivariant cigar reduction principle of section I-25.2 in this example, as a warm up for its application in Section 19.

9.2.3. *The three dimensional A model.* In Section 12, we give an overview of the construction of the three dimensional A model to $Y$ as a factorization $\mathbb{E}_1$ algebra $\mathcal{A}(Y) \in \mathrm{Alg}_{\mathbb{E}_1, \mathrm{un}}^{\mathrm{fact}}(X)$ on a smooth curve $X$, following [BFN18] and references therein. In this case, the category of line operators is given by the genuine factorization category $D(Y_{\mathcal{K}}) \in \mathrm{Cat}_{\mathrm{un}}^{\mathrm{fact}}(X)$. We recall several expected results including an internal lift of the construction and the existence of a $\mathbb{G}_a \rtimes \mathbb{G}_m$ equivariant structure on $\mathcal{A}(Y)$, inducing a $\mathbb{BD}_1^u$ algebra, as in the preceding section for the three dimensional B model.

**9.2.4.** *The Braverman-Finkelberg-Nakajima construction.* In Section 13, we recall the complete construction of the three dimensional A model $\mathcal{A}(G, N) \in \mathrm{Alg}^{\mathrm{fact}}_{\mathbb{E}_1, \mathrm{un}}(X)$ in the case $Y = N/G$, carried out in [BFN18]. In this case, the constructions and results of the preceding section can be formulated carefully in terms of sheaf theory on infinite dimensional varieties and stacks; we give an account of this construction reconciling that in [BFN18] with the theory of $D$ modules on infinite dimensional varieties constructed in [Ras15b] and references therein.

**9.2.5.** *Chiral differential operators and the three dimensional A model.* In Section 14, we give an overview of the geometric construction of the factorization algebra over $X$ of chiral differential operators $\mathcal{D}^{\mathrm{ch}}(Y) \in \mathrm{Alg}^{\mathrm{fact}}_{\mathrm{un}}(X)$ on $Y$, and outline a geometric construction of an action of the factorization $\mathbb{E}_1$ algebra $\mathcal{A}(Y) \in \mathrm{Alg}^{\mathrm{fact}}_{\mathbb{E}_1, \mathrm{un}}(X)$ of local observables of the three dimensional $A$ model on the algebra of chiral differential operators $\mathcal{D}^{\mathrm{ch}}(Y) \in \mathcal{A}(Y)\text{-}\mathrm{Mod}(\mathrm{Alg}^{\mathrm{fact}}_{\mathrm{un}}(X))$. This module structure was predicted in [CG18], motivated by the fact that three dimensional $A$ model to $Y$ is expected to admit a chiral boundary condition such that the boundary local observables are given by the chiral differential operators to $Y$.

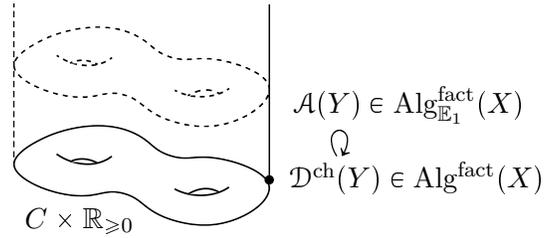

**9.2.6.** *Chiral differential operators on schemes.* In Section 15, we give a careful account of the geometric construction of the algebra of chiral differential operators $\mathcal{D}^{\mathrm{ch}}(Y) \in \mathrm{Alg}^{\mathrm{fact}}_{\mathrm{un}}(X)$ in the case that $Y$ is a well-behaved smooth, affine scheme, following [BD04] and [KV06]. In particular, we reconcile the constructions of *loc. cit.* with the theories of indcoherent sheaves in infinite type and their relations with theories of $D$ modules, as outlined in [Ras20b] following [Gai13],[Ras15b] and references therein. In particular, in these terms, we explain the requirement of a trivialization of the determinant gerbe of $Y$ for the construction of the global sections functor necessary to define $\mathcal{D}^{\mathrm{ch}}(Y)$ following [KV06].

**9.2.7.** *Chiral differential operators on quotient stacks and semi-infinite cohomology.* In Section 16, we study the analogous geometric construction of chiral differential operators in the case that $Y = N/G$. We identify the analogue of the requisite trivialization of the determinant gerbe of $Y$ with the requirement of a lift of the $G_{\mathcal{O}}$ action on $\mathcal{D}^{\mathrm{ch}}(N)$ to an action of $\hat{\mathfrak{g}}$ at level $-$Tate, and given a choice of such we identify the resulting chiral algebra $\mathcal{D}^{\mathrm{ch}}(N/G) \cong \mathrm{C}^{\frac{\infty}{2}}(\hat{\mathfrak{g}}, \mathfrak{g}_{\mathcal{O}}, G_{\mathcal{O}}; \mathcal{D}^{\mathrm{ch}}(N))$ with the semi-infinite cohomology of $\mathcal{D}^{\mathrm{ch}}(N)$ with respect to $\hat{\mathfrak{g}}$, following the results of [Ras20b].

**9.2.8.** *Action of the three dimensional A model on chiral differential operators.* In Section 17, we show that under the hypotheses identified in the preceding section, the chiral differential operators admit a canonical module structure $\mathcal{D}^{\mathrm{ch}}(N/G) \in \mathcal{A}(G, N)\text{-}\mathrm{Mod}(\mathrm{Alg}^{\mathrm{fact}}_{\mathrm{un}}(X))$ over the Coulomb branch factorization $\mathbb{E}_1$ algebra $\mathcal{A}(G, N) \in \mathrm{Alg}^{\mathrm{fact}}_{\mathbb{E}_1, \mathrm{un}}(X)$ the factorization $\mathbb{E}_1$ algebra constructed in [BFN18] in the case $Y = N/G$, as outlined in Section 14.

**9.2.9.** *The three dimensional holomorphic-B model and its deformations.* In Section 18, we outline a new construction of a factorization $\mathbb{E}_1$ algebra $\mathcal{C}(Y) \in \mathrm{Alg}^{\mathrm{fact}}_{\mathbb{E}_1, \mathrm{un}}(X)$ corresponding to the three dimensional mixed holomorphic-B model. Moreover, we show that this factorization $\mathbb{E}_1$ algebra admits a deformation $\mathcal{C}(Y)^{\hbar} \in \mathrm{Alg}^{\mathrm{fact}}_{\mathbb{E}_1, \mathrm{un}}(X)_{/\mathbb{K}[\hbar]}$ over $\mathbb{A}^1_{\hbar}$ such that the generic fibre $\mathcal{C}(Y)^{\hbar}|_{\{\hbar=1\}} \cong \mathcal{A}(Y)$ is equivalent to the three dimensional $A$ model; this realizes the prediction that the holomorphic-B twist of three dimensional $\mathcal{N} = 4$ gauge theory admits a deformation to the A twist.

Further, we show that the action of $\mathcal{A}(Y)$ on $\mathcal{D}^{\mathrm{ch}}(Y)$ extends to an action on the family of factorization algebras $\mathcal{D}^{\mathrm{ch}}(Y)_\hbar \in \mathrm{Alg}^{\mathrm{fact}}_{\mathrm{un}}(X)_{/\mathbb{K}[\hbar]}$ given by the Rees algebra of chiral differential operators.

9.2.10. *The four dimensional holmorphic-B model and chiral quantization.* In Section 19, we outline a new construction of a factorization $\mathbb{E}_2$ algebra $\mathcal{F}(Y) \in \mathrm{Alg}^{\mathrm{fact}}_{\mathbb{E}_2,\mathrm{un}}(X)$ corresponding to the four dimensional mixed holomorphic-B model. Moreover, we sketch an arguement that $\mathcal{F}(Y)$ admits an $S^1$ equivariant structure in the topological direction if and only if $Y$ admits a Tate structure, and that in this case the factorization $\mathbb{BD}_0^u$ algebra $\mathcal{F}(Y)_u \in \mathrm{Alg}^{\mathrm{fact}}_{\mathbb{BD}_0^u,\mathrm{un}}(X)$ corresponding to $\mathcal{F}(Y)$ under the equivalence of Proposition I-25.3.1 is given by the (two-periodic) Rees algebra of chiral differential operators $\mathcal{D}^{\mathrm{ch}}(Y)_u$. This gives a mathematical formulation of a variant of the main physical construction of [BLL+15], which associates a chiral algebra to each four dimensional $\mathcal{N} = 2$ superconformal gauge theory.

Further, we outline an identification $\mathrm{CC}^-_\bullet(\mathcal{F}(Y)) \cong \mathcal{C}(Y)^u \in \mathrm{Alg}^{\mathrm{fact}}_{\mathbb{E}_1,\mathrm{un}}(X)_{/\mathbb{K}[u]}$ of the negative cyclic chains on $\mathcal{F}(Y)$ with the (two-periodic variant of) the deformation of the holomorphic-B to A twist, such that the canonical action of $\mathrm{CC}^-_\bullet(\mathcal{F}(Y))$ on $\mathcal{F}(Y)_u$ induced by the equivariant cigar reduction principle of Proposition I-25.4.3 identifies with (the two-periodic variant of) the module structure constructed in the preceding section. Figure 3 below summarizes these statements, and should be compared with the general case of the equivariant cigar reduction principle illustrated in Figure I-4. This explains the relationship between our account of the predictions of [BLL+15] in the preceding paragraph and our formulation of the predictions of [CG18] summarized above.

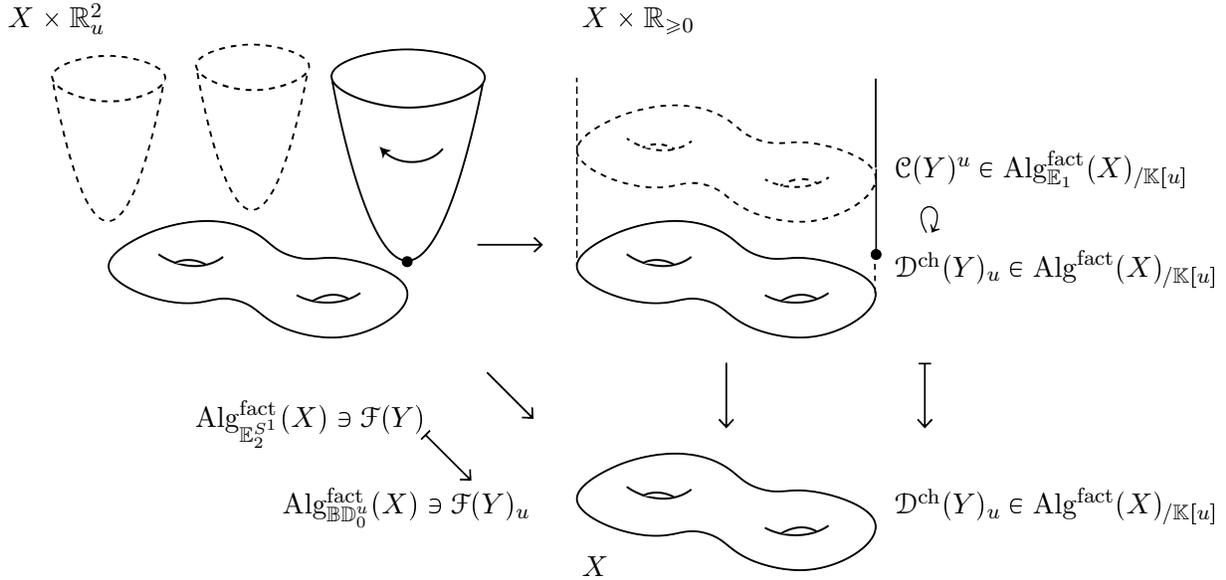

FIGURE 3. The equivariant cigar reduction principle for the four dimensional holomorphic-B model

## 10. The two dimensional B model

In this section, we recall the well-known construction of the two dimensional B model with target space $Y$ from the general perspective outlined in Section 9.1. The two dimensional B model is a topological field theory and thus described by an $\mathbb{E}_2$ algebra $\mathcal{V}(Y) \in \mathrm{Alg}_{\mathbb{E}_2}(\mathrm{Vect})$, which is computed in terms of the category of line operators $\mathrm{QCoh}(Y \times Y) \in \mathrm{Alg}_{\mathbb{E}_1}(\mathrm{DGCat})$, in keeping with the general format described in *loc. cit.*, by

$$\mathcal{V}(Y) = \mathrm{Hom}_{\mathrm{QCoh}(Y \times Y)}(\Delta_\bullet \mathcal{O}_Y, \Delta_\bullet \mathcal{O}_Y) \cong \Gamma(Y, \mathrm{PV}_Y^\bullet) \ .$$

This is precisely its definition as the Hochschild cochains $\mathrm{CH}_\bullet$ of the boundary condition category $\mathrm{QCoh}(Y)$, as suggested in Kontsevich's initial proposal of Homological Mirror Symmetry [Kon95]. As we explain below, the monoid structure on $Y \times Y$ is the analogue of the factorization space structure in this example, and the construction is summarized as in Equation 8.1.4 by the diagrams (10.0.1)

10.0.1. *Summary.* In Section 10.1, we construct the two dimensional B model $\mathbb{E}_2$ algebra, and in Section 10.2 we give the internal variant of the construction. In Section 10.3, we discuss the existence of $S^1$ equivariant structures on the two dimensional B model.

*Warning* 10.0.1. In keeping with Warning 9.1.1, we do not formulate specific hypotheses on the space $Y$ used in this section, so that the results stated throughout are only an outline of the general expectations. Since this section is primarily motivational, we will not give a more careful treatment, but the concerned reader may assume $Y$ is a smooth, finite type variety, for example.

### 10.1. The two dimensional B model factorization algebra.

*Example* 10.1.1. Recall that the space $Y \times Y$ is a monoid in prestacks under correspondences, with multiplication and unit structure maps given by the correspondences

*Remark* 10.1.2. The monoid structure on $Y \times Y$ is used here analogously to the factorization structures in the 'higher dimensional' examples of the following sections. In particular, we present this section following the exposition of sections 5 and 8 under this analogy, as a warm up for the more involved examples that follow:

We construct an $\mathbb{E}_1$-algebra in the category of $\mathbb{E}_1$-algebras (or equivalently, an $\mathbb{E}_2$ algebra) as endomorphisms of an object in an $\mathbb{E}_1$-monoidal category of sheaves on a monoid in spaces, in analogy with the construction of an $\mathbb{E}_1$-algebra in the category of factorization algebras (or equivalently, a factorization $\mathbb{E}_1$-algebra) as endomorphisms of an object in a factorization category of sheaves on a factorization space, as explained in *loc. cit.*, and carried out concretely in the sections following this one.

*Proposition* 10.1.3. The category $\mathrm{QCoh}(Y^{\times 2})^\star \in \mathrm{Alg}_{\mathbb{E}_1}(\mathrm{DGCat}_{\mathrm{cont}})$ is naturally an $\mathbb{E}_1$-monoidal category with respect to the convolution tensor product $(\cdot) \star (\cdot) : \mathrm{QCoh}(Y^{\times 2})^{\otimes 2} \to \mathrm{QCoh}(Y^{\times 2})$ defined by the composition

$$\mathrm{QCoh}(Y^{\times 2}) \otimes \mathrm{QCoh}(Y^{\times 2}) \xrightarrow{\pi_{12}^\bullet \boxtimes \pi_{23}^\bullet} \mathrm{QCoh}(Y^{\times 3})^{\otimes 2} \xrightarrow{\Delta^\bullet} \mathrm{QCoh}(Y^{\times 3}) \xrightarrow{\pi_{13\bullet}} \mathrm{QCoh}(Y^{\times 2}) \ .$$

*Example* 10.1.4. The unit $\mathrm{u} \in \mathrm{QCoh}(Y^{\times 2})$ for the convolution monoidal tensor product is given by

$$\mathrm{unit}_{Y^{\times 2}} = \Delta_\bullet \mathrm{p}_Y^\bullet \mathcal{O}_{\mathrm{pt}} \cong \Delta_\bullet \mathcal{O}_Y \ \in \mathrm{QCoh}(Y^{\times 2}) \ .$$

In particular, $\mathrm{u} \in \mathrm{Alg}_{\mathbb{E}_1}(\mathrm{QCoh}(Y^{\times 2})^\star)$ canonically lifts to an $\mathbb{E}_1$-algebra object in $\mathrm{QCoh}(Y^{\times 2})^\star$.

*Definition* 10.1.5. The two dimensional B model $\mathcal{V}(Y) \in \mathrm{Alg}_{\mathbb{E}_2}(\mathrm{Vect})$ is the $\mathbb{E}_2$ algebra defined by

$$\mathcal{V}(Y) = \mathrm{End}_{\mathrm{QCoh}(Y^{\times 2})}(\mathrm{unit}_{Y^{\times 2}}) = \mathrm{Hom}_{\mathrm{QCoh}(Y^{\times 2})}(\Delta_\bullet \mathcal{O}_Y, \Delta_\bullet \mathcal{O}_Y) \ .$$

*Remark* 10.1.6. The two dimensional B model $\mathcal{V}(Y) = CH^\bullet(\mathcal{O}_Y) \in \mathrm{Alg}_{\mathbb{E}_2}(\mathrm{Vect})$ is by definition the Hochschild cochains on $\mathcal{O}_Y$ and thus naturally an $\mathbb{E}_2$ algebra, by Deligne's conjecture. Concretely, $\mathcal{V}(Y)$ has an $\mathbb{E}_1$ algebra structure given by composition of endomorphisms, and a second compatible $\mathbb{E}_1$ algebra structure as a space of maps between $\mathbb{E}_1$ objects internal to the category $\mathrm{QCoh}(Y^{\times 2})^\star \in \mathrm{Alg}_{\mathbb{E}_1}(\mathrm{DGCat}_{\mathrm{cont}})$.

*Remark* 10.1.7. The vector space underlying $\mathcal{V}(Y)$ is given by

$$\begin{aligned}
\mathcal{V}(Y) &= \mathrm{Hom}_{\mathrm{QCoh}(Y^{\times 2})}(\Delta_\bullet \mathcal{O}_Y, \Delta_\bullet \mathcal{O}_Y) \\
&\cong \mathrm{Hom}_{\mathrm{QCoh}(Y)}(\Delta^\bullet \Delta_\bullet \mathcal{O}_Y, \mathcal{O}_Y) \\
&\cong \mathrm{Hom}_{\mathrm{QCoh}(Y)}(\mathrm{Sym}_{\mathcal{O}_Y}^\bullet(\Omega_Y^1[1]), \mathcal{O}_Y) \\
&\cong \Gamma(X, \mathrm{PV}_Y^\bullet)
\end{aligned}$$

where $\mathrm{PV}_Y^\bullet = \mathrm{Sym}_{\mathcal{O}_Y}^\bullet(\Theta_Y[-1]) \in \mathrm{QCoh}(Y)$ is the sheaf of polyvector fields.

*Remark* 10.1.8. There is a natural identification $H^\bullet(\mathcal{V}(Y)) = \mathcal{O}(T^\vee[1]Y) \in \mathrm{Alg}_{\mathbb{P}_2}(\mathrm{Vect}_{\mathbb{K}})$ of the cohomology $\mathbb{P}_2$ algebra of $\mathcal{V}(Y)$ with the space of global functions on the 1-shifted cotangent bundle to $Y$, which is naturally 1-shifted symplectic and thus also defines a $\mathbb{P}_2$ algebra; under the isomorphism with the space of polyvector fields above, the commutative multiplication identifies with the wedge product of polyvector fields and the shifted Poisson structure identifies with the Schouten bracket.

*Example* 10.1.9. For $Y = \mathbb{C}^n$ with global coordinates $y_i$, the algebra is given by

$$\mathcal{V}(Y) = \mathrm{Sym}^\bullet(\mathbb{K}_{y_i}^n \oplus \mathbb{K}_{\partial_{y_i}}^n[-1]) \qquad \text{or less formally} \qquad \mathcal{V}(Y) \cong \mathbb{C}[y_i, \partial_{y_i}] \ ,$$

with $\partial_{y_i}$ of cohomological degree 1. In this notation, the Poisson bracket is determined by the relation $\{y_i, \partial_{y_j}\} = \delta_{ij}$.

## 10.2. Internal construction of the two dimensional B model.

Alternatively, under the analogy of Remark 10.1.2, we give an internal variant of the preceding construction, following Subsection 8.1.

*Example* 10.2.1. The topological free loop space

$$\mathcal{L}Y := Y \times_{Y \times Y} Y \in \mathrm{PreStk}$$

is a monoid in prestacks under correspondences, with multiplication and unit structure maps given by the correspondences

$$\mathcal{L}Y_{(3)} := Y \times_{Y \times Y} Y \times_{Y \times Y} Y \xrightarrow{\pi_{13}} Y \times_{Y \times Y} Y$$

with $\pi_{12}$ down to $Y \times_{Y \times Y} Y$ and $\pi_{23}$ to $Y \times_{Y \times Y} Y$, and

$$Y$$

with $\mathrm{p}_Y$ to $\mathrm{pt}$ and $\mathrm{u}$ to $Y \times_{Y \times Y} Y$.

Moreover, this monoid structure admits a natural $\mathbb{E}_2$ enhancement, in analogy with Remark 10.1.6, coming from the compatibility of this construction with the underlying monoid structure on $Y \times Y$ of Example 10.1.1.

*Remark* 10.2.2. For $Y \in \mathrm{DGSch}_{\mathrm{ft}}$ a finite type derived scheme, the topological loop space

$$\mathcal{L}Y \cong T[-1]Y \in \mathrm{DGSch}_{\mathrm{ft}}$$

is equivalent to the total space of the $(-1)$-shifted tangent bundle to $Y$. In particular, we have

$$\pi_{Y \bullet} \mathcal{O}_{\mathcal{L}Y} \cong \mathrm{Sym}^\bullet_{\mathcal{O}_Y}(\Omega^1_Y[1]) \in \mathrm{QCoh}(Y) \ .$$

Following Proposition 8.1.5 under the analogy of Remark 10.1.2, we have:

*Proposition* 10.2.3. The category $\mathrm{QCoh}(\mathcal{L}Y) \in \mathrm{Alg}_{\mathbb{E}_2}(\mathrm{DGCat}_{\mathrm{cont}})$ is naturally an $\mathbb{E}_2$-monoidal category with respect to the convolution tensor product $(\cdot) \star (\cdot) : \mathrm{QCoh}(\mathcal{L}Y)^{\otimes 2} \to \mathrm{QCoh}(\mathcal{L}Y)$ defined by the composition

$$\mathrm{QCoh}(\mathcal{L}Y) \otimes \mathrm{QCoh}(\mathcal{L}Y) \xrightarrow{\pi_{12}^* \boxtimes \pi_{23}^*} \mathrm{QCoh}(\mathcal{L}Y_{(3)})^{\otimes 2} \xrightarrow{\Delta^\bullet} \mathrm{QCoh}(\mathcal{L}Y_{(3)}) \xrightarrow{\pi_{13 \bullet}} \mathrm{QCoh}(\mathcal{L}Y) \ .$$

*Remark* 10.2.4. The space $\mathcal{L}Y_{(3)} \cong \mathcal{L}Y \times_Y \mathcal{L}Y \in \mathrm{PreStk}$ is interpreted as the space of topological maps to $Y$ from the pushout $S^1 \cup_{\mathrm{pt}} S^1$, the latter of which is homotopy equivalent to $\mathbb{P}^1 \backslash \{0, 1, \infty\}$ and should be thought of as the 'pair of pants' cobordism from $S^1 \sqcup S^1$ to $S^1$.

In analogy with Example 8.1.15, and in turn Example 8.1.14, there is an internal variant of the two dimensional B model $\mathbb{E}_2$ algebra:

*Example* 10.2.5. There is a natural lifted internal Hom object $\tilde{\mathcal{V}}(Y) \in \mathrm{Alg}_{\mathbb{E}_2}(\mathrm{QCoh}(\mathcal{L}Y))$ defined by

$$\tilde{\mathcal{V}}(Y) = \underline{\tilde{\mathrm{Hom}}}_{\mathrm{QCoh}(Y^{\times 2})}(\Delta_\bullet \mathcal{O}_Y, \Delta_\bullet \mathcal{O}_Y) = \mathcal{O}_Y \underset{Y^{\times 2}}{\tilde{\boxtimes}} \omega_{Y/Y^{\times 2}} \cong \omega_{\mathcal{L}Y/Y} \ .$$

There is a canonical equivalence $\mathrm{p}_{\mathcal{L}Y \bullet} \tilde{\mathcal{V}}(Y) \cong \mathcal{V}(Y) \in \mathrm{Alg}_{\mathbb{E}_2}(\mathrm{Vect})$ given by

$$\begin{aligned}
\mathrm{p}_{\mathcal{L}Y \bullet} \tilde{\mathcal{V}}(Y) &\cong \mathrm{p}_{Y \bullet} \pi_{Y \bullet}(\mathcal{O}_Y \underset{Y^{\times 2}}{\tilde{\boxtimes}} \omega_{Y/Y^{\times 2}}) \\
&\cong \mathrm{p}_{Y \bullet}(\mathrm{Sym}^\bullet_{\mathcal{O}_Y}(\Omega^1_Y[1]) \otimes \omega_Y^{-1}) \\
&\cong \Gamma(X, \mathrm{PV}^\bullet_Y) \\
&\cong \mathcal{V}(Y)
\end{aligned}$$

Thus, we also write

$$(10.2.1) \qquad \tilde{\mathcal{V}}(Y) = \tilde{\mathrm{PV}}^\bullet_Y := \omega_{\mathcal{L}Y/Y} \in \mathrm{QCoh}(\mathcal{L}Y) \qquad \text{noting} \qquad \pi_{Y \bullet} \tilde{\mathrm{PV}}^\bullet_Y \cong \mathrm{PV}^\bullet_Y \in \mathrm{QCoh}(Y) \ .$$

*Remark* 10.2.6. The preceding description of $\tilde{\mathcal{V}}(Y) \in \mathrm{Alg}_{\mathbb{E}_2}(\mathrm{QCoh}(\mathcal{L}Y))$ should be interpreted as an algebraic analogue of the space of distributions on the loop space, with (singular) support conditions prescribed in such a way that these distributions naturally form an $\mathbb{E}_2$ algebra with respect to the convolution monoidal structure.

*Remark* 10.2.7. Following the preceding remark, there is an expected 'de Rham' (as opposed to 'Betti') variant of the preceding constructions given by

$$\mathcal{V}(Y) = \Gamma(\mathrm{Maps}(\mathbb{D}_{\mathrm{dR}}^\circ, Y), \tilde{\mathcal{V}}(Y)) \in \mathrm{Alg}^{\mathrm{fact}}(C) ,$$

where $\tilde{\mathcal{V}}(Y) \in \mathrm{Alg}^{\mathrm{fact}}(\mathrm{QCoh}_{\mathcal{J}_{\mathrm{dR}}^{\mathrm{mer}}(Y)})$ is the analogous relative dualizing sheaf on the de Rham locally constant maps from the algebraic formal punctured disk to $Y$.

## 10.3. Equivariant structures on the two dimensional B model.
Throughout this subsection, let $Y \in \mathrm{Sch}_{\mathrm{ft}}$ be smooth algebraic variety of dimension $d_Y$.

*Definition* 10.3.1. A Calabi-Yau structure on $Y$ is a non-vanishing section $\eta \in \Gamma(X, \Omega_Y^{d_Y})$ of the top exterior power of the cotangent sheaf.

A variety $Y$ as above equipped with a Calabi-Yau structure is called a Calabi-Yau variety.

*Remark* 10.3.2. Equivalently, we can view $\eta : \mathcal{O}_Y \xrightarrow{\cong} \omega_Y$ as a trivialization of the canonical bundle, and in particular the existence of a Calabi-Yau structure on $Y$ is equivalent to the condition $c_1(Y) = 0$ of the vanishing of the first Chern class of $Y$.

*Remark* 10.3.3. Let $Y$ be a Calabi-Yau variety and recall the natural isomorphism

$$\mathrm{Sym}^\bullet_{\mathcal{O}_Y}(\Omega_Y^1[1]) \otimes \omega_Y^{-1} \xrightarrow{\cong} \mathrm{PV}_Y^\bullet .$$

The Calabi-Yau structure $\eta$ on $Y$ induces a trivialization of $\Omega_Y^{d_Y}$ and thus an isomorphism

$$(10.3.1) \qquad \eta : \mathrm{Sym}^\bullet_{\mathcal{O}_Y}(\Omega_Y^1[1])[-n] \xrightarrow{\cong} \mathrm{PV}_Y^\bullet .$$

*Example* 10.3.4. Recall the two periodic de Rham complex is defined by

$$\Omega_{Y,u}^{-\bullet} := \mathrm{Sym}^\bullet_{\mathcal{O}_Y}(\Omega_Y^1[1]) \otimes_{\mathbb{K}} \mathbb{K}[u] \qquad \text{with differential} \qquad d = d_{\mathrm{dR}} \otimes m_u ,$$

where $d_{\mathrm{dR}} : \Omega_Y^\bullet \to \Omega_Y^\bullet[1]$ is the de Rham differential, $\mathbb{K}[u] := H_{S^1}^\bullet(\mathrm{pt}; \mathbb{K})$ and $m_u : \mathbb{K}[u] \to \mathbb{K}[u]$ denotes the multiplication by $u$ map. Note that the differential is of total cohomological degree $+1$ since the de Rham differential is of degree $-1$ in the given grading and $m_u$ is of degree $+2$.

*Definition* 10.3.5. The divergence complex is defined by

$$\mathrm{Div}_{Y,u}^\bullet = \mathrm{PV}_Y^\bullet \otimes_{\mathbb{K}} \mathbb{K}[u] \qquad \text{with differential} \qquad d_\eta = \mathrm{Div}_\eta \otimes m_u ,$$

where $\mathrm{Div}_\eta := \eta \circ d_{\mathrm{dR}} \circ \eta^{-1} : \mathrm{PV}_Y^\bullet \to \mathrm{PV}_Y^\bullet[-1]$ is the divergence operator, defined in terms of the de Rham differential by transport of structure via the isomorphism in Equation 10.3.1 above.

*Proposition* 10.3.6. The divergence operator $\mathrm{Div}_\eta : \mathrm{PV}_Y^\bullet \to \mathrm{PV}_Y^\bullet[-1]$ lifts the $\mathbb{P}_2$ structure on $\Gamma(Y, \mathrm{PV}_Y^\bullet) \in \mathrm{Alg}_{\mathbb{P}_2}(\mathrm{Vect})$ of Remark 10.1.8 to $\Gamma(Y, \mathrm{PV}_Y^\bullet) \in \mathrm{Alg}_{\mathbb{P}_2}(\mathrm{Vect})$ an algebra over the operad $\mathrm{BV} \in \mathrm{Op}(\mathrm{Vect}_{\mathbb{Z}})$ of Example I-22.0.5.

*Proof.* It suffices to check the BV relation of Equation I-22.0.1 on generators $f \in \Gamma(Y, \mathcal{O}_Y)$ and $\theta \in \Gamma(Y, T_Y[-1])$, where it is given by

$$\begin{aligned}
\mathrm{Div}_\eta(f\theta) &= \eta^{-1} d_{\mathrm{dR}}(f\eta(\theta)) \\
&= \eta^{-1}(d_{\mathrm{dR}} f \wedge \eta(\theta) + f d_{\mathrm{dR}}\eta(\theta)) \\
&= \theta(f) + f\mathrm{Div}_\eta(\theta) \\
&= [f, \theta] + f\mathrm{Div}_\eta(\theta) + \mathrm{Div}_\eta(f)\theta
\end{aligned},$$

noting $\mathrm{Div}_\eta(f) = 0$, and recalling from Remark 10.1.8 that the multiplication and Poisson bracket defining the initial $\mathbb{P}_2$ structure are given by wedge product and Schouten bracket, respectively. $\square$

*Corollary* 10.3.7. The divergence complex $\Omega_{Y,u}^{-\bullet} \in \mathrm{Alg}_{\mathbb{B}\mathbb{D}_0^u}(\mathrm{D}_{\mathrm{fg}}^b(\mathbb{K}[u])$ defines an algebra over the operad $\mathbb{B}\mathbb{D}_0^u \in \mathrm{Op}(\mathrm{D}_{\mathrm{fg}}^b(\mathbb{K}[u]))$ of Example I-23.0.4.

*Proof.* This follows directly from Proposition I-23.0.1. $\square$

*Proposition* 10.3.8. A choice of $S^1$ equivariant structure $\mathcal{V}(Y) \in \mathrm{Alg}_{\mathbb{E}_2^{S^1}}$ on the two dimensional B model to $Y$ is equivalent to a Calabi-Yau structure on $Y$, and the corresponding $\mathbb{B}\mathbb{D}_0^u$ algebra

$$\mathcal{V}(Y)_u = \mathrm{Div}_{Y,u}^\bullet \in \mathrm{Alg}_{\mathbb{B}\mathbb{D}_0^u}(\mathrm{D}(\mathbb{K}[u]))$$

is given by the divergence complex of Definition 10.3.5 above.

*Proof.* $\square$

The preceding description of $S^1$ equivariant structures on the two dimensional B model can be interpreted in terms of the internal construction of Subsection 10.2, as we explain in Example 10.3.14 below, after the following preliminaries:

*Example* 10.3.9. Let $\mathbb{G}_m$ act on $\mathbb{A}_\hbar^1$ with coordinate $\hbar$ of weight $+1$, and

$$\{0\} = \mathrm{pt} \hookrightarrow \mathbb{A}^1/\mathbb{G}_m \qquad \text{and} \qquad \{1\} = \mathrm{pt} = (\mathbb{A}^1\backslash\{0\})/\mathbb{G}_m \hookrightarrow \mathbb{A}^1/\mathbb{G}_m$$

denote the inclusion of the fixed point and the complementary embedding of the open generic point. There is a canonical family of prestacks

$$Y_{\mathrm{Hdg}}^\hbar \in \mathrm{PreStk}_{/(\mathbb{A}^1/\mathbb{G}_m)}$$

over $\mathbb{A}_\hbar^1/\mathbb{G}_m$, called the Hodge stack of $Y$, such that

$$Y_{\mathrm{Hdg}}^\hbar \times_{\mathbb{A}^1/\mathbb{G}_m} \{1\} \cong Y_{\mathrm{dR}} \quad \text{and} \quad Y_{\mathrm{Hdg}}^\hbar \times_{\mathbb{A}^1/\mathbb{G}_m} \{0\} = Y_{\mathrm{Dol}}$$

where the Dolbeault stack

$$Y_{\mathrm{Dol}} := (T[1]Y)_Y^\wedge = \mathrm{B}(TY)_Y^\wedge \in \mathrm{PreStk}$$

is defined as the classifying stack of the formal tangent bundle considered as a (formal) additive group scheme over $Y$. This is a special case of the main construction of Chapter 9 in [GR17b], as explained therein for example, and was introduced by Simpson [Sim09] in the context of non-abelian Hodge theory.

*Example* 10.3.10. In particular, the category of quasicoherent sheaves on the Hodge stack

$$D^\hbar(Y) := \mathrm{QCoh}(Y_{\mathrm{Hdg}}^\hbar) \in \mathrm{ShvCat}(\mathbb{A}_\hbar^1/\mathbb{G}_m)$$

defines a quasicoherent sheaf of categories over $\mathbb{A}^1/\mathbb{G}_m$ (or 'filtered DG category' in the sense of e.g. Appendix A in [Ras16]) such that

$$D^h(Y)|_{\{1\}} \cong \mathrm{QCoh}(Y_{\mathrm{dR}}) = D(Y) \quad \text{and} \quad D^h(Y)|_{\{0\}} = \mathrm{QCoh}(Y_{\mathrm{Dol}}) \ .$$

Following Example 10.3.9 above, $D^h(Y)$ is equivalent to the category of modules over the Rees differential operators, in the sense of Definition 11.3.1.

*Remark* 10.3.11. Concretely, following remarks A.5.1 and A.5.6, the category $\mathrm{QCoh}(Y_{\mathrm{Dol}})$ is modelled by

$$\mathrm{Sym}^\bullet_{\mathcal{O}_Y}(\Omega^1_Y[-1])\text{-Mod} \qquad \text{noting} \qquad \mathrm{Sym}^\bullet_{\mathcal{O}_Y}(\Omega^1_Y[-1]) = C^\bullet_{\mathrm{CE}}(\Theta^{\mathrm{ab}}_X)$$

is given by the Chevalley-Eilenberg cochains on the abelian Lie algebroid defined by the tangent sheaf equipped with the trivial Lie bracket.

Moreover, the full family $D^h(Y)$ over $\mathbb{A}^1/\mathbb{G}_m$ is modelled by $\Omega^\bullet_{Y,\hbar}$-Mod where

$$\Omega^\bullet_{Y,\hbar} = \mathrm{Sym}^\bullet_{\mathcal{O}_Y}(\Omega^1_Y[-1]) \otimes_\mathbb{K} \mathbb{K}[\hbar] \qquad \text{with differential} \qquad d = d_{\mathrm{dR}} \otimes m_\hbar \ .$$

*Remark* 10.3.12. Let $\mathbb{K}[u] = H^\bullet_{S^1}(\mathrm{pt})$ as above and $\mathbb{A}^1_u = \mathrm{Spec}\ \mathbb{K}[u]$. The results of [BZN12] can be interpreted as proving that the topological loop space $\mathcal{L}Y \cong T[-1]Y \in \mathrm{PreStk}$ admits a 2-periodic variant of the Hodge stack deformation

$$\mathcal{L}^u Y \in \mathrm{PreStk}_{/\mathbb{A}^1_u/\mathbb{G}_m}$$

over $\mathbb{A}^1_u/\mathbb{G}_m$ such that

$$\mathcal{L}^u Y \times_{\mathbb{A}^1/\mathbb{G}_m} \{1\} \cong Y_{\mathrm{dR}} \quad \text{and} \quad \mathcal{L}^u Y \times_{\mathbb{A}^1/\mathbb{G}_m} \{0\} = \mathcal{L}Y \ .$$

Correspondingly, there is a quasicoherent sheaf of categories

$$D^u(Y) := \mathrm{QCoh}(\mathcal{L}^u Y) \in \mathrm{ShvCat}(\mathbb{A}^1_u/\mathbb{G}_m)$$

over $\mathbb{A}^1_u/\mathbb{G}_m$ such that

$$D^u(Y)|_{\{1\}} \cong \mathrm{QCoh}(Y_{\mathrm{dR}}) = D(Y) \quad \text{and} \quad D^u(Y)|_{\{0\}} = \mathrm{QCoh}(\mathcal{L}Y) \ .$$

*Remark* 10.3.13. Concretely, following Remark 10.3.11 above, the family of categories $D^u(Y)$ over $\mathbb{A}^1_u/\mathbb{G}_m$ is modeled by $\Omega^{-\bullet}_{Y,u}$-Mod where

$$\Omega^{-\bullet}_{Y,u} := \mathrm{Sym}^\bullet_{\mathcal{O}_Y}(\Omega^1_Y[1]) \otimes_\mathbb{K} \mathbb{K}[u] \qquad \text{with differential} \qquad d = d_{\mathrm{dR}} \otimes m_u \ ,$$

as in Example 10.3.4.

*Remark* 10.3.14. The requirement of a Calabi-Yau structre on $Y$ to define an $S^1$ equivariant structure can be understood from this perspective as follows: The object $\tilde{\mathrm{PV}}^\bullet_Y \in \mathrm{QCoh}(\mathcal{L}Y)$ corresponds to $\mathrm{PV}^\bullet_Y$ viewed as a module over $\mathrm{Sym}^\bullet_{\mathcal{O}_Y}(\Omega^1_Y[1])$. In order to deform $\mathrm{PV}^\bullet_Y$ to a DG module over $\Omega^{-\bullet}_{Y,u}$, it is necessary to identify

$$\mathrm{PV}^\bullet_Y \xrightarrow{\cong} \mathrm{PV}^\bullet_Y \otimes \omega_Y[-n] \cong \mathrm{Sym}^\bullet_{\mathcal{O}_Y}(\Omega^1_Y[1])[-n]$$

as in Remark 10.3.3, the latter of which deforms to the module $\Omega^{-\bullet}_{Y,u}[-n]$ over $\Omega^{-\bullet}_{Y,u}$ as desired.

## 11. The three dimensional B model

In this section, we recall a construction of the three dimensional B model to $Y$, also known as Rozansky-Witten theory to $T^\vee Y$, which is the B type topological twist of the three dimensional $\mathcal{N} = 4$ sigma model with target $T^\vee Y$ [RW97]. The three dimensional B model is a topological field theory and thus described by an $\mathbb{E}_3$ algebra $\mathcal{B}(Y) \in \mathrm{Alg}^{\mathrm{fact}}_{\mathrm{un},\mathbb{E}_1}(X)$, which is computed in terms of the category of line operators $\mathrm{QCoh}(\mathcal{L}Y) \in \mathrm{Alg}_{\mathbb{E}_2}^\bullet(\mathrm{DGCat})$ by

$$\mathcal{B}(Y) = \mathrm{Hom}_{\mathrm{QCoh}(\mathcal{L}Y)}(\mathrm{u}_\bullet \mathcal{O}_Y, \mathrm{u}_\bullet \mathcal{O}_Y) \cong \Gamma(Y, \mathrm{Sym}^\bullet_{\mathcal{O}_Y}(\Theta_Y[-2])) \ ,$$

following the discussion in [BZ14] and in keeping with the general format outlined in Section 9.1. The $\mathbb{E}_2$ monoid structure on $\mathcal{L}Y$ is the analogue of the factorization space structure in this example, and the construction is summarized as in Equation 8.1.4 by the diagrams:

(11.0.1)

11.0.1. *Summary.* In Section 11.1, we construct the three dimensional B model $\mathbb{E}_3$ algebra, and in Section 11.2 we give the internal variant of the construction. In Section 11.3 we describe the canonical $S^1$ equivariant structure on the three dimensional B model and corresponding $\mathbb{BD}_1^u$ algebra, and in section 11.4 we explain the results of the equivariant cigar reduction principle of Section I-25.2 in this example.

*Warning* 11.0.1. In keeping with Warning 9.1.1, we do not formulate specific hypotheses on the space $Y$ used in this section, so that the results stated throughout are only an outline of the general expectations. Since this section is primarily motivational, we will not give a more careful treatment, but the concerned reader may assume $Y$ is a smooth, finite type variety, for example.

11.1. **The three dimensional B model factorization algebra.** Recall from Example 10.2.1 that $\mathcal{L}Y$ is naturally a unital $\mathbb{E}_2$ monoid in prestacks, and correspondingly from Proposition 10.2.3 that $\mathrm{QCoh}(\mathcal{L}Y) \in \mathrm{Alg}_{\mathbb{E}_2}(\mathrm{DGCat}_{\mathrm{cont}})$ is naturally a unital $\mathbb{E}_2$ category.

*Remark* 11.1.1. Similarly to Remark 10.1.2, this $\mathbb{E}_2$ structure will be analogous to the factorization structure of the later examples. In fact, it is expected that $\mathbb{E}_2$ categories are a special case of the notion of factorization categories, but this construction has only been written down carefully in the $\mathbb{E}_\infty$ case, which appears in Section 7 of [Ras15a].

*Example* 11.1.2. Recall that the unit correspondence for $\mathcal{L}Y$ is given by and thus $\quad \mathrm{unit}_{\mathcal{L}Y} = \mathrm{u}_\bullet \mathcal{O}_Y \ \in \mathrm{Alg}_{\mathbb{E}_2}(\mathrm{QCoh}(\mathcal{L}Y))$

is the $\mathbb{E}_2$ unit object of $\mathrm{QCoh}(\mathcal{L}Y)$.

*Definition* 11.1.3. The three dimensional B model $\mathcal{B}(Y) \in \mathrm{Alg}_{\mathbb{E}_3}(\mathrm{Vect})$ is the $\mathbb{E}_3$ algebra defined by

$$\mathcal{B}(Y) = \mathrm{End}_{\mathrm{QCoh}(\mathcal{L}Y)}(\mathrm{unit}_{\mathcal{L}Y}) = \mathrm{Hom}_{\mathrm{QCoh}(\mathcal{L}Y)}(\mathrm{u}_\bullet \mathcal{O}_Y, \mathrm{u}_\bullet \mathcal{O}_Y) \ .$$

*Remark* 11.1.4. The $\mathbb{E}_3$ structure is determined by the $\mathbb{E}_1$ structure coming from composition of endomorphisms compatible with the underlying $\mathbb{E}_2$ structure on the object $\mathrm{unit}_{\mathcal{L}Y} \in \mathrm{Alg}_{\mathbb{E}_2}(\mathrm{QCoh}(\mathcal{L}Y))$ internal to $\mathrm{QCoh}(\mathcal{L}Y)$, similarly to Remark 10.1.6.

*Remark* 11.1.5. The vector space underlying $\mathcal{B}(Y)$ is given by

$$\begin{aligned}
\mathcal{B}(Y) &= \mathrm{Hom}_{\mathrm{QCoh}(\mathcal{L}Y)}(\mathrm{u}_\bullet \mathcal{O}_Y, \mathrm{u}_\bullet \mathcal{O}_Y) \\
&\cong \mathrm{Hom}_{\mathrm{QCoh}(Y)}(\mathrm{u}^\bullet \mathrm{u}_\bullet \mathcal{O}_Y, \mathcal{O}_Y) \\
&\cong \mathrm{Hom}_{\mathrm{QCoh}(Y)}(\mathrm{Sym}^\bullet_{\mathcal{O}_Y}(\Omega_Y[2]), \mathcal{O}_Y) \\
&\cong \Gamma(Y, \mathrm{SV}^\bullet_Y)
\end{aligned}$$

where $\mathrm{SV}^\bullet_Y = \mathrm{Sym}^\bullet_{\mathcal{O}_Y}(\Theta_Y[-2]) \in \mathrm{QCoh}(Y)$ is the space of symmetric polyvector fields on $Y$ generated in degree 2.

*Remark* 11.1.6. There is a natural identification $H^\bullet(\mathcal{B}(Y)) = \mathcal{O}(T^\vee[2]Y) \in \mathrm{Alg}_{\mathbb{P}_3}(\mathrm{Vect}_{\mathbb{K}})$ of the cohomology $\mathbb{P}_3$ algebra of $\mathcal{B}(Y)$ with the space of global functions on the 2-shifted cotangent bundle to $Y$, which is naturally 2-shifted symplectic and thus also defines a $\mathbb{P}_3$ algebra; under the isomorphism with the space of polyvector fields above, the commutative multiplication identifies with the symmetric product of symmetric polyvector fields and the shifted Poisson structure identifies with the Schouten bracket.

*Example* 11.1.7. For $Y = \mathbb{C}^n$ with global coordinates $y_i$, the algebra is given by

$$\mathcal{B}(Y) = \mathrm{Sym}^\bullet(\mathbb{K}^n_{y_i} \oplus \mathbb{K}^n_{\partial_{y_i}}[-2]) \qquad \text{or less formally} \qquad \mathcal{V}(Y) \cong \mathbb{C}[y_i, \partial_{y_i}] \ ,$$

with $\partial_{y_i}$ of cohomological degree 2. In this notation, the Poisson bracket is determined by the relation $\{y_i, \partial_{y_j}\} = \delta_{ij}$.

*Remark* 11.1.8. Analogously to Remark 10.2.7, there is a de Rham (as opposed to Betti) variant of the above construction, given by

$$\mathcal{B}(Y) = \mathrm{Hom}_{\mathrm{QCoh}(\mathrm{Maps}(\mathbb{D}^\circ_{\mathrm{dR}}, Y))}(\mathrm{u}_\bullet \mathcal{O}_{\mathrm{Maps}(\mathbb{D}_{\mathrm{dR}}, Y)}, \mathrm{u}_\bullet \mathcal{O}_{\mathrm{Maps}(\mathbb{D}_{\mathrm{dR}}, Y)}) \ .$$

## 11.2. Internal construction of the three dimensional B model.
There is again an internal variant of the preceding construction, which we now give following Subsection 8.1:

*Example* 11.2.1. The topological sphere space

$$\mathcal{S}Y := Y \times_{\mathcal{L}Y} Y \ \in \mathrm{PreStk}$$

is a monoid in prestacks under correspondences, with multiplication and unit structure maps given by the correspondences

$$\mathcal{S}Y_{(3)} := Y \times_{\mathcal{L}Y} Y \times_{\mathcal{L}Y} Y \xrightarrow{\pi_{13}} Y \times_{\mathcal{L}Y} Y \qquad\qquad Y$$

with $\pi_{12}$, $\pi_{23}$ and $\mathrm{pt}$, $Y \times_{\mathcal{L}Y} Y$ as shown, "and"

Moreover, this monoid structure admits a natural $\mathbb{E}_3$ enhancement, in analogy with the $\mathbb{E}_2$ enhancement in Example 10.2.1, coming from the compatibility of this construction with the underlying $\mathbb{E}_2$ monoid structure on $\mathcal{L}Y$ given in *loc. cit.*.

Following Proposition 8.1.5 under the analogy of Remark 11.1.1, we have:

*Proposition* 11.2.2. The category $\mathrm{QCoh}(\mathcal{S}Y) \in \mathrm{Alg}_{\mathbb{E}_3}(\mathrm{DGCat}_{\mathrm{cont}})$ is naturally an $\mathbb{E}_3$-monoidal category with respect to the convolution tensor product $(\cdot) \star (\cdot) : \mathrm{QCoh}(\mathcal{S}Y)^{\otimes 2} \to \mathrm{QCoh}(\mathcal{S}Y)$ defined by the composition

$$\mathrm{QCoh}(\mathcal{S}Y) \otimes \mathrm{QCoh}(\mathcal{S}Y) \xrightarrow{\pi_{12}^! \boxtimes \pi_{23}^\bullet} \mathrm{QCoh}(\mathcal{S}Y_{(3)})^{\otimes 2} \xrightarrow{\Delta^\bullet} \mathrm{QCoh}(\mathcal{S}Y_{(3)}) \xrightarrow{\pi_{13\bullet}} \mathrm{QCoh}(\mathcal{S}Y) .$$

Similarly, in analogy with Example 8.1.15, and in turn Example 8.1.14, there is an internal variant of the three dimensional B model $\mathbb{E}_3$ algebra:

*Example* 11.2.3. There is a natural lifted internal Hom object $\tilde{\mathcal{B}}(Y) \in \mathrm{Alg}_{\mathbb{E}_3}(\mathrm{QCoh}(\mathcal{S}Y))$ defined by

$$\tilde{\mathcal{B}}(Y) = \underline{\mathrm{Hom}}_{\mathrm{QCoh}(\mathcal{L}Y)}(\mathrm{u}_\bullet \mathcal{O}_Y, \mathrm{u}_\bullet \mathcal{O}_Y) = \mathcal{O}_Y \underset{\mathcal{L}Y}{\boxtimes} \omega_{Y/\mathcal{L}Y} \cong \omega_{\mathcal{S}Y/Y} .$$

There is a canonical equivalence $\mathrm{p}_{\mathcal{S}Y\bullet}\tilde{\mathcal{B}}(Y) \cong \mathcal{B}(Y) \in \mathrm{Alg}_{\mathbb{E}_2}(\mathrm{Vect})$ given by

$$\begin{aligned}
\mathrm{p}_{\mathcal{S}Y\bullet}\tilde{\mathcal{B}}(Y) &\cong \mathrm{p}_{Y\bullet}\pi_{Y\bullet}(\mathcal{O}_Y \underset{\mathcal{L}Y}{\boxtimes} \omega_{Y/\mathcal{L}Y}) \\
&\cong \mathrm{p}_{Y\bullet}(\mathrm{Sym}^\bullet_{\mathcal{O}_Y}(\Omega^1_Y[2]) \otimes \omega_Y^{-1}[-n]) \\
&\cong \Gamma(X, \mathrm{S}\mathrm{V}^\bullet_Y) \\
&\cong \mathcal{B}(Y)
\end{aligned}.$$

*Remark* 11.2.4. The preceding description of $\tilde{\mathcal{B}}(Y) \in \mathrm{Alg}_{\mathbb{E}_3}(\mathrm{QCoh}(\mathcal{S}Y))$ should be interpreted as an algebraic analogue of the space of distributions on the sphere space, with (singular) support conditions prescribed in such a way that these distributions naturally form an $\mathbb{E}_3$ algebra with respect to the convolution monoidal structure.

## 11.3. **Equivariant structures on the three dimensional B model.**

*Definition* 11.3.1. The Rees differential operators $\mathcal{D}_{Y,\hbar} \in \mathrm{Alg}_{\mathbb{B}\mathbb{D}_1^\hbar}(\mathrm{D}^b(Y) \otimes \mathrm{D}^b_{\mathrm{fg}}(\mathbb{K}[\hbar]))$ are defined to be the (sheaf of) algebras

$$\mathcal{D}_{Y,\hbar} \cong \otimes^\bullet_{\mathcal{O}_Y}(\Theta_Y) \otimes_{\mathbb{K}} \mathbb{K}[\hbar]/(\theta \cdot \chi - \chi \cdot \theta - \hbar[\theta, \chi], \ \theta \cdot f - f \cdot \theta - \hbar\theta(f)) ,$$

over the operad $\mathbb{B}\mathbb{D}_1^\hbar \in \mathrm{Op}(\mathrm{D}^b_{\mathrm{fg}}(\mathbb{K}[\hbar]))$ of Definition I-C.6.5, where the relations are imposed on all $\theta, \chi \in \Gamma(Y, \Theta_Y)$ and $f \in \Gamma(Y, \mathcal{O}_Y)$.

*Remark* 11.3.2. In analogy with the discussion in remarks 10.3.12 and 10.3.13, there is a natural two-periodic analogue $\mathcal{D}_{Y,u} \in \mathrm{Alg}_{\mathbb{B}\mathbb{D}_1^u}(\mathrm{D}^b(Y) \otimes \mathrm{D}^b_{\mathrm{fg}}(\mathbb{K}[u]))$ of the Rees differential operators given by the sheaf of algebras

$$\mathcal{D}_{Y,u} \cong \otimes^\bullet_{\mathcal{O}_Y}(\Theta_Y[-2]) \otimes_{\mathbb{K}} \mathbb{K}[u]/(\theta \cdot \chi - \chi \cdot \theta - u[\theta, \chi], \ \theta \cdot f - f \cdot \theta - u\theta(f)) ,$$

over the operad $\mathbb{B}\mathbb{D}_1^u \in \mathrm{Op}(\mathrm{D}^b(\mathbb{K}[u]))$ of Definition I-23.0.5, where the relations are imposed on all $\theta, \chi \in \Gamma(Y, \Theta_Y[-2])$ and $f \in \Gamma(Y, \mathcal{O}_Y)$.

*Proposition* 11.3.3. There is a canonical $S^1$ equivariant structure $\mathcal{B}(Y) \in \mathrm{Alg}_{\mathbb{E}_3^{S^1}}(\mathrm{Vect})$ such that the corresponding $\mathbb{B}\mathbb{D}_1^u$ algebra

$$\mathcal{B}(Y)^u = \Gamma(Y, \mathcal{D}_{Y,u}) \ \in \mathrm{Alg}_{\mathbb{B}\mathbb{D}_1^u}(\mathrm{D}^b_{\mathrm{fg}}(\mathbb{K}[u])) ,$$

is given by the global sections of the 2-periodic Rees differential operators of Remark 11.3.2 above.

*Proof.* The object $u_\bullet \mathcal{O}_Y \in \mathrm{QCoh}(\mathcal{L}Y)$ admits a natural deformation

$$u_* \mathcal{O}_Y \in \mathrm{QCoh}(\mathcal{L}^u Y) \qquad \text{such that} \qquad u_* \mathcal{O}_Y \cong \mathcal{D}_{Y,u} \in \mathrm{QCoh}(\mathcal{L}^u Y) = D^u(Y) \ .$$

Thus, the endomorphism algebra naturally deforms to

$$\mathcal{B}(Y)^u = \mathrm{Hom}_{\mathrm{QCoh}(\mathcal{L}^u Y)}(u_* \mathcal{O}_Y, u_* \mathcal{O}_Y) = \mathrm{Hom}_{D^u(Y)}(\mathcal{D}_{Y,u}, \mathcal{D}_{Y,u}) \ ,$$

the latter of which is computed by the induction adjunction, giving

$$\mathcal{B}(Y)^u = \mathrm{Hom}_{D^u(Y)}(\mathcal{D}_{Y,u}, \mathcal{D}_{Y,u}) \cong \Gamma(Y, \mathcal{D}_{Y,u}) \ \in \mathrm{Alg}_{\mathbb{BD}_1^u}(\mathrm{D}_{\mathrm{fg}}^b(\mathbb{K}[u]))$$

so that the claim follows by Proposition I-25.2.3. $\qquad\square$

As in the discussion of equivariant structures on the two dimensional B model, there is an interpretation of the preceding proposition in terms of the internal construction of the three dimensional B model given in Subsection 11.2. To begin, following Remark 10.3.12, we have:

*Remark* 11.3.4. The topological sphere space $\mathcal{S}Y = Y \times_{\mathcal{L}Y} Y \in \mathrm{PreStk}$ admits a deformation over $\mathbb{A}_u^1/\mathbb{G}_m$ defined by

$$\mathcal{S}^u Y := Y^u \times_{\mathcal{L}^u Y} Y^u \ \in \mathrm{PreStk}_{/(\mathbb{A}_u^1/\mathbb{G}_m)} \qquad \text{where} \qquad Y^u := Y \times \mathbb{A}_u^1/\mathbb{G}_m \ ,$$

and $\mathcal{L}^u Y$ is as defined in Remark 10.3.12. In particular, we have

$$\mathcal{S}^u Y \times_{\mathbb{A}^1/\mathbb{G}_m} \{1\} \cong Y \times_{Y_{\mathrm{dR}}} Y \quad \text{and} \quad \mathcal{S}^u Y \times_{\mathbb{A}^1/\mathbb{G}_m} \{0\} = Y \times_{\mathcal{L}Y} Y = \mathcal{S}Y$$

Correspondingly, there is a quasicoherent sheaf of categories

$$\mathrm{QCoh}(\mathcal{S}^u Y) \in \mathrm{ShvCat}(\mathbb{A}_u^1/\mathbb{G}_m)$$

over $\mathbb{A}_u^1/\mathbb{G}_m$ such that

$$\mathrm{QCoh}(\mathcal{S}^u Y)|_{\{1\}} \cong \mathrm{QCoh}(Y \times_{Y_{\mathrm{dR}}} Y) \quad \text{and} \quad \mathrm{QCoh}(\mathcal{S}^u Y)|_{\{0\}} = \mathrm{QCoh}(\mathcal{S}Y) \ .$$

*Example* 11.3.5. The internal variant of the three dimensional B model $\tilde{\mathcal{B}}(Y) \cong \pi_Y^! \mathcal{O}_Y \in \mathrm{Alg}_{\mathbb{E}_3}(\mathrm{QCoh}(\mathcal{S}Y))$ defined in Example 11.2.3 admits a natural deformation

$$\tilde{\mathcal{B}}(Y)^u := \mathcal{O}_{Y^u} \underset{\mathcal{L}^u Y}{\boxtimes} \omega_{Y^u/\mathcal{L}^u Y} \cong \omega_{\mathcal{S}^u Y/Y^u} = \pi_{Y^u}^! \mathcal{O}_{Y^u} \ \in \mathrm{QCoh}(\mathcal{S}^u Y)$$

with fibre over the generic point given by

$$\pi_{Y^u}^! \mathcal{O}_{Y^u}|_{\{1\}} = \pi_Y^! \mathcal{O}_Y \in \mathrm{QCoh}(Y \times_{Y_{\mathrm{dR}}} Y) \qquad \text{so that} \qquad \pi_{Y^u *} \pi_{Y^u}^! \mathcal{O}_Y|_{\{1\}} \cong \mathcal{D}_Y \in \mathrm{QCoh}(Y)$$

by the quasicoherent variant of Proposition 5.1.3 in [GR14a], which follows from the indcoherent variant stated there by applying Proposition 5.2.7(i) in *loc. cit.*. In particular, we have

$$\mathrm{p}_{\mathcal{S}^u Y *} \tilde{\mathcal{B}}(Y)^u \cong \Gamma(Y, \mathcal{D}_{Y,u}) \cong B(Y)^u \ .$$

*Remark* 11.3.6. Concretely, the restriction to the generic fibre of the internal Hom object $\tilde{\mathcal{B}}(Y)^u \in \mathrm{QCoh}(\mathcal{S}^u Y)$ is given by

$$\tilde{\mathcal{B}}(Y)^u|_{\{1\}} = \mathcal{O}_Y \boxtimes_{Y_{\mathrm{dR}}} \omega_Y \in \mathrm{QCoh}(Y \times_{Y_{\mathrm{dR}}} Y) \ ,$$

which agrees with the above computation of $\mathcal{B}(Y)^u$ by Proposition 5.2.4(a) in [GR14a].

11.4. **Equivariant cigar compactification for the three dimensional B model.** In this subsection, we apply the results of Subsection I-25.2 to the example of the three dimensional B model.

*Proposition* 11.4.1. There is an equivalence of prestacks over $\mathbb{A}^1/\mathbb{G}_m$

$$\mathcal{S}^u Y = Y^u \times_{\mathcal{L}^u Y} Y^u \cong \mathcal{L}^u Y \times_{(\mathcal{L}^u Y)^{\times 2}} Y^{u \times 2}$$

*Proof.* Each of the spaces admits the required maps such that, by the universal property of fibre products, it maps to the other; the induced maps are evidently inverse equivalences. □

*Corollary* 11.4.2. There is an equivalence of sheaves of categories over $\mathbb{A}^1/\mathbb{G}_m$

$$\mathrm{QCoh}(\mathcal{S}^u Y) = \mathrm{QCoh}(Y^u \times_{\mathcal{L}^u Y} Y^u) \cong \mathrm{QCoh}(\mathcal{L}^u Y \times_{(\mathcal{L}^u Y)^{\times 2}} Y^{u \times 2}) \ .$$

*Example* 11.4.3. Under the equivalence of the preceding corollary, the internal variant of the three dimensional B model is given by

$$\tilde{B}(Y)^u = \mathcal{O}_{Y^u} \underset{\mathcal{L}^u Y}{\boxtimes} \omega_{Y^u/\mathcal{L}^u Y} \ \mapsto \ \Omega^{-\bullet}_{Y,u} \underset{(\mathcal{L}^u Y)^{\times 2}}{\boxtimes} \left( \mathcal{O}_{Y^u} \boxtimes \omega_{Y^u/\mathcal{L}^u Y} \right)$$

where $\Omega^{-\bullet}_{Y,u} \in \mathrm{QCoh}(\mathcal{L}^u Y)$ is as defined in Equation 10.3.13.

Towards the statement of the main result of this subsection, we now give a description of the negative cyclic chains on $\mathcal{B}(Y)$:

*Proposition* 11.4.4. There is a natural equivalence

$$\mathrm{CC}^-_\bullet(\mathcal{B}(Y)) \cong \mathrm{Hom}_{\mathrm{QCoh}((\mathcal{L}^u Y)^{\times 2})}(\Delta_\bullet \mathcal{O}_{\mathcal{L}^u Y}, \Delta_\bullet \mathcal{O}_{\mathcal{L}^u Y}) \ \in \mathrm{Alg}_{\mathbb{E}_2}(\mathrm{D}^b_{\mathrm{fg}}(\mathbb{K}[u])) \ .$$

*Proof.* There is a natural commutative diagram
(11.4.1)

$$\begin{array}{ccc}
\mathrm{Alg}_{\mathbb{E}_2}(\mathrm{DGCat}) & \xrightarrow{\mathrm{End}_{(\cdot)}(\mathrm{unit})} & \mathrm{Alg}_{\mathbb{E}_3}(\mathrm{Vect}) \\
\downarrow{\mathrm{CC}^-_\bullet} & & \downarrow{\mathrm{CC}^-_\bullet} \\
\mathrm{Alg}_{\mathbb{E}_1}(\mathrm{ShvCat}(\mathbb{A}^1_u/\mathbb{G}_m)) & \xrightarrow{\mathrm{End}_{(\cdot)}(\mathrm{unit})} & \mathrm{Alg}_{\mathbb{E}_2}(\mathrm{D}(\mathbb{K}[u]))
\end{array}$$

under which

$$\begin{array}{ccc}
\mathrm{QCoh}(\mathcal{L}Y) & \longmapsto & \mathcal{B}(Y) \\
\downarrow & & \downarrow \\
\mathrm{QCoh}((\mathcal{L}^u Y)^{\times 2}) & \longmapsto & \mathrm{End}_{\mathrm{QCoh}((\mathcal{L}^u Y)^{\times 2})}(\Delta_\bullet \mathcal{O}_{\mathcal{L}^u Y})
\end{array}$$

where we have assumed an extension of Theorem 5.3 of [BZFN10] (in fact, its variant for hochschild chains, also outlined in *loc. cit.*) to the homotopy $S^1$ invariants to compute $\mathrm{CC}^-_\bullet(\mathrm{QCoh}(\mathcal{L}Y)) \cong \mathrm{QCoh}((\mathcal{L}^u Y)^{\times 2})$. □

We can now apply the equivariant cigar compactification principle in the present example:

*Proposition* 11.4.5. The $S^1$ equivariant structure $\mathcal{B}(Y) \in \mathrm{Alg}_{\mathbb{E}^{S^1}_3}(\mathrm{Perf}_\mathbb{K})$ of Proposition 11.3.3 corresponds, under the equivalence of Proposition I-25.2.3, to the module structure

$$\mathcal{B}(Y)^u \cong \mathrm{Hom}_{\mathrm{QCoh}((\mathcal{L}^u Y)^{\times 2})}(\Delta_\bullet \mathcal{O}_{\mathcal{L}^u Y}, \mathcal{D}_{Y,u} \boxtimes \mathcal{D}_{Y,u}) \in \mathrm{CC}^-_\bullet(\mathcal{B}(Y))\text{-}\mathrm{Mod}(\mathrm{Alg}_{\mathbb{E}_1}(\mathrm{D}^b_{\mathrm{fg}}(\mathbb{K}[u]))) \ ,$$

given by precomposition by endomorphisms under the equivalence of Proposition 11.4.4 above.

*Proof.* Passing to global sections in the computation in Example 11.4.3 identifies $\mathcal{B}(Y)^u$ with the Hom space as claimed; the module structure corresponds to that given by precomposition with endomorphisms, by functoriality. □

*Remark* 11.4.6. Concretely, the restriction to the generic fibre is given by the central action of the de Rham cochains

$$C^\bullet_{\mathrm{dR}}(Y) = \mathrm{Hom}_{D(Y^{\times 2})}(\Delta_* \omega_Y, \Delta_* \omega_Y) \ \in \mathrm{Alg}_{\mathbb{E}_2}(\mathrm{Perf}_{\mathbb{K}})$$

on the algebra of global differential operators

$$\Gamma(Y, \mathcal{D}_Y) = \mathrm{Hom}_{D(Y^{\times 2})}(\Delta_* \omega_Y, \mathcal{D}_Y \boxtimes \mathcal{D}_Y) \ \in C^\bullet_{\mathrm{dR}}(Y)\text{-Mod}(\mathrm{Alg}_{\mathbb{E}_1}(\mathrm{Perf}_{\mathbb{K}})) \ .$$

## 12. The three dimensional A model: Overview

In this section, we recall the construction of the three dimensional A model to $Y$, also known as twisted Rozansky-Witten theory to $T^\vee Y$, which is the A type topological twist of the three dimensional $\mathcal{N} = 4$ sigma model with target $T^\vee Y$; the local observables of the A type twist are precisely the Coulomb branch chiral ring, which was studied extensively in the physics literature, for example in [BDG17] and references therein, but a general definition of the algebra in non-abelian gauge theories remained ellusive. A major breakthrough in this area was the work of Braverman, Finkelberg, and Nakajima [BFN18], which gave a geometric construction of the Coulomb branch algera for a general gauge theory with complexified gauge group $G$ and matter representation $N$. In this section, we recall the main construction schematically in terms of a general space $Y$, postponing the more careful construction given in *loc. cit.* in the case $Y = N/G$ to the following section 13.

The Coulomb branch construction of [BFN18] can naturally be interpreted as defining a factorization $\mathbb{E}_1$ algebra $\mathcal{A}(Y) \in \mathrm{Alg}^{\mathrm{fact}}_{\mathbb{E}_1, \mathrm{un}}(X)$ on any smooth algebraic curve $X$, defined in keeping with the general format outlined in Section 9.1, by

$$\mathcal{A}(Y)_x := \mathcal{H}om_{D(Y_\mathcal{K})}(\iota_* \omega_{Y_\mathcal{O}}, \iota_* \omega_{Y_\mathcal{O}}) \cong C^{\mathrm{BM}}_\bullet(Y_\mathcal{O} \times_{Y_\mathcal{K}} Y_\mathcal{O})[-2 \dim Y_\mathcal{O}] \ ,$$

where the line operator category $D(Y_\mathcal{K}) \in \mathrm{Cat}^{\mathrm{fact}}_{\mathrm{un}}(X)$ is given by the factorization category $D$ modules on the meromorphic jet scheme to $Y$; the construction is summarized as in Equation 8.1.4 by the following diagrams in factorization spaces and categories, which for simplicity we denote by their fibre over a fixed point $x \in X$:

(12.0.1)

In fact, the three dimensional A model factorization $\mathbb{E}_1$ algebra admits a canonical $\mathbb{G}_a \rtimes \mathbb{G}_m$ equivariant structure $\mathcal{A}(Y) \in \mathrm{Alg}^{\mathrm{fact}}_{\mathbb{E}_1}(\mathbb{A}^1)^{\mathbb{G}_a \rtimes \mathbb{G}_m}$, which corresponds to the fact that it is actually a three dimensional topological field theory and admits an $\Omega$-background deformation; the $\mathbb{BD}^u_1$ algebra $\mathcal{A}(Y)_u \in \mathrm{Alg}_{\mathbb{BD}^u_1}(\mathrm{Vect})$ corresponding to $\mathcal{A}$ under Proposition I-25.1.1 is the filtered quantization which is the more explicit object of study in [BFN18].

12.0.1. *Summary.* In Section 12.1, we give an overview of the construction of the three dimensional A model as a factorization $\mathbb{E}_1$ algebra, in Section 12.2 we describe the internal variant of the construction, and in Section 12.3 we recall the existence of its canonical $\mathbb{G}_a \rtimes \mathbb{G}_m$ equivariant structure.

*Warning* 12.0.1. In keeping with Warning 9.1.1, we do not formulate specific hypotheses on the space $Y$ used in this section, so that the results stated throughout are only an outline of the general expectations. In Section 13, we restrict to the case $Y = N/G$ and give careful statements and proofs of these results, again following [BFN18].

12.1. **The three dimensional A-model factorization algebra.** The starting point for the construction of the three dimensional A model is the de Rham stack of the loop space of $Y$, and its associated factorization category:

*Example* 12.1.1. The de Rham stack of the loop space of $Y$ is the (unital) factorization space $\mathcal{J}^{\mathrm{mer}}(Y)_{\mathrm{dR}} \in \mathrm{PreStk}^{\mathrm{fact}}_{\mathrm{un}}(X)$ defined by considering the de Rham stack, as in Remark 5.2.10, of the (unital) factorization space $\mathcal{J}^{\mathrm{mer}}(Y)$ of Example 4.2.4.

*Remark* 12.1.2. Concretely, the prestacks $\mathcal{J}^{\mathrm{mer}}(Y)_{\mathrm{dR},I} \in \mathrm{PreStk}_{/X^I}$ assigned to each $I \in \mathrm{fSet}$ and the fibre $\mathcal{J}^{\mathrm{mer}}(Y)_{\mathrm{dR},x} \in \mathrm{PreStk}$ over $x \in X$ are given by

$$\mathcal{J}^{\mathrm{mer}}(Y)_{\mathrm{dR},I} = \{x \in X_{\mathrm{dR}}^I, \ a : \mathbb{D}_x^\circ \to Y \text{ a map over } X\}_{\mathrm{dR}} \qquad \text{and} \qquad \mathcal{J}^{\mathrm{mer}}(Y)_{\mathrm{dR},x} = Y_{\mathcal{K}_x,\mathrm{dR}} \ .$$

Similarly, the unit data corresponding to $\pi : \emptyset \to I$ and its restriction to the point $x \in X$ are given by

*Example* 12.1.3. The category of $D$ modules on the loop space of $Y$ is the (unital) factorization category

$$D_{\mathcal{J}^{\mathrm{mer}}(Y)} := \mathrm{QCoh}_{\mathcal{J}^{\mathrm{mer}}(Y)_{\mathrm{dR}}} \in \mathrm{Cat}^{\mathrm{fact}}_{\mathrm{un}}(X)$$

associated to $\mathcal{J}^{\mathrm{mer}}(Y)_{\mathrm{dR}} \in \mathrm{PreStk}^{\mathrm{fact}}_{\mathrm{un}}(X)$, following Definition 5.2.7.

*Remark* 12.1.4. Concretely, following Remark 5.2.8, the sheaf of categories assigned to each $I \in \mathrm{fSet}$, its sections category, and the fibre category over $x \in X$, are given by

$$\mathrm{p}_{I,\mathrm{dR}} {}_* \mathrm{QCoh}_{\mathcal{J}^{\mathrm{mer}}(Y)_{I,\mathrm{dR}}} \in \mathrm{ShvCat}(X_{\mathrm{dR}}^I) \qquad D(\mathcal{J}^{\mathrm{mer}}(Y)_I) \in \mathrm{DGCat} \qquad \text{and} \qquad D(Y_{\mathcal{K}_x}) \in \mathrm{DGCat} \ .$$

*Warning* 12.1.5. For simplicity, the preceding Example is stated in terms of 'left' $D$ modules, in the sense that we use $\mathrm{QCoh}(X_{\mathrm{dR}})$ rather than $\mathrm{IndCoh}(X_{\mathrm{dR}})$ as our model for $D(X)$. In the following subsections, we give careful statements of the results outlined below using the theory of $D$ modules on infinite dimensional varieties recalled in subappendices B.7 and B.8, following [Ras15b] and references therein, the definition of which is somewhat agnostic about the relation to coherent sheaves. However, in Section 14 and its successors we will see that the indcoherent setting is in fact the correct variant for our applications of interest.

Following Example 5.4.3, we have:

*Example* 12.1.6. The factorization unit object defines a factorization algebra

$$\mathrm{unit}_{D_{\mathcal{J}^{\mathrm{mer}}(Y)}} \in \mathrm{Alg}^{\mathrm{fact}}_{\mathrm{un}}(D_{\mathcal{J}^{\mathrm{mer}}(Y)}) \qquad \mathrm{unit}_{D_{\mathcal{J}^{\mathrm{mer}}(Y)}} = \iota_{\mathcal{J}(Y)*}\mathrm{p}^!_{\mathcal{J}(Y)}\omega_{\mathrm{Ran}_{X,\mathrm{un}}} = \iota_{\mathcal{J}(Y)*}\omega_{\mathcal{J}(Y)}$$

internal to $D_{\mathcal{J}^{\mathrm{mer}}(Y)}$, where

$$\mathrm{p}_{\mathcal{J}(Y)} : \mathcal{J}(Y) \to \mathrm{Ran}_{X_{\mathrm{dR}},\mathrm{un}} \qquad \text{and} \qquad \iota_{\mathcal{J}(Y)} : \mathcal{J}(Y) \to \mathcal{J}^{\mathrm{mer}}(Y)$$

are the factorization space structure map for $\mathcal{J}(Y) = \mathrm{unit}_{\mathcal{J}^{\mathrm{mer}}(Y)}$ and the map of unital factorization spaces given by the inclusion of arcs into loops, respectfully, and

$$\mathrm{p}^!_{\mathcal{J}(Y)} : D_{\mathrm{Ran}_{X,\mathrm{un}}} \to D_{\mathcal{J}(Y)} \qquad \text{and} \qquad \iota_{\mathcal{J}(Y)*} : D_{\mathcal{J}(Y)} \to D_{\mathcal{J}^{\mathrm{mer}}(Y)}$$

are the induced unital factorization functors.

*Remark* 12.1.7. Concretely, the object of the sections category assigned to each $I \in \mathrm{fSet}$, and the object of the fibre category over $x \in X$, are given by

$$\mathrm{unit}_{D_{\mathcal{J}^{\mathrm{mer}}(Y)}, I} = \iota_{\mathcal{J}(Y)_I*} \omega_{\mathcal{J}(Y)_I} \in D(\mathcal{J}^{\mathrm{mer}}(Y)_I) \qquad \text{and} \qquad \mathrm{unit}_{D_{\mathcal{J}^{\mathrm{mer}}(Y)}, x} = \iota_{x*} \omega_{Y_{\mathcal{O}_x}} \in D(Y_{\mathcal{K}_x}) \;,$$

where $\iota_x : Y_{\mathcal{O}_x} \to Y_{\mathcal{K}_x}$ is the map corresponding to the inclusion of arcs into loops and $\omega_{Y_{\mathcal{O}_x}}$ is the dualizing sheaf of the arc space.

Now, suppose that $\mathrm{unit}_{D_{\mathcal{J}^{\mathrm{mer}}(Y)}} \in \mathrm{Alg}^{\mathrm{fact}}_{\mathrm{un}}(D_{\mathcal{J}^{\mathrm{mer}}(Y)})$ admits internal Hom objects over $X$, in the sense of Definition 8.1.9. Then, following Example 7.1.7, we give the following tentative definition of the factorization $\mathbb{E}_1$ algebra describing the three dimensional A model to $Y$:

*Tentative Definition* 12.1.8. The three dimensional A model to $Y$ is the unital factorization $\mathbb{E}_1$-algebra

$$\mathcal{A}(Y) = \mathcal{H}\mathrm{om}_{D(\mathcal{J}^{\mathrm{mer}}(Y))}(\mathrm{unit}_{D_{\mathcal{J}^{\mathrm{mer}}(Y)}}, \mathrm{unit}_{D_{\mathcal{J}^{\mathrm{mer}}(Y)}}) \; \in \mathrm{Alg}^{\mathrm{fact}}_{\mathbb{E}_1, \mathrm{un}}(X) \;.$$

*Remark* 12.1.9. Concretely, following Remark 8.1.10, the $\mathbb{E}_1$ algebra internal to $D(X^I)$ assigned to each $I \in \mathrm{fSet}_\emptyset$, and the $\mathbb{E}_1$ algebra in Vect over each $x \in X$, are given by

$$\mathcal{A}(Y)_I = \mathrm{p}_{\mathcal{J}^{\mathrm{mer}}(Y)_I*} \mathcal{H}\mathrm{om}_{D(\mathcal{J}^{\mathrm{mer}}(Y)_I)}(\iota_{\mathcal{J}(Y)_I*} \omega_{\mathcal{J}(Y)_I}, \iota_{\mathcal{J}(Y)_I*} \omega_{\mathcal{J}(Y)_I}) \qquad \in \mathrm{Alg}_{\mathbb{E}_1}(D(X^I)) \;, \text{ and}$$

$$\mathcal{A}(Y)_x = \mathrm{Hom}_{D(Y_{\mathcal{K}_x})}(\iota_{x*} \omega_{Y_{\mathcal{O}_x}}, \iota_{x*} \omega_{Y_{\mathcal{O}_x}}) \qquad \in \mathrm{Alg}_{\mathbb{E}_1}(\mathrm{Vect}) \;.$$

Further, following Example A.6.7, these are expected to be calculated by the 'renormalized' Borel-Moore homology groups:

$$\mathcal{A}(Y)_x = \mathrm{Hom}_{D(Y_{\mathcal{K}_x})}(\iota_{x*} \omega_{Y_{\mathcal{O}_x}}, \iota_{x*} \omega_{Y_{\mathcal{O}_x}}) \cong C^{\mathrm{BM}}_\bullet(Y_{\mathcal{O}_x} \times_{Y_{\mathcal{K}_x}} Y_{\mathcal{O}_x})[-2 \dim Y_{\mathcal{O}_x}] \;,$$

and similarly for $\mathcal{A}(Y)_I$, though this can not be formulated carefully at the current level of generality, as we explain following Remark.

*Remark* 12.1.10. The appearence of the ill-defined infinite cohomological degree shift in the above expression makes apparent that a more careful treatment of the theory of $D$ modules in the preceding (and following) exposition of this subsection is necessary. In the following subsections, we will restrict the class of prestacks under consideration to $Y = N/G$, the quotient by a reductive, affine algebraic group $G$ of a finite dimensional $G$ representation $N$. We then recall an explicit construction of the above putative factorization $\mathbb{E}_1$ algebra in this case, following [BFN18].

## 12.2. Internal construction of the three dimensional A model.
There is also an internal analogue of the above construction, following a variant of Example 8.1.15 as in Remark 8.1.16, as well as [BFN19b] and references therein:

*Example* 12.2.1. Following Proposition 8.1.1, the self fibre product

$$\mathcal{Z}(Y)_{\mathrm{dR}} = \mathcal{J}(Y)_{\mathrm{dR}} \times_{\mathcal{J}^{\mathrm{mer}}(Y)_{\mathrm{dR}}} \mathcal{J}(Y)_{\mathrm{dR}} \in \mathrm{PreStk}^{\mathrm{fact}}_{\mathrm{un}}(X)$$

defines a unital factorization space over $X$.

*Remark* 12.2.2. Concretely, the prestacks $\mathcal{Z}(Y)_{\mathrm{dR},I} \in \mathrm{PreStk}_{/X^I}$ over $X^I$ assigned to each $I \in \mathrm{fSet}$ and the fibre $\mathcal{Z}_{\mathrm{dR},x} \in \mathrm{PreStk}$ over $x \in X$ are given by

$$\mathcal{Z}(Y)_{\mathrm{dR},I} = \mathcal{J}(Y)_{I,\mathrm{dR}} \times_{\mathcal{J}^{\mathrm{mer}}(Y)_{I,\mathrm{dR}}} \mathcal{J}(Y)_{I,\mathrm{dR}} \qquad \text{and} \qquad \mathcal{Z}(Y)_{\mathrm{dR},x} = Y_{\mathcal{O}_x,\mathrm{dR}} \times_{Y_{\mathcal{K}_x,\mathrm{dR}}} Y_{\mathcal{O}_x\mathrm{dR}} \ .$$

*Remark* 12.2.3. More generally, the iterated fibre products

$$\mathcal{Z}(Y)_{(n),\mathrm{dR}} := \mathcal{J}(Y)_{\mathrm{dR}} \times_{\mathcal{J}^{\mathrm{mer}}(Y)_{\mathrm{dR}}} \times \ldots \times_{\mathcal{J}^{\mathrm{mer}}(Y)_{\mathrm{dR}}} \mathcal{J}(Y)_{\mathrm{dR}} \in \mathrm{PreStk}_{\mathrm{un}}^{\mathrm{fact}}(X)$$

and projections $\pi_{ij} : \mathcal{Z}(Y)_{(n),\mathrm{dR}} \to \mathcal{Z}(Y)_{\mathrm{dR}}$, define factorization spaces and maps of such for each $n \in \mathbb{N}$ and $i,j \in \{1,...,n\}$.

Following Proposition 8.1.5, we have:

*Proposition* 12.2.4. The factorization category

$$D_{\mathcal{Z}(Y)}^{\star} := \mathrm{QCoh}_{\mathcal{Z}(Y)_{\mathrm{dR}}}^{\star} \in \mathrm{Cat}_{\mathbb{E}_1,\mathrm{un}}^{\mathrm{fact}}(X_{\mathrm{dR}})$$

is naturally a (unital) $\mathbb{E}_1$-factorization category with respect to the convolution monoidal structure defined by the composition

$$D_{\mathcal{Z}(Y)} \otimes^{\star} D_{\mathcal{Z}(Y)} \xrightarrow{\pi_{12}^{*} \boxtimes \pi_{23}^{*}} D_{\mathcal{Z}_{(3)}(Y)^{\times 2}} \xrightarrow{\Delta^{*}} D_{\mathcal{Z}(Y)_{(3)}} \xrightarrow{\pi_{13,*}} D_{\mathcal{Z}(Y)} \ .$$

Further, the pushforward functor $\mathrm{p}_{\mathcal{Z}(Y)*} : D_{\mathcal{Z}(Y)}^{\star} \to D_{\mathrm{Ran}_{X,\mathrm{un}}}^{\otimes^!}$ is a unital, lax $\mathbb{E}_1$-monoidal factorization functor.

Further, following Corollary 8.1.6, we have:

*Corollary* 12.2.5. The pushforward functor $\mathrm{p}_{\mathcal{Z}*} : D_{\mathcal{Z}(Y)}^{\star} \to D(\mathrm{Ran}_{X,\mathrm{un}})$ induces a functor

(12.2.1)                    $$\mathrm{p}_* : \mathrm{Alg}_{\mathbb{E}_1,\mathrm{un}}^{\mathrm{fact}}(D_{\mathcal{Z}(Y)}) \to \mathrm{Alg}_{\mathbb{E}_1,\mathrm{un}}^{\mathrm{fact}}(X) \ .$$

*Example* 12.2.6. Let $G$ be a reductive algebraic group and $Y = BG$ be its classifying stack. Then the factorization space

$$\mathcal{Z}(Y)_{\mathrm{dR}} \cong \mathcal{J}(G)_{\mathrm{dR}} \backslash \mathrm{Gr}_{G,\mathrm{Ran}_{X,\mathrm{un}},\mathrm{dR}}$$

is isomorphic to the quotient of the de Rham stack of the Beilinson-Drinfeld Grassmannian, as in Example 4.1.9, by the factorization group stack $\mathcal{J}(G)_{\mathrm{dR}}$ defined by the de Rham stack of the space of jets to $G$, as in Example 4.1.11. The resulting $\mathbb{E}_1$ factorization category

$$D_{\mathcal{Z}(Y)} \cong D_{\mathcal{J}(G) \backslash \mathrm{Gr}_{G,\mathrm{Ran}_{X,\mathrm{un}}}} \in \mathrm{Cat}_{\mathbb{E}_1,\mathrm{un}}^{\mathrm{fact}}(X)$$

is (a DG enhancement of) the derived Satake category considered in [BF08]; concretely, its fibre category over $x \in X$ is given by

$$D_{\mathcal{Z}(Y),x} = D_{G_{\mathcal{O}_x}}(\mathrm{Gr}_{G,x}) \ .$$

Following Example 8.1.15 and in turn Example A.6.7, we define the internal variant of the three dimensional A model factorization $\mathbb{E}_1$ algebra:

*Tentative Definition* 12.2.7. The internal variant of the three dimensional A model to $Y$ is the factorization $\mathbb{E}_1$ algebra

$$\tilde{\mathcal{A}}(Y) = \mathcal{H}\tilde{\mathrm{om}}_{D_{\mathcal{J}^{\mathrm{mer}}(Y)}}(\iota_{\mathcal{J}(Y)*}\omega_{\mathcal{J}(Y)}, \iota_{\mathcal{J}(Y)*}\omega_{\mathcal{J}(Y)}) \cong \omega_{\mathcal{Z}(Y)}[-2\dim \mathcal{J}(Y)] \in \mathrm{Alg}_{\mathbb{E}_1,\mathrm{un}}^{\mathrm{fact}}(D_{\mathcal{Z}(Y)}^{\star})$$

internal to $D_{\mathcal{Z}(Y)}^{\star} \in \mathrm{Cat}_{\mathbb{E}_1,\mathrm{un}}^{\mathrm{fact}}(X_{\mathrm{dR}})$.

*Remark* 12.2.8. Concretely, in terms of the description of Remark 8.1.7, $\tilde{\mathcal{A}}(Y) \in \mathrm{Alg}^{\mathrm{fact}}_{\mathbb{E}_1, \mathrm{un}}(D_{\mathcal{Z}(Y)})$ has underlying $\mathbb{E}_1$ algebra internal to $D(\mathcal{Z}(Y)_I)$ for each $I \in \mathrm{fSet}_{\varnothing}$, and internal to $D(\mathcal{Z}(Y)_x)$ for each $x \in X$, given by

$$\tilde{\mathcal{A}}(Y)_I = \tilde{\mathrm{Hom}}_{D(\mathcal{J}^{\mathrm{mer}}(Y)_I)}(\iota_{\mathcal{J}(Y)_I *}\omega_{\mathcal{J}(Y)_I}, \iota_{\mathcal{J}(Y)_I *}\omega_{\mathcal{J}(Y)_I}) \cong \omega_{\mathcal{Z}(Y)_I}[-2\dim \mathcal{J}(Y)_I] \quad \in \mathrm{Alg}_{\mathbb{E}_1}(D(\mathcal{Z}(Y)_I)) \ , \text{ and}$$

$$\tilde{\mathcal{A}}(Y)_x = \tilde{\mathrm{Hom}}_{D(Y_{\mathcal{K}_x})}(\iota_{x *}\omega_{Y_{\mathcal{O}_x}}, \iota_{x *}\omega_{Y_{\mathcal{O}_x}}) \cong \omega_{\mathcal{Z}(Y)_x}[-2\dim Y_{\mathcal{O}_x}] \qquad\qquad \in \mathrm{Alg}_{\mathbb{E}_1}(D(\mathcal{Z}(Y)_x)).$$

Following Example 8.1.15 and Corollary 8.1.6, we have

*Proposition* 12.2.9. The image of $\tilde{\mathcal{A}}(Y) \in \mathrm{Alg}^{\mathrm{fact}}_{\mathbb{E}_1, \mathrm{un}}(D_{\mathcal{Z}(Y)})$ under the functor of Equation 12.2.1 above is canonically equivalent to $\mathcal{A}(Y) \in \mathrm{Alg}^{\mathrm{fact}}_{\mathbb{E}_1, \mathrm{un}}(X)$.

## 12.3. **Equivariant structures on the three dimensional A model.** In this subsection, we outline the application of the results of Section I-25 to the three dimensional A model.

*Proposition* 12.3.1. The (unital) factorization $\mathbb{E}_1$ algebra of the three dimensional A model to $Y$ over $X = \mathbb{A}^1$

$$\mathcal{A}(Y) \in \mathrm{Alg}^{\mathrm{fact}}_{\mathbb{E}_1, \mathrm{un}}(\mathbb{A}^1)^{\mathbb{G}_a \rtimes \mathbb{G}_m}$$

admits a canonical $\mathbb{G}_a \rtimes \mathbb{G}_m$ equivariant structure.

*Remark* 12.3.2. This is proved in Subsection 13.5 in the generality $Y = N/G$ as in Remark 12.1.10, again following [BFN18].

Thus, applying the equivalence of Proposition I-24.0.1, and following the discussion in section I-25 and in particular Proposition I-25.1.1 and Example I-25.1.5, we have:

*Corollary* 12.3.3. The three dimensional A model $\mathcal{A}(Y) \in \mathrm{Alg}^{S^1}_{\mathbb{E}_3}(\mathrm{Vect}_{\mathbb{K}})$ to $Y$ over $X = \mathbb{A}^1$ defines

$$\mathcal{A}_u \in \mathrm{Alg}_{\mathbb{BD}^u_1}(\mathrm{D}^b(\mathbb{K}[u])) \ ,$$

an algebra over the operad $\mathbb{BD}^u_1 \in \mathrm{Op}(\mathrm{D}^b(\mathbb{K}[u]))$ of Definition I-23.0.5.

*Remark* 12.3.4. Concretely, $\mathcal{A}_u \in \mathrm{Alg}_{\mathbb{BD}^u_1}(\mathrm{D}^b(\mathbb{K}[u]))$ defines a two-periodic filtered quantization of the homology $\mathbb{P}_3$ algebra $H_\bullet(\mathcal{A}) \in \mathrm{Alg}_{\mathbb{P}_3}(\mathrm{Perf}_{\mathbb{K}})$ to an associative (or $\mathbb{E}_1$) algebra $\mathcal{A}_u|_{\{1\}} \in \mathrm{Alg}_{\mathbb{E}_1}(\mathrm{Perf}_{\mathbb{K}})$. This is the data which more evidently matches the description of the construction given in [BFN18].

## 13. The Braverman-Finkelberg-Nakajima construction

In this section, we give a careful account of the construction of the three dimensional A model in the preceding section 12 in the special case when the space $Y = N/G$ is a quotient stack, again following [BFN18] throughout; this corresponds physically to the case of three dimensional $\mathcal{N} = 4$ A-twisted gauge theory.

In Section 13.1, we recall a particularly well behaved presentation of the relevant infinite dimensional stack $Y_{\mathcal{O}} \times_{Y_{\mathcal{K}}} Y_{\mathcal{O}}$, and in Section 13.2 we use this presentation to define a theory of $D$ modules on this space, following [Ras15b], and in particular its renormalized Borel-Moore homology groups. In Section 13.3, we explain some concrete calculations of the underlying vector spaces of these Borel-Moore homology groups, and in Section 13.4 we explain that they naturally carry a (factorization) associative (or $\mathbb{E}_1$) algebra structure coming from the convolution monoid structure on the underlying space. Finally, in Section 13.5, we recall that this factorization algebra admits a canonical $\mathbb{G}_a \rtimes \mathbb{G}_m$ equivariant structure, and thus defines a filtered quantization by the results of Section I-25, recovering the more familiar description of the Coulomb branch algebra from [BFN18].

13.1. **The variety of triples.** Throughout the remainder of this section, let $G$ be a finite-type, affine, reductive algebraic group, $N$ a finite dimensional $G$ representation, and $Y = N/G$ be the quotient stack. In this subsection, we recall a particular presentation of the factorization space $\mathcal{Z}(Y)$ (and in turn $\mathcal{Z}(Y)_{\mathrm{dR}}$) in this generality, towards accomplishing the narrative of the preceding subsection, as explained in Remark 12.1.10. Again, we follow [BFN18] and [BFN19b] throughout.

To begin, we define a factorization space $\mathcal{T}(G, N)$, which is a straightforward simeltaneous generalization of the Beilinson-Drinfeld Grassmannian $\mathrm{Gr}_{G, \mathrm{Ran}_{X,\mathrm{un}}}$ of Example 4.1.9 and the jet space $\mathcal{J}(Y)$ of Example 4.1.11:

*Definition* 13.1.1. The variety of pre-triples is the factorization space $\mathcal{T}(G, N) \in \mathrm{PreStk}^{\mathrm{fact}}_{\mathrm{un}}(X)$ defined by

$$\mathcal{T}(G, N)_I(S) = \{x : S \to X^I, P \in \mathrm{Bun}_G(S \times X), \varphi : P|_{S \times X \backslash \Gamma_x} \xrightarrow{\cong} P^{\mathrm{triv}}, \eta \in \Gamma(\mathbb{D}_x, P \times_G N)\}$$

with the evident factorization structure maps inherited from $\mathrm{Gr}_{G, \mathrm{Ran}_{X,\mathrm{un}}}$ and $\mathcal{J}(N)$, and unital structure defined, for simplicity for $\pi : I \hookrightarrow J$ injective, by the correspondence

$$\mathrm{unit}^\pi_{\mathcal{T}(G,N)}(S) = \{x = (x_1, x_2) : S \to X^J, P \in \mathrm{Bun}_G(S \times X), \varphi : P|_{S \times X \backslash \Gamma_{x_2}} \xrightarrow{\cong} P^{\mathrm{triv}}, a \in \Gamma(\mathbb{D}_x, P \times_G N)\}$$

together with the maps

$$\mathrm{unit}^\pi_{\mathcal{T}(G,N)}$$

$$X^{I_\pi} \times \mathcal{T}(G, N)_I \qquad\qquad\qquad\qquad \mathcal{T}(G, N)_J$$

defined by

$$(x = (x_1, x_2), P, \varphi : P|_{S \times X \backslash \Gamma_{x_2}} \xrightarrow{\cong} P^{\mathrm{triv}}, a : \mathbb{D}_x \to Y) \mapsto (x_1, (x_2, P, \varphi, a|_{\mathbb{D}_{x_2}} : \mathbb{D}_{x_2} \to Y)) \qquad , \text{and}$$

$$(x = (x_1, x_2), P, \varphi : P|_{S \times X \backslash \Gamma_{x_2}} \xrightarrow{\cong} P^{\mathrm{triv}}, a : \mathbb{D}_x \to Y) \mapsto (x, P, \varphi|_{S \times X \backslash \Gamma_x} : P|_{S \times X \backslash \Gamma_x} \xrightarrow{\cong} P^{\mathrm{triv}}, a) \qquad ,$$

respectfully, where the notation is as in Remark 4.1.8 and examples 4.1.9 and 4.1.11.

*Remark* 13.1.2. Concretely, the fibre of the factorization space $\mathcal{T}(G, N) \in \mathrm{PreStk}^{\mathrm{fact}}_{\mathrm{un}}(X)$ over each $x \in X$ is given by

$$\mathcal{T}(G, N)_x = (G_{\mathcal{K}_x} \times N_{\mathcal{O}_x})/G_{\mathcal{O}_x} .$$

*Warning* 13.1.3. We use the notation $\mathcal{T}(G, N)$ to denote the factorization space, in contrast with the conventions of [BFN18] where this notation is used for the fibre space $\mathcal{T}(G, N)_x$.

*Proposition* 13.1.4. There are natural maps of factorization spaces

$$\mathcal{T}(G, N) \to \mathrm{Gr}_{G, \mathrm{Ran}_{X,\mathrm{un}}} \qquad \mathcal{T}(G, N) \to \mathcal{J}(N/G) \qquad \text{and} \qquad \mathcal{T}(G, N) \to \mathcal{J}^{\mathrm{mer}}(N)$$

given by forgetting the section $\eta$, forgetting the trivialization $\varphi$ and restricting the bundle $P$ to $\mathbb{D}_x$, and applying the trivialization $\varphi$ to the section $\eta$, respectively.

*Remark* 13.1.5. Concretely, the maps on prestacks over $X^I$, and on the fibre over each $x \in X$, are given by

$$\mathcal{T}(G, N)_I \to \mathrm{Gr}_{G, I} \qquad \text{and} \quad (G_{\mathcal{K}} \times N_{\mathcal{O}})/G_{\mathcal{O}} \to G_{\mathcal{K}}/G_{\mathcal{O}} \qquad \text{defined by} \quad (P, \varphi, \eta) \mapsto (P, \varphi)$$

$$\mathcal{T}(G, N)_I \to \mathcal{J}(N/G)_I \quad \text{and} \quad (G_{\mathcal{K}} \times N_{\mathcal{O}})/G_{\mathcal{O}} \to N_{\mathcal{O}}/G_{\mathcal{O}} \qquad \text{defined by} \quad (P, \varphi, \eta) \mapsto (P|_{\mathbb{D}}, \eta)$$

$$\mathcal{T}(G, N)_I \to \mathcal{J}^{\mathrm{mer}}(N)_I \quad \text{and} \quad (G_{\mathcal{K}} \times N_{\mathcal{O}})/G_{\mathcal{O}} \to N_{\mathcal{K}} \qquad \text{defined by} \quad (P, \varphi, \eta) \mapsto \varphi \cdot \eta .$$

We now describe the primary factorization space of interest

*Definition* 13.1.6. The variety of triples, or BFN space, is the factorization space

$$\mathcal{R}(G,N) \in \mathrm{PreStk}_{\mathrm{un}}^{\mathrm{fact}}(X) \qquad \text{defined by} \qquad \mathcal{R}(G,N) := \mathcal{T}(G,N) \times_{\mathcal{J}^{\mathrm{mer}}(N)} \mathcal{J}(N) \ ,$$

as in the construction of fibre product factorization spaces recalled in Example 8.2.1.

*Remark* 13.1.7. Concretely, the prestacks over $X^I$ for each $I \in \mathrm{fSet}_\emptyset$, and the fibre over each $x \in X$, are given by

$$\mathcal{R}(G,N)_I = \mathcal{T}(G,N)_I \times_{\mathcal{J}^{\mathrm{mer}}(N)_I} \mathcal{J}(N)_I \cong \{(x,P,\varphi,\eta) \in \mathcal{T}(G,N)_{X^I} \mid \varphi \cdot \eta \in \mathcal{J}(N)_I \subset \mathcal{J}^{\mathrm{mer}}(N)_I\} \ , \text{ and}$$

$$\mathcal{R}(G,N)_x = (G_{\mathcal{K}_x} \times N_{\mathcal{O}_x})/G_{\mathcal{O}_x} \times_{N_{\mathcal{K}_x}} N_{\mathcal{O}_x} \cong \{(P,\varphi,\eta) \in \mathcal{T}(G,N)_x \mid \varphi \cdot \eta \in N_{\mathcal{O}_x} \subset N_{\mathcal{K}_x}\} \ .$$

*Example* 13.1.8. The unital factorization space $\mathcal{Z}(G,N)_{\mathrm{dR}} \in \mathrm{PreStk}_{\mathrm{un}}^{\mathrm{fact}}(X)$ is defined by $\mathcal{Z}(G,N)_{\mathrm{dR}} = \mathcal{Z}(N/G)_{\mathrm{dR}}$ where the latter is as defined in Example 12.2.1 in the case $Y = N/G$ is a quotient stack, and similarly for $\mathcal{Z}(G,N) \in \mathrm{PreStk}_{\mathrm{un}}^{\mathrm{fact}}(X)$ without passing to de Rham stacks.

The motivation for the preceding construction is the following:

*Proposition* 13.1.9. There are equivalences of factorization spaces

$$\mathcal{Z}(G,N) = \mathcal{J}(N/G) \times_{\mathcal{J}^{\mathrm{mer}}(N/G)} \mathcal{J}(N/G) \cong (\mathcal{T}(G,N) \times_{\mathcal{J}^{\mathrm{mer}}(N)} \mathcal{J}(N))/\mathcal{J}(G) = \mathcal{R}(G,N)/\mathcal{J}(G) \ .$$

*Proof.* The statement over each point $x \in X$ given in the following remark is immidiate, and the compatibility of the factorization data is a direct check after unpacking the definitions, following the proof of Proposition 8.1.1, Remark 8.1.2, and Example 8.2.1. $\square$

*Remark* 13.1.10. Concretely, over each point $x \in X$ the preceding equivalences are given by

$$\mathcal{Z}(G,N)_x = N_{\mathcal{O}}/G_{\mathcal{O}} \underset{N_{\mathcal{K}}/G_{\mathcal{K}}}{\times} N_{\mathcal{O}}/G_{\mathcal{O}} \cong ([(N_{\mathcal{O}} \times G_{\mathcal{K}})/G_{\mathcal{O}}] \times_{N_{\mathcal{K}}} N_{\mathcal{O}})/G_{\mathcal{O}} = \mathcal{R}(G,N)_x/G_{\mathcal{O}}$$

The main application of this description will be to give a careful construction of the three dimensional A model in this setting, as follows:

*Remark* 13.1.11. Following Remark 12.1.9, our tentative definition of the three dimensional A model is given by

$$\mathcal{A}(G,N) = C_\bullet^{\mathrm{BM}}(\mathcal{Z}(G,N))[-2\dim N_{\mathcal{O}}/G_{\mathcal{O}}] \cong C_\bullet^{G_{\mathcal{O}},\mathrm{BM}}(\mathcal{R}(G,N))[-2\dim N_{\mathcal{O}}] \quad \in \mathrm{Alg}_{\mathbb{E}_1}^{\mathrm{fact}}(X) \ .$$

Thus, we would like to give a careful definition of the latter Borel-Moore chains, and prove that is carries a (factorization) $\mathbb{E}_1$-algebra structure.

*Remark* 13.1.12. More generally, following Proposition 12.2.4 and Corollary 13.4.10, and Definition 12.2.7 and Remark 12.2.8, we would like to prove that

$$D_{\mathcal{Z}(G,N)} = D_{\mathcal{R}(G,N)/\mathcal{J}(G)} \in \mathrm{Cat}_{\mathbb{E}_1,\mathrm{un}}^{\mathrm{fact}}(X) \qquad \text{and} \qquad \mathrm{p}_{\mathcal{Z}*} : D_{\mathcal{R}(G,N)/\mathcal{J}(G)} \to D_{\mathrm{Ran}_{X,\mathrm{un}}}^{\otimes !}$$

define a factorization $\mathbb{E}_1$-category and a lax $\mathbb{E}_1$-monoidal factorization functor, such that

$$\tilde{\mathcal{A}}(G,N) := \omega_{\mathcal{Z}(G,N)}[-2\dim N_{\mathcal{O}}/G_{\mathcal{O}}] \in \mathrm{Alg}_{\mathbb{E}_1}^{\mathrm{fact}}(D_{\mathcal{Z}(G,N)}) \qquad \text{and} \qquad \mathrm{p}_{\mathcal{Z}*}\tilde{\mathcal{A}}(G,N) = \mathcal{A}(G,N) \quad \in \mathrm{Alg}_{\mathbb{E}_1}^{\mathrm{fact}}(X) \ .$$

13.2. **Borel-Moore Homology of the BFN space.** Following the final remarks of the preceding subsection, we now proceed to define the renormalized, $G_{\mathcal{O}}$-equivariant Borel-Moore homology of the BFN space $\mathcal{R}(G, N)$, following [BFN18]. For simplicity, we will restrict our attention to the space $\mathcal{R} = \mathcal{R}(G, N)_x$, having outlined the factorization compatibility of the constructions in Subsection 12.1.

*Example* 13.2.1. Recall that $\mathrm{Gr}_G = G_{\mathcal{K}}/G_{\mathcal{O}}$ admits a left action of $G_{\mathcal{O}}$, which induces a stratification

$$\mathrm{Gr}_G = \sqcup_{\lambda \in \Lambda_G^+} \mathrm{Gr}_G^\lambda \qquad \text{where} \qquad \mathrm{Gr}_G^\lambda = G_{\mathcal{O}} \cdot t^\lambda \,,$$

for $\Lambda_G^+$ is the set of dominant coweights of $G$ and $t^\lambda \in \mathrm{Gr}_G(\mathbb{C})$ defined by the restriction of the corresponding cocharacter $\lambda : \mathbb{C}^\times \to G$ to the formal neighbourhood of $\mathbf{0} \in \mathbb{C} \supset \mathbb{C}^\times$.

Further, each $\mathrm{Gr}_G^\lambda$ is a finite type scheme, and the $\mathbb{C}_h^\times$ action on $\mathrm{Gr}_G$ by 'loop rotation' contracts $\mathrm{Gr}_G^\lambda \twoheadrightarrow G \cdot t^\lambda = G/P_\mu$ where $P_\mu \subset G$ is the corresponding parabolic subgroup, with Levi quotient $\mathrm{Stab}_G(\lambda)$. Thus, $\mathrm{Gr}_G^\lambda$ is a finite dimensional vector bundle over $G/P^\lambda$, as the attracting set fibres are equidimensional affine spaces.

Moreover, recall that

$$\mathrm{Gr}_G^{\leqslant \lambda} := \overline{\mathrm{Gr}}_G^\lambda = \sqcup_{\mu \leqslant \lambda} \mathrm{Gr}_G^\mu \qquad \text{so that} \qquad \mathrm{Gr}_G = \mathrm{colim}_{\lambda \in \Lambda_G^+} \mathrm{Gr}_G^{\leqslant \lambda} \,.$$

Further, each $\mathrm{Gr}_G^{\leqslant \lambda}$ is a projective scheme of dimension $\langle 2\rho^\vee, \lambda \rangle$, so that $\mathrm{Gr}_G$ is an ind-projective, ind-finite type indscheme.

*Definition* 13.2.2. Let $\mathcal{T} = \mathcal{T}(G, N)_x$ and $\mathcal{R} = \mathcal{R}(G, N)_x$ be the fibres of the factorization spaces of pre-triples and triples as in Definitions 13.1.1 and 13.1.6, respectively, and let $\pi : \mathcal{T} \to \mathrm{Gr}_G$ denote the map from Remark 13.1.4 and $\pi' : \mathcal{R} \to \mathrm{Gr}_G$ its restriction to $\mathcal{R}(G, N)$. Let $\mathcal{T}^{\leqslant \lambda}$ and $\mathcal{R}^{\leqslant \lambda}$ denote the pullbacks

$$
\begin{array}{ccc}
\mathcal{T}^{\leqslant \lambda} & \longrightarrow & \mathcal{T} \\
\downarrow & & \downarrow{\scriptstyle \pi} \\
\mathrm{Gr}_G^{\leqslant \lambda} & \longrightarrow & \mathrm{Gr}_G
\end{array}
\qquad \text{and} \qquad
\begin{array}{ccc}
\mathcal{R}^{\leqslant \lambda} & \longrightarrow & \mathcal{R} \\
\downarrow & & \downarrow{\scriptstyle \pi'} \\
\mathrm{Gr}_G^{\leqslant \lambda} & \longrightarrow & \mathrm{Gr}_G
\end{array}
\,,
$$

for each $\lambda \in \Lambda_G^+$ a dominant coweight, and similarly $\mathcal{T}^\lambda, \mathcal{R}^\lambda$ and $\mathcal{T}^{<\lambda}, \mathcal{R}^{<\lambda}$ as the analogous fibre products with $\mathrm{Gr}_G^\lambda$ and $\mathrm{Gr}_G^{<\lambda} = \sqcup_{\mu < \lambda} \mathrm{Gr}_G^\mu$, respectively.

*Example* 13.2.3. The space $\mathcal{T}^{\leqslant \lambda} \twoheadrightarrow \mathrm{Gr}_G^{\leqslant \lambda}$ is an infinite dimensional vector bundle with fibre $N_{\mathcal{O}}$, a pro-finite type scheme presented by $N_{\mathcal{O}} = \lim_n N_{\mathcal{O}}/(z^n \cdot N_{\mathcal{O}})$ as in Example B.2.7. Thus, $\mathcal{T}^{\leqslant \lambda}$ is also pro-finite type scheme, presented by the limit of the corresponding finite dimensional vector bundles over the finite type scheme $\mathrm{Gr}_G^{\leqslant \lambda}$:

$$\mathcal{T}^{\leqslant \lambda} = \lim_{i \in \mathcal{I}} \mathcal{T}_i^{\leqslant \lambda} = \lim \left[ \ldots \to \mathcal{T}_i^{\leqslant \lambda} \xrightarrow{\varphi_{ij}^\lambda} \mathcal{T}_j^{\leqslant \lambda} \to \ldots \to \mathcal{T}_0^{\leqslant \lambda} = \mathrm{Gr}_G^{\leqslant \lambda} \right] \qquad \mathcal{T}_i^{\leqslant \lambda} = ((N_{\mathcal{O}}/z^i N_{\mathcal{O}}) \times G_{\mathcal{K}}^{\leqslant \lambda})/G_{\mathcal{O}} \,.$$

Similarly, $\mathcal{R}^{\leqslant \lambda} \hookrightarrow \mathcal{T}^{\leqslant \lambda}$ is a pro-finite type subscheme of finite codimension, and presented analogously as

$$\mathcal{R}^{\leqslant \lambda} = \lim_{i \in \mathcal{I}} \mathcal{R}_i^{\leqslant \lambda} = \lim \left[ \ldots \to \mathcal{R}_i^{\leqslant \lambda} \xrightarrow{\varphi_{ij}^\lambda} \mathcal{R}_j^{\leqslant \lambda} \to \ldots \to \mathcal{R}_0^{\leqslant \lambda} = \mathrm{Gr}_G^{\leqslant \lambda} \right] \qquad \mathcal{R}_i^{\leqslant \lambda} = (((N_{\mathcal{O}} \cap z^\lambda N_{\mathcal{O}})/z^i N_{\mathcal{O}}) \times G_{\mathcal{K}}^{\leqslant \lambda})/G_{\mathcal{O}} \,.$$

Note that for $i < d_\lambda = \mathrm{rank}(\mathcal{T}^\lambda/\mathcal{R}^\lambda)$, the spaces $\mathcal{R}_i^{\leqslant \lambda} = \mathcal{R}_0^{\leqslant \lambda} = \mathrm{Gr}_G^{\leqslant \lambda}$ are equivalent, and the transition maps are all given by the identity.

Moreover, the maps $\varphi_{ij}^\lambda$ are each smooth with affine space fibres, and thus define placid presentations of $\mathcal{R}^{\leqslant\lambda}$ and $\mathcal{T}^{\leqslant\lambda}$, in the sense if Definition B.7.9.

*Example* 13.2.4. Following Remark B.7.3 and Proposition B.7.11, the category of $D$ modules $D^*(\mathcal{R}^{\leqslant\lambda}) \in$ DGCat is presented as

$$\mathcal{D}^*(\mathcal{R}^{\leqslant\lambda}) = \lim_{i\in\mathcal{I}} D(\mathcal{R}_i^{\leqslant\lambda}) = \lim \left[ \ldots \to D(\mathcal{R}_i^{\leqslant\lambda}) \xrightarrow{\varphi_{ij*}^\lambda} D(\mathcal{R}_j^{\leqslant\lambda}) \to \ldots \to D(\mathrm{Gr}_G^{\leqslant\lambda}) \right] \text{ , and}$$

$$\cong \operatorname*{colim}_{i\in\mathcal{I}^{\mathrm{op}}} D(\mathcal{R}_i^{\leqslant\lambda}) = \operatorname{colim} \left[ \ldots \leftarrow D(\mathcal{R}_i^{\leqslant\lambda}) \xleftarrow{\varphi_{ij}^{\lambda*}} D(\mathcal{R}_j^{\leqslant\lambda}) \leftarrow \ldots \leftarrow D(\mathrm{Gr}_G^{\leqslant\lambda}) \right] ,$$

where the latter description follows from passing to left adjoints, which by placidity exist and are given by

$$\varphi_{ij}^{\lambda*} \cong \varphi_{ij}^{\lambda!}[-2d_{\mathcal{R}_i^{\leqslant\lambda}/\mathcal{R}_j^{\leqslant\lambda}}] : D(\mathcal{R}_j^{\leqslant\lambda}) \to D(\mathcal{R}_i^{\leqslant\lambda}) .$$

*Remark* 13.2.5. Concretely, from the limit description, a $D$ module $M^\lambda \in \mathcal{D}^*(\mathcal{R}^{\leqslant\lambda})$ is given by an assignment

$$i \mapsto \left( M_i^\lambda \in \mathcal{D}^*(\mathcal{R}_i^{\leqslant\lambda}) \right) \qquad [i \to j] \mapsto \left[ \varphi_{ij*}^\lambda M_i^\lambda \xrightarrow{\cong} M_j^\lambda \right] .$$

*Remark* 13.2.6. Alternatively, from the colimit description, there are objects $M^\lambda \in \mathcal{D}^*(\mathcal{R}^{\leqslant\lambda})$ given by $M^\lambda = \varphi_{i_0}^* M_{i_0}^\lambda$ where $\varphi_{i_0}^* : \mathcal{D}^*(\mathcal{R}_{i_0}^{\leqslant\lambda}) \to \mathcal{D}^*(\mathcal{R}^{\leqslant\lambda})$ is the canonical functor of the colimit description, for some $i_0 \in \mathcal{I}$. In terms of the limit presentation above, the resulting object is defined by
(13.2.1)
$$M_{j_0}^\lambda = \operatorname*{colim}_{i\in\mathcal{I}^{\mathrm{op}}_{/j_0}} \varphi_{ij_0*}^\lambda \varphi_{ii_0}^{\lambda*} M_{i_0}^\lambda \qquad \text{under} \qquad \varphi_{jj_0*}^\lambda \varphi_{ji_0}^{\lambda*} M_{i_0}^\lambda \to \varphi_{jj_0*}^\lambda \varphi_{ij*}^\lambda \varphi_{ij}^{\lambda*} \varphi_{ji_0}^{\lambda*} M_{i_0}^\lambda \cong \varphi_{ij_0*}^\lambda \varphi_{ii_0}^{\lambda*} M_{i_0}^\lambda .$$

*Example* 13.2.7. Recall from Definition B.7.12 that there is a canonical renormalized dualizing sheaf $\omega_{\mathcal{R}^{\leqslant\lambda}}^{\mathrm{ren}} \in D^*(\mathcal{R}^{\leqslant\lambda})$ defined by

$$\omega_{\mathcal{R}^{\leqslant\lambda}}^{\mathrm{ren}} = \varphi_i^* \omega_{\mathcal{R}_i^{\leqslant\lambda}}[-2d_{\mathcal{R}_i^{\leqslant\lambda}}] \qquad \text{noting} \qquad \varphi_i^{\lambda*} \omega_{\mathcal{R}_i^{\leqslant\lambda}}[-2d_{\mathcal{R}_i^{\leqslant\lambda}}] = \varphi_i^{\lambda*} \varphi_{ij}^{\lambda*} \omega_{\mathcal{R}_j^{\leqslant\lambda}}[-2d_{\mathcal{R}_j^{\leqslant\lambda}}] = \varphi_j^{\lambda*} \omega_{\mathcal{R}_j^{\leqslant\lambda}}[-2d_{\mathcal{R}_j^{\leqslant\lambda}}] .$$

*Remark* 13.2.8. Heuristically, the renormalized dualizing sheaf can be written

$$\omega_{\mathcal{R}^{\leqslant\lambda}}^{\mathrm{ren}} = \omega_{\mathcal{R}^{\leqslant\lambda}}[-2\dim\mathcal{R}^{\leqslant\lambda}] = \omega_{\mathcal{R}^{\leqslant\lambda}}[-2(\dim N_{\mathcal{O}} + \dim\mathrm{Gr}^{\leqslant\lambda} - d_\lambda)] \qquad d_\lambda = \mathrm{rank}(\mathcal{T}^\lambda/\mathcal{R}^\lambda) : \mathcal{R} \to \mathbb{Z}.$$

*Example* 13.2.9. There is a canonical pushforward functor

$$\mathrm{p}_{\mathcal{R}^{\leqslant\lambda}*} : D^*(\mathcal{R}^{\leqslant\lambda}) \to \mathrm{Vect} \qquad \text{defined by} \qquad M^\lambda = (M_i^\lambda)_{i\in\mathcal{I}} \mapsto \mathrm{p}_{\mathcal{R}_i^\lambda*} M_i^\lambda ;$$

the result is independent of $i$ since

$$\mathrm{p}_{\mathcal{R}_i^{\leqslant\lambda}*} M_i^\lambda \cong \mathrm{p}_{\mathcal{R}_j^{\leqslant\lambda}*} \varphi_{ij*}^\lambda M_i^\lambda \cong \mathrm{p}_{\mathcal{R}_j^{\leqslant\lambda}*} M_j^\lambda .$$

*Example* 13.2.10. Following the preceding three examples, the renormalized Borel-Moore chains are given by

$$C_\bullet^{\mathrm{BM}}(\mathcal{R}^{\leqslant\lambda}) = C^\bullet(\mathcal{R}^{\leqslant\lambda}, \omega_{\mathcal{R}^{\leqslant\lambda}}^{\mathrm{ren}}) = \mathrm{p}_{\mathcal{R}^{\leqslant\lambda}*} \omega_{\mathcal{R}^{\leqslant\lambda}}^{\mathrm{ren}} = \mathrm{p}_{\mathcal{R}_j^{\leqslant\lambda}*} \mathcal{M}_j^\lambda \cong \operatorname*{colim}_{i\in\mathcal{I}^{\mathrm{op}}_{/j}} \mathrm{p}_{\mathcal{R}_i^{\leqslant\lambda}*} \varphi_{ij}^{\lambda*} \omega_{\mathcal{R}_j^{\leqslant\lambda}}[-2d_{\mathcal{R}_j^{\leqslant\lambda}}]$$

for any $j \in \mathcal{I}$.

*Example* 13.2.11. The spaces $\mathfrak{T}, \mathfrak{R} \in \mathrm{IndSch}$ are presented by

$$\mathfrak{T} = \operatorname*{colim}_{\lambda \in \Lambda_G^+} \mathfrak{T}^{\leqslant \lambda} = \operatorname{colim} \left[ \ldots \leftarrow \mathfrak{T}^{\leqslant \mu} \xleftarrow{\iota^{\lambda\mu}} \mathfrak{T}^{\leqslant \lambda} \leftarrow \ldots \leftarrow \mathfrak{T}^{\leqslant 0} = N_{\mathcal{O}} \right] \qquad , \text{ and}$$

$$\mathfrak{R} = \operatorname*{colim}_{\lambda \in \Lambda_G^+} \mathfrak{R}^{\leqslant \lambda} = \operatorname{colim} \left[ \ldots \leftarrow \mathfrak{R}^{\leqslant \mu} \xleftarrow{\iota^{\lambda\mu}} \mathfrak{R}^{\leqslant \lambda} \leftarrow \ldots \leftarrow \mathfrak{R}^{\leqslant 0} = N_{\mathcal{O}} \right] \qquad ,$$

exhibiting $\mathfrak{T}$ and $\mathfrak{R}$ as indschemes. Moreover, $\mathfrak{T}^{\leqslant \lambda}$ and $\mathfrak{R}^{\leqslant \lambda}$ are reasonable subschemes of $\mathfrak{T}$ and $\mathfrak{R}$, and thus present $\mathfrak{T}$ and $\mathfrak{R}$ as placid indschemes, in the sense of definitions B.3.7 and B.8.3, respectfully.

*Example* 13.2.12. Following propositions B.8.1 and B.8.2, the category of $D$ modules $D^*(\mathfrak{R})$ is presented as

$$D^*(\mathfrak{R}) = \operatorname*{colim}_{\lambda \in \Lambda_G^+} D^*(\mathfrak{R}^{\leqslant \lambda}) = \operatorname{colim} \left[ \ldots \leftarrow D^*(\mathfrak{R}^{\leqslant \mu}) \xleftarrow{\iota_*^{\lambda\mu}} D^*(\mathfrak{R}^{\leqslant \lambda}) \leftarrow \ldots \leftarrow D^*(N_{\mathcal{O}}) \right] \qquad , \text{ and}$$

$$\cong \operatorname*{lim}_{\lambda \in \Lambda_G^{+,\mathrm{op}}} D^*(\mathfrak{R}^{\leqslant \lambda}) = \operatorname{lim} \left[ \ldots \to D^*(\mathfrak{R}^{\leqslant \mu}) \xrightarrow{\iota^{\lambda\mu(!,*)}} D^*(\mathfrak{R}^{\leqslant \lambda}) \to \ldots \to D^*(N_{\mathcal{O}}) \right]$$

where the latter description follows from passing to right adjoints, which exist since $\iota^{\lambda\mu}$ are closed embeddings of finite presentation, as inclusions of reasonable subschemes.

*Remark* 13.2.13. Concretely, by the latter description, an object $M \in D^*(\mathfrak{R})$ is given by an assignment

$$\lambda \mapsto \left( M^\lambda \in D^*(\mathfrak{R}^{\leqslant \lambda}) \right) \qquad [\lambda \to \mu] \mapsto \left[ \iota^{\lambda\mu(!,*)} M^\mu \xrightarrow{\cong} M^\lambda \right] .$$

*Example* 13.2.14. There is a canonical pushforward functor

$$\mathrm{p}_{\mathfrak{R}*} : D^*(\mathfrak{R}) \to \mathrm{Vect} \qquad \text{defined by} \qquad M = (M^\lambda)_{\lambda \in \Lambda_G^+} \mapsto \operatorname*{colim}_{\lambda \in \Lambda_G^+} \mathrm{p}_{\mathfrak{R}^{\leqslant \lambda}*} M^\lambda$$

where the structure maps in the colimit diagram are defined for each $k < l$ by

$$\mathrm{p}_{\mathfrak{R}^{\leqslant \lambda}*} M^\lambda \cong \mathrm{p}_{\mathfrak{R}^{\leqslant \mu}*} \iota_*^{\lambda\mu} M^\lambda \cong \mathrm{p}_{\mathfrak{R}^{\leqslant \mu}*} \iota_*^{\lambda\mu} \iota^{\lambda\mu(!,*)} M^\mu \to \mathrm{p}_{\mathfrak{R}^{\leqslant \mu}*} M^\mu ,$$

where the map is given by the counit of the $(\iota_*^{\lambda\mu}, \iota^{\lambda\mu(!,*)})$ adjunction.

*Remark* 13.2.15. The placid scheme structure on each $\mathfrak{R}^{\leqslant \lambda}$ equips it with a renormalized dualizing sheaf $\omega_{\mathfrak{R}^{\leqslant \lambda}}^{\mathrm{ren}} \in D^*(\mathfrak{R}^{\leqslant \lambda})$ as in Example 13.2.7. However, as in Remark B.8.5, these objects do not glue together to define an object of $D(\mathfrak{R})$ as the cohomological shifts in the renormalized dualizing sheaves are not intertwined under $\iota^{\lambda\mu(!,*)}$ as required in the preceding remark; to arrange for this we must choose a dimension theory on $\mathfrak{R}^{\leqslant \lambda}$, in the sense of Definition B.8.6.

*Remark* 13.2.16. The affine Grassmannian $\mathrm{Gr}_G$ is ind-finite type, so that it has a canonical dimension theory defined by $\tau_{\mathrm{Gr}_G^{\leqslant \lambda}} = \dim \mathrm{Gr}_G^{\leqslant \lambda}$, as in Example B.8.8. The jet scheme $N_{\mathcal{O}}$ is placid, so carries a canonical dimension theory such that $\tau_{N_{\mathcal{O}}} = 0$, as in Example B.8.7.

As dimension theories are etale local, the product dimension theory construction of Example B.8.9 induces a dimension theory on the vector bundle $\mathfrak{T}$. Further, $\mathfrak{R}$ inherits the dimension theory of $\mathfrak{T}$, as every reasonable subscheme of the finite codimension sub indscheme $\mathfrak{R}$ is reasonable in $\mathfrak{T}$.

*Definition* 13.2.17. The canonical dimension theory $\tau$ on $\mathcal{R}$ is the dimension theory it inherits as a finite codimension sub indscheme of $\mathcal{T}$ shifted by the locally constant function $d_\lambda : \mathcal{R} \to \mathbb{Z}$ given by $d_\lambda = \operatorname{rank}(\mathcal{T}^\lambda/\mathcal{R}^\lambda)$.

The renormalized dualizing sheaf $\omega_\mathcal{R}^{\mathrm{ren}} := \omega_\mathcal{R}^{\tau\text{-ren}} \in D^*(\mathcal{R})$ of the variety of triples is defined as that corresponding to the canonical dimension theory $\tau$.

The Borel-Moore chains of $\mathcal{R}$, and its $G_\mathcal{O} \rtimes \mathbb{C}^\times$ equivariant analogue, are defined as

$$C_\bullet^{\mathrm{BM}}(\mathcal{R}) = C_{\mathrm{dR}}^\bullet(\mathcal{R}, \omega_\mathcal{R}^{\mathrm{ren}}) = \mathrm{p}_{\mathcal{R}*}\omega_\mathcal{R}^{\mathrm{ren}} \qquad \text{and} \qquad C_\bullet^{G_\mathcal{O} \rtimes \mathbb{C}_\hbar^\times}(\mathcal{R}) = C_{G_\mathcal{O} \rtimes \mathbb{C}_\hbar^\times}^\bullet(\mathcal{R}, \omega_\mathcal{R}^{\mathrm{ren}}) \ ,$$

noting that the dimension theory on $\mathcal{R}$ is $G_\mathcal{O} \rtimes \mathbb{C}_\hbar^\times$ equivariant, and thus the canonical equivariant structure on $\omega_\mathcal{R} \in D_{G_\mathcal{O} \rtimes \mathbb{C}_\hbar^\times}^!(\mathcal{R})$ induces an equivariant structure on $\omega_\mathcal{R}^{\mathrm{ren}} \in D_{G_\mathcal{O} \rtimes \mathbb{C}^\times}^*(\mathcal{R})$.

*Remark* 13.2.18. For notational simplicity, we will state many of the remaining results in terms of $C_\bullet^{G_\mathcal{O}}(\mathcal{R})$ with the deformation to $C_\bullet^{G_\mathcal{O} \rtimes \mathbb{C}_\hbar^\times}(\mathcal{R})$ left implicit in cases where it does not play an essential role.

*Remark* 13.2.19. Concretely, the canonical dimension theory on $\mathcal{R}$ assigns to each reasonable subscheme $\mathcal{R}^{\leqslant\lambda}$ the function $\tau_{\leqslant\lambda} = \dim \mathrm{Gr}^\lambda - d_\lambda$; note that this is evidently $G_\mathcal{O} \rtimes \mathbb{C}^\times$ equivariant. Thus, the renormalized dualizing sheaf $\omega_\mathcal{R}^{\mathrm{ren}} \in D^*(\mathcal{R})$ is defined by the condition that

$$\iota_{\leqslant\lambda}^{!*}\omega_\mathcal{R}^{\mathrm{ren}} = \omega_{\mathcal{R}^{\leqslant\lambda}}^{\mathrm{ren}}[2(\dim \mathrm{Gr}_G^\lambda - d_\lambda)] \qquad \text{which is given heuristically by} \qquad \text{``}\omega_{\mathcal{R}^{\leqslant\lambda}}[-2\dim N_\mathcal{O}]\text{''}$$

for each $\lambda \in \Lambda_G^+$, following Example 13.2.7 and Remark 13.2.8. These objects are intertwined under $\iota_{\lambda,\mu}^{!*}$ by construction, as in Definition B.8.10.

In particular, the (equivariant) Borel-Moore chains are computed by

$$C_\bullet^{G_\mathcal{O} \rtimes \mathbb{C}_\hbar^\times}(\mathcal{R}) = \operatorname*{colim}_{\lambda \in \Lambda_G^+} C_{G_\mathcal{O} \rtimes \mathbb{C}_\hbar^\times}^\bullet(\mathcal{R}^{\leqslant\lambda}, \omega_{\mathcal{R}^{\leqslant\lambda}}^{\mathrm{ren}})[2(\dim \mathrm{Gr}_G^\lambda - d_\lambda)]$$

$$= \operatorname*{colim}_{\lambda \in \Lambda_G^+} \operatorname*{colim}_{i \in \mathcal{I}_{/j}^{\mathrm{op}}} \mathrm{p}_{\mathcal{R}_i^{\leqslant\lambda}*}\varphi_{ij}^{\lambda*}\omega_{\mathcal{R}_j^{\leqslant\lambda}}[2(\dim \mathrm{Gr}_G^\lambda - d_\lambda) - 2d_{\mathcal{R}_j^{\leqslant\lambda}}] \quad \text{, for any } j \in \mathcal{I},$$

where the structure maps in the first line are given by the counit of the $(\iota_*^{\lambda\mu}, \iota^{\lambda\mu(!,*)})$ adjunction, as in Example 13.2.14 above, and the second line follows by Example 13.2.10.

*Remark* 13.2.20. Heuristically, the renormalized dualizing sheaf on $\mathcal{R}$ is well-defined because the shift $[-2\dim N_\mathcal{O}]$ is independent of $\lambda, \mu$ and $\iota_{\lambda,\mu}^!\omega_{\mathcal{R}^{\leqslant\lambda}} \cong \omega_{\mathcal{R}^{\leqslant\mu}}$, so that the induced limit object can be written as simply

$$\omega_\mathcal{R}^{\mathrm{ren}} = \omega_\mathcal{R}[-2\dim N_\mathcal{O}] \in D^*(\mathcal{R}) \ .$$

Moreover, the (equivariant) Borel-Moore homology is thus computed heuristically by

$$C_\bullet^{G_\mathcal{O} \rtimes \mathbb{C}_\hbar^\times}(\mathcal{R}) = \text{``} \operatorname*{colim}_{\lambda \in \Lambda_G^+} C_{G_\mathcal{O} \rtimes \mathbb{C}_\hbar^\times}^\bullet(\mathcal{R}^{\leqslant\lambda}, \omega_{\mathcal{R}^{\leqslant\lambda}})[-2\dim N_\mathcal{O}]\text{''} \ .$$

## 13.3. **Concrete calculations.**

In this section, we recall some concrete calculations of the resulting algebras, again following [BFN18].

*Proposition* 13.3.1. For each $\lambda \in \Lambda_+^G$, we have

$$H_\bullet^{G_\mathcal{O}}(\mathcal{R}^\lambda) = H_{\mathrm{Stab}_G(\lambda)}^\bullet(\mathrm{pt})[-2\dim \mathrm{Gr}_G^\lambda - 2d_\lambda] = \mathrm{Sym}^\bullet(\mathfrak{t}^\vee[-2])^{W_\lambda}[4\langle\rho, \lambda\rangle - 2d_\lambda]$$

where $\mathrm{Stab}_G(\lambda)$ denote the stabilizer of $\lambda$ in $G$, $W_\lambda$ denotes its Weyl group, and $\mathfrak{t}$ is a Cartan of $\mathfrak{g} = \mathrm{Lie}(G)$.

*Proof.* The smooth adjunction for the projection $\pi : \mathcal{R}^\lambda \to \mathrm{Gr}_G^\lambda$ induces an isomorphism

$$H_\bullet^{G_\mathcal{O}}(\mathcal{R}^\lambda) \xrightarrow{\cong} H_\bullet^{G_\mathcal{O}}(\mathrm{Gr}_G^\lambda)[-2d_\lambda] \ .$$

Similarly, the projection $\pi : \mathrm{Gr}_G^\lambda \to G \cdot t^\lambda \cong G/P_\lambda$ induces an isomorphism

$$H_\bullet^{G_\mathcal{O}}(\mathrm{Gr}_G^\lambda) \xrightarrow{\cong} H_\bullet^{G_\mathcal{O}}(G/P_\lambda)[2(\dim \mathrm{Gr}_G^\lambda - \dim G/P_\lambda)] \cong H_G^\bullet(G/P_\lambda)[2\dim \mathrm{Gr}_G^\lambda] \ ,$$

where $P_\lambda$ is the parabolic corresponding to $\lambda$ with Levi quotient $\mathrm{Stab}_G(\lambda)$. The latter is a standard computation

$$H_G^\bullet(G/P_\lambda) \cong H_{\mathrm{Stab}_G(\lambda)}^\bullet(\mathrm{pt}) = \mathrm{Sym}^\bullet(\mathfrak{t}^\vee[-2])^{W_\lambda} \ .$$

$\square$

*Corollary* 13.3.2. The equivariant Poincare polynomial of $\mathcal{R}^\lambda$ is given by

$$P_t^{G_\mathcal{O}}(\mathcal{R}^\lambda) = t^{2d_\lambda - 4\langle \rho, \lambda \rangle} \prod_i (1 - t^{2d_i})^{-1}$$

where the product runs over the degrees $d_i$ of the generators of the $W_\lambda$ invariant polynomials on $\mathfrak{t}[2]$.

*Proposition* 13.3.3. The complex $H_\bullet^{G_\mathcal{O}}(\mathcal{R}^{\leq \lambda})$ is concentrated in even degree. In particular, the long exact sequence in cohomology induced by the exact triangle of the complimentary closed and open embeddings

$$\mathcal{R}^{<\lambda} \xrightarrow{\iota} \mathcal{R}^{\leq \lambda} \xleftarrow{j} \mathcal{R}^\lambda \qquad \text{splits as} \qquad H_\bullet^{G_\mathcal{O}}(\mathcal{R}^{<\lambda}) \hookrightarrow H_\bullet^{G_\mathcal{O}}(\mathcal{R}^{\leq \lambda}) \twoheadrightarrow H_\bullet^{G_\mathcal{O}}(\mathcal{R}^\lambda) \ .$$

*Proof.* The proof is by induction on $\lambda$. The base case is when $\lambda$ is miniscule, and $\mathcal{R}^{\leq \lambda} = \mathcal{R}^\lambda$, in which case the result follows from the previous proposition. For the inductive step, we know additionally that $H_\bullet^{G_\mathcal{O}}(\mathcal{R}^{<\lambda})$ is concentrated in even degree, so that the long exact sequence implies the same for $H_\bullet^{G_\mathcal{O}}(\mathcal{R}^{\leq \lambda})$. $\square$

*Corollary* 13.3.4. As a vector space, the Borel-Moore homology is given by

$$H_\bullet^{G_\mathcal{O}}(\mathcal{R}) \cong \bigoplus_{\lambda \in \Lambda_G^+} \mathrm{Sym}^\bullet(\mathfrak{t}^\vee[-2])^{W_\lambda}[4\langle \rho, \lambda \rangle - 2d_\lambda] \ .$$

*Example* 13.3.5. Suppose $G = T = (\mathbb{C}^\times)^n$ is a rank $n$ torus and $N$ be arbitrary, so that $\Lambda_G^+ = \Lambda_G = \mathbb{Z}^n$, $\mathcal{R}^\lambda(\mathbb{C}) = \mathrm{Gr}_G^\lambda(\mathbb{C}) = \mathrm{pt}$ for each $\lambda \in \Lambda_G^+$, and thus $\mathrm{Gr}_G(\mathbb{C}) = \mathbb{Z}^n$. Then the above computation gives

$$H_\bullet^{G_\mathcal{O}}(\mathcal{R}) = \bigoplus_{\lambda \in \mathbb{Z}^n} \mathrm{Sym}^\bullet(\mathfrak{t}^\vee[-2])[-2d_\lambda] \cong \bigoplus_{\lambda \in \mathbb{Z}^n} \mathbb{C}[\xi_i]_{i=1}^n \cdot r^\lambda = \mathbb{C}[\xi_i, r_i^{\pm 1}]_{i=1}^n \ ,$$

where $\xi_i \in H_T^\bullet(\mathrm{pt})[-2]$ are the generators of the $T$ equivariant cohomology of a point and $r^\lambda := [\mathcal{R}^\lambda] \in H_{2d_\lambda}^{G_\mathcal{O}}(\mathcal{R}^\lambda)$ is the fundamental class of $\mathcal{R}^\lambda$.

For $N$ the zero representation $N = \{0\}$, this lifts to an isomorphism of algebras

$$H_\bullet^{G_\mathcal{O}}(\mathcal{R}) = \mathbb{C}[\xi_i, r_i^{\pm 1}]_{i=1}^n = \mathcal{O}(\mathfrak{t} \times T^\vee) \cong \mathcal{O}(T^* T^\vee)$$

with the coordinate ring of the cotangent bundle of the dual torus $T^\vee$.

In general, the product structure depends on the data of the representation; see Proposition 13.5.6 below, which again follows [BFN18].

*Example* 13.3.6. Suppose $G = \mathrm{PGL}_2$, so that $\Lambda_G^+ \cong \mathbb{N}$ and for each $\lambda = n \in \mathbb{N}_{\geqslant 0}$ we have $\langle \lambda, \rho \rangle = \frac{1}{2}n$. Then $\mathrm{Gr}_G^0 = \mathrm{pt}$, and for each $n \geqslant 0$ there is a fibre sequence

$$\mathbb{A}^{n-1} \hookrightarrow \mathrm{Gr}_G^n \twoheadrightarrow G/B = \mathbb{P}^1 \qquad \text{inducing} \qquad H_\bullet^{G_\mathcal{O}}(\mathrm{Gr}_G^n) \cong H_B^\bullet(\mathrm{pt})[2n] \cong \mathbb{C}[\xi] \cdot \eta^n$$

where $\eta := [\mathrm{Gr}_G^1] \in H_{-2}^{G_\mathcal{O}}(\mathrm{Gr}_G^1)$, more generally $\eta^n := [\mathrm{Gr}_G^n] \in H_{-2n}^{G_\mathcal{O}}(\mathrm{Gr}_G^n)$, and $\xi \in H_\bullet^T(\mathrm{pt})$ is the standard degree 2 generator. In summary, we find

$$H_\bullet^{G_\mathcal{O}}(\mathcal{R}) = H_\bullet^\bullet(\mathrm{pt}) \oplus \bigoplus_{n \geqslant 1} H_B^\bullet(\mathrm{pt})[2n] \cong \mathbb{C}[\delta] \oplus \bigoplus_{n \in \mathbb{N}} \mathbb{C}[\xi] \cdot \eta^n \cong \mathbb{C}[\delta, \eta, \tilde{\xi}]/(\tilde{\xi}^2 - \delta\eta^2)$$

where $\delta = \xi^2 \in H_G^\bullet(\mathrm{pt})[-4]$ is the degree 4 generator, $\tilde{\xi} := \xi \cdot \eta$, and the relation in the last expression follows from the redundancy in parameterization $(\tilde{\xi})^2 = \xi^2\eta^2 = \delta\eta^2$.

This lifts to an isomorphism of algebras

$$H_\bullet^{G_\mathcal{O}}(\mathrm{Gr}_G) = \mathbb{C}[\delta, \eta, \tilde{\xi}]/(\tilde{\xi}^2 - \delta\eta^2 - 1) = \mathcal{O}(\mathfrak{z}(\mathfrak{sl}_2))$$

with the coordinate ring of the group scheme $\mathfrak{z}(\mathfrak{sl}_2)$ of regular centralizers in $\mathfrak{sl}_2$; see [BFM05].

*Example* 13.3.7. Let $G = \mathrm{SL}_2$, and identify $\Lambda_{\mathrm{SL}_2}^+ \cong 2\mathbb{N} \subset \mathbb{N} \cong \Lambda_{\mathrm{PGL}_2}^+$, with the conventions as above. Then $\mathrm{Gr}_{\mathrm{SL}_2} \hookrightarrow \mathrm{Gr}_{\mathrm{PGL}_2}$ includes as the first of the two connected components of $\mathrm{Gr}_{\mathrm{PGL}_2}$, given by

$$\mathrm{Gr}_{\mathrm{SL}_2} = \sqcup_{n \in \mathbb{N}} \mathrm{Gr}_{\mathrm{PGL}_2}^{2n} \qquad \text{inducing} \qquad H_\bullet^{\mathrm{SL}_2, \mathcal{O}}(\mathrm{Gr}_{\mathrm{SL}_2}) = H_\bullet^{\mathrm{PGL}_2, \mathcal{O}}(\mathrm{Gr}_{\mathrm{PGL}_2})^{\mathbb{Z}/2\mathbb{Z}},$$

where the $\mathbb{Z}/2\mathbb{Z}$ action is of weight $-1$ on components of degree $4n$ and weight $1$ on those of degree $4n - 2$.

This lifts to an isomorphism of algebras

$$H_\bullet^{G_\mathcal{O}}(\mathrm{Gr}_G) = \mathbb{C}[\delta, \eta, \tilde{\xi}]^{\mathbb{Z}/2\mathbb{Z}}/(\tilde{\xi}^2 - \delta\eta^2 - 1) = \mathcal{O}(\mathfrak{z}(\mathfrak{sl}_2)/\mathbb{Z}/2\mathbb{Z}) \cong \mathcal{O}(\mathfrak{z}(\mathfrak{pgl}_2))$$

with the coordinate ring of the group scheme $\mathfrak{z}(\mathfrak{pgl}_2)$ of regular centralizers in $\mathfrak{pgl}_2$; see [BFM05].

13.4. **The convolution diagram for the variety of triples.** Following Remark 12.1.9, we would like to define the renormalized Borel-Moore chains

$$\mathcal{A}(G, N)_x = C_\bullet^{\mathrm{BM}}(\mathcal{Z}(G, N)_x) \in \mathrm{Alg}_{\mathbb{E}_1}(\mathrm{Vect})$$

on $\mathcal{Z}(G, N)_x$ and show that they form an $\mathbb{E}_1$-algebra under convolution. Moreover, we would like to show that

$$\mathcal{A}(G, N) = C_\bullet^{\mathrm{BM}}(\mathcal{Z}(G, N)) \in \mathrm{Alg}_{\mathbb{E}_1, \mathrm{un}}^{\mathrm{fact}}(X)$$

defines a factorization $\mathbb{E}_1$-algebra, though this will follow from the factorization compatibility of the construction, as outlined in Subsection 12.1, and thus for simplicity we present the arguments over a fixed point $x \in X$. Again, we follow [BFN18] and [BFN19b] throughout.

Recall from Proposition 13.1.9 and Remark 13.1.10 that we have an isomorphism of stacks

$$(13.4.1) \quad \mathcal{Z}(G, N)_x := N_\mathcal{O}/G_\mathcal{O} \underset{N_\mathcal{K}/G_\mathcal{K}}{\times} N_\mathcal{O}/G_\mathcal{O} \cong ([(N_\mathcal{O} \times G_\mathcal{K})/G_\mathcal{O}] \times_{N_\mathcal{K}} N_\mathcal{O})/G_\mathcal{O} = \mathcal{R}(G, N)_x/G_\mathcal{O} .$$

Thus, we have

$$(13.4.2) \quad D^*(\mathcal{Z}(G, N)_x) \cong D_{G_\mathcal{O}}^*(\mathcal{R}(G, N)_x) \in \mathrm{DGCat} ,$$

and following Remark 13.1.11, the results of Subsection 13.2 allow us to define the underlying vector space of our putative (factorization) $\mathbb{E}_1$-algebra:

*Definition* 13.4.1. The renormalized dualizing sheaf of $\mathcal{Z}(G,N)_x$ is defined by

$$\omega^{\mathrm{ren}}_{\mathcal{Z}(G,N)_x} := \omega^{\mathrm{ren}}_{\mathcal{R}(G,N)_x} \in D^*(\mathcal{Z}(G,N)_x) \ ,$$

where the latter is as in Definition 13.2.17.

*Remark* 13.4.2. In particular, we have

$$C^{\mathrm{BM}}_{\bullet}(\mathcal{Z}(G,N)_x) := \mathrm{p}_{\mathcal{Z}(G,N)_x *}\omega^{\mathrm{ren}}_{\mathcal{Z}(G,N)_x} \cong C^{G_{\mathcal{O}}}_{\bullet}(\mathcal{R}(G,N)_x) \ ,$$

where the latter is as in Definition 13.2.17.

*Remark* 13.4.3. Concretely, the object $\omega^{\mathrm{ren}}_{\mathcal{Z}(G,N)_x} \in D^*(\mathcal{Z}(G,N)_x)$ and its image $\omega^{\mathrm{ren}}_{\mathcal{R}(G,N)_x} \in D^*_{G_{\mathcal{O}}}(\mathcal{R}(G,N)_x)$ under the equivalence of Equation 13.4.2 above, are given by

$$D^*(\mathcal{Z}(G,N)_x) = D^*(N_{\mathcal{O}}/G_{\mathcal{O}} \underset{N_{\mathcal{K}}/G_{\mathcal{K}}}{\times} N_{\mathcal{O}}/G_{\mathcal{O}}) \qquad \overset{\cong}{\longrightarrow} \qquad D^*_{G_{\mathcal{O}}}([(N_{\mathcal{O}} \times G_{\mathcal{K}})/G_{\mathcal{O}}] \times_{N_{\mathcal{K}}} N_{\mathcal{O}})$$

$$\omega^{\mathrm{ren}}_{\mathcal{Z}(G,N)_x} = \quad \omega_{N_{\mathcal{O}}/G_{\mathcal{O}}} \underset{N_{\mathcal{K}}/G_{\mathcal{K}}}{\boxtimes} \mathbb{K}_{N_{\mathcal{O}}/G_{\mathcal{O}}} \qquad \mapsto \qquad \omega_{N_{\mathcal{O}}} \tilde{\boxtimes} \, \omega_{\mathrm{Gr}_G} \underset{N_{\mathcal{K}}}{\boxtimes} \mathbb{K}_{N_{\mathcal{O}}}$$

where we recall that

$$\mathbb{K}_{N_{\mathcal{O}}/G_{\mathcal{O}}} = \omega^{\mathrm{ren}}_{N_{\mathcal{O}}/G_{\mathcal{O}}} \qquad \text{and} \qquad \mathbb{K}_{N_{\mathcal{O}}} = \omega^{\mathrm{ren}}_{N_{\mathcal{O}}}$$

are as in Definition B.7.4, and the stated equalities follow from Example B.7.14. Note that the former description is in keeping with the internal Hom interpretation of Remark 12.2.8, and following Example A.6.7.

*Warning* 13.4.4. Note there is an automorphism of the prequotient of $\mathcal{Z}(G,N)_x$ by $G_{\mathcal{O}}$ given by

$$[(N_{\mathcal{O}} \times G_{\mathcal{K}})/G_{\mathcal{O}}] \times_{N_{\mathcal{K}}} N_{\mathcal{O}} \cong N_{\mathcal{O}} \times_{N_{\mathcal{K}}} [(N_{\mathcal{O}} \times G_{\mathcal{K}})/G_{\mathcal{O}}] \qquad \text{under which}$$

$$\omega_{N_{\mathcal{O}}} \tilde{\boxtimes} \, \omega_{\mathrm{Gr}_G} \underset{N_{\mathcal{K}}}{\boxtimes} \mathbb{K}_{N_{\mathcal{O}}} \mapsto \omega_{N_{\mathcal{O}}} \underset{N_{\mathcal{K}}}{\boxtimes} \mathbb{K}_{N_{\mathcal{O}}} \tilde{\boxtimes} \, \omega_{\mathrm{Gr}_G} \qquad .$$

The presentation of $\omega^{\mathrm{ren}}_{\mathcal{Z}(G,N)_x} \in D^*(\mathcal{Z}(G,N)_x)$ resulting from the latter description is evidently identified with the renormalized dualizing sheaf of Definition 13.2.17, but we will use the former description throughout the remainder of the text for notational convenience.

*Remark* 13.4.5. Following remarks 13.2.8 and 13.2.20, the identifications of the preceding remark are given heuristically by:

$$\omega_{\mathcal{Z}(G,N)_x}[-2\dim N_{\mathcal{O}}/G_{\mathcal{O}}] = \omega_{N_{\mathcal{O}}/G_{\mathcal{O}}} \underset{N_{\mathcal{K}}/G_{\mathcal{K}}}{\boxtimes} \omega_{N_{\mathcal{O}}/G_{\mathcal{O}}}[-2\dim N_{\mathcal{O}}/G_{\mathcal{O}}] \mapsto \omega_{N_{\mathcal{O}}} \tilde{\boxtimes} \, \omega_{\mathrm{Gr}_G} \underset{N_{\mathcal{K}}/G_{\mathcal{K}}}{\boxtimes} \omega_{N_{\mathcal{O}}}[-2\dim N_{\mathcal{O}}] \ .$$

Towards establishing Proposition 12.2.4 in the present context, we recall the analogue of Remark 18.3.2:

*Example* 13.4.6. The convolution diagram for $\mathcal{Z}(G,N)_x$ is given by

$$\mathcal{Z}_{(3)}(G,N)_x = N_{\mathcal{O}}/G_{\mathcal{O}} \underset{N_{\mathcal{K}}/G_{\mathcal{K}}}{\times} N_{\mathcal{O}}/G_{\mathcal{O}} \underset{N_{\mathcal{K}}/G_{\mathcal{K}}}{\times} N_{\mathcal{O}}/G_{\mathcal{O}}$$

together with the canonical projections $\pi_{ij} : \mathcal{Z}_{(3)}(G,N)_x \to \mathcal{Z}(G,N)_x$. There is a natural analogue for $\mathcal{Z}(G,N)_{(3)}$ of the isomorphism recalled in equation 13.4.1 above, given by

$$(13.4.3) \qquad \mathcal{Z}_{(3)}(G,N)_x \cong ([(N_{\mathcal{O}} \times G_{\mathcal{K}})/G_{\mathcal{O}}] \times_{N_{\mathcal{K}}} [(N_{\mathcal{O}} \times G_{\mathcal{K}})/G_{\mathcal{O}}] \times_{N_{\mathcal{K}}} N_{\mathcal{O}})/G_{\mathcal{O}} \ .$$

To describe the convolution diagram, we pass to the prequotient by $G_{\mathcal{O}}$ as usual, defining

$$\tilde{\mathcal{Z}}_{(3)}(G,N)_x = [(N_{\mathcal{O}} \times G_{\mathcal{K}})/G_{\mathcal{O}}] \times_{N_{\mathcal{K}}} [(N_{\mathcal{O}} \times G_{\mathcal{K}})/G_{\mathcal{O}}] \times_{N_{\mathcal{K}}} N_{\mathcal{O}} \cong ((G_{\mathcal{K}} \times N_{\mathcal{O}})/G_{\mathcal{O}})^{\times 2} \times_{N_{\mathcal{K}}^{\times 2}} N_{\mathcal{O}} \ ,$$

so that the prequotient by $G_{\mathbb{O}}$ of the convolution diagram becomes

$$
\begin{array}{ccc}
((G_{\mathcal{K}} \times N_{\mathbb{O}})/G_{\mathbb{O}})^{\times 2} \times_{N_{\mathcal{K}}^{\times 2}} N_{\mathbb{O}} & \xrightarrow{\pi_{13}} & (G_{\mathcal{K}} \times N_{\mathbb{O}})/G_{\mathbb{O}} \times_{N_{\mathcal{K}}} N_{\mathbb{O}} \\
\downarrow{\scriptstyle \pi_{12}} & {\scriptstyle \pi_{23}} \nearrow & \\
(G_{\mathcal{K}} \times N_{\mathbb{O}})/G_{\mathbb{O}} \times_{N_{\mathcal{K}}} N_{\mathbb{O}} & & G_{\mathbb{O}}\backslash(G_{\mathcal{K}} \times N_{\mathbb{O}})/G_{\mathbb{O}} \times_{N_{\mathcal{K}}} N_{\mathbb{O}}
\end{array}
\qquad
\begin{array}{ccc}
(g_1, g_2, s) & \longmapsto & (g_1 g_2, s) \\
\downarrow & & \\
(g_1, g_2 \cdot s) & & (g_2, s)
\end{array}
.
$$

The indscheme $\tilde{\mathcal{Z}}_{(3)}(G, N)_x$ is also denoted "$qp^{-1}(\mathcal{R} \times \mathcal{R})$" in [BFN18].

*Remark* 13.4.7. The convolution diagram for $\mathcal{Z}(G, N)_x$ fits into a larger commutative diagram covering that for the usual geometric Satake:

$$
\begin{array}{ccccccc}
\mathcal{R} \times \mathcal{R} & \longleftarrow & p^{-1}(\mathcal{R} \times \mathcal{R}) & \longrightarrow & q(p^{-1}(\mathcal{R} \times \mathcal{R})) & \longrightarrow & \mathcal{R} \\
\downarrow & & \downarrow & & \downarrow & & \| \\
\mathcal{T} \times \mathcal{R} & \xleftarrow{p} & G_{\mathcal{K}} \times \mathcal{R} & \xrightarrow{q} & (G_{\mathcal{K}} \times \mathcal{R})/G_{\mathbb{O}} & \longrightarrow & \mathcal{R} \\
\downarrow & & \downarrow & & \downarrow & & \downarrow \\
\mathrm{Gr}_G \times \mathrm{Gr}_G & \xleftarrow{\tilde{p}} & G_{\mathcal{K}} \times \mathrm{Gr}_G & \xrightarrow{\tilde{q}} & (G_{\mathcal{K}} \times \mathrm{Gr}_G)/G_{\mathbb{O}} & \longrightarrow & \mathrm{Gr}_G
\end{array}
.
$$

Following Proposition 12.2.4 and in turn Proposition 8.1.5, we have:

*Proposition* 13.4.8. The factorization category $(D^*_{\mathcal{Z}(G,N)})^\star \in \mathrm{Cat}^{\mathrm{fact}}_{\mathbb{E}_1, \mathrm{un}}(X_{\mathrm{dR}})$ is naturally a (unital) $\mathbb{E}_1$-factorization category with respect to the convolution monoidal structure defined by the composition

$$
D^*_{\mathcal{Z}(G,N)} \otimes^* D^*_{\mathcal{Z}(G,N)} \xrightarrow{\pi^*_{12} \boxtimes \pi^*_{23}} D^*_{\mathcal{Z}_{(3)}(G,N)^{\times 2}} \xrightarrow{\Delta^*} D^*_{\mathcal{Z}_{(3)}(G,N)} \xrightarrow{\pi_{13,*}} D^*_{\mathcal{Z}(G,N)} .
$$

Further, the pushforward functor $\mathrm{p}_{\mathcal{Z}*} : (D^*_{\mathcal{Z}(G,N)})^\star \to D^{\otimes^!}_{\mathrm{Ran}_{X,\mathrm{un}}}$ is a unital, lax $\mathbb{E}_1$-monoidal factorization functor.

*Warning* 13.4.9. The functors $\pi^*_{12}$ and $\pi^*_{23}$ of the preceding proposition are only partially defined in general. For example, they are each defined on the full subcategories of holonomic $D$ modules, in the sense of Definition B.8.11, which suffices for the application in Theorem 13.4.12 below.

Moreover, the composition

$$
(\pi_{12} \times \pi_{23})^* := \Delta^* \circ (\pi^*_{12} \boxtimes \pi^*_{23}) : D^*(\mathcal{Z}(G,N)_x) \otimes D^*(\mathcal{Z}(G,N)_x) \to D^*(\mathcal{Z}_{(3)}(G,N))
$$

used in the preceding proposition is in fact well-defined on the entire categories, so that the proposition holds as stated. This follows from the fact that the map $\pi_{12} \times \pi_{23}$ is pro-smooth on the prequotient, even though $\pi_{12}$ and $\pi_{23}$ individually are not, so that the (equivariant) $D$ module pullback functor $(\pi_{12} \times \pi_{23})^*$ is well-defined on the entire category.

Further, following Corollary 8.1.6, we have:

*Corollary* 13.4.10. The pushforward functor $\mathrm{p}_{\mathcal{Z}(G,N)*} : D^\star_{\mathcal{Z}(G,N)} \to D_{\mathrm{Ran}_{X,\mathrm{un}}}$ induces a functor

$$
(13.4.4) \qquad \mathrm{p}_{\mathcal{Z}(G,N)*} : \mathrm{Alg}^{\mathrm{fact}}_{\mathbb{E}_1, \mathrm{un}}(D^\star_{\mathcal{Z}(G,N)}) \to \mathrm{Alg}^{\mathrm{fact}}_{\mathbb{E}_1, \mathrm{un}}(X) .
$$

*Remark* 13.4.11. Note the maps $\pi_{ij}$ are ind-proper (and in particular ind-finitely presented and reasonable), as they are closed sub fibre bundles of $\mathrm{Gr}_G$ bundles; thus, the $(\pi_{13,*}, \pi_{13}^{!*})$ adjunction yields the natural transformation

$$\pi_{13,*}\pi_{13}^{!*} \to \mathbb{1} \quad \text{in} \quad D^*_{G_{\mathcal{O}} \times \mathbb{C}_\hbar^\times}(\mathcal{R}) .$$

Following Remark 8.1.8 and in particular Example 8.1.15, we obtain the analogue of Proposition A.6.2 and in particular Example A.6.17 via Proposition 13.4.8 above:

*Theorem* 13.4.12. The renormalized dualizing sheaf of Definition 13.4.1 defines

$$\tilde{\mathcal{A}}(G, N) := \omega^{\mathrm{ren}}_{\mathcal{Z}(G,N)} \in \mathrm{Alg}^{\mathrm{fact}}_{\mathbb{E}_1, \mathrm{un}}(D^\star_{\mathcal{Z}(G,N)})$$

a factorization $\mathbb{E}_1$ algebra internal to $D^\star_{\mathcal{Z}(G,N)} \in \mathrm{Cat}^{\mathrm{fact}}_{\mathbb{E}_1, \mathrm{un}}(X)$ as in Proposition 13.4.8 above. In particular, the renormalized Borel-Moore chains

$$\mathcal{A}(G, N) := \mathrm{p}_{\mathcal{Z}(G,N)*}\omega^{\mathrm{ren}}_{\mathcal{Z}(G,N)} = C^{\mathrm{BM}}_\bullet(\mathcal{Z}(G,N)) \in \mathrm{Alg}^{\mathrm{fact}}_{\mathbb{E}_1, \mathrm{un}}(X)$$

defines a factorization $\mathbb{E}_1$ algebra.

*Proof.* We follow the proof of Proposition A.6.2, and begin by establishing the formal setup in our present setting. The ind-proper adjunction of Remark 13.4.11 yields the desired map

$$\pi_{13,*}\pi_{13}^{!*}\omega^{\mathrm{ren}}_{\mathcal{R}} \to \omega^{\mathrm{ren}}_{\mathcal{R}} \quad \text{in} \quad D^*_{G_{\mathcal{O}} \times \mathbb{C}_\hbar^\times}(\mathcal{R}) .$$

while the $(\pi^*, \pi_*)$ adjunction for holonomic $D$ modules on $\mathcal{R}$ yields the desired maps

$$\omega^{\mathrm{ren}}_{\mathcal{R}} \to \pi_{12,*}\pi^{12,*}\omega^{\mathrm{ren}}_{\mathcal{R}} \qquad \text{and} \qquad \omega^{\mathrm{ren}}_{\mathcal{R}} \to \pi_{23,*}\pi^{23,*}\omega^{\mathrm{ren}}_{\mathcal{R}} \quad \text{in} \quad D^*_{G_{\mathcal{O}} \times \mathbb{C}_\hbar^\times}(\mathcal{R}) ,$$

as well as

$$\pi^{12,*}\omega^{\mathrm{ren}}_{\mathcal{R}}\boxtimes\pi^{23,*}\omega^{\mathrm{ren}}_{\mathcal{R}} \to \Delta_*\Delta^*(\pi^{12,*}\omega^{\mathrm{ren}}_{\mathcal{R}}\boxtimes\pi^{23,*}\omega^{\mathrm{ren}}_{\mathcal{R}}) \cong \Delta_*(\pi^{12,*}\omega^{\mathrm{ren}}_{\mathcal{R}}\otimes^*\pi^{23,*}\omega^{\mathrm{ren}}_{\mathcal{R}}) \quad \text{in} \quad D^*_{G_{\mathcal{O}} \times \mathbb{C}_\hbar^\times}(\mathcal{Z}^{\times 2}_{(3)}) .$$

Thus, it suffices to construct the multiplication map, as in Equation 8.1.2 of Remark 8.1.7, following Proposition A.6.2 and Example A.6.6. Following the description of $\omega^{\mathrm{ren}}_{\mathcal{Z}(G,N)}$ in Remark 13.4.3, under the identification induced by Equation 13.4.3, we compute

$$D^*(\mathcal{Z}_{(3)}(G, N)_x) \cong D^*_{G_{\mathcal{O}}}([(N_{\mathcal{O}} \times G_{\mathcal{K}})/G_{\mathcal{O}}] \times_{N_{\mathcal{K}}} [(N_{\mathcal{O}} \times G_{\mathcal{K}})/G_{\mathcal{O}}] \times_{N_{\mathcal{K}}} N_{\mathcal{O}})$$

$$\pi^*_{12}\omega^{\mathrm{ren}}_{\mathcal{Z}(G,N)} \mapsto \quad \omega_{N_{\mathcal{O}}}\tilde{\boxtimes} \, \omega_{\mathrm{Gr}_G} \boxtimes_{N_{\mathcal{K}}} \underline{\mathbb{K}}_{N_{\mathcal{O}}}\tilde{\boxtimes} \, \underline{\mathbb{K}}_{\mathrm{Gr}_G} \boxtimes_{N_{\mathcal{K}}} \underline{\mathbb{K}}_{N_{\mathcal{O}}}$$

$$\pi^*_{23}\omega^{\mathrm{ren}}_{\mathcal{Z}(G,N)} \mapsto \quad \underline{\mathbb{K}}_{N_{\mathcal{O}}}\tilde{\boxtimes} \, \underline{\mathbb{K}}_{\mathrm{Gr}_G} \boxtimes_{N_{\mathcal{K}}} \omega_{N_{\mathcal{O}}}\tilde{\boxtimes} \, \omega_{\mathrm{Gr}_G} \boxtimes_{N_{\mathcal{K}}} \underline{\mathbb{K}}_{N_{\mathcal{O}}}$$

$$\pi^!_{13}\omega^{\mathrm{ren}}_{\mathcal{Z}(G,N)} \mapsto \quad \omega_{N_{\mathcal{O}}}\tilde{\boxtimes} \, \omega_{\mathrm{Gr}_G} \boxtimes_{N_{\mathcal{K}}} \omega_{N_{\mathcal{O}}}\tilde{\boxtimes} \, \omega_{\mathrm{Gr}_G} \boxtimes_{N_{\mathcal{K}}} \underline{\mathbb{K}}_{N_{\mathcal{O}}}$$

.

Thus, recalling $\underline{\mathbb{K}}$ is the unit of the $\otimes^*$ tensor structure, we have a canonical map

$$\pi^*_{12}\omega^{\mathrm{ren}}_{\mathcal{Z}(G,N)}\otimes^*\pi^*_{23}\omega^{\mathrm{ren}}_{\mathcal{Z}(G,N)} \xrightarrow{\cong} \pi^!_{13}\omega^{\mathrm{ren}}_{\mathcal{Z}(G,N)} \qquad \text{and thus} \qquad \pi_{13*}\left(\pi^*_{12}\omega^{\mathrm{ren}}_{\mathcal{Z}(G,N)} \otimes^* \pi^*_{23}\omega^{\mathrm{ren}}_{\mathcal{Z}(G,N)}\right) \to \omega^{\mathrm{ren}}_{\mathcal{Z}(G,N)}$$

by Remark 13.4.11, as desired. $\square$

13.5. **Equivariant structures on the three dimensional A model.** In this subsection, we explain that the factorization algebra underlying the three dimensional A model to $Y = N/G$ on $X = \mathbb{A}^1$ is canonically $\mathbb{G}_a \rtimes \mathbb{G}_m$ equivariant, and thus defines a 2-periodic filtered quantization by the results of Section I-25, and in particular Example I-25.1.5. Again, these results were essentially obtained in [BFN18]. We also recall several resulting concrete computations, following *loc. cit.*.

To begin, following Proposition 12.3.1, we have:

*Proposition* 13.5.1. The (unital) factorization $\mathbb{E}_1$ algebra of the three dimensional A model to $N/G$ over $X = \mathbb{A}^1$

$$\mathcal{A}(G, N) \in \mathrm{Alg}^{\mathrm{fact}}_{\mathbb{E}_1, \mathrm{un}}(\mathbb{A}^1)^{\mathbb{G}_a \rtimes \mathbb{G}_m}$$

admits a canonical $\mathbb{G}_a \rtimes \mathbb{G}_m$ equivariant structure.

Similarly, following Corollary 12.3.3, we have:

*Corollary* 13.5.2. The three dimensional A model $\mathcal{A}(G, N) \in \mathrm{Alg}^{S^1}_{\mathbb{E}_3}(\mathrm{Vect}_{\mathbb{K}})$ to $N/G$ over $X = \mathbb{A}^1$ defines

$$\mathcal{A}_u(G, N) \in \mathrm{Alg}_{\mathbb{BD}^u_1}(\mathrm{D}^b(\mathbb{K}[u]))\ ,$$

an algebra over the operad $\mathbb{BD}^u_1 \in \mathrm{Op}(\mathrm{D}^b(\mathbb{K}[u]))$ of Definition I-23.0.5.

*Remark* 13.5.3. Concretely, $\mathcal{A}_u(G, N) \in \mathrm{Alg}_{\mathbb{BD}^u_1}(\mathrm{D}^b(\mathbb{K}[u]))$ defines a two-periodic filtered quantization of the homology $\mathbb{P}_3$ algebra $H_\bullet(\mathcal{A}) \in \mathrm{Alg}_{\mathbb{P}_3}(\mathrm{Perf}_{\mathbb{K}})$ to an $\mathbb{E}_1$ algebra, or equivalently a usual associative algebra $\mathcal{A}_u|_{\{1\}} \in \mathrm{Alg}_{\mathbb{E}_1}(\mathrm{Perf}_{\mathbb{K}})$. The homology $\mathbb{P}_3$ algebra is denoted by $\mathcal{A}(G, N)$ in [BFN18], while the $\mathbb{BD}^u_1$ algebra is denoted $\mathcal{A}_\hbar(G, N)$ in *loc. cit.*. We use the notation $u$ in place of $\hbar$ as a reminder that the generator of the polynomial algebra is of cohomological degree $+2$.

*Example* 13.5.4. Let $G = T = (\mathbb{C}^\times)^n$ be a rank $n$ torus and $N$ the zero representation, as in Example 13.3.5. Then we have seen $\mathcal{R}(\mathbb{C}) = \mathrm{Gr}_G(\mathbb{C}) \cong \mathbb{Z}^n$ and that

$$H^{G_\mathcal{O}}_\bullet(\mathcal{R}) = \bigoplus_{\lambda \in \mathbb{Z}^n} \mathrm{Sym}^\bullet(\mathfrak{t}^\vee[-2]) \cong \bigoplus_{\lambda \in \mathbb{Z}^n} \mathbb{C}[\xi_i]^n_{i=1} \cdot r^\lambda = \mathbb{C}[\xi_i, r^{\pm 1}_i]^n_{i=1}\ ,$$

where $r^\lambda = [\mathrm{Gr}^\lambda_T] \in H^{T_\mathcal{O}}_0(\mathcal{R})$. The convolution diagram is given by

$$\begin{array}{ccc} (\mathbb{Z}^n)^{\times 2}/T \longrightarrow \mathbb{Z}^n/T & \quad & (\lambda_1, \lambda_2) \longmapsto \lambda_1 + \lambda_2 \\ \downarrow & & \downarrow \\ \mathbb{Z}^n/T \qquad \mathbb{Z}^n/T & & \lambda_1 \qquad\qquad \lambda_2 \end{array}$$

so that we find $\pi_{13,*}(\pi^*_{12} r^\lambda \otimes^* \pi^*_{23} r^\mu) = r^{\lambda + \mu}$ as desired. Thus, we obtain

$$H^{T_\mathcal{O}}_\bullet(\mathrm{Gr}_T) = \mathbb{C}[\xi_i, r^{\pm 1}_i]^n_{i=1} = \mathcal{O}(\mathfrak{t} \times T^\vee) \cong \mathcal{O}(T^* T^\vee)$$

is given by the coordinate ring of the cotangent bundle of the dual torus $T^\vee$. The filtered quantization is given by

$$\mathcal{A}_u(T, \{0\}) = H^{T_\mathcal{O} \rtimes \mathbb{C}^\times_u}_\bullet(\mathrm{Gr}_T) = \Gamma(T^\vee, \mathcal{D}_{T^\vee, u})$$

the algebra of global Reese differential operators on $T^\vee$.

*Example* 13.5.5. More generally, for $G = T$ as above and $N$ an arbirary representation, we have

$$H^{G_\mathcal{O}}_\bullet(\mathcal{R}) = \bigoplus_{\lambda \in \mathbb{Z}^n} \mathrm{Sym}^\bullet(\mathfrak{t}^\vee[-2])[-2d_\lambda] \cong \bigoplus_{\lambda \in \mathbb{Z}^n} \mathbb{C}[\xi_i]^n_{i=1} \cdot r^\lambda = \mathbb{C}[\xi_i, r^{\pm 1}_i]^n_{i=1}\ ,$$

where $\xi_i \in H_T^\bullet(\mathrm{pt})[-2]$ are the generators of the $T$ equivariant cohomology of a point and $r^\lambda := [\mathcal{R}^\lambda] \in H_{2d_\lambda}^{G_\mathcal{O}}(\mathcal{R}^\lambda)$ is the fundamental class of $\mathcal{R}^\lambda$.

The main result facilitating more general computations in the abelian case is the following:

*Proposition* 13.5.6. The inclusion of the zero section $\mathrm{Gr}_T \hookrightarrow \mathcal{R}$ induces an inclusion of algebras

$$\mathcal{A}_\hbar(T, N) \hookrightarrow \mathcal{A}_\hbar(T, \{0\}) \qquad \text{defined by} \qquad r^\lambda \mapsto \tilde{r}^\lambda \cdot \mathrm{eu}(z^\lambda N_\mathcal{O}/(N_\mathcal{O} \cap z^\lambda N_\mathcal{O}))$$

where $r^\lambda := [\mathcal{R}^\lambda] \in H_{2d_\lambda}^{T_\mathcal{O}}(\mathcal{R}^\lambda)$ is the fundamental class of $\mathcal{R}^\lambda$, as above, and $\tilde{r}^\lambda = [\mathrm{Gr}_T^\lambda] \in H_0^{T_\mathcal{O}}(\mathrm{Gr}_T)$.

*Remark* 13.5.7. Concretely, the above map is given by $r^\lambda \mapsto \begin{cases} \xi^{-\xi(\lambda)} r^\lambda & \text{for } \xi(\lambda) < 0 \\ r^\lambda & \text{otherwise} \end{cases}$.

*Example* 13.5.8. Let $T = \mathbb{C}^\times$, $N = \mathbb{C}^n$ under the diagonal action, so that the corresponding stable quotient variety is $\mathbb{P}^{n-1}$. Letting $x = r, y = r^{-1}$, we find

$$\mathcal{A} = \mathbb{C}[\xi, x, y]/(xy - \xi^n) \qquad \text{or} \qquad \mathcal{A}_\hbar = \mathbb{C}\langle x, y, \xi, \hbar \mid xy = (\xi + \hbar/2)^n, yx = (\xi - \hbar/2)^n, [x, \xi] = \hbar x, [y, \xi] = -\hbar y \rangle$$

after passing to $\mathbb{C}^\times$ equivariant cohomology with respect to loop rotation. This follows from the above proposition, noting $xy = r \cdot r^{-1} = \xi^n \tilde{r} \cdot \tilde{r}^{-1} = \xi^n$.

Note this is an affine filtered quantization of the $A_{n-1}$ singularity, which is the expected symplectic dual to $\mathbb{P}^{n-1}$, as desired. In particular, we have $\mathcal{A}_\hbar = \mathcal{U}_\hbar(\mathfrak{pgl}_2)$ for $n = 2$.

# 14. Chiral differential operators and the three dimensional A model: Overview

In this section, we outline the geometric construction of the factorization algebra of chiral differential operators $\mathcal{D}^{\mathrm{ch}}(Y) \in \mathrm{Alg}_{\mathrm{un}}^{\mathrm{fact}}(X)$ on $Y$ over $X$ a smooth algebraic curve, following [KV06], and in turn [BD04] and the series of papers [GMS00], [GMS04], and [GMS03]; the construction is recalled in terms of the general format outlined in Section 9.1, in the variant described in Section 8.2. Moreover, we explain that the module structure given by the variant of Example 8.2.8 in this setting defines a module structure $\mathcal{D}^{\mathrm{ch}}(Y) \in \mathcal{A}(Y)\text{-Mod}(\mathrm{Alg}_{\mathrm{un}}^{\mathrm{fact}}(X))$ over the factorization $\mathbb{E}_1$ algebra $\mathcal{A}(Y) \in \mathrm{Alg}_{\mathbb{E}_1, \mathrm{un}}^{\mathrm{fact}}(X)$ describing the three dimensional A model, as constructed in the preceding sections following [BFN18]; this establishes the conjecture of Costello-Gaiotto [CG18] which predicted that there was a chiral boundary condition for the three dimensional A model with boundary observables given by the algebra of chiral differential operators.

The geometric definition of the factorization algebra $\mathcal{D}^{\mathrm{ch}}(Y) \in \mathrm{Alg}_{\mathrm{un}}^{\mathrm{fact}}(X)$ is given by

$$\mathcal{D}^{\mathrm{ch}}(Y)_x = \mathcal{H}\mathrm{om}_{D(Y_\mathcal{K})}(\mathcal{D}_{Y_\mathcal{K}}, \iota_* \omega_{Y_\mathcal{O}}) \cong \Gamma(Y_\mathcal{K}, \iota_* \omega_{Y_\mathcal{O}}) ,$$

where the line operator category $D(Y_\mathcal{K}) \in \mathrm{Cat}_{\mathrm{un}}^{\mathrm{fact}}(X)$ is given by the factorization category $D$ modules on the meromorphic jet scheme to $Y$, as in Section 12; the construction is summarized by the following diagrams in factorization spaces and categories, which for simplicity we denote by their fibre over a fixed point $x \in X$:

(14.0.1)

$$
\begin{array}{ccc}
& (Y_{\mathcal{K}})^{\wedge}_{Y_{\mathcal{O}}} & \\
{\scriptstyle \pi_{Y_{\mathcal{O},\mathrm{dR}}}}\swarrow & & \searrow{\scriptstyle \pi_{Y_{\mathcal{K}}}} \\
Y_{\mathcal{O},\mathrm{dR}} & & Y_{\mathcal{K}} \\
{\scriptstyle \mathrm{p}_{Y_{\mathcal{O},\mathrm{dR}}}}\downarrow \quad {\scriptstyle \iota_{\mathrm{dR}}}\searrow \quad \swarrow{\scriptstyle \mathrm{q}} & & \downarrow{\scriptstyle \mathrm{p}_{Y_{\mathcal{K}}}} \\
\mathrm{pt} \qquad Y_{\mathcal{K},\mathrm{dR}} & & \mathrm{pt}
\end{array}
\qquad\text{and}\qquad
\begin{array}{ccc}
& \mathrm{IndCoh}((Y_{\mathcal{K}})^{\wedge}_{Y_{\mathcal{O}}}) & \\
{\scriptstyle \pi^{!}_{Y_{\mathcal{O},\mathrm{dR}}}}\nearrow & & \nwarrow{\scriptstyle \pi_{Y_{\mathcal{K}},\bullet}} \\
D(Y_{\mathcal{O}}) & & \mathrm{IndCoh}(Y_{\mathcal{K}}) \\
{\scriptstyle \mathrm{p}^{!}_{Y_{\mathcal{O}}}}\uparrow \quad {\scriptstyle \iota_{*}}\nwarrow \quad \nearrow{\scriptstyle \mathrm{q}^{!}} & & \uparrow{\scriptstyle \mathrm{p}_{Y_{\mathcal{K}},\bullet}} \\
\mathrm{Vect} \qquad D(Y_{\mathcal{K}}) & & \mathrm{Vect}
\end{array}
.
$$

14.0.1. *Summary.* In Section 14.3 we outline the construction of the factorization algebra of chiral differential operators $\mathcal{D}^{\mathrm{ch}}(Y)$ on a space $Y$, in Section 14.2 we explain the schematic construction of the expected module structure on $\mathcal{D}^{\mathrm{ch}}(Y)$ over the factorization $\mathbb{E}_1$ algebra $\mathcal{A}(Y) \in \mathrm{Alg}^{\mathrm{fact}}_{\mathbb{E}_1,\mathrm{un}}(X)$ describing the three dimensional A model, and in Section 14.1 we outline the internal variant of the construction of $\mathcal{D}^{\mathrm{ch}}(Y)$.

*Warning* 14.0.1. In keeping with Warning 9.1.1, we do not formulate specific hypotheses on the space $Y$ used in this section, so that the results stated throughout are only an outline of the general expectations. However, in the following several sections, we will recall the precise geometric construction of chiral differential operators in the case that $Y$ is a well-behaved smooth, affine scheme, generalize this construction to the case that $Y = N/G$ is given by a quotient stack as in Section 13, and give a careful construction of the canonical module structure on $\mathcal{D}^{\mathrm{ch}}(N/G)$ over $\mathcal{A}(G, N)$ in this setting.

## 14.1. **The factorization algebra of chiral differential operators.**

*Example* 14.1.1. Recall the (unital) factorization spaces $\mathcal{J}^{\mathrm{mer}}(Y)$, $\mathcal{J}^{\mathrm{mer}}(Y)_{\mathrm{dR}} \in \mathrm{PreStk}^{\mathrm{fact}}_{\mathrm{un}}(X)$ of meromorphic jets to $Y$ and its de Rham stack, from Examples 4.2.4 and 12.1.1, respectively. The projection to the de Rham stack defines a map of factorization spaces

$$\mathrm{q}_{\mathcal{J}^{\mathrm{mer}}(Y)} : \mathcal{J}^{\mathrm{mer}}(Y) \to \mathcal{J}^{\mathrm{mer}}(Y)_{\mathrm{dR}} .$$

*Remark* 14.1.2. Concretely, the maps on prestacks over each $X^I$, and over $x \in X$, are given by

$$\mathrm{q}_{\mathcal{J}^{\mathrm{mer}}(Y)_I} : \mathcal{J}^{\mathrm{mer}}(Y)_I \to \mathcal{J}^{\mathrm{mer}}(Y)_{I,\mathrm{dR}} \qquad \text{and} \qquad \mathrm{q}_{Y_{\mathcal{K}_x}} : Y_{\mathcal{K}_x} \to Y_{\mathcal{K}_x,\mathrm{dR}} .$$

Towards defining the factorization functor corresponding to pullback along the maps of factorization spaces of the preceding example, following Definition A.5.16, we introduce the following factorization categories:

*Example* 14.1.3. Consider the (putative, see Remark 14.1.5 below) indcoherent variants of the factorization categories of examples 5.2.9 and 12.1.3: the category of indcoherent sheaves, and that of 'right' $D$ modules, on the loop space of $Y$ are the unital factorization categories defined by

$$\mathrm{IndCoh}_{\mathcal{J}^{\mathrm{mer}}(Y)} \in \mathrm{Cat}^{\mathrm{fact}}_{\mathrm{un}}(X) \qquad \text{and} \qquad D^r_{\mathcal{J}^{\mathrm{mer}}(Y)} := \mathrm{IndCoh}_{\mathcal{J}^{\mathrm{mer}}(Y)_{\mathrm{dR}}} \in \mathrm{Cat}^{\mathrm{fact}}_{\mathrm{un}}(X) .$$

*Remark* 14.1.4. Concretely, the sheaf of categories assigned to each $I \in \mathrm{fSet}$, its sections category, and the fibre category over $x \in X$, are given for $\mathrm{IndCoh}_{\mathcal{J}^{\mathrm{mer}}(Y)}$ and $D^r_{\mathcal{J}^{\mathrm{mer}}(Y)}$ by

$$\mathrm{p}_{I*}\,\mathrm{IndCoh}_{\mathcal{J}^{\mathrm{mer}}(Y)_I} \in \mathrm{ShvCat}(X^I_{\mathrm{dR}}) \quad \mathrm{IndCoh}(\mathcal{J}^{\mathrm{mer}}(Y)_I) \in \mathrm{DGCat} \quad \text{and} \quad \mathrm{IndCoh}(Y_{\mathcal{K}_x}) \in \mathrm{DGCat} , \text{ and}$$

$$\mathrm{p}_{I,\mathrm{dR}*}\,\mathrm{IndCoh}_{\mathcal{J}^{\mathrm{mer}}(Y)_{I,\mathrm{dR}}} \in \mathrm{ShvCat}(X^I_{\mathrm{dR}}) \qquad D^r(\mathcal{J}^{\mathrm{mer}}(Y)_I) \in \mathrm{DGCat} \quad \text{and} \qquad D^r(Y_{\mathcal{K}_x}) \in \mathrm{DGCat} .$$

*Remark* 14.1.5. The relevant (sheaves of) categories of indcoherent sheaves are not defined in this level of generality, so the above example and the remaining discussion in the present subsection is only given as a heuristic overview of the contents of the remainder of this section.

In the following subsections, we will make sense of the constructions below under more restrictive assumptions on the input space $Y$, using the results recalled in Appendix B, following [KV06], [BD04], [Ras15a], [Ras20b] and references therein. For the remainder of the present section, we summarize the narrative structure of the results, implicitly assuming $Y$ is well-behaved enough that the various constructions can be defined.

The above subtlety is not just a technical obstacle: careful consideration of this issue leads to the appearance of the 'determinant anomaly' and the requirement of a choice of Tate structure on $\mathcal{J}^{\mathrm{mer}}(Y)$, as we explain in Remark 15.1.22 in the setting where $Y$ is a scheme, and in Remark 16.1.16 in the case $Y = N/G$ is a quotient stack.

*Example* 14.1.6. Recall from Example 12.1.6 that the factorization unit object

$$\mathrm{unit}_{D_{\mathcal{J}^{\mathrm{mer}}(Y)}} \in \mathrm{Alg}_{\mathrm{un}}^{\mathrm{fact}}(D_{\mathcal{J}^{\mathrm{mer}}(Y)}) \qquad \text{is defined by} \qquad \mathrm{unit}_{D_{\mathcal{J}^{\mathrm{mer}}(Y)}} = \iota_{\mathcal{J}(Y),*} \mathrm{p}_{\mathcal{J}(Y)}^{!} \omega_{\mathrm{Ran}_{X,\mathrm{un}}} = \iota_{\mathcal{J}(Y)*} \omega_{\mathcal{J}(Y)}$$

where $\mathrm{p}_{\mathcal{J}(Y)} : \mathcal{J}(Y) \to \mathrm{Ran}_{X,\mathrm{un}}$ is the factorization space structure map and $\iota_{\mathcal{J}(Y)} : \mathcal{J}(Y) \to \mathcal{J}^{\mathrm{mer}}(Y)$ is the map of factorization spaces defined by the unit section.

*Example* 14.1.7. The map $\mathrm{q}_{\mathcal{J}^{\mathrm{mer}}(Y)} : \mathcal{J}^{\mathrm{mer}}(Y) \to \mathcal{J}^{\mathrm{mer}}(Y)_{\mathrm{dR}}$ induces a unital factorization functor

$$o_{\mathcal{J}^{\mathrm{mer}}(Y)}^{r} := \mathrm{q}_{\mathcal{J}^{\mathrm{mer}}(Y)}^{!} : D_{\mathcal{J}^{\mathrm{mer}}(Y)} \to \mathrm{IndCoh}_{\mathcal{J}^{\mathrm{mer}}(Y)} \ ,$$

following Example 5.3.2, which is interpreted as the forgetful functor from right $D$ modules to (ind)coherent sheaves, following Proposition A.5.16.

*Example* 14.1.8. The map $\mathrm{p}_{\mathcal{J}^{\mathrm{mer}}(Y)} : \mathcal{J}^{\mathrm{mer}}(Y) \to \mathrm{Ran}_{X_{\mathrm{dR}},\mathrm{un}}$ induces a unital factorization functor

$$\mathrm{p}_{\mathcal{J}^{\mathrm{mer}}(Y)\bullet} : \mathrm{IndCoh}_{\mathcal{J}^{\mathrm{mer}}(Y)} \to D_{\mathrm{Ran}_{X,\mathrm{un}}} \ ,$$

defined following Example 5.3.1; we have used the fact that $\mathcal{J}^{\mathrm{mer}}(Y)$ is naturally a factorization space over $X_{\mathrm{dR}}$, as explained in Remark 4.2.5.

*Remark* 14.1.9. Concretely, the functor on sections categories over each $X^I$, and over $x \in X$, are given by

$$\mathrm{p}_{\mathcal{J}^{\mathrm{mer}}(Y)_I\bullet} : \mathrm{IndCoh}(\mathcal{J}^{\mathrm{mer}}(Y)_I) \to \mathrm{IndCoh}(X_{\mathrm{dR}}^I) \qquad \text{and} \qquad \mathrm{p}_{Y_{\mathcal{K}x}\bullet} : \mathrm{IndCoh}(Y_{\mathcal{K}}) \to \mathrm{Vect} \ .$$

*Tentative Definition* 14.1.10. The chiral differential operators on $Y$ is the unital factorization algebra

$$\mathcal{D}^{\mathrm{ch}}(Y) = \mathrm{p}_{\mathcal{J}^{\mathrm{mer}}(Y)\bullet} \circ \mathrm{q}_{\mathcal{J}^{\mathrm{mer}}(Y)}^{!}(\mathrm{unit}_{D_{\mathcal{J}^{\mathrm{mer}}(Y)}}) \ \in \mathrm{Alg}_{\mathrm{un}}^{\mathrm{fact}}(X) \ .$$

*Remark* 14.1.11. Concretely, the $D$ module $\mathcal{D}^{\mathrm{ch}}(Y)_I \in D(X^I)$ assigned to each $I \in \mathrm{fSet}_{\varnothing}$, and the fibre vector space $\mathcal{D}^{\mathrm{ch}}(Y)_x \in \mathrm{Vect}$ over each $x \in X$, are given by

$$\mathcal{D}^{\mathrm{ch}}(Y)_I = \mathrm{p}_{\mathcal{J}^{\mathrm{mer}}(Y)_I\bullet} \mathrm{q}_{\mathcal{J}^{\mathrm{mer}}(Y)_I}^{!} \iota_{\mathcal{J}(Y)_I*} \omega_{\mathcal{J}(Y)_I} \qquad \text{and} \qquad \mathcal{D}^{\mathrm{ch}}(Y)_x = \mathrm{p}_{Y_{\mathcal{K}x}\bullet} \mathrm{q}_{Y_{\mathcal{K}_x}}^{!} \iota_{x*} \omega_{Y_{\mathcal{O}_x}} = \Gamma(Y_{\mathcal{K}_x}, \iota_{x*} \omega_{Y_{\mathcal{O}_x}}) \ ,$$

where in the latter expression $\Gamma$ denotes the space of sections as an $\mathcal{O}$ module.

*Remark* 14.1.12. In examples 15.1.34 and 16.1.25, we give concrete calculations following the above definition in the cases $Y = N$ and $Y = N/G$, respectively.

14.2. **Action of the three dimensional A model on chiral differential operators.** Applying Example 7.1.10 to $\mathcal{A}(Y) \in \mathrm{Alg}^{\mathrm{fact}}_{\mathbb{E}_1,\mathrm{un}}(X)$ of Definition 12.1.8 together with the factorization functor $\mathrm{p}_{\mathcal{J}^{\mathrm{mer}}(Y)\bullet} \circ \mathrm{q}^!_{\mathcal{J}^{\mathrm{mer}}(Y)} : D_{\mathcal{J}^{\mathrm{mer}}(Y)} \to D_{\mathrm{Ran}_{X,\mathrm{un}}}$, we obtain:

*Theorem* 14.2.1. The chiral differential operators $\mathcal{D}^{\mathrm{ch}}(Y)$ on $Y$ admit a canonical factorization $\mathbb{E}_1$-module structure

$$\mathcal{D}^{\mathrm{ch}}(Y) \in \mathcal{A}(Y)\text{-Mod}(\mathrm{Alg}^{\mathrm{fact}}_{\mathrm{un}}(X))$$

over the three dimensional A model factorization $\mathbb{E}_1$-algebra $\mathcal{A}(Y) \in \mathrm{Alg}^{\mathrm{fact}}_{\mathbb{E}_1,\mathrm{un}}(X)$.

*Remark* 14.2.2. In keeping with Warning 14.0.1, the preceding Theorem is only an outline of the result. We give a proof in the case $Y = N/G$ in Theorem 17.0.15 below, defining an action of the $\mathbb{E}_1$ factorization algebra $\mathcal{A}(G, N) \in \mathrm{Alg}^{\mathrm{fact}}_{\mathbb{E}_1,\mathrm{un}}(X)$ constructed in Theorem 13.4.12, following [BFN18] and [BFN19b], on $\mathcal{D}^{\mathrm{ch}}(N/G) \in \mathrm{Alg}^{\mathrm{fact}}_{\mathrm{un}}(X)$ which we study in Section 16 below.

14.3. **Internal construction of chiral differential operators and the action of the three dimensional A model.** There is also an internal variant of the above constructions, following Subsection 8.2 and in particular Example 8.2.9.

*Example* 14.3.1. Following example 8.2.1, the fibre product

$$\mathcal{J}^{\mathrm{mer}}(Y)^{\wedge}_{\mathcal{J}(Y)} := \mathcal{J}(Y)_{\mathrm{dR}} \times_{\mathcal{J}^{\mathrm{mer}}(Y)_{\mathrm{dR}}} \mathcal{J}^{\mathrm{mer}}(Y) \in \mathrm{PreStk}^{\mathrm{fact}}_{\mathrm{un}}(X)$$

defines a unital factorization space over $X$.

*Remark* 14.3.2. Concretely, the prestacks $\mathcal{J}^{\mathrm{mer}}(Y)^{\wedge}_{\mathcal{J}(Y),I} \in \mathrm{PreStk}_{/X^I}$ over $X^I$ assigned to each $I \in \mathrm{fSet}_{\varnothing}$ and the fibre $\mathcal{J}^{\mathrm{mer}}(Y)^{\wedge}_{\mathcal{J}(Y),x} \in \mathrm{PreStk}$ over each $x \in X$ are given by

$$\mathcal{J}^{\mathrm{mer}}(Y)^{\wedge}_{\mathcal{J}(Y),I} = \mathcal{J}(Y)_{I,\mathrm{dR}} \times_{\mathcal{J}^{\mathrm{mer}}(Y)_{I,\mathrm{dR}}} \mathcal{J}^{\mathrm{mer}}(Y)_I = (\mathcal{J}^{\mathrm{mer}}(Y)_I)^{\wedge}_{\mathcal{J}(Y)_I} \qquad \text{, and}$$

$$\mathcal{J}^{\mathrm{mer}}(Y)^{\wedge}_{\mathcal{J}(Y),x} = Y_{\mathcal{O}_X,\mathrm{dR}} \times_{Y_{\mathcal{K}_x,\mathrm{dR}}} Y_{\mathcal{K}_x} = (Y_{\mathcal{K}_x})^{\wedge}_{Y_{\mathcal{O}_x}}$$

*Remark* 14.3.3. More generally, the iterated fibre products

$$\mathcal{J}^{\mathrm{mer}}(Y)^{\wedge}_{(n)} := \mathcal{J}(Y)_{\mathrm{dR}} \times_{\mathcal{J}^{\mathrm{mer}}(Y)_{\mathrm{dR}}} \mathcal{J}(Y)_{\mathrm{dR}} \times_{\mathcal{J}^{\mathrm{mer}}(Y)_{\mathrm{dR}}} \times \ldots \times_{\mathcal{J}^{\mathrm{mer}}(Y)_{\mathrm{dR}}} \mathcal{J}^{\mathrm{mer}}(Y) \in \mathrm{PreStk}^{\mathrm{fact}}_{\mathrm{un}}(X)$$

and projections

$$\pi_{ij} : \mathcal{J}^{\mathrm{mer}}(Y)^{\wedge}_{(n)} \to \mathcal{J}^{\mathrm{mer}}(Y)^{\wedge}_{\mathcal{J}(Y)} \qquad \text{for } i \in \{1, ..., n-1\} \text{ and } j = n, \text{ and}$$

$$\pi_{ij} : \mathcal{J}^{\mathrm{mer}}(Y)^{\wedge}_{(n)} \to \mathcal{Z}(Y)_{\mathrm{dR}} \qquad \text{for } i, j \in \{1, ..., n-1\},$$

define factorization spaces and maps of such for each $n \in \mathbb{N}$.

Following Proposition 8.2.4, analogously to the statement of Proposition 12.2.4 following Proposition 8.1.5, we have:

*Proposition* 14.3.4. The factorization category

$$\mathrm{IndCoh}_{\mathcal{J}^{\mathrm{mer}}(Y)^{\wedge}_{\mathcal{J}(Y)}} \in D^{\star}_{\mathcal{Z}(Y)}\text{-Mod}(\mathrm{Cat}^{\mathrm{fact}}_{\mathrm{un}}(X))$$

is naturally a factorization $\mathbb{E}_1$-module category over $D^{\star}_{\mathcal{Z}(Y)} \in \mathrm{Cat}^{\mathrm{fact}}_{\mathbb{E}_1,\mathrm{un}}(X)$, with respect to the convolution module structure $(\cdot) \star (\cdot) : D_{\mathcal{Z}(Y)} \otimes \mathrm{IndCoh}_{\mathcal{J}^{\mathrm{mer}}(Y)^{\wedge}_{\mathcal{J}(Y)}} \to \mathrm{IndCoh}_{\mathcal{J}^{\mathrm{mer}}(Y)^{\wedge}_{\mathcal{J}(Y)}}$ defined by the composition

$$D_{\mathcal{Z}(Y)} \otimes^{\star} \mathrm{IndCoh}_{\mathcal{J}^{\mathrm{mer}}(Y)^{\wedge}_{\mathcal{J}(Y)}} \xrightarrow{\pi^{\bullet}_{12}\boxtimes\pi^{\bullet}_{23}} \mathrm{IndCoh}_{(\mathcal{J}^{\mathrm{mer}}(Y)^{\wedge}_{(3)})^{\times 2}} \xrightarrow{\Delta^{\bullet}} \mathrm{IndCoh}_{\mathcal{J}^{\mathrm{mer}}(Y)^{\wedge}_{(3)}} \xrightarrow{\pi_{13\bullet}} \mathrm{IndCoh}_{\mathcal{J}^{\mathrm{mer}}(Y)^{\wedge}_{\mathcal{J}(Y)}} \ .$$

Further, the pushforward functors $p_{\mathcal{Z}*} : D_{\mathcal{Z}}^{\star} \to D_{\mathrm{Ran}_{X,\mathrm{un}}}^{\otimes^!}$ as defined in Proposition 12.2.4 and $p_{\mathcal{J}^{\mathrm{mer}}(Y)_{\hat{\mathcal{J}}(Y)}^{\wedge} \bullet} : \mathrm{IndCoh}_{\mathcal{J}^{\mathrm{mer}}(Y)_{\hat{\mathcal{J}}(Y)}^{\wedge}} \to D_{\mathrm{Ran}_{X,\mathrm{un}}}$ define unital, lax compatible factorization functors with respect to the above module structure.

Further, following Corollary 8.2.5, we have:

*Corollary* 14.3.5. The pushforward functor $p_{\mathcal{J}^{\mathrm{mer}}(Y)_{\hat{\mathcal{J}}(Y)}^{\wedge} \bullet} : \mathrm{IndCoh}_{\mathcal{J}^{\mathrm{mer}}(Y)_{\hat{\mathcal{J}}(Y)}^{\wedge}} \to D_{\mathrm{Ran}_{X,\mathrm{un}}}$ induces a functor

$$(14.3.1) \qquad \tilde{\mathcal{A}}(Y)\text{-Mod}(\mathrm{Alg}_{\mathrm{un}}^{\mathrm{fact}}(\mathrm{IndCoh}_{\mathcal{J}^{\mathrm{mer}}(Y)_{\hat{\mathcal{J}}(Y)}^{\wedge}})) \to \mathcal{A}(Y)\text{-Mod}(\mathrm{Alg}_{\mathrm{un}}^{\mathrm{fact}}(X)) \ .$$

*Example* 14.3.6. The map of factorization spaces $\pi_{\hat{\mathcal{J}}(Y)_{\mathrm{dR}}} : \mathcal{J}^{\mathrm{mer}}(Y)_{\hat{\mathcal{J}}(Y)}^{\wedge} \to \mathcal{J}(Y)_{\mathrm{dR}}$ defines a factorization functor

$$\pi^!_{\hat{\mathcal{J}}(Y)_{\mathrm{dR}}} : D_{\mathcal{J}(Y)} \to \mathrm{IndCoh}_{\mathcal{J}^{\mathrm{mer}}(Y)_{\hat{\mathcal{J}}(Y)}^{\wedge}} \ .$$

*Definition* 14.3.7. The internal variant $\tilde{\mathcal{D}}^{\mathrm{ch}}(Y) \in \mathrm{Alg}_{\mathrm{un}}^{\mathrm{fact}}(\mathrm{IndCoh}_{\mathcal{J}^{\mathrm{mer}}(Y)_{\hat{\mathcal{J}}(Y)}^{\wedge}})$ of chiral differential operators is defined by

$$\tilde{\mathcal{D}}^{\mathrm{ch}}(Y) = \pi^!_{\hat{\mathcal{J}}(Y)_{\mathrm{dR}}} \omega_{\mathcal{J}(Y)} = p^!_{\mathcal{J}^{\mathrm{mer}}(Y)_{\hat{\mathcal{J}}(Y)}^{\wedge}} \omega_{\mathrm{Ran}_{X,\mathrm{un}}}$$

*Remark* 14.3.8. Concretely, in terms of the description of Remark 3.2.6, $\tilde{\mathcal{D}}^{\mathrm{ch}}(Y) \in \mathrm{Alg}_{\mathrm{un}}^{\mathrm{fact}}(\mathrm{IndCoh}(\mathcal{J}^{\mathrm{mer}}(Y)_{\hat{\mathcal{J}}(Y)}^{\wedge}))$ has underlying object of $\mathrm{IndCoh}(\mathcal{J}^{\mathrm{mer}}(Y)_{\hat{\mathcal{J}}(Y),I}^{\wedge})$ for each $I \in \mathrm{fSet}_{\varnothing}$, and object of $\mathrm{IndCoh}(\mathcal{J}^{\mathrm{mer}}(Y)_{\hat{\mathcal{J}}(Y),x}^{\wedge})$ for each $x \in X$, given by

$$\tilde{\mathcal{D}}^{\mathrm{ch}}(Y)_I = \pi^!_{\hat{\mathcal{J}}(Y)_{\mathrm{dR},I}} \omega_{\mathcal{J}(Y),I} = \omega_{\mathcal{J}^{\mathrm{mer}}(Y)_{\hat{\mathcal{J}}(Y),I}^{\wedge}} \qquad \in \mathrm{IndCoh}(\mathcal{J}^{\mathrm{mer}}(Y)_{\hat{\mathcal{J}}(Y),I}^{\wedge}), \text{ and}$$

$$\tilde{\mathcal{D}}^{\mathrm{ch}}(Y)_x = \pi^!_{\hat{\mathcal{J}}(Y)_{\mathrm{dR},x}} \omega_{\mathcal{J}(Y)_x} = \omega_{\mathcal{J}^{\mathrm{mer}}(Y)_{\hat{\mathcal{J}}(Y),x}^{\wedge}} \qquad \in \mathrm{IndCoh}(\mathcal{J}^{\mathrm{mer}}(Y)_{\hat{\mathcal{J}}(Y),x}^{\wedge}).$$

Following Example 8.2.9 and Corollary 8.2.5, we have:

*Proposition* 14.3.9. The image of $\tilde{\mathcal{D}}^{\mathrm{ch}}(Y) \in \tilde{\mathcal{A}}(Y)\text{-Mod}(\mathrm{Alg}_{\mathrm{un}}^{\mathrm{fact}}(\mathrm{IndCoh}_{\mathcal{J}^{\mathrm{mer}}(Y)_{\hat{\mathcal{J}}(Y)}^{\wedge}}))$ under the functor of Equation 17.0.1 above is canonically equivalent to $\mathcal{D}^{\mathrm{ch}}(Y) \in A(Y)\text{-Mod}(\mathrm{Alg}_{\mathrm{un}}^{\mathrm{fact}}(X))$.

## 15. Chiral differential operators on schemes

Towards establishing the results outlined in the preceding section in some examples, we recall the construction of chiral differential operators in the case of $Y$ a smooth, finite type variety, following [KV06], and in turn [BD04] and the series of papers [GMS00], [GMS04], and [GMS03]. We employ the theory of $D$ modules in infinite type following [Ras15b] and references therein, as recalled in subappendices B.7 and B.8, and the theory of indcoherent sheaves in infinite type following Section 6 of [Ras20b], and in turn Section 2 of [Gai15], as recalled in Subappendix B.5.

For simplicity, we will by default assume $Y$ is affine throughout, though in Remark 15.1.22 we discuss the obstruction that arises in considering the general case. Further, we will focus on the details of the construction on the fibre over a point $x \in X$, and discuss only a few technical points regarding factorization compatibility, having outlined the general approach in the preceding section and references therein.

15.0.1. *Summary.* In Section 15.1 we recall the construction of the factorization algebra of chiral differential operators on a scheme $Y$, and in Section 15.2, we give the internal variant of the construction.

15.1. **The factorization algebra of chiral differential operators on a scheme.** To begin, we explain that the space $\mathcal{Y}_x = Y_{\mathcal{K}_x}$ defines a reasonable indscheme, in the sense of Subappendix B.3.1.

*Example* 15.1.1. Let $Y = \mathbb{A}^d$, and $\mathcal{Y} = \mathcal{J}^{\mathrm{mer}}(Y) \in \mathrm{PreStk}^{\mathrm{fact}}_{\mathrm{un}}(X)$ be the meromorphic jet space of $Y$, as in the preceding subsection. Following example B.3.5, there is a natural presentation of the fibre space $\mathcal{Y}_x = \mathbb{A}^d_{\mathcal{K}_x} \in \mathrm{PreStk}$ of $\mathcal{Y}$ over $x \in X$ as an indscheme

$$\mathbb{A}^d_{\mathcal{K}} = \operatorname*{colim}_{k \in \mathcal{L}} \lim_{i \in \mathcal{I}} \mathrm{Maps}(\mathbb{D}^\circ_i, \mathbb{A}^d)^{<k} \qquad \text{where}$$

$$\begin{array}{ccc}
\mathrm{Maps}(\mathbb{D}^\circ_j, \mathbb{A}^d)^{<l} & \xrightarrow{\mathrm{rest}^l_{ij}} & \mathrm{Maps}(\mathbb{D}^\circ_i, \mathbb{A}^d)^{<l} \\
{\scriptstyle \iota^{kl}_j} \downarrow & & \downarrow {\scriptstyle \iota^{kl}_i} \\
\mathrm{Maps}(\mathbb{D}^\circ_j, \mathbb{A}^d)^{<k} & \xrightarrow{\mathrm{rest}^k_{ij}} & \mathrm{Maps}(\mathbb{D}^\circ_i, \mathbb{A}^d)^{<k}
\end{array}$$

is the general component of the bidiagram, $\mathbb{K}[x^{\pm 1}]^{<k} \in \mathrm{Sch}_{\mathrm{pft}}$ denotes the subspace of laurent polynomials with poles bounded by $k$, and $\mathrm{Maps}(\mathbb{D}^\circ_i, \mathbb{A}^d)^{<k} := (\mathbb{K}[x^{\pm 1}]^{<k}/x^i)^{\times d} \in \mathrm{Sch}_{\mathrm{ft}}$, for each $i \in \mathcal{I} = \mathbb{N}$ and $j \in \mathcal{K} = \mathbb{N}$.

*Example* 15.1.2. More generally, let $Y$ be a finite type affine scheme, and recall from Example B.3.6 that choosing a closed embedding $Y \hookrightarrow \mathbb{A}^d$ of $Y$ into affine space induces an embedding $Y_{\mathcal{K}} \hookrightarrow \mathbb{A}^d_{\mathcal{K}}$. Thus, choosing a uniformizer, the fibre $\mathcal{Y}_x = Y_{\mathcal{K}_x} \in \mathrm{PreStk}$ of the factorization space $\mathcal{Y} = \mathcal{J}^{\mathrm{mer}}(Y) \in \mathrm{PreStk}^{\mathrm{fact}}_{\mathrm{un}}(X)$ is presented as an ind scheme:

(15.1.1)

$$\mathcal{Y}_x = \operatorname*{colim}_{k \in \mathcal{L}} \lim_{i \in \mathcal{I}} \mathcal{Y}^k_i \qquad \text{where} \qquad \mathcal{Y}^k_i = Y_{\mathcal{K}} \times_{\mathbb{A}^d_{\mathcal{K}}} \mathrm{Maps}(\mathbb{D}^\circ_i, \mathbb{A}^d)^{<k} \in \mathrm{Sch}_{\mathrm{ft}} \ ,$$

$$\begin{array}{ccc}
\mathcal{Y}^l_j & \xrightarrow{\varphi^l_{ij}} & \mathcal{Y}^l_i \\
{\scriptstyle \iota^{kl}_j} \downarrow & & \downarrow {\scriptstyle \iota^{kl}_i} \\
\mathcal{Y}^k_j & \xrightarrow{\varphi^k_{ij}} & \mathcal{Y}^k_i
\end{array}$$

is the general component of the bidiagram, $\varphi^k_{ij}$ is affine, and $\iota^{kl}_i$ is a closed embedding, for each combination of $i, j \in \mathcal{I}$ and $k, l \in \mathcal{L}$. In particular, we have presentations

$$\mathcal{Y}_x = \operatorname*{colim}_{k \in \mathcal{L}} \mathcal{Y}^k = \mathrm{colim} \left[ \dots \leftarrow \mathcal{Y}^l \xleftarrow{\iota^{kl}} \mathcal{Y}^k \leftarrow \dots \right] \qquad \text{and} \qquad \mathcal{Y}^k = \lim_{i \in \mathcal{I}} \mathcal{Y}^k_i = \lim \left[ \dots \to \mathcal{Y}^k_i \xrightarrow{\varphi^k_{ij}} \mathcal{Y}^k_j \to \dots \right] \ .$$

Further, each $\mathcal{Y}^k$ is a reasonable subscheme, and thus the above presents $\mathcal{Y}_x \in \mathrm{IndSch}_{\mathrm{reas}}$ as a reasonable indscheme, in the sense of Definition B.3.7. Further, in favourable situations (e.g. for $Y = N$ a vector space, or $Y = G$ an affine algebraic group), the above defines a placid presentation of each $\mathcal{Y}^k$, in the sense of Definition B.7.9, and thus a placid presentation of $\mathcal{Y}_x$, in the sense of Definition B.8.3.

Next, following subappendices B.7 and B.8, we recall the definition of the category of $D$ modules $D(\mathcal{Y}_x)$ in the present context, and construct the unit object $\mathrm{unit}_{D_{\mathcal{J}^{\mathrm{mer}}(Y)}, x} \in D(\mathcal{Y}_x)$.

*Example* 15.1.3. The category of $D$ modules $D^!(\mathcal{Y}_x) \in \mathrm{DGCat}$ on $\mathcal{Y}_x \in \mathrm{IndSch}_{\mathrm{reas}}$ is presented as

$$D^!(\mathcal{Y}_x) = \lim_{k \in \mathcal{L}} D^!(\mathcal{Y}^k) = \lim \left[ \dots \to D^!(\mathcal{Y}^l) \xrightarrow{\iota^{kl!}} D^!(\mathcal{Y}^k) \to \dots \right] \ ,$$

by Proposition B.8.1. Further, for each $k \in \mathcal{L}$ the category $D^!(\mathcal{Y}^k) \in \mathrm{DGCat}$ of $D$ modules on $\mathcal{Y}^k \in \mathrm{Sch}_{\mathrm{pft}}$ is presented by

$$D^!(\mathcal{Y}^k) = \operatorname*{colim}_{i \in \mathcal{I}} D(\mathcal{Y}_i^k) = \operatorname{colim}\left[\ldots \leftarrow D(\mathcal{Y}_i^k) \xleftarrow{\varphi_{ij}^{k!}} D(\mathcal{Y}_j^k) \leftarrow \ldots\right]$$

by Remark B.7.3.

*Remark* 15.1.4. Concretely, an object $M \in D^!(\mathcal{Y}_x)$ is specified by an assignment

$$k \mapsto \left(M^k \in D^!(\mathcal{Y}^k)\right) \qquad [k \to l] \mapsto \left[\iota^{kl!}M^l \xrightarrow{\cong} M^k\right]$$

defined for each $k, l \in \mathcal{L}$ and map $k \to l$.

Similarly, a typical object $M^k \in D^!(\mathcal{Y}^k)$ is given by a pullback $M^k = \varphi_{i_0}^{k!}M_{i_0}^k$ where $\varphi_{i_0}^{k!} : D(\mathcal{Y}_{i_0}^k) \to D^!(\mathcal{Y}^k)$ is the canonical functor corresponding to pullback along $\varphi_{i_0}^k : \mathcal{Y}^k \to \mathcal{Y}_{i_0}^k$.

*Example* 15.1.5. Following Proposition B.7.11, if the presentation of $\mathcal{Y}^k$ in Example 15.1.2 defines a placid structure, then there is a presentation

$$D^!(\mathcal{Y}^k) = \lim_{i \in \mathcal{I}} D^!(\mathcal{Y}_i^k) = \lim\left[\ldots \to D^!(\mathcal{Y}_i^l) \xrightarrow{\varphi_{ij*}^k[-2d_{\mathcal{Y}_i^k/\mathcal{Y}_j^k}]} D^!(\mathcal{Y}^k) \to \ldots\right].$$

*Remark* 15.1.6. Concretely, if $\mathcal{Y}^k$ admits a placid structure, then an object $M^k \in D^!(\mathcal{Y}^k)$ is given by an assignment

$$i \mapsto (M_i^k \in D(\mathcal{Y}_i^k)) \qquad [i \to j] \mapsto \left[\varphi_{ij*}^k[-2d_{\mathcal{Y}_i^k/\mathcal{Y}_j^k}]M_i^k \xrightarrow{\cong} M_j^k\right].$$

In these terms, the object $\varphi_{i_0}^{k!}M_{i_0}^k \in D^!(\mathcal{Y}^k)$ is given by the assignment

$$(15.1.2) \qquad M_{j_0}^k = \operatorname*{colim}_{i \in \mathcal{I}_{/j_0}^{\mathrm{op}}} \varphi_{ij_0(*,\mathrm{ren})}^k \varphi_{ii_0}^{k!}M_{i_0}^k \qquad \text{under} \qquad \varphi_{jj_0(*,\mathrm{ren})}^k \varphi_{ji_0}^{k!}M_{i_0}^k \to \varphi_{ij_0(*,\mathrm{ren})}^k \varphi_{ii_0}^{k!}M_{i_0}^k,$$

where the structure maps are defined analogously to those in Equation 13.2.1.

*Example* 15.1.7. For each $k$, there is a canonical dualizing sheaf object $\omega_{\mathcal{Y}^k} \in D^!(\mathcal{Y}^k)$ defined by $\omega_{\mathcal{Y}^k} = \varphi_i^{k!}\omega_{\mathcal{Y}_i^k}$ for any $i \in \mathcal{I}$; the definition is independent of $i$ since

$$\varphi_i^{k!}\omega_{\mathcal{Y}_i^k} \cong \varphi_i^{k!}\varphi_{ij}^{k!}\omega_{\mathcal{Y}_j^k} \cong \varphi_j^{k!}\omega_{\mathcal{Y}_j^k}.$$

*Example* 15.1.8. The map $\iota^0 : \mathcal{Y}^0 \to \mathcal{Y}_x$, which is the inclusion $\iota_x : Y_{\mathcal{O}_x} \hookrightarrow Y_{\mathcal{K}_x}$, is ind-proper and thus by Proposition B.7.1 it admits a left adjoint

$$\iota_{*,!}^0 : D^!(\mathcal{Y}^0) \to D^!(\mathcal{Y}_x) \qquad \text{given by} \qquad \iota_{x*,!} : D^!(Y_{\mathcal{O}_x}) \to D^!(Y_{\mathcal{K}_x}).$$

We can now give the formal definition of the unit object of Example 14.1.6 in the present context:

*Definition* 15.1.9. The factorization unit object $\mathrm{unit}_{D_{\mathfrak{J}^{\mathrm{mer}}(Y)}, x} \in D(\mathcal{Y}_x)$ is defined by

$$\mathrm{unit}_{D_{\mathfrak{J}^{\mathrm{mer}}(Y)}, x} = \iota_{x*,!}\omega_{Y_{\mathcal{O}_x}} \in D^!(Y_{\mathcal{K}_x}).$$

*Remark* 15.1.10. Concretely, in terms of the description of Remark 15.1.4, the factorization unit object $\mathrm{unit}_{D_{\mathfrak{J}^{\mathrm{mer}}(Y)}, x} \in D(\mathcal{Y}_x)$ is given by the assignment

$$k \mapsto \iota_{*,!}^{0k}\omega_{\mathcal{Y}^0} \in D^!(\mathcal{Y}^k) \qquad [k \to l] \mapsto \iota^{kl!}\iota_{*,!}^{0l}\omega_{\mathcal{Y}^0} \cong \iota^{kl!}\iota_{*,!}^{kl}\iota_{*,!}^{0k}\omega_{\mathcal{Y}^0} \cong \iota_{*,!}^{0k}\omega_{\mathcal{Y}^0},$$

where the first equivalence is given by functoriality data for the $\iota_{*,!}$ construction and the second by Kashiwara's lemma.

Next, we recall the construction of the category $\mathrm{IndCoh}(\mathcal{Y}_x) \in \mathrm{DGCat}$ of incoherent sheaves on $\mathcal{Y}_x \in \mathrm{IndSch}_{\mathrm{reas}}$ and the forgetful functor $D(\mathcal{Y}_x) \to \mathrm{IndCoh}(\mathcal{Y}_x)$ of Example 14.1.7 in the present context:

*Example* 15.1.11. The category of indcoherent sheaves $\mathrm{IndCoh}^!(\mathcal{Y}_x)$ is presented as

$$\mathrm{IndCoh}^!(\mathcal{Y}_x) = \lim_{k \in \mathcal{L}} \mathrm{IndCoh}^!(\mathcal{Y}^k) = \lim \left[ \ldots \to \mathrm{IndCoh}^!(\mathcal{Y}^l) \xrightarrow{\iota^{kl!}} \mathrm{IndCoh}^!(\mathcal{Y}^k) \to \ldots \right],$$

by Remark B.5.7. Further, for each $k \in \mathcal{L}$ the category $\mathrm{IndCoh}^!(\mathcal{Y}^k) \in \mathrm{DGCat}$ of indcoherent sheaves on $\mathcal{Y}^k \in \mathrm{Sch}_{\mathrm{pft}}$ is presented by

$$\mathrm{IndCoh}^!(\mathcal{Y}^k) = \operatorname*{colim}_{i \in \mathcal{I}} \mathrm{IndCoh}(\mathcal{Y}_i^k) = \operatorname*{colim} \left[ \ldots \leftarrow \mathrm{IndCoh}(\mathcal{Y}_i^k) \xleftarrow{\varphi_{ij}^{k!}} \mathrm{IndCoh}(\mathcal{Y}_j^k) \leftarrow \ldots \right],$$

by Remark B.5.5.

As in Remark 15.1.4, we have:

*Remark* 15.1.12. Concretely, an object $\mathcal{F} \in \mathrm{IndCoh}^!(\mathcal{Y}_x)$ is specified by an assignment

$$k \mapsto \left( \mathcal{F}^k \in \mathrm{IndCoh}^!(\mathcal{Y}^k) \right) \qquad [k \to l] \mapsto \left[ \iota^{kl!} \mathcal{F}^l \xrightarrow{\cong} \mathcal{F}^k \right]$$

defined for each $k, l \in \mathcal{L}$ and map $k \to l$.

Further, a typical object $\mathcal{F}^k \in \mathrm{IndCoh}^!(\mathcal{Y}^k)$ is given by $\mathcal{F}^k = \varphi_i^{k!} \mathcal{F}_i^k$ under $\varphi_i^k : \mathcal{Y}^k \to \mathcal{Y}_i^k$.

*Proposition* 15.1.13. There is a natural forgetful functor

$$\mathrm{oblv}_{Y_{\mathcal{K}}}^r := \mathrm{q}_{Y_{\mathcal{K}}}^! : D^!(\mathcal{Y}_x) \to \mathrm{IndCoh}^!(\mathcal{Y}_x)$$

which restricts to the usual forgetful functor of Definition A.5.16 on each finite type approximation $\mathcal{Y}_i^k \in \mathrm{Sch}_{\mathrm{ft}}$.

*Proof.* Recall from Proposition A.5.23 that the forgetful functor $\mathrm{oblv}^r := \mathrm{q}_X^! : D(X) \to \mathrm{IndCoh}(X)$ intertwines the $D$ module and incoherent pullback functors $f^!$ for $f : X \to Y$ a map of finite type schemes. Thus, by comparing the descriptions of Examples 15.1.3 and 15.1.11 above, this functor induces a forgetful functor on the colimit categories $\mathrm{oblv}_{\mathcal{Y}^k}^r := \mathrm{q}_{\mathcal{Y}_k}^! : D^!(\mathcal{Y}^k) \to \mathrm{IndCoh}^!(\mathcal{Y}^k)$ and in turn on the limit categories as desired. $\qquad\square$

*Remark* 15.1.14. Concretely, in terms of the descriptions in remarks 15.1.4 and 15.1.12, the functor $\mathrm{oblv}_{Y_{\mathcal{K}}}^r := \mathrm{q}_{Y_{\mathcal{K}}}^! : D^!(\mathcal{Y}_x) \to \mathrm{IndCoh}^!(\mathcal{Y}_x)$ is defined on objects by

$$\left( M^k \in D^!(\mathcal{Y}^k) \right)_{k \in \mathcal{L}} \mapsto \left( \mathcal{F}^k := \mathrm{oblv}_{\mathcal{Y}^k}^r (M^k) \in \mathrm{IndCoh}(\mathcal{Y}^k) \right)_{k \in \mathcal{L}}.$$

To define the factorization algebra $\mathcal{D}_Y^{\mathrm{ch}} \in \mathrm{Alg}_{\mathrm{un}}^{\mathrm{fact}}(X)$ of chiral differential operators following the outline in the preceding subsection, it remains to define the pushforward functor $\mathrm{p}_{Y_{\mathcal{K}}} \bullet : \mathrm{IndCoh}(\mathcal{Y}_{\mathcal{K}}) \to \mathrm{Vect}$. Following subsection B.5, there is a second, in general distinct, variant of the category of incoherent sheaves on $Y_{\mathcal{K}}$, which canonically admits such a pushforward functor:

*Example* 15.1.15. The category $\mathrm{IndCoh}^\bullet(\mathcal{Y}_x) \in \mathrm{DGCat}$ of indcoherent sheaves on $Y_{\mathcal{K}_x} \in \mathrm{IndSch}_{\mathrm{reas}}$ is presented as

(15.1.3)
$$\mathrm{IndCoh}^\bullet(\mathcal{Y}_x) = \operatorname*{colim}_{k \in \mathcal{L}} \mathrm{IndCoh}^\bullet(\mathcal{Y}^k) = \mathrm{colim} \left[ \ldots \leftarrow \mathrm{IndCoh}^\bullet(\mathcal{Y}^k) \xleftarrow{\iota^{kl}_\bullet} \mathrm{IndCoh}^\bullet(\mathcal{Y}^l) \leftarrow \ldots \right] \text{, or}$$

(15.1.4)
$$\mathrm{IndCoh}^\bullet(\mathcal{Y}_x) = \lim_{k \in \mathcal{L}} \mathrm{IndCoh}^\bullet(\mathcal{Y}^k) = \lim \left[ \ldots \to \mathrm{IndCoh}^\bullet(\mathcal{Y}^k) \xrightarrow{\iota^{kl,!}} \mathrm{IndCoh}^\bullet(\mathcal{Y}^l) \to \ldots \right] .$$

where the former description follows from Remark B.5.7, and the latter is equivalent by passing to right adjoints as in Remark B.5.9, which exist by Proposition B.5.8.

*Remark* 15.1.16. Concretely, an object $\mathcal{F} \in \mathrm{IndCoh}^\bullet(\mathcal{Y}_x)$ is given by an assignment
$$k \mapsto \left( \mathcal{F}^k \in \mathrm{IndCoh}^\bullet(\mathcal{Y}^k) \right) \qquad [l \to k] \mapsto \left[ \iota^{kl,!}\mathcal{F}^l \xrightarrow{\cong} \mathcal{F}^k \right]$$
for each $k, l \in \mathcal{L}$ and arrow $l \to k$ in $\mathcal{L}$.

*Example* 15.1.17. Further, for each $\mathcal{Y}^k \in \mathrm{Sch}_{\mathrm{pft}}$, the presentation $\mathcal{Y}^k = \lim_{i \in \mathcal{I}} \mathcal{Y}^k_i$ as a limit of finite type schems $\mathcal{Y}^k_i \in \mathrm{Sch}_{\mathrm{ft}}$ induces presentations

(15.1.5)
$$\mathrm{IndCoh}^\bullet(\mathcal{Y}^k) = \lim_{i \in \mathcal{I}} \mathrm{IndCoh}^\bullet(\mathcal{Y}^k_i) = \lim \left[ \ldots \to \mathrm{IndCoh}(\mathcal{Y}^k_j) \xrightarrow{\varphi^k_{ij\bullet}} \mathrm{IndCoh}(\mathcal{Y}^k_i) \to \ldots \right] \text{, and}$$

(15.1.6)
$$\mathrm{IndCoh}^\bullet(\mathcal{Y}^k) = \operatorname*{colim}_{i \in \mathcal{I}} \mathrm{IndCoh}^\bullet(\mathcal{Y}^k_i) = \mathrm{colim} \left[ \ldots \leftarrow \mathrm{IndCoh}(\mathcal{Y}^k_j) \xleftarrow{\varphi^{k\bullet}_{ij}} \mathrm{IndCoh}(\mathcal{Y}^k_i) \leftarrow \ldots \right] ,$$

where the latter follows by passing to left adjoints, which exist by Proposition 3.1.6 of [GR17a].

*Remark* 15.1.18. Concretely, an object $\mathcal{F}^k \in \mathrm{IndCoh}^\bullet(\mathcal{Y}^k)$ is given by an assignment
$$i \mapsto \left( \mathcal{F}^k_i \in \mathrm{IndCoh}(\mathcal{Y}^k_i) \right) \qquad [j \to i] \mapsto [\varphi^k_{ij\bullet}\mathcal{F}^k_i \xrightarrow{\cong} \mathcal{F}^k_j]$$
for each $i, j \in \mathcal{I}$ and arrow $j \to i$ in $\mathcal{I}$.

*Example* 15.1.19. There is a canonical pushforward functor
$$\mathrm{p}_{\mathcal{Y}^k\bullet} : \mathrm{IndCoh}^\bullet(\mathcal{Y}^k) \to \mathrm{Vect} \qquad \text{defined by} \qquad \mathcal{F}^k = (\mathcal{F}^k_i)_{i \in \mathcal{I}} \mapsto \mathrm{p}_{\mathcal{Y}^k_i\bullet}\mathcal{F}^k_i \text{ ;}$$
the result is independent of $i$ since
$$\mathrm{p}_{\mathcal{Y}^k_i\bullet}\mathcal{F}^k_i \cong \mathrm{p}_{\mathcal{Y}^k_j\bullet}\varphi^k_{ij\bullet}\mathcal{F}^k_i \cong \mathrm{p}_{\mathcal{Y}^k_j\bullet}\mathcal{F}^k_j \text{ .}$$

Moreover, we have:

*Example* 15.1.20. There is a pushforward functor
$$\mathrm{p}_{\mathcal{Y}_x\bullet} : \mathrm{IndCoh}^\bullet(\mathcal{Y}_x) \to \mathrm{Vect} \qquad \text{defined by} \qquad \mathcal{F} = (\mathcal{F}^k)_{k \in \mathcal{J}} \mapsto \operatorname*{colim}_{k \in \mathcal{L}} \mathrm{p}_{\mathcal{Y}^k\bullet}\mathcal{F}^k$$
where the structure maps in the colimit diagram are defined for each $k < l$ by

(15.1.7)
$$\mathrm{p}_{\mathcal{Y}^k\bullet}\mathcal{F}^k \cong \mathrm{p}_{\mathcal{Y}^l\bullet}\iota^{kl}_\bullet\mathcal{F}^k \cong \mathrm{p}_{\mathcal{Y}^l\bullet}\iota^{kl}_\bullet\iota^{kl!}\mathcal{F}^l \to \mathrm{p}_{\mathcal{Y}^l\bullet}\mathcal{F}^l \text{ ,}$$

where the map is given by the counit of the $(\iota^{kl}_\bullet, \iota^{kl!})$ adjunction.

To conclude the construction of $\mathcal{D}_Y^{\mathrm{ch}}$, it remains to construct an equivalence $\mathrm{IndCoh}^!(Y_{\mathcal{K}}) \cong \mathrm{IndCoh}^\bullet(Y_{\mathcal{K}})$ so that we can apply the pushforward functor of the preceding Example 15.1.20 to the object of Example 15.1.9; we now explain a putative construction of such an equivalence and the potential obstruction to its existence:

*Remark* 15.1.21. Comparing the limit descriptions of examples 15.1.11 and 15.1.15, it suffices to construct equivalences $\mathrm{IndCoh}^!(\mathcal{Y}^k) \cong \mathrm{IndCoh}^\bullet(\mathcal{Y}^k)$ compatible with $\iota^{kl!}$ for each $k \in \mathcal{L}$.

Applying commutativity of the diagram in Equation 5.23 of [Gai15] together with the discussion in Section 3.1.5 of [GR17a], for each $X, Y \in \mathrm{DGSch}_{\mathrm{ft}}$ and $f : X \to Y$ eventually coconnective, we have the commutative diagram

(15.1.8)

$$
\begin{array}{ccc}
\mathrm{IndCoh}(Y) & \xrightarrow{(\cdot)\otimes_{\mathcal{O}_Y}\omega_Y} & \mathrm{IndCoh}(Y) \\
\downarrow{\scriptstyle f^\bullet} & & \downarrow{\scriptstyle f^!} \\
\mathrm{IndCoh}(X) & \xrightarrow{(\cdot)\otimes_{\mathcal{O}_X}\omega_X} & \mathrm{IndCoh}(X)
\end{array}
\qquad \text{where} \qquad (\cdot) \otimes_{\mathcal{O}_X} \omega_X := \Upsilon_X \circ \Psi_X : \mathrm{IndCoh}(X) \to \mathrm{IndCoh}(X) \ .
$$

This is essentially just the standard identification $f^! \cong f^\bullet \otimes_{\mathcal{O}_X} \omega_{X/Y}$. Thus, we obtain identifications

(15.1.9)

$$
\begin{array}{ccc}
\mathrm{IndCoh}^!(\mathcal{Y}^k) & & \\
\downarrow{\scriptstyle \omega_{\mathcal{Y}^k}^{-1}} & \text{induced by} & \\
\mathrm{IndCoh}^\bullet(\mathcal{Y}^k) & &
\end{array}
\qquad
\begin{array}{ccccc}
\ldots \longleftarrow & \mathrm{IndCoh}(\mathcal{Y}_i^k) & \xleftarrow{\varphi_{ij}^{k!}} & \mathrm{IndCoh}(\mathcal{Y}_j^k) & \longleftarrow \ldots \\
& \downarrow{\scriptstyle \omega_{\mathcal{Y}_i^k}^{-1}} & & \downarrow{\scriptstyle \omega_{\mathcal{Y}_j^k}^{-1}} & \\
\ldots \longleftarrow & \mathrm{IndCoh}(\mathcal{Y}_j^k) & \xleftarrow{\varphi_{ij}^{k\bullet}} & \mathrm{IndCoh}(\mathcal{Y}_i^k) & \longleftarrow \ldots
\end{array} \ .
$$

However, these are not all compatible with the functors $\iota^{kl!}$, in the sense that

(15.1.10)

$$
\begin{array}{ccc}
\ldots \longrightarrow \mathrm{IndCoh}^!(\mathcal{Y}^k) & \xrightarrow{\iota^{kl!}} & \mathrm{IndCoh}^!(\mathcal{Y}^l) \longrightarrow \ldots \\
\downarrow{\scriptstyle \omega_{\mathcal{Y}^k}^{-1}} & & \downarrow{\scriptstyle \omega_{\mathcal{Y}^l}^{-1}} \\
\ldots \longrightarrow \mathrm{IndCoh}^\bullet(\mathcal{Y}^k) & \xrightarrow{\iota^{kl!}} & \mathrm{IndCoh}^\bullet(\mathcal{Y}^l) \longrightarrow \ldots
\end{array}
\qquad \text{does } not \text{ commute.}
$$

Instead, we use the equivalence of Equation 15.1.9 to identify $\mathrm{IndCoh}^!(\mathcal{Y}^0) \cong \mathrm{IndCoh}^\bullet(\mathcal{Y}^0)$, and proceed by attempting to extend in a way which ensures compatibility:

*Remark* 15.1.22. In sections 1.7 and 1.8 of [KV06], it is explained that the pushforward functor of Example 15.1.20 above can be modified to define a functor $\mathrm{p}_{\mathcal{Y}_x, \tilde{\varsigma}} : \mathrm{IndCoh}^!(\mathcal{Y}_x) \to \mathrm{Vect}$ as desired, given a choice of global object of the $\mathcal{O}_{\mathcal{Y}_x}^\times$ gerbe $\mathrm{Det}_Y$ defined in *loc. cit.*; we recall this construction in Remark 15.1.29 below. Moreover, it is explained in loc. cit. that this construction induces an equivalence of (factorization) gerbes $\mathrm{Det}_{Y,\mathrm{fact}} \xrightarrow{\cong} \mathrm{CDO}_Y$ with that of locally trivial sheaves of chiral differential operators on $Y$.

Alternatively, in 3.9.20 of [BD04], it is explained that there is an equivalence $\mathrm{CDO}(Y) \xrightarrow{\cong} \mathrm{Tate}(Y)$ of the (sheaf of) groupoid(s) of sheaves of chiral differential operators on $Y$ with that of Tate structures on $Y$, which we now define:

First, we introduce some notation, following 1.4.16 and 2.5.16 in [BD04]:

*Example* 15.1.23. Let $\mathcal{Y} \in D\text{-Sch}(X)$ be a $D$ scheme on $X$, $V$ a vector $\mathcal{D}_X$ bundle on $\mathcal{Y}$, and let $\mathcal{E}(V)$ be the universal Lie* algebroid over $\mathcal{Y}$ acting on $V$. Concretely, the underlying relative Lie* algebra $\mathcal{E}(V) \in \mathrm{Lie}^*(\mathcal{Y})$ over $\mathcal{Y}$ is an extension

$$\mathfrak{gl}(V) \hookrightarrow \mathcal{E}(V) \xrightarrow{\tau} \Theta_{\mathcal{Y}}$$

where $\tau$ is the anchor map, and $\mathfrak{gl}(V) \in \mathrm{Lie}^*(\mathcal{Y})$ is the relative analogue of $\mathfrak{gl}(V) \in \mathrm{Lie}^*(X)$ as in Definition I-13.0.5.

Further, for a Lie* algebroid $\mathcal{L}$ over $\mathcal{Y}$, we let $\mathcal{L}^{\diamond} = \ker(\tau_L) \in \mathrm{Lie}^*(\mathcal{Y})$ denote the relative Lie* algebra given by the kernel of the anchor map, so that in particular $\mathcal{E}(V)^{\diamond} = \mathfrak{gl}(V)$.

Now, following 2.8.1 in [BD04], we recall the definition of Tate structure:

*Definition* 15.1.24. A Tate structure on $V$ over $\mathcal{Y}$ is given by:

- an extension of Lie* algebroids on $\mathcal{Y}$

$$\omega_{\mathcal{Y}} \hookrightarrow \mathcal{E}(V)^{\flat} \twoheadrightarrow \mathcal{E}(V) \qquad \text{, and}$$

- an isomorphism of $\mathrm{Lie}^*(\mathcal{Y})$ algebras $\mathcal{E}(V)^{\flat,\diamond} \cong \mathfrak{gl}(V)^{\flat}$, lifting the natural identification $\mathcal{E}(V)^{\diamond} \cong \mathfrak{gl}(V)$ and compatible with the natural adjoint actions of $\mathcal{E}(V)^{\flat}$ on the two $\mathrm{Lie}^*(\mathcal{Y})$ algebras, where we $\mathfrak{gl}(V)^{\flat}$ denotes the Tate extension of $\mathfrak{gl}(V)$ as in Definition I-13.0.5.

Let $\mathrm{Tate}(\mathcal{Y}, V)$ denote the sheaf of groupoids on $\mathcal{Y}$ defined by the space of Tate structures on restrictions of $V$. In particular, for $\mathcal{Y} \in D\text{-}\mathrm{Sch}(X)$ a smooth $D$ scheme on $X$, $\Theta_{\mathcal{Y}}$ defines a vector $\mathcal{D}_X$ bundle on $\mathcal{Y}$, and we let $\mathrm{Tate}(\mathcal{Y}) = \mathrm{Tate}(\mathcal{Y}; \Theta_{\mathcal{Y}})$ denote the sheaf of groupoids of Tate structures on $\Theta_{\mathcal{Y}}$.

*Warning* 15.1.25. Throughout the remainder of this subsection, we assume that $\mathcal{Y} = \mathcal{J}(Y) \in D\text{-}\mathrm{Sch}(X)$ is equipped with a Tate structure. However, we will supress the dependence on the choice of Tate structure from the notation throughout, for example letting $\mathcal{D}^{\mathrm{ch}}(Y) \in \mathrm{Alg}^{\mathrm{fact}}_{\mathrm{un}}(X)$ denote the resulting factorization algebra of chiral differential operators.

For concreteness, following [KV06], we now explain the construction of the resulting functor

$$\mathrm{p}_{\mathcal{Y}_x \bullet} : \mathrm{IndCoh}^!(\mathcal{Y}_x) \to \mathrm{Vect}$$

analogous to that of Example 15.1.20, since this is all that is required to define the factorization algebra $\mathcal{D}^{\mathrm{ch}}(Y) \in \mathrm{Alg}^{\mathrm{fact}}_{\mathrm{un}}(X)$.

*Example* 15.1.26. For each $k \in \mathcal{L}$, the equivalence of Equation 15.1.9 above induces an equivalent presentation

$$\mathrm{IndCoh}^!(\mathcal{Y}^k) = \lim \left[ \dots \longrightarrow \mathrm{IndCoh}(\mathcal{Y}^k_i) \xrightarrow{\underline{\mathrm{Hom}}(\omega_{\mathcal{Y}^k_i/\mathcal{Y}^k_j}, \cdot)} \mathrm{IndCoh}(\mathcal{Y}^k_j) \longrightarrow \dots \right]$$

$$\downarrow{\scriptstyle (\cdot)\otimes_{\mathcal{O}_{\mathcal{Y}^k}} \omega_{\mathcal{Y}^k}^{-1}} \qquad \downarrow{\scriptstyle (\cdot)\otimes_{\mathcal{O}_{\mathcal{Y}^k_i}} \omega_{\mathcal{Y}^k_i}^{-1}} \qquad \downarrow{\scriptstyle (\cdot)\otimes_{\mathcal{O}_{\mathcal{Y}^k_j}} \omega_{\mathcal{Y}^k_j}^{-1}}$$

$$\mathrm{IndCoh}^{\bullet}(\mathcal{Y}^k) = \lim \left[ \dots \longrightarrow \mathrm{IndCoh}(\mathcal{Y}^k_j) \xrightarrow{\varphi^k_{ij\bullet}} \mathrm{IndCoh}(\mathcal{Y}^k_i) \longrightarrow \dots \right]$$

where

$$\underline{\mathrm{Hom}}(\omega_{\mathcal{Y}^k_i/\mathcal{Y}^k_j}, \cdot) = \varphi^k_{ij\bullet} \circ \left( (\cdot) \otimes_{\mathcal{O}_{\mathcal{Y}^k_i}} \omega^{-1}_{\mathcal{Y}^k_i/\mathcal{Y}^k_j} \right) : \mathrm{IndCoh}(\mathcal{Y}^k_i) \to \mathrm{IndCoh}(\mathcal{Y}^k_j) \ .$$

Following Remark 15.1.18, we have:

*Remark* 15.1.27. Concretely, an object $\mathcal{F}^k \in \mathrm{IndCoh}^!(\mathcal{Y}^k)$ is given by an assignment

$$i \mapsto \Big( \mathcal{F}_i^k \in \mathrm{IndCoh}(\mathcal{Y}_i^k) \Big) \qquad [j \to i] \mapsto \Big[ \underline{\mathrm{Hom}}(\omega_{\mathcal{Y}_i^k / \mathcal{Y}_j^k}, \mathcal{F}_i^k) \xrightarrow{\cong} \mathcal{F}_j^k \Big]$$

for each $i, j \in \mathcal{I}$ and arrow $j \to i$ in $\mathcal{I}$.

Further, following Example 15.1.19, we have:

*Example* 15.1.28. In particular, the presentation of Example 15.1.26 above induces a functor

$$\mathrm{Hom}(\omega_{\omega_{\mathcal{Y}^k}}, \cdot) := \mathrm{p}_{\mathcal{Y}^k \bullet} \circ \Big( (\cdot) \otimes_{\mathcal{O}_{\mathcal{Y}_k}} \omega_{\mathcal{Y}_k}^{-1} \Big) : \mathrm{IndCoh}^!(\mathcal{Y}^k) \to \mathrm{Vect} \qquad \text{defined by} \qquad \mathcal{F}^k = (\mathcal{F}_i^k)_{i \in \mathcal{I}} \mapsto \mathrm{Hom}(\omega_{\mathcal{Y}_i^k}, \mathcal{F}_i^k) \, ;$$

the result in independent of $i$ since

$$\mathrm{Hom}(\omega_{\mathcal{Y}_i^k}, \mathcal{F}_i^k) \cong \mathrm{Hom}(\omega_{\mathcal{Y}_j^k}, \underline{\mathrm{Hom}}(\omega_{\mathcal{Y}_i^k / \mathcal{Y}_j^k}, \mathcal{F}_i^k)) \cong \mathrm{Hom}(\omega_{\mathcal{Y}_j^k}, \mathcal{F}_j^k) \ .$$

*Remark* 15.1.29. The analogue of the non-commutativity of the diagram in Equation 15.1.10 in this context is that the functors of the preceding example do *not* satisfy the compatibility condition $\mathrm{Hom}(\omega_{\mathcal{Y}^l}, \iota_\bullet^{kl}(\cdot)) \cong \mathrm{Hom}(\omega_{\mathcal{Y}^k}, \cdot)$ and thus the construction of Example 15.1.20 can not be repeated without some additional data.

Following sections 1.7 and 1.8 of [KV06], as in Remark 15.1.22, a choice of global object $\mathcal{E} \in \mathrm{Det}_{\mathcal{Y}_x}$ defines the desired pushforward functor as follows: given $\mathcal{E} = (\mathcal{E}^k)_{k \in \mathcal{L}}$ where

$$\mathcal{E}^k = (\mathcal{E}_i^k)_{i \in \mathcal{I}} \qquad \text{such that} \qquad \iota_i^{kl \bullet} \mathcal{E}_i^l \cong \omega_{\mathcal{Y}_i^k / \mathcal{Y}_i^l}^{-1} \otimes \mathcal{E}_i^k$$

we define the functors

$$\mathrm{p}_{\mathcal{Y}^k \bullet} := \mathrm{Hom}(\omega_{\mathcal{Y}^k}, \mathcal{E}^k \otimes_{\mathcal{O}_{\mathcal{Y}^k}} (\cdot)) : \mathrm{IndCoh}^!(\mathcal{Y}^k) \to \mathrm{Vect} \ .$$

These satisfy the desired relation

$$\begin{aligned}
\mathrm{Hom}(\omega_{\mathcal{Y}^l}, \mathcal{E}^l \otimes \iota_\bullet^{kl} \mathcal{F}^k) &\cong \mathrm{Hom}(\omega_{\mathcal{Y}^l}, \iota_\bullet^{kl}(\iota^{kl \bullet} \mathcal{E}^l \otimes \mathcal{F}^k)) \\
&\cong \mathrm{Hom}(\omega_{\mathcal{Y}^l}, \iota_\bullet^{kl}(\omega_{\mathcal{Y}^k / \mathcal{Y}^l}^{-1} \otimes \mathcal{E}^k \otimes \mathcal{F}^k)) \\
&\cong \mathrm{Hom}(\omega_{\mathcal{Y}^k}, \mathcal{E}^k \otimes \mathcal{F}^k)
\end{aligned}$$

and thus for $\mathcal{F} = (\mathcal{F}^k)_{k \in \mathcal{L}} \in \mathrm{IndCoh}^!(\mathcal{Y}_x)$ we have canonical maps

$$\mathrm{Hom}(\omega_{\mathcal{Y}^k}, \mathcal{E}^k \otimes_{\mathcal{O}_{\mathcal{Y}^k}} \mathcal{F}^k) \cong \mathrm{Hom}(\omega_{\mathcal{Y}^l}, \mathcal{E}^l \otimes_{\mathcal{O}_{\mathcal{Y}^l}} \iota_\bullet^{kl} \iota^{kl!} \mathcal{F}^l) \to \mathrm{Hom}(\omega_{\mathcal{Y}^l}, \mathcal{E}^l \otimes_{\mathcal{O}_{\mathcal{Y}^l}} \mathcal{F}^l)$$

in analogy with those of Equation 15.1.7.

Thus, we make the following definition:

*Definition* 15.1.30. The global sections functor

$$\mathrm{p}_{\mathcal{Y}_x, \bullet} : \mathrm{IndCoh}^!(\mathcal{Y}_x) \to \mathrm{Vect} \qquad \text{is defined by} \qquad \mathcal{F} = (\mathcal{F}^k)_{k \in \mathcal{L}} \mapsto \mathrm{colim}_{k \in \mathcal{L}} \mathrm{p}_{\mathcal{Y}^k \bullet} \mathcal{F}^k \ ,$$

where the structure maps are as in the preceding Remark.

*Remark* 15.1.31. We encourage the reader to compare the preceding construction with that of the semi-infinite cohomology functor in Example 16.2.14 below.

Moreover, we can now define the factorization algebra of chiral differential operators, following the tentative Definition 14.1.10:

*Definition* 15.1.32. The chiral differential operators $\mathcal{D}^{\mathrm{ch}}(Y) \in \mathrm{Alg}_{\mathrm{un}}^{\mathrm{fact}}(X)$ are defined by

$$\mathcal{D}^{\mathrm{ch}}(Y) = \mathrm{p}_{\mathcal{J}^{\mathrm{mer}}(Y) \bullet} \, \mathrm{oblv}_{\mathcal{J}^{\mathrm{mer}}(Y)}^r \, \iota_{\mathcal{J}(Y) *} \omega_{\mathcal{J}(Y)} \ \in \mathrm{Alg}_{\mathrm{un}}^{\mathrm{fact}}(X) \ .$$

*Remark* 15.1.33. Concretely, the underlying vector space $\mathcal{D}^{\mathrm{ch}}(Y)_x \in \mathrm{Vect}$ is given by

$$\mathcal{D}^{\mathrm{ch}}(Y)_x = \operatorname*{colim}_{k \in \mathcal{L}} \mathrm{p}_{\mathrm{y}^k{}_\bullet} \circ \mathrm{oblv}^r_{\mathrm{y}^k} \circ \iota^{0k}_{*!} \omega_{\mathrm{y}^0}$$

$$\cong \operatorname*{colim}_{k \in \mathcal{L}} \mathrm{Hom}(\omega_{\mathrm{y}^k_j}, \mathcal{E}^k_j \otimes (\cdot)) \circ \mathrm{oblv}^r_{\mathrm{y}^k_j} \circ \iota^{0k}_{j*} (\omega_{\mathrm{y}^0})_j \qquad \text{, for any } j \in \mathcal{J} \text{ by Example 15.1.28,}$$

$$\cong \operatorname*{colim}_{k \in \mathcal{L}, i \in \mathcal{J}^{\mathrm{op}}} \mathrm{Hom}(\omega_{\mathrm{y}^k_j}, \mathcal{E}^k_j \otimes (\cdot)) \circ \mathrm{oblv}^r_{\mathrm{y}^k_j} \circ \iota^{0k}_{j*} \circ \varphi^0_{ij*} \omega_{\mathrm{y}^0_i} \text{, by Equation 15.1.2.}$$

*Example* 15.1.34. In the case $Y = N$ is given by a vector space, the underlying vector space is calculated explicitly as

$$\mathcal{D}^{\mathrm{ch}}(N)_x = \operatorname*{colim}_{k \in \mathcal{L}, i \in \mathcal{J}^{\mathrm{op}}} \mathrm{Sym}^\bullet_{\mathbb{K}}\left(\left[N \otimes_{\mathbb{K}} z^{-1}\mathbb{K}[z^{-1}]/(z^{-k})\right] \oplus \left[N^\vee \otimes_{\mathbb{K}} z^{-1}\mathbb{K}[z^{-1}]/(z^{-i}) \cdot dz\right]\right)$$

$$\cong \mathrm{Sym}^\bullet_{\mathbb{K}}(\left[N \otimes_{\mathbb{K}} z^{-1}\mathbb{K}[z^{-1}]\right] \oplus \left[N^\vee \otimes_{\mathbb{K}} z^{-1}\mathbb{K}[z^{-1}]\right] \cdot dz)$$

just as in Equation I-12.1.1. Indeed, the chiral differential operators $\mathcal{D}^{\mathrm{ch}}(N) \in \mathrm{Alg}^{\mathrm{fact}}_{\mathrm{un}}(X)$ in the case $X = \mathbb{A}^1$ admit a canonical weak $\mathbb{G}_a$ equivariant structure, and under the equivalence of Theorem I-8.0.5 the corresponding vertex algebra is given by the chiral Weyl algebra of Example I-12.1.10.

15.2. **Internal construction of chiral differential operators on schemes.** Alternatively, we now give the internal variant of the above construction, as outlined in Subsection 14.3. As above, we give the construction over a fixed point $x \in X$, having outlined the factorization compatible analogues in *loc. cit.*.

*Example* 15.2.1. The factorization space $\mathcal{J}^{\mathrm{mer}}(Y)^{\wedge}_{\mathcal{J}(Y)} \in \mathrm{PreStk}^{\mathrm{fact}}_{\mathrm{un}}(X)$ of Example 14.3.1 has fibre at each $x \in X$ given by

$$\mathcal{J}^{\mathrm{mer}}(Y)^{\wedge}_{\mathcal{J}(Y)x} = (\mathcal{Y}_x)^{\wedge}_{\mathrm{y}^0} = \mathcal{Y}^0_{\mathrm{dR}} \times_{\mathrm{y}_{x,\mathrm{dR}}} \mathcal{Y}_x \in \mathrm{PreStk}$$

which admits a presentation

(15.2.1)

$$(\mathcal{Y}_x)^{\wedge}_{\mathrm{y}^0} = \operatorname*{colim}_{k \in \mathcal{L}} (\mathcal{Y}^k)^{\wedge}_{\mathrm{y}^0} = \mathrm{colim}\left[\ldots \leftarrow (\mathcal{Y}^l)^{\wedge}_{\mathrm{y}^0} \xleftarrow{(\iota^{kl})^{\wedge}_{\mathrm{y}^0}} (\mathcal{Y}^k)^{\wedge}_{\mathrm{y}^0} \leftarrow \ldots\right] \quad \text{under} \quad (\iota^{kl})^{\wedge}_{\mathrm{y}^0} : (\mathcal{Y}^k)^{\wedge}_{\mathrm{y}^0} \hookrightarrow (\mathcal{Y}^l)^{\wedge}_{\mathrm{y}^0} \text{, with}$$

(15.2.2)

$$(\mathcal{Y}^k)^{\wedge}_{\mathrm{y}^0} = \operatorname*{lim}_{i \in \mathcal{J}} (\mathcal{Y}^k_i)^{\wedge}_{\mathrm{y}^0_i} = \lim\left[\ldots \to (\mathcal{Y}^k_i)^{\wedge}_{\mathrm{y}^0_i} \xrightarrow{(\varphi^k_{ij})^{\wedge}_{\mathrm{y}^0}} (\mathcal{Y}^k_j)^{\wedge}_{\mathrm{y}^0_j} \to \ldots\right] \quad \text{under} \quad (\varphi^k_{ij})^{\wedge}_{\mathrm{y}^0} : (\mathcal{Y}^k_i)^{\wedge}_{\mathrm{y}^0_i} \to (\mathcal{Y}^k_j)^{\wedge}_{\mathrm{y}^0_j},$$

by Example 15.1.2, so that $(\mathcal{Y}_x)^{\wedge}_{\mathrm{y}^0}$ defines an ind-pro-finite type formal scheme.

Following Example 15.1.11, we have:

*Example* 15.2.2. The category $\mathrm{IndCoh}^!((\mathcal{Y}_x)^{\wedge}_{\mathrm{y}^0}) \in \mathrm{DGCat}$ is presented by

(15.2.3)

$$\mathrm{IndCoh}^!((\mathcal{Y}_x)^{\wedge}_{\mathrm{y}^0}) = \operatorname*{lim}_{k \in \mathcal{L}} \mathrm{IndCoh}^!((\mathcal{Y}^k)^{\wedge}_{\mathrm{y}^0}) = \lim\left[\ldots \to \mathrm{IndCoh}^!((\mathcal{Y}^l)^{\wedge}_{\mathrm{y}^0}) \xrightarrow{(\iota^{kl})^{\wedge,!}_{\mathrm{y}^0}} \mathrm{IndCoh}^!((\mathcal{Y}^k)^{\wedge}_{\mathrm{y}^0}) \to \ldots\right],$$

where for each $k \in \mathcal{L}$ the category $\mathrm{IndCoh}^!((\mathcal{Y}^k)^{\wedge}_{\mathrm{y}^0}) \in \mathrm{DGCat}$ is presented by

(15.2.4)

$$\mathrm{IndCoh}^!((\mathcal{Y}^k)^{\wedge}_{\mathrm{y}^0}) = \operatorname*{colim}_{i \in \mathcal{J}} \mathrm{IndCoh}^!((\mathcal{Y}^k_i)^{\wedge}_{\mathrm{y}^0_i}) = \mathrm{colim}\left[\ldots \leftarrow \mathrm{IndCoh}^!((\mathcal{Y}^k_i)^{\wedge}_{\mathrm{y}^0_i}) \xleftarrow{(\varphi^k_{ij})^{\wedge,!}_{\mathrm{y}^0}} \mathrm{IndCoh}^!((\mathcal{Y}^k_j)^{\wedge}_{\mathrm{y}^0_j}) \leftarrow \ldots\right].$$

Similarly, following Remark 15.1.12, we have:

*Remark* 15.2.3. Concretely, an object $\mathcal{F} \in \mathrm{IndCoh}^!((\mathcal{Y}_x)^{\wedge}_{\mathcal{Y}^0})$ is specified by an assignment

$$k \mapsto \left( \mathcal{F}^k \in \mathrm{IndCoh}^!((\mathcal{Y}^k)^{\wedge}_{\mathcal{Y}^0}) \right) \qquad [k \to l] \mapsto \left[ (\iota^{kl})^{\wedge !}_{\mathcal{Y}^0} \mathcal{F}^l \xrightarrow{\cong} \mathcal{F}^k \right]$$

defined for each $k, l \in \mathcal{L}$ and map $k \to l$.

Further, a typical object $\mathcal{F}^k \in \mathrm{IndCoh}^!((\mathcal{Y}^k)^{\wedge}_{\mathcal{Y}^0})$ is given by $\mathcal{F}^k = (\varphi^k_i)^{\wedge !}_{\mathcal{Y}^0} \mathcal{F}^k_i$ for $\mathcal{F}^k_i \in \mathrm{IndCoh}^!((\mathcal{Y}^k_i)^{\wedge}_{\mathcal{Y}^0_i})$.

*Example* 15.2.4. There is a canonical functor

$$\mathrm{p}^!_{(\mathcal{Y}_x)^{\wedge}_{\mathcal{Y}^0}} : \mathrm{Vect} \to \mathrm{IndCoh}^!((\mathcal{Y}_x)^{\wedge}_{\mathcal{Y}^0}) \qquad \text{induced by} \qquad \mathrm{p}^!_{(\mathcal{Y}^k_i)^{\wedge}_{\mathcal{Y}^0_i}} : \mathrm{Vect} \to \mathrm{IndCoh}^!((\mathcal{Y}^k_i)^{\wedge}_{\mathcal{Y}^0_i}) \;,$$

where $\mathrm{p}_{(\mathcal{Y}^k_i)^{\wedge}_{\mathcal{Y}^0_i}} : (\mathcal{Y}^k_i)^{\wedge}_{\mathcal{Y}^0_i} \to \mathrm{pt}$. In particular, there is a canonical object $\omega_{(\mathcal{Y}_x)^{\wedge}_{\mathcal{Y}^0}} \in \mathrm{IndCoh}^!((\mathcal{Y}_x)^{\wedge}_{\mathcal{Y}^0})$.

Following Example 14.3.6, we have:

*Example* 15.2.5. There is a canonical functor

$$\pi^!_{\mathcal{Y}^0_{\mathrm{dR}}} : D^!(\mathcal{Y}^0) \to \mathrm{IndCoh}^!((\mathcal{Y}_x)^{\wedge}_{\mathcal{Y}^0}) \qquad \text{induced by} \qquad \pi^{k!}_{\mathcal{Y}^0_{i,\mathrm{dR}}} : D(\mathcal{Y}^0_i) \to \mathrm{IndCoh}^!((\mathcal{Y}^k_i)^{\wedge}_{\mathcal{Y}^0_i})$$

where $\pi^k_{\mathcal{Y}^0_{i,\mathrm{dR}}} : (\mathcal{Y}^k_i)^{\wedge}_{\mathcal{Y}^0_i} \to \mathcal{Y}^0_{i,\mathrm{dR}}$ is the canonical projection.

Following Definition 14.3.7, we have:

*Definition* 15.2.6. The internal variant $\tilde{\mathcal{D}}^{\mathrm{ch}}(Y) \in \mathrm{Alg}^{\mathrm{fact}}_{\mathrm{un}}(\mathrm{IndCoh}^!_{\mathcal{J}^{\mathrm{mer}}(Y)^{\wedge}_{\mathcal{J}(Y)}})$ of chiral differential operators is defined for $Y \in \mathrm{Sch}_{\mathrm{ft}}$ a finite-type scheme satisfying the hypotheses of this subsection by

$$\tilde{\mathcal{D}}^{\mathrm{ch}}(Y) = \pi^!_{\mathcal{J}(Y)_{\mathrm{dR}}} \omega_{\mathcal{J}(Y)} = \omega_{\mathcal{J}^{\mathrm{mer}}(Y)^{\wedge}_{\mathcal{J}(Y)}} \quad \in \mathrm{Alg}^{\mathrm{fact}}_{\mathrm{un}}(\mathrm{IndCoh}^!_{\mathcal{J}^{\mathrm{mer}}(Y)^{\wedge}_{\mathcal{J}(Y)}}) \;.$$

*Remark* 15.2.7. Concretely, the underlying object of the fibre category $\tilde{\mathcal{D}}^{\mathrm{ch}}(Y)_x \in \mathrm{IndCoh}^!((\mathcal{Y}_x)^{\wedge}_{\mathcal{Y}^0})$ at $x \in X$ is given by

$$\tilde{\mathcal{D}}^{\mathrm{ch}}(Y)_x = \pi^!_{\mathcal{Y}^0_{\mathrm{dR}}} \omega_{\mathcal{Y}^0} = \omega_{(\mathcal{Y}_x)^{\wedge}_{\mathcal{Y}^0}} \in \mathrm{IndCoh}^!((\mathcal{Y}_x)^{\wedge}_{\mathcal{Y}^0}) \;.$$

Following Example 15.1.26, we have:

*Example* 15.2.8. There is a presentation

$$\mathrm{IndCoh}^!((\mathcal{Y}^k)^{\wedge}_{\mathcal{Y}^0}) = \lim \left[ \ldots \to \mathrm{IndCoh}^!((\mathcal{Y}^k_i)^{\wedge}_{\mathcal{Y}^0_i}) \xrightarrow{\underline{\mathrm{Hom}}(\omega_{(\mathcal{Y}^k_i)^{\wedge}_{\mathcal{Y}^0_i}/(\mathcal{Y}^k_j)^{\wedge}_{\mathcal{Y}^0_j}}, \cdot)} \mathrm{IndCoh}^!((\mathcal{Y}^k_j)^{\wedge}_{\mathcal{Y}^0_j}) \to \ldots \right]$$

Following Remark 15.1.27, we have:

*Remark* 15.2.9. Concretely, an object $\mathcal{F}^k \in \mathrm{IndCoh}^!((\mathcal{Y}^k)^{\wedge}_{\mathcal{Y}^0})$ is given by an assignment

$$i \mapsto \left( \mathcal{F}^k_i \in \mathrm{IndCoh}((\mathcal{Y}^k_i)^{\wedge}_{\mathcal{Y}^0_i}) \right) \qquad [j \to i] \mapsto \left[ \underline{\mathrm{Hom}}(\omega_{(\mathcal{Y}^k_i)^{\wedge}_{\mathcal{Y}^0_i}/(\mathcal{Y}^k_j)^{\wedge}_{\mathcal{Y}^0_j}}, \mathcal{F}^k_i) \xrightarrow{\cong} \mathcal{F}^k_j \right]$$

for each $i, j \in \mathcal{I}$ and arrow $j \to i$ in $\mathcal{I}$.

Further, following Example 15.1.28, we have:

*Example* 15.2.10. In particular, the presentation of Example 15.2.8 above induces a functor

$$\mathrm{Hom}(\omega_{(\mathcal{Y}^k)^\wedge_{y0}}, \cdot) : \mathrm{IndCoh}^!((\mathcal{Y}^k)^\wedge_{y0}) \to \mathrm{Vect} \qquad \text{defined by} \qquad \mathcal{F}^k = (\mathcal{F}^k_i)_{i \in \mathcal{J}} \mapsto \mathrm{Hom}(\omega_{\mathcal{Y}^k_i}, \mathcal{F}^k_i) \ .$$

Following Remark 15.1.29 and Definition 15.1.30, we have:

*Example* 15.2.11. The functors of the preceding Example do *not* satisfy the compatibility condition $\mathrm{Hom}(\omega_{(\mathcal{Y}^l)^\wedge_{y0}}, (\iota^{kl})^\wedge_{y0\bullet}(\cdot)) \cong \mathrm{Hom}(\omega_{(\mathcal{Y}^k)^\wedge_{y0}}, \cdot)$, and thus do not canonically extend to the limit category $\mathrm{IndCoh}^!((\mathcal{Y}_x)^\wedge_{y0})$. However, given a choice of global object $\mathcal{E} \in \mathrm{Det}_Y$, the functors

$$\mathrm{p}_{(\mathcal{Y}^k)^\wedge_{y0}\bullet} := \mathrm{Hom}(\omega_{(\mathcal{Y}^k)^\wedge_{y0}}, \mathcal{E}^k \otimes_{\mathcal{O}_{(\mathcal{Y}^k)^\wedge_{y0}}} (\cdot)) : \mathrm{IndCoh}^!((\mathcal{Y}^k)^\wedge_{y0}) \to \mathrm{Vect}$$

admit natural structure maps along $k \in \mathcal{L}$, and thus define

$$\mathrm{p}_{(\mathcal{Y}_x)^\wedge_{y0}\bullet} : \mathrm{IndCoh}^!((\mathcal{Y}_x)^\wedge_{y0}) \to \mathrm{Vect} \qquad \text{is defined by} \qquad \mathcal{F} = (\mathcal{F}^k)_{k \in \mathcal{L}} \mapsto \mathrm{colim}_{k \in \mathcal{L}} \mathrm{p}_{(\mathcal{Y}^k)^\wedge_{y0}\bullet} \mathcal{F}^k \ .$$

*Remark* 15.2.12. There are canonical functors

$$\pi_{\mathcal{Y}^k_\bullet} : \mathrm{IndCoh}^!((\mathcal{Y}^k)^\wedge_{y0}) \to \mathrm{IndCoh}^!(\mathcal{Y}^k) \qquad \text{induced by} \qquad \pi_{\mathcal{Y}^i_\bullet} : \mathrm{IndCoh}^!((\mathcal{Y}^k_i)^\wedge_{y0}) \to \mathrm{IndCoh}^!(\mathcal{Y}^i_k)$$

where $\pi_{\mathcal{Y}^i_{k\bullet}} : (\mathcal{Y}^k_i)^\wedge_{y0} \to \mathcal{Y}^k_i$ is the canonical projection. Moreover, the functors of the preceding example factor as

Following Proposition 14.3.9, we have:

*Proposition* 15.2.13. There is a canonical equivalence

$$\mathcal{D}^{\mathrm{ch}}(Y) = \mathrm{p}_{\mathcal{J}^{\mathrm{mer}}(Y)^\wedge_{\mathcal{J}(Y)}\bullet} \tilde{\mathcal{D}}^{\mathrm{ch}}(Y) \ \in \mathrm{Alg}^{\mathrm{fact}}_{\mathrm{un}}(X) \ .$$

*Proof.* Over each point $x \in X$ the equivalence is given in terms of the colimit description

$$\begin{aligned}
\mathcal{D}^{\mathrm{ch}}(Y)_x &\cong \mathrm{colim}_{k \in \mathcal{L}} \mathrm{Hom}(\omega_{\mathcal{Y}^k_j}, \mathcal{E}^k_j \otimes (\cdot)) \circ \mathrm{oblv}^r_{\mathcal{Y}^k_j} \circ \iota^{0k}_{j*}(\omega_{\mathcal{Y}0})_j \\
&\cong \mathrm{colim}_{k \in \mathcal{L}} \mathrm{Hom}(\omega_{\mathcal{Y}^k_j}, \mathcal{E}^k_j \otimes (\cdot)) \circ \pi_{\mathcal{Y}^k_\bullet} \circ \pi^!_{\mathcal{Y}0_{j,\mathrm{dR}}}(\omega_{\mathcal{Y}0})_j \\
&\cong \mathrm{p}_{(\mathcal{Y}^k)^\wedge_{y0}\bullet} \circ \pi^!_{\mathcal{Y}0_{\mathrm{dR}}} \omega_{\mathcal{Y}_0}
\end{aligned}$$

where the first isomorphism follows by definition, as in Remark 15.1.33, the second follows by base change along and the third follows from the preceding remark and the definition of $\tilde{\mathcal{D}}^{\mathrm{ch}}(Y)$.  □

*Remark* 15.2.14. The compatibility of the preceding Proposition with the module structure over $\tilde{\mathcal{A}}(Y) \in \mathrm{Alg}_{\mathbb{E}_1,\mathrm{un}}^{\mathrm{fact}}(D_{\mathcal{Z}(Y)})$ is proved in the case $Y = N/G$ in Subsection 17.

## 16. Chiral differential operators on quotient stacks and semi-infinite cohomology

In this section, we extend the geometric construction of chiral differential operators to the case $Y = N/G$ is a quotient stack of $N$ a finite dimensional $G$ representation by $G$ a finite-type, affine, reductive algebraic group, as considered in Section 13 as the input data for the construction of the three dimensional A model in [BFN18]. Moreover, we identify the resulting factorization algebra with the $G_\mathcal{O}$ equivariant semi-infinite cohomology $C^{\frac{\infty}{2}}(\hat{\mathfrak{g}}, \mathfrak{g}_\mathcal{O}, G_\mathcal{O}; (\mathcal{D}^{\mathrm{ch}}(N)))$ (or BRST reduction, in the sense of Section I-13) of the chiral differential operators $\mathcal{D}^{\mathrm{ch}}(N)$ on $N$ with respect to $\hat{\mathfrak{g}}$. The goal is to establish the results outlined in Subsection 14 in this generality, culminating in the following section 17 with the proof of Theorem 14.2.1 in this setting.

16.0.1. *Summary.* In Section 16.1, we extend the geometric construction of $\mathcal{D}^{\mathrm{ch}}(Y)$ in the case $Y = N/G$, and in Section 16.2, we recall a geometric construction of the semi-infinite cohomology functor which is manifestly equivalent to the given construction.

16.1. **Chiral differential operators on quotient stacks.** In this subsection, we explain the definition of $\mathcal{D}^{\mathrm{ch}}(G, N) := \mathcal{D}^{\mathrm{ch}}(N/G) \in \mathrm{Alg}_{\mathrm{un}}^{\mathrm{fact}}(X)$, and show that it is calculated as the BRST reduction (alias semi-infinite cohomology) of $\mathcal{D}^{\mathrm{ch}}(N)$ in the sense of Section I-13; in particular, we explain that the analogue of the determinant anomaly of Remark 15.1.22 in this context is given by the requirement that Lie* algebra $\Theta_{\mathcal{J}(G)} = \mathfrak{g} \otimes \mathcal{D}_X$ acts on $\mathcal{D}^{\mathrm{ch}}(N)$ at level $-\kappa_{\mathrm{Tate}}$.

*Example* 16.1.1. Following Example 14.3.1, we consider the factorization space

$$\mathcal{J}^{\mathrm{mer}}(G, N)_{\mathcal{J}(G,N)}^\wedge = \mathcal{J}(G,N)_{\mathrm{dR}} \times_{\mathcal{J}^{\mathrm{mer}}(G,N)_{\mathrm{dR}}} \mathcal{J}^{\mathrm{mer}}(G,N) \in \mathrm{PreStk}_{\mathrm{un}}^{\mathrm{fact}}(X)$$

where $\mathcal{J}^{(\mathrm{mer})}(G, N) = \mathcal{J}^{(\mathrm{mer})}(N/G) \in \mathrm{PreStk}_{\mathrm{un}}^{\mathrm{fact}}(X)$ are the spaces of (meromorphic) jets to the quotient stack $Y = N/G$.

*Remark* 16.1.2. Concretely, the prestacks $\mathcal{J}^{\mathrm{mer}}(G, N)_{\mathcal{J}(G,N),I}^\wedge \in \mathrm{PreStk}_{/X^I}$ over $X^I$ assigned to each $I \in \mathrm{fSet}_\emptyset$ and the fibre $\mathcal{J}^{\mathrm{mer}}(G, N)_{\mathcal{J}(G,N),x}^\wedge \in \mathrm{PreStk}$ over each $x \in X$ are given by

$$\mathcal{J}^{\mathrm{mer}}(G, N)_{\mathcal{J}(G,N),I}^\wedge = \mathcal{J}(G,N)_{I,\mathrm{dR}} \times_{\mathcal{J}^{\mathrm{mer}}(G,N)_{I,\mathrm{dR}}} \mathcal{J}^{\mathrm{mer}}(G,N)_I = (\mathcal{J}^{\mathrm{mer}}(G,N)_I)_{\mathcal{J}(G,N)_I}^\wedge \quad , \text{and}$$

$$\mathcal{J}^{\mathrm{mer}}(G, N)_{\mathcal{J}(G,N),x}^\wedge = [N/G]_{\mathcal{O}_x,\mathrm{dR}} \times_{[N/G]_{\mathcal{K}_x,\mathrm{dR}}} [N/G]_{\mathcal{K}_x} = ([N/G]_{\mathcal{K}_x})_{[N/G]_{\mathcal{O}_x}}^\wedge$$

Towards defining the factorization category $\mathrm{IndCoh}_{\mathcal{J}^{\mathrm{mer}}(G,N)_{\mathcal{J}(G,N)}^\wedge}^! \in \mathrm{Cat}_{\mathrm{un}}^{\mathrm{fact}}(X)$, we introduce a convenient presentation of the underlying factorization space:

*Proposition* 16.1.3. There is an equivalence of factorization spaces

$$\mathcal{J}^{\mathrm{mer}}(G, N)_{\mathcal{J}(G,N)}^\wedge \cong (\mathcal{J}^{\mathrm{mer}}(N)_{\mathcal{J}(N)}^\wedge \times \mathcal{J}^{\mathrm{mer}}(G)_{\mathrm{dR}})/(\mathcal{J}(G)_{\mathrm{dR}} \times \mathcal{J}^{\mathrm{mer}}(G)) \ .$$

*Remark* 16.1.4. Concretely, over each point $x \in X$ the preceding equivalence is given by

$$\mathcal{J}^{\mathrm{mer}}(G, N)_{\mathcal{J}(G,N)_x}^\wedge = [N/G]_{\mathcal{O},\mathrm{dR}} \times_{[N/G]_{\mathcal{K},\mathrm{dR}}} [N/G]_{\mathcal{K}} \cong \left((N_\mathcal{O} \times_{N_{\mathcal{K},\mathrm{dR}}} N_\mathcal{K}) \times G_{\mathcal{K},\mathrm{dR}}\right)/(G_{\mathcal{O},\mathrm{dR}} \times G_\mathcal{K}) \ .$$

*Remark* 16.1.5. We define $\mathrm{IndCoh}_{\mathcal{J}^{\mathrm{mer}}(N)_{\mathcal{J}(N)}^\wedge}^! \in \mathrm{Cat}_{\mathrm{un}}^{\mathrm{fact}}(X)$ following Remark 5.2.8, with fibre category over each $x \in X$ given by $\mathrm{IndCoh}^!((N_\mathcal{K})_{N_\mathcal{O}}^\wedge)$ as defined in Example 15.2.2.

Following the proceeding proposition and remarks, we make the following definition:

*Definition* 16.1.6. The factorization category $\mathrm{IndCoh}^!_{\mathcal{J}^{\mathrm{mer}}(G,N)^{\wedge}_{\mathcal{J}(G,N)}} \in \mathrm{Cat}^{\mathrm{fact}}_{\mathrm{un}}(X)$ is defined by

$$\mathrm{IndCoh}^!_{\mathcal{J}^{\mathrm{mer}}(G,N)^{\wedge}_{\mathcal{J}(G,N)}} = \left(\mathrm{IndCoh}^!_{\mathcal{J}^{\mathrm{mer}}(N)^{\wedge}_{\mathcal{J}(G,N)}} \otimes^* D^!_{\mathcal{J}^{\mathrm{mer}}(G)}\right)^{\mathcal{J}(G)\times(\mathcal{J}^{\mathrm{mer}}(G),w)} .$$

*Remark* 16.1.7. Concretely, the fibre category $\mathrm{IndCoh}^!(\mathcal{J}^{\mathrm{mer}}(G,N)^{\wedge}_{\mathcal{J}(G,N)_x}) \in \mathrm{DGCat}$ over each $x \in X$ is given by

$$\mathrm{IndCoh}^!(\mathcal{J}^{\mathrm{mer}}(G,N)^{\wedge}_{\mathcal{J}(G,N)x}) = \left(\mathrm{IndCoh}^!((N_{\mathcal{K}})^{\wedge}_{N_{\mathcal{O}}}) \otimes D^!(G_{\mathcal{K}})\right)^{G_{\mathcal{O}}\times(G_{\mathcal{K}},w)} .$$

where the superscript denotes simeltaneously taking the category strong $G_{\mathcal{O}}$ invariants and weak $G_{\mathcal{K}}$ invariants, the latter of which is as defined in Section 7 of [Ras20b].

*Example* 16.1.8. In the case $N = \{0\}$ is the zero representation, the factorization space above is given by

$$\mathcal{J}^{\mathrm{mer}}(G,\{0\})^{\wedge}_{\mathcal{J}(G,\{0\})} = \mathcal{J}^{\mathrm{mer}}(BG)^{\wedge}_{\mathcal{J}(BG)} = \mathcal{J}(G)_{\mathrm{dR}}\backslash \mathcal{J}^{\mathrm{mer}}(G)_{\mathrm{dR}}/\mathcal{J}^{\mathrm{mer}}(G) .$$

Moreover, there are canonical identifications of factorization categories

$$\mathrm{IndCoh}^!_{\mathcal{J}^{\mathrm{mer}}(BG)^{\wedge}_{\mathcal{J}(BG)}} := \mathrm{IndCoh}^!_{\mathcal{J}^{\mathrm{mer}}(G,\{0\})^{\wedge}_{\mathcal{J}(G,\{0\})}} = (D^!_{\mathcal{J}^{\mathrm{mer}}(G)})^{\mathcal{J}(G)\times(\mathcal{J}^{\mathrm{mer}}(G),w)} = \hat{\mathfrak{g}}_0\text{-}\mathrm{Mod}^{\mathcal{J}(G)}_{\mathrm{Ran}_{X,\mathrm{un}}} ,$$

where the latter is the factorization enhancement of the Kazhdan-Lusztig category of $G_{\mathcal{O}}$ integrable $\hat{\mathfrak{g}}$ modules at level $\kappa = 0$, recalling that

$$\mathrm{IndCoh}^!_{\mathcal{J}^{\mathrm{mer}}(BG)^{\wedge}_{\sigma_{\{e\}}}} := (D^!_{\mathcal{J}^{\mathrm{mer}}(G)})^{\mathcal{J}^{\mathrm{mer}}(G),w} = \hat{\mathfrak{g}}_0\text{-}\mathrm{Mod}_{\mathrm{Ran}_{X,\mathrm{un}}} \in \mathrm{Cat}^{\mathrm{fact}}_{\mathrm{un}}(X) ,$$

where $\sigma_{\{e\}} : \mathrm{Ran}_{X_{\mathrm{dR}},\mathrm{un}} \to \mathcal{J}^{\mathrm{mer}}(BG)$ is the factorizable unit section, as in Example 16.1.10 below.

*Remark* 16.1.9. Concretely, over each point $x \in X$ the preceding identifications are given by

$$\mathrm{IndCoh}^!((BG_{\mathcal{K}})^{\wedge}_{BG_{\mathcal{O}}}) := \mathrm{IndCoh}^!(\mathcal{J}^{\mathrm{mer}}(G,\{0\})^{\wedge}_{\mathcal{J}(G,\{0\})},x) = D^!(G_{\mathcal{K}})^{G_{\mathcal{O}}\times(G_{\mathcal{K}},w)} = \hat{\mathfrak{g}}_0\text{-}\mathrm{Mod}^{G_{\mathcal{O}}}$$

recalling that

$$\mathrm{IndCoh}^!((BG_{\mathcal{K}})^{\wedge}_{\{e\}}) := D^!(G_{\mathcal{K}})^{G_{\mathcal{K}},w} = \hat{\mathfrak{g}}_0\text{-}\mathrm{Mod} \in \mathrm{DGCat}$$

as in e.g. 9.12 of [Ras20b], where $\{e\} : \mathrm{pt} \to BG$ is the canonical map.

*Example* 16.1.10. There is a canonical strict factorizable section $\sigma_{\{e\}} : \mathrm{Ran}_{X_{\mathrm{dR}},\mathrm{un}} \to \mathcal{J}^{\mathrm{mer}}(BG)$ given fibrewise by the canonical map $\mathrm{pt} \to BG_{\mathcal{K}}$. The factorization functor induced by pullback, as in Example 5.3.2, defines a natural forgetful functor

$$\hat{\mathfrak{g}}_0\text{-}\mathrm{Mod}_{\mathrm{Ran}_{X,\mathrm{un}}} \to D_{\mathrm{Ran}_{X,\mathrm{un}}} \qquad \text{and thus} \qquad \mathrm{Alg}^{\mathrm{fact}}_{\mathrm{un}}(\hat{\mathfrak{g}}_0\text{-}\mathrm{Mod}_{\mathrm{Ran}_{X,\mathrm{un}}}) \to \mathrm{Alg}^{\mathrm{fact}}_{\mathrm{un}}(X) .$$

The latter is the canonical forgetful functor from the category of factorization algebras admitting an action of the Lie* algebra $\Theta_{\mathcal{J}(G)} = \mathfrak{g} \otimes \mathcal{D}_X \in \mathrm{Lie}^*(X)$ of Example I-10.2.4.

Following definitions 14.3.7 and 15.2.6 and Remark 15.2.7, we make the following definition:

*Definition* 16.1.11. The internal variant

$$\tilde{\mathcal{D}}^{\mathrm{ch}}(G,N) \in \mathrm{Alg}^{\mathrm{fact}}_{\mathrm{un}}(\mathrm{IndCoh}^!_{\mathcal{J}^{\mathrm{mer}}(G,N)^{\wedge}_{\mathcal{J}(G,N)}})$$

of chiral differential operators to $N/G$ is defined as

$$\tilde{\mathcal{D}}^{\mathrm{ch}}(G,N) = \omega_{\mathcal{J}^{\mathrm{mer}}(N)^{\wedge}_{\mathcal{J}(N)}} \boxtimes \omega_{\mathcal{J}^{\mathrm{mer}}(G)} \in \mathrm{Alg}^{\mathrm{fact}}_{\mathrm{un}}\left(\mathrm{IndCoh}^!_{\mathcal{J}^{\mathrm{mer}}(N)^{\wedge}_{\mathcal{J}(N)}} \otimes^* D^!_{\mathcal{J}^{\mathrm{mer}}(G)}\right)^{\mathcal{J}(G)\times(\mathcal{J}^{\mathrm{mer}}(G),w)} .$$

*Remark* 16.1.12. Concretely, the underlying object $\tilde{\mathcal{D}}^{\mathrm{ch}}(G, N)_x \in \mathrm{IndCoh}^!(\mathcal{J}^{\mathrm{mer}}(G, N)^{\wedge}_{\mathcal{J}(G,N)_x})$ of the fibre category over each $x \in X$ is given by

$$\tilde{\mathcal{D}}^{\mathrm{ch}}(G, N)_x = \omega_{(N_{\mathcal{K}})^{\wedge}_{N_{\mathcal{O}}}} \boxtimes \omega_{G_{\mathcal{K}}} \in \left( \mathrm{IndCoh}^!((N_{\mathcal{K}})^{\wedge}_{N_{\mathcal{O}}}) \otimes D^!(G_{\mathcal{K}}) \right)^{G_{\mathcal{O}} \times (G_{\mathcal{K}}, w)} \ .$$

*Remark* 16.1.13. Following Example 15.2.11, we now wish to define an analogous factorization functor

$$\mathrm{p}_{\mathcal{J}^{\mathrm{mer}}(G,N)^{\wedge}_{\mathcal{J}(G,N)}, \bullet} : \mathrm{IndCoh}^!_{\mathcal{J}^{\mathrm{mer}}(G,N)^{\wedge}_{\mathcal{J}(G,N)}} \to D_{\mathrm{Ran}_{X, \mathrm{un}}} \ ,$$

so that, following Proposition 15.2.13, we can define

$$\mathcal{D}^{\mathrm{ch}}(G, N) = \mathrm{p}_{\mathcal{J}^{\mathrm{mer}}(G,N)^{\wedge}_{\mathcal{J}(G,N)}, \bullet} \tilde{\mathcal{D}}^{\mathrm{ch}}(G, N) \ \in \mathrm{Alg}^{\mathrm{fact}}_{\mathrm{un}}(X) \ .$$

The analogue of the requirement of a Tate structure on $Y$ will be that of a $G$ equivariant Tate structure on $N$ at level $-\mathrm{Tate}$, and we will identify the resulting factorization algebra with the BRST reduction of $\mathcal{D}^{\mathrm{ch}}(N) \in \mathrm{Alg}^{\mathrm{fact}}_{\mathrm{un}}(X)$, in the sense of Section I-13.

We now proceed with the construction outlined in the preceding remark:

*Proposition* 16.1.14. There are natural maps of factorization spaces

$$\mathcal{J}^{\mathrm{mer}}(N)^{\wedge}_{\mathcal{J}(N)} \xrightarrow{\iota_{\mathcal{J}^{\mathrm{mer}}(N)^{\wedge}_{\mathcal{J}(N)}}} \mathcal{J}^{\mathrm{mer}}(G, N)^{\wedge}_{\mathcal{J}(G,N)} \xrightarrow{\mathrm{p}^{G}_{\mathcal{J}^{\mathrm{mer}}(N)^{\wedge}_{\mathcal{J}(N)}}} \mathcal{J}^{\mathrm{mer}}(BG)^{\wedge}_{\mathcal{J}(BG)} \ ,$$

which define a fibre sequence of factorization spaces.

*Remark* 16.1.15. Concretely, over each point $x \in X$, the fibre sequence in prestacks is given by

$$(16.1.1) \qquad N_{\mathcal{O}, \mathrm{dR}} \times_{N_{\mathcal{K}, \mathrm{dR}}} N_{\mathcal{K}} \hookrightarrow [N/G]_{\mathcal{O}, \mathrm{dR}} \times_{[N/G]_{\mathcal{K}, \mathrm{dR}}} [N/G]_{\mathcal{K}} \twoheadrightarrow BG_{\mathcal{O}, \mathrm{dR}} \times_{BG_{\mathcal{K}, \mathrm{dR}}} BG_{\mathcal{K}}$$

or under the equivalence of Proposition 16.1.3 above,

$$(16.1.2) \qquad (N_{\mathcal{K}})^{\wedge}_{N_{\mathcal{O}}} \xrightarrow{\iota_{(N_{\mathcal{K}})^{\wedge}_{N_{\mathcal{O}}}}} \left( (N_{\mathcal{K}})^{\wedge}_{N_{\mathcal{O}}} \times G_{\mathcal{K}, \mathrm{dR}} \right) / (G_{\mathcal{O}, \mathrm{dR}} \times G_{\mathcal{K}}) \xrightarrow{\mathrm{p}^{G}_{(N_{\mathcal{K}})^{\wedge}_{N_{\mathcal{O}}}}} G_{\mathcal{O}, \mathrm{dR}} \backslash G_{\mathcal{K}, \mathrm{dR}} / G_{\mathcal{K}} \ .$$

*Remark* 16.1.16. The desired factorization functor of Remark 16.1.13 will be constructed as a composition of factorization functors, corresponding to pushforward along the maps

$$\mathcal{J}^{\mathrm{mer}}(G, N)^{\wedge}_{\mathcal{J}(G,N)} \xrightarrow{\mathrm{p}^{G}_{\mathcal{J}^{\mathrm{mer}}(N)^{\wedge}_{\mathcal{J}(N)}}} \mathcal{J}^{\mathrm{mer}}(BG)^{\wedge}_{\mathcal{J}(BG)} \xrightarrow{\mathrm{p}_{\mathcal{J}^{\mathrm{mer}}(BG)^{\wedge}_{\mathcal{J}(BG)}}} \mathrm{Ran}_{X_{\mathrm{dR}}, \mathrm{un}}.$$

Restricting to a point $x \in X$ for concreteness, we now give a more detailed overview before proceeding with the construction. An overview of the construction of the semi-infinite cohomology functors used below is also recalled in the following subsection 16.2.

The map $\mathrm{p}^{G}_{(N_{\mathcal{K}})^{\wedge}_{N_{\mathcal{O}}}}$ is a fibration with each fibre equivalent to $(N_{\mathcal{K}})^{\wedge}_{N_{\mathcal{O}}}$, so that a Tate structure on $N$ which is strongly $G$ equivariant induces a functor

$$\mathrm{p}^{G}_{(N_{\mathcal{K}})^{\wedge}_{N_{\mathcal{O}}}, \bullet} : \mathrm{IndCoh}^!(\mathcal{J}^{\mathrm{mer}}(G, N)^{\wedge}_{\mathcal{J}(G,N)_x}) \to \mathrm{IndCoh}^!((BG_{\mathcal{K}})^{\wedge}_{BG_{\mathcal{O}}}) \ ,$$

by an equivariant enhancement of the construction of $\mathrm{p}_{(N_{\mathcal{K}})^{\wedge}_{N_{\mathcal{O}}}, \bullet}$ in Example 15.2.11, as we explain in Definition 16.1.20 below.

However, the desired pushforward functor

$$\mathrm{IndCoh}^!((BG_{\mathcal{K}})^{\wedge}_{BG_{\mathcal{O}}}) = D^!(G_{\mathcal{K}})^{G_{\mathcal{O}} \times (G_{\mathcal{K}}, w)} \to \mathrm{Vect} \qquad \text{is } not \text{ defined.}$$

Instead, following Section 8 of [Ras20b], there is a canonical functor

$$\mathrm{p}_{(\mathrm{B}G_{\mathcal{K}})^{\wedge}_{\mathrm{B}G_{\mathcal{O}}},\bullet}:\mathrm{IndCoh}^{\bullet}((\mathrm{B}G_{\mathcal{K}})^{\wedge}_{\mathrm{B}G_{\mathcal{O}}}):=D^{*}(G_{\mathcal{K}})^{G_{\mathcal{O}}}_{G_{\mathcal{K}},w}\to\mathrm{Vect}\ ,$$

and the analogue of the identification $\mathrm{IndCoh}^{!}\cong\mathrm{IndCoh}^{*}$ is given by

$$\mathrm{IndCoh}^{!}_{\chi-\mathrm{Tate}}((\mathrm{B}G_{\mathcal{K}})^{\wedge}_{\mathrm{B}G_{\mathcal{O}}}):=\left(\mathrm{Vect}_{\chi-\mathrm{Tate}}\otimes D^{!}(G_{\mathcal{K}})\right)^{G_{\mathcal{O}}\times(G_{\mathcal{K}},w)}\cong D^{*}(G_{\mathcal{K}})^{G_{\mathcal{O}}}_{G_{\mathcal{K}},w}=\mathrm{IndCoh}^{\bullet}((\mathrm{B}G_{\mathcal{K}})^{\wedge}_{\mathrm{B}G_{\mathcal{O}}})\ ,$$

where $\mathrm{Vect}_{\chi_{\kappa}}\in G_{\mathcal{K}}\text{-}\mathrm{Mod}_{\mathrm{weak}}$ is the weak $G_{\mathcal{K}}$ category with underlying DG category $\mathrm{Vect}\in\mathrm{DGCat}$ determined by the Kac-Moody central extension of $G_{\mathcal{K}}$ at level $\kappa$.

Moreover, it is shown in *loc. cit.* that the resulting functor
(16.1.3)
$$\mathrm{p}_{(\mathrm{B}G_{\mathcal{K}})^{\wedge}_{\mathrm{B}G_{\mathcal{O}}},\bullet}:\mathrm{IndCoh}^{!}_{\chi-\mathrm{Tate}}((\mathrm{B}G_{\mathcal{K}})^{\wedge}_{\mathrm{B}G_{\mathcal{O}}})\to\mathrm{Vect}\qquad\text{identifies with}\qquad C^{\frac{\infty}{2}}(\hat{\mathfrak{g}},\mathfrak{g}_{\mathcal{O}},G_{\mathcal{O}};(\cdot)):\hat{\mathfrak{g}}_{-\mathrm{Tate}}\text{-}\mathrm{Mod}^{G_{\mathcal{O}}}\to\mathrm{Vect}\ ,$$

under the equivalence of Remark 16.1.9, where the latter is the functor of $G_{\mathcal{O}}$-equivariant semi-infinite cohomology relative to $\mathfrak{g}_{\mathcal{O}}$.

Thus, in order to construct the desired total pushforward functor, we need to require that the Tate structure on $N$ is strongly $G$ equivariant up to a twist by the level $\kappa=-\mathrm{Tate}$. This corresponds to the requirement that the induced Lie* algebra action of $\Theta_{\mathfrak{d}(G)}=\mathfrak{g}\otimes\mathcal{D}_X$ on $D^{\mathrm{ch}}(N)$ is at level $\kappa=-\mathrm{Tate}$, so that the BRST reduction of $D^{\mathrm{ch}}(N)$ by $\mathfrak{g}\otimes\mathcal{D}_X$ is well defined, following Section I-13.

*Remark* 16.1.17. Following Example 15.2.11, given a choice of Tate structure on $N$, we obtain a functor

$$\mathrm{p}_{(N_{\mathcal{K}})^{\wedge}_{N_{\mathcal{O}}},\bullet}:\mathrm{IndCoh}^{!}((N_{\mathcal{K}})^{\wedge}_{N_{\mathcal{O}}})\to\mathrm{Vect}\qquad\text{and thus}\qquad\mathrm{IndCoh}^{!}((N_{\mathcal{K}})^{\wedge}_{N_{\mathcal{O}}})\otimes D^{!}(G_{\mathcal{K}})\to D^{!}(G_{\mathcal{K}})\ ,$$

Note that these are functors of plain DG categories, but do not in general define functors of $G_{\mathcal{O}}\times(G_{\mathcal{K}},w)$ categories. Thus, we make the following definition:

*Definition* 16.1.18. A strongly $G$ equivariant Tate structure on $N$ at level $\kappa$ is an enhancement of the functor of the preceding Remark to a functor

$$\mathrm{p}^{G_{\mathcal{O}}}_{(N_{\mathcal{K}})^{\wedge}_{N_{\mathcal{O}}},\bullet}:\left(\mathrm{IndCoh}^{!}((N_{\mathcal{K}})^{\wedge}_{N_{\mathcal{O}}})\otimes D^{!}(G_{\mathcal{K}})\right)^{G_{\mathcal{O}}}\to\mathrm{Vect}_{\chi_{\kappa}}\otimes D^{!}(G_{\mathcal{K}})^{G_{\mathcal{O}}}$$

of weak $G_{\mathcal{K}}$ categories, where $\mathrm{Vect}_{\chi_{\kappa}}\in G_{\mathcal{K}}\text{-}\mathrm{Mod}_{\mathrm{weak}}$ is as in Remark 16.1.16.

*Remark* 16.1.19. Concretely, following Section 2.8.16 of [BD04], a strongly $G$-equivariant Tate structure on $N$ at level $\kappa$ in particular induces an action of the Kac-Moody extension at level $\kappa$ (as in Example I-12.1.6) of the associated Lie* algebra $\Theta_{\mathfrak{d}(G)}=\mathfrak{g}\otimes\mathcal{D}_X\in\mathrm{Lie}^*(X)$ on $D^{\mathrm{ch}}(N)\in\mathrm{Alg}^{\mathrm{fact}}_{\mathrm{un}}(X)$.

*Definition* 16.1.20. Let $N$ be equipped with a strongly $G$ equivariant Tate structure at level $\kappa$. The induced pushforward factorization functor

$$\mathrm{p}^{G,\kappa}_{\mathfrak{d}^{\mathrm{mer}}(N)^{\wedge}_{\mathfrak{d}(N)},\bullet}:=\mathrm{p}^{\mathfrak{d}(G)\times(\mathfrak{d}^{\mathrm{mer}}(G),w)}_{\mathfrak{d}^{\mathrm{mer}}(N)^{\wedge}_{\mathfrak{d}(N)},\bullet}:\mathrm{IndCoh}^{!}_{\mathfrak{d}^{\mathrm{mer}}(G,N)^{\wedge}_{\mathfrak{d}(G,N)}}\to\hat{\mathfrak{g}}_{\kappa}\text{-}\mathrm{Mod}^{\mathfrak{d}(G)}_{\mathrm{Ran}_{X,\mathrm{un}}}$$

is defined by taking weak $G_{\mathcal{K}}$ invariants of the factorization enhancement of that in Definition 16.1.18 above.

*Remark* 16.1.21. Concretely, the functor on fibre categories over each $x\in X$ is given by

$$\mathrm{p}^{G,\kappa}_{(N_{\mathcal{K}})^{\wedge}_{N_{\mathcal{O}}},\bullet}:=\mathrm{p}^{G_{\mathcal{O}}\times(G_{\mathcal{K}},w)}_{(N_{\mathcal{K}})^{\wedge}_{N_{\mathcal{O}}},\bullet}:\left(\mathrm{IndCoh}^{!}((N_{\mathcal{K}})^{\wedge}_{N_{\mathcal{O}}})\otimes D^{!}(G_{\mathcal{K}})\right)^{G_{\mathcal{O}}\times(G_{\mathcal{K}},w)}\to\left(\mathrm{Vect}_{\chi_{\kappa}}\otimes D^{!}(G_{\mathcal{K}})\right)^{G_{\mathcal{O}}\times(G_{\mathcal{K}},w)}\ ,$$

the functor induced by taking weak $G_{\mathcal{K}}$ invariants of that in Definition 16.1.18.

*Example* 16.1.22. The image of $\tilde{\mathcal{D}}^{\mathrm{ch}}(G,N) \in \mathrm{Alg}^{\mathrm{fact}}_{\mathrm{un}}(\mathrm{IndCoh}^!_{\mathcal{J}^{\mathrm{mer}}(G,N)^{\wedge}_{\mathcal{J}(G,N)}})$ under the induced functor on internal factorization algebras is given by

$$\tilde{\mathcal{D}}^{\mathrm{ch}}(N)^{G,\kappa} \in \mathrm{Alg}^{\mathrm{fact}}_{\mathrm{un}}(\hat{\mathfrak{g}}_{\kappa}\text{-}\mathrm{Mod}^{\mathcal{J}(G)}_{\mathrm{Ran}_{X,\mathrm{un}}}) \ ,$$

the object whose image under the forgetful functor of Example 16.1.10 is given by $D^{\mathrm{ch}}(N) \in \mathrm{Alg}^{\mathrm{fact}}_{\mathrm{un}}(X)$, lifted by the factorizable, $G_{\mathcal{O}}$-integrable, $\hat{\mathfrak{g}}$ action at level $\kappa$ of Remark 16.1.19.

*Example* 16.1.23. There is a natural factorization enhancement of the $G_{\mathcal{O}}$ equivariant semi-infinite cohomology functor of Equation 16.1.3, which we denote by

$$\mathrm{p}_{\mathcal{J}^{\mathrm{mer}}(BG)^{\wedge}_{\mathcal{J}(BG)}}\text{-} : \hat{\mathfrak{g}}_{-\mathrm{Tate}}\text{-}\mathrm{Mod}^{\mathcal{J}(G)}_{\mathrm{Ran}_{X,\mathrm{un}}} \to D_{\mathrm{Ran}_{X,\mathrm{un}}} \ .$$

In the following subsection, we recall a geometric construction of the semi-infinite cohomology functor, emphasizing the analogy with the internal construction of chiral differential operators on schemes.

We now define our desired factorization algebra:

*Definition* 16.1.24. Let $N$ be equipped with a strongly $G$ equivariant Tate structure at level $\kappa = -\mathrm{Tate}$. The chiral differential operators to $Y = N/G$ are defined by

$$\mathcal{D}^{\mathrm{ch}}(G,N) = \mathrm{p}_{\mathcal{J}^{\mathrm{mer}}(G,N)^{\wedge}_{\mathcal{J}(G,N)}}\text{-}\tilde{\mathcal{D}}^{\mathrm{ch}}(G,N) \ \in \mathrm{Alg}^{\mathrm{fact}}_{\mathrm{un}}(X) \ ,$$

where:

- $\tilde{\mathcal{D}}^{\mathrm{ch}}(G,N) \in \mathrm{Alg}^{\mathrm{fact}}_{\mathrm{un}}(\mathrm{IndCoh}^!_{\mathcal{J}^{\mathrm{mer}}(G,N)^{\wedge}_{\mathcal{J}(G,N)}})$ is as in Example 16.1.11,
- $\mathrm{p}_{\mathcal{J}^{\mathrm{mer}}(G,N)^{\wedge}_{\mathcal{J}(G,N)}}\text{-} := \mathrm{p}_{\mathcal{J}^{\mathrm{mer}}(BG)^{\wedge}_{\mathcal{J}(BG)}}\text{-} \circ \mathrm{p}^{G,-\mathrm{Tate}}_{\mathcal{J}^{\mathrm{mer}}(N)^{\wedge}_{\mathcal{J}(N)}}\text{-} : \mathrm{IndCoh}^!_{\mathcal{J}^{\mathrm{mer}}(G,N)^{\wedge}_{\mathcal{J}(G,N)}} \to D_{\mathrm{Ran}_{X,\mathrm{un}}}$,
- $\mathrm{p}^{G,-\mathrm{Tate}}_{\mathcal{J}^{\mathrm{mer}}(N)^{\wedge}_{\mathcal{J}(N)}}\text{-} : \mathrm{IndCoh}^!_{\mathcal{J}^{\mathrm{mer}}(G,N)^{\wedge}_{\mathcal{J}(G,N)}} \to \hat{\mathfrak{g}}_{-\mathrm{Tate}}\text{-}\mathrm{Mod}^{\mathcal{J}(G)}_{\mathrm{Ran}_{X,\mathrm{un}}}$ is as in Example 16.1.20, and
- $\mathrm{p}_{\mathcal{J}^{\mathrm{mer}}(BG)^{\wedge}_{\mathcal{J}(BG)}}\text{-} : \hat{\mathfrak{g}}_{-\mathrm{Tate}}\text{-}\mathrm{Mod}^{\mathcal{J}(G)}_{\mathrm{Ran}_{X,\mathrm{un}}} \to D_{\mathrm{Ran}_{X,\mathrm{un}}}$ is as in Example 16.1.23 ,

*Example* 16.1.25. Vertex algebra

## 16.2. **A geometric construction of semi-infinite cohomology.** In an attempt to make the constructions of the preceding subsection more palatable, we digress to recall a geometric construction of the $G_{\mathcal{O}}$ equivariant semi-infinite cohomology functor

$$C^{\frac{\infty}{2}}(\hat{\mathfrak{g}}, \mathfrak{g}_{\mathcal{O}}, G_{\mathcal{O}}; (\cdot)) : \hat{\mathfrak{g}}_{-\mathrm{Tate}}\text{-}\mathrm{Mod}^{G_{\mathcal{O}}} \to \mathrm{Vect} \qquad \text{as} \qquad \mathrm{p}_{(BG_{\mathcal{K}})^{\wedge}_{BG_{\mathcal{O}}}}\text{-} : \mathrm{IndCoh}^!_{\chi_{-\mathrm{Tate}}}((BG_{\mathcal{K}})^{\wedge}_{BG_{\mathcal{O}}}) \to \mathrm{Vect},$$

as in Example 16.1.23 and Equation 16.1.3, following Appendix A of [Ras16] and sections 8 and 9 of [Ras20b]. In particular, we emphasize as much as possible the relationship with the construction of the pushforward functor

$$\mathrm{p}_{(\mathcal{Y}_x)^{\wedge}_{\mathcal{Y}_0}}\text{-} : \mathrm{IndCoh}^!((\mathcal{Y}_x)^{\wedge}_{\mathcal{Y}_0}) \to \mathrm{Vect}$$

used to define chiral differential operators in Subsection 15.2, under the analogy $Y = BG$.

In place of $G_{\mathcal{K}}$ and $G_{\mathcal{O}}$, we consider a general polarizeable Tate group indscheme $H$ and compact open subgroup $H^0 \subset H$, in the sense of 7.17 in [Ras20b]. Further, we let $\mathfrak{h} = \mathrm{Lie}(H)$ be the corresponding Tate Lie algebra, and $\mathfrak{h}^0 = \mathrm{Lie}(H^0)$ the corresponding compact open subalgebra. The basic premise of the analogy is that

$$B\exp(\mathfrak{h}) := BH^{\wedge}_{\{e\}} = H_{\mathrm{dR}}/H \qquad \text{and thus} \qquad \mathfrak{h}\text{-}\mathrm{Mod} := \mathrm{IndCoh}^!(BH^{\wedge}_{\{e\}}) = D^!(H)^{H,w} \ .$$

As outlined in Remark 16.1.16 following sections 8 and 9 of [Ras20b], in analogy with the case that $Y$ is a scheme as explained in Example 15.2.11 following Definition 15.1.30, Remark 15.1.22 and ultimately [KV06], the desired pushforward functor is constructed as follows:

There is a canonical pushforward functor

$$\mathrm{p}_{\mathrm{B}H^\wedge_{\{e\}},\bullet} : \mathrm{IndCoh}^\bullet(\mathrm{B}H^\wedge_{\{e\}}) := D^*(H)_{H,w} \to \mathrm{Vect} ,$$

as well as a canonical (up to the choice of compact open subgroup $H^0$) identification

$$\mathrm{IndCoh}^!_{\chi_{-\mathrm{Tate}}}(\mathrm{B}H^\wedge_{\{e\}}) := \left(\mathrm{Vect}_{\chi_{-\mathrm{Tate}}} \otimes D^!(H)\right)^{H,w} \cong D^*(H)_{H,w} = \mathrm{IndCoh}^\bullet(\mathrm{B}H^\wedge_{\{e\}}) ,$$

so that we obtain an induced pushforward functor

$$\mathrm{p}_{\mathrm{B}H^\wedge_{\{e\}},\bullet} : \mathrm{IndCoh}^!_{\chi_{-\mathrm{Tate}}}(\mathrm{B}H^\wedge_{\{e\}}) \to \mathrm{Vect} \qquad \text{which identifies with} \qquad C^{\frac{\infty}{2}}(\mathfrak{h}, \mathfrak{h}^0; (\cdot)) : \mathfrak{h}_{-\mathrm{Tate}}\text{-Mod} \to \mathrm{Vect}$$

under the equivalence of Remark 16.1.9, where the latter is the functor of classical semi-infinite cohomology relative to $\mathfrak{h}_0$. The $H^0$-equivariant analogue is constructed similarly, as outlined in Remark 16.1.16, and explained below.

We now proceed to define the semi-infinite cohomology functor in a series of increasingly general situations, emphasizing the analogy with the construction of $\mathrm{p}_{(\mathcal{Y}_x)^\wedge_{y_0},\bullet}$ in Subsection 15.2.

### 16.2.1. *Chevalley-Eilenberg cohomology for pro-finite type Lie algebras.* We begin by recalling the geometric construction of the usual Chevalley-Eilenberg Lie algebra cochains functor on pro-finite type Lie algebras.

*Example* 16.2.1. Let $H$ be a pro-finite type group scheme, presented as a limit of finite type quotient groups

$$H = \lim_{i \in \mathcal{I}} H/H_i = \lim \left[\ldots \to H/H_i \to H/H_j \to \ldots\right] \qquad , \text{for} \qquad \ldots \subset H_i \subset H_j \subset \ldots \subset H$$

a sequence of normal subgroups of finite codimension in $H$. There is a corresponding presentation

$$\mathrm{B}H^\wedge_{\{e\}} = \lim_{i \in \mathcal{I}} \mathrm{B}(H/H_i)^\wedge_{\{e\}} = \lim \left[\ldots \to \mathrm{B}(H/H_i)^\wedge_{\{e\}} \xrightarrow{\varphi_{ij}} \mathrm{B}(H/H_j)^\wedge_{\{e\}} \to \ldots\right] .$$

The presentation of $\mathrm{B}H^\wedge_{\{e\}}$ in the preceding example is analogous to that of Equation 15.2.2, and the analogue of the resulting presentation of $\mathrm{IndCoh}^!$ in Equation 15.2.4 is given as follows:

*Example* 16.2.2. The category $\mathfrak{h}\text{-Mod} := \mathrm{IndCoh}^!(\mathrm{B}H^\wedge_{\{e\}}) \in \mathrm{DGCat}$ is presented as

(16.2.1)
$$\mathrm{IndCoh}^!(\mathrm{B}H^\wedge_{\{e\}}) = \underset{i \in \mathcal{I}}{\mathrm{colim}} \, \mathrm{IndCoh}^!(\mathrm{B}(H/H_i)^\wedge_{\{e\}}) = \mathrm{colim} \left[\ldots \leftarrow \mathrm{IndCoh}^!(\mathrm{B}(H/H_i)^\wedge_{\{e\}}) \xleftarrow{\varphi^!_{ij}} \mathrm{IndCoh}^!(\mathrm{B}(H/H_j)^\wedge_{\{e\}}) \leftarrow \ldots\right]$$

or equivalently, identifying $\mathfrak{h}/\mathfrak{h}_i\text{-Mod} = \mathrm{IndCoh}^!(\mathrm{B}(H/H_i)^\wedge_{\{e\}})$, as

(16.2.2)
$$\mathfrak{h}\text{-Mod} = \underset{i \in \mathcal{I}}{\mathrm{colim}} \, \mathfrak{h}/\mathfrak{h}_i\text{-Mod} = \mathrm{colim} \left[\ldots \leftarrow \mathfrak{h}/\mathfrak{h}_i\text{-Mod} \xleftarrow{\mathrm{oblv}_{ij}} \mathfrak{h}/\mathfrak{h}_j\text{-Mod} \leftarrow \ldots\right] ,$$

where

(16.2.3)
$$\varphi^!_{ij} = \mathrm{oblv}_{ij} : \mathfrak{h}/\mathfrak{h}_j\text{-Mod} \to \mathfrak{h}/\mathfrak{h}_i\text{-Mod}$$

is the forgetful functor on Lie algebra modules induced by the map of Lie algebras $\mathfrak{h}/\mathfrak{h}_i \to \mathfrak{h}/\mathfrak{h}_j$.

*Example* 16.2.3. More generally, for $\mathfrak{h}$ a profinite type Lie algebra, presented as a limit $\mathfrak{h} = \lim_i \mathfrak{h}/\mathfrak{h}_i$ of quotients by finite type subalgebras $\mathfrak{h}_i \subset \mathfrak{h}$, we have a presentation of $\mathfrak{h}$-Mod $\in$ DGCat as in Equation 16.2.11 above.

Now, following Example 15.2.8 and in turn Example 15.1.26, we have:

*Example* 16.2.4. The category $\mathfrak{h}$-Mod $\in$ DGCat admits an equivalent presentation
$$(16.2.4)$$
$$\mathrm{IndCoh}^!(\mathrm{B}(H)^\wedge_{\{e\}}) = \lim \left[ \ldots \to \mathrm{IndCoh}^!(\mathrm{B}(H/H_i)^\wedge_{\{e\}}) \xrightarrow{\underline{\mathrm{Hom}}(\omega_{B(H_j/H_i)^\wedge_{\{e\}}}, \cdot)} \mathrm{IndCoh}^!(\mathrm{B}(H/H_j)^\wedge_{\{e\}}) \to \ldots \right]$$

or equivalently

$$(16.2.5) \qquad \mathfrak{h}\text{-Mod} = \lim_{i \in \mathcal{I}} \mathfrak{h}/\mathfrak{h}_i\text{-Mod} = \lim \left[ \ldots \to \mathfrak{h}/\mathfrak{h}_i\text{-Mod} \xrightarrow{C^\bullet(\mathfrak{h}_j/\mathfrak{h}_i; (\cdot))} \mathfrak{h}/\mathfrak{h}_j\text{-Mod} \to \ldots \right] \, ,$$

where

$$\underline{\mathrm{Hom}}(\omega_{B(H_j/H_i)^\wedge_{\{e\}}}, \cdot) = \varphi_{ij\bullet} \circ \left( (\cdot) \otimes_{\mathcal{O}_{B(H_j/H_i)^\wedge_{\{e\}}}} \omega^{-1}_{B(H_j/H_i)^\wedge_{\{e\}}} \right) = C^\bullet(\mathfrak{h}_j/\mathfrak{h}_i; (\cdot)) : \mathfrak{h}/\mathfrak{h}_i\text{-Mod} \to \mathfrak{h}/\mathfrak{h}_j\text{-Mod}$$

denotes the functor of internal Hom out of the relative dualizing sheaf as defined in Example 15.1.26, which identifies algebraically with the the relative Chevalley-Eilenberg cochains functor, as we now explain:

The usual IndCoh pushforward functor identifies with the Chevalley-Eilenberg chains functor

$$\varphi_{ij\bullet} = C_\bullet(\mathfrak{h}_j/\mathfrak{h}_i; (\cdot)) : \mathfrak{h}/\mathfrak{h}_i\text{-Mod} \to \mathfrak{h}/\mathfrak{h}_j\text{-Mod}$$

and thus the preceding identification follows from the standard equivalence

$$C^\bullet(\mathfrak{h}_j/\mathfrak{h}_i; (\cdot)) = C_\bullet(\mathfrak{h}_j/\mathfrak{h}_i; (\cdot) \otimes \det(\mathfrak{h}_j/\mathfrak{h}_i)^{-1}[-\dim \mathfrak{h}_j/\mathfrak{h}_i]) \, .$$

In particular, the relative Chevalley-Eilenberg cochains functor is right adjoint to the forgetful functor of Equation 16.2.3 above, and thus the limit presentation of $\mathfrak{h}$-Mod $\in$ DGCat claimed above follows by passing to right adjoints.

Thus, we have the analogue of remarks 15.1.27 and 16.2.5 in the present setting:

*Remark* 16.2.5. Concretely, an object $M \in \mathfrak{h}$-Mod is given by an assignment

$$i \mapsto (M_i \in \mathfrak{h}/\mathfrak{h}_i\text{-Mod}) \qquad [j \to i] \mapsto \left[ C^\bullet(\mathfrak{h}_j/\mathfrak{h}_i; M_i) \xrightarrow{\cong} M_j \right]$$

for each $i, j \in \mathcal{I}$ and arrow $j \to i$ in $\mathcal{I}$.

Moreover, following examples 15.1.28 and 15.2.10 we have:

*Example* 16.2.6. The presentation of Example 16.2.4 above induces a functor

$$\underline{\mathrm{Hom}}(\omega_{BH^\wedge_{\{e\}}}, \cdot) = C^\bullet(\mathfrak{h}; (\cdot)) : \mathfrak{h}\text{-Mod} \to \mathrm{Vect} \qquad \text{defined by} \qquad M = (M_i)_{i \in \mathcal{I}} \mapsto C^\bullet(\mathfrak{h}/\mathfrak{h}_i; M_i) \, ,$$

which is independent of the choice of $i \in \mathcal{I}$ as in *loc. cit.*.

**16.2.2.** *Semi-infinite cohomology for Tate Lie algebras filtered by compact open subalgebras.* Now, we consider the case of a (polarizeable) Tate group indscheme, with the additional hypothesis that $H$ is filtered by compact open subgroups in the following sense:

*Example* 16.2.7. Let $H$ be a polarizable Tate group indscheme presented as a colimit

$$H = \operatorname{colim}_{k \in \mathcal{L}} H^k = \operatorname{colim}\left[\ldots \leftarrow H^l \leftarrow H^k \leftarrow \ldots\right]$$

where each $H^k$ is a pro-finite type group scheme as above. There is a corresponding presentation

$$(\mathrm{B}H)^{\wedge}_{\{e\}} = \operatorname{colim}_{k \in \mathcal{L}}(\mathrm{B}H^k)^{\wedge}_{\{e\}} = \operatorname{colim}\left[\ldots \leftarrow (\mathrm{B}H^l)^{\wedge}_{\{e\}} \xleftarrow{\iota^{kl}} (\mathrm{B}H^k)^{\wedge}_{\{e\}} \leftarrow \ldots\right] \ .$$

*Warning* 16.2.8. This assumption does not hold in our main example of interest $H = G_{\mathcal{K}}$ for $G$ a finite type, reductive, affine algebraic group. However, it is required to maintain the analogy with the construction in the case that $Y$ is a scheme.

The presentation of $\mathrm{B}(H)^{\wedge}_{\{e\}}$ in the preceding example is analogous to that of Equation 15.2.1, and the analogue of the resulting presentation of $\mathrm{IndCoh}^!$ in Equation 15.2.3 is given as follows:

*Example* 16.2.9. The category $\mathfrak{h}\text{-Mod} = \mathrm{IndCoh}^!((\mathrm{B}H)^{\wedge}_{\{e\}}) \in \mathrm{DGCat}$ is presented as

$$\mathrm{IndCoh}^!((\mathrm{B}H)^{\wedge}_{\{e\}}) = \lim_{k \in \mathcal{L}} \mathrm{IndCoh}^!((\mathrm{B}H^k)^{\wedge}_{\{e\}}) = \lim\left[\ldots \to \mathrm{IndCoh}^!((\mathrm{B}H^l)^{\wedge}_{\{e\}}) \xrightarrow{\iota^{kl!}} \mathrm{IndCoh}^!((\mathrm{B}H^k)^{\wedge}_{\{e\}}) \to \ldots\right] \ ,$$

or equivalently

$$(16.2.6) \qquad \mathfrak{h}\text{-Mod} = \lim_{k \in \mathcal{L}} \mathfrak{h}^k\text{-Mod} = \lim\left[\ldots \to \mathfrak{h}^l\text{-Mod} \xrightarrow{\mathrm{oblv}^{kl}} \mathfrak{h}^k\text{-Mod} \to \ldots\right] \ ,$$

where

$$(16.2.7) \qquad \iota^{kl!} = \mathrm{oblv}^{kl} : \mathfrak{h}^l\text{-Mod} \to \mathfrak{h}^k\text{-Mod}$$

is the forgetful functor on Lie algebra modules induced by the map of Lie algebras $\mathfrak{h}^k \hookrightarrow \mathfrak{h}^l$.

In particular, in analogy with Remark 15.2.3, we have:

*Remark* 16.2.10. Concretely, an object $M \in \mathfrak{h}\text{-Mod}$ is specified by an assignment

$$k \mapsto \left(M^k \in \mathfrak{h}^k\text{-Mod}\right) \qquad [k \to l] \mapsto \left[\mathrm{oblv}^{kl} M^l \xrightarrow{\cong} M^k\right]$$

defined for each $k, l \in \mathcal{L}$ and map $k \to l$.

Further, following Example 15.1.15, we have:

*Example* 16.2.11. There is an equivalent presentation

$$\mathrm{IndCoh}^!((\mathrm{B}H)^{\wedge}_{\{e\}}) = \operatorname*{colim}_{k \in \mathcal{L}} \mathrm{IndCoh}^!((\mathrm{B}H^k)^{\wedge}_{\{e\}}) = \operatorname{colim}\left[\ldots \leftarrow \mathrm{IndCoh}^!((\mathrm{B}H^l)^{\wedge}_{\{e\}}) \xleftarrow{\iota^{kl}_{\bullet}} \mathrm{IndCoh}^!((\mathrm{B}H^k)^{\wedge}_{\{e\}}) \leftarrow \ldots\right] \ ,$$

or equivalently

$$(16.2.8) \qquad \mathfrak{h}\text{-Mod} = \operatorname*{colim}_{k \in \mathcal{L}} \mathfrak{h}^k\text{-Mod} = \operatorname{colim}\left[\ldots \leftarrow \mathfrak{h}^l\text{-Mod} \xleftarrow{\mathrm{Ind}^{\mathfrak{h}^l}_{\mathfrak{h}^k}} \mathfrak{h}^k\text{-Mod} \leftarrow \ldots\right] \ ,$$

where

$$\iota^{kl}_{\bullet} = \mathrm{Ind}^{\mathfrak{h}^l}_{\mathfrak{h}^k} : \mathfrak{h}^k\text{-Mod} \to \mathfrak{h}^l\text{-Mod}$$

is the induction functor on Lie algebra modules induced by the map of Lie algebras $\mathfrak{h}^k \hookrightarrow \mathfrak{h}^l$. In particular, the induction functor is left adjoint to the forgetful functor of Equation 16.2.7, and thus the colimit presentation of $\mathfrak{h}$-Mod $\in$ DGCat claimed above follows by passing to left adjoints.

*Example* 16.2.12. More generally, for $\mathfrak{h}$ a Tate Lie algebra presented as a colimit $\mathfrak{h} = \operatorname{colim}_{k \in \mathcal{L}} \mathfrak{h}^k$ of compact open subalgebras $\mathfrak{h}^k \subset \mathfrak{h}$, we have presentations of $\mathfrak{h}$-Mod $\in$ DGCat as in equations 16.2.6 and 16.2.8 above.

*Remark* 16.2.13. For each $k \in \mathcal{L}$, the category $\mathfrak{h}^k$-Mod admits presentations as in examples 16.2.2 and 16.2.4. In particular, following Example 16.2.6, there are natural Chevalley-Eilenberg cochains functors

$$\underline{\operatorname{Hom}}(\omega_{(\mathrm{B}H^k)^{\wedge}_{\{e\}}}, \cdot) = C^{\bullet}(\mathfrak{h}^k; (\cdot)) : \mathfrak{h}\text{-Mod} \to \mathrm{Vect} \ .$$

Following Example 15.2.11, and in turn Remark 15.1.29 and Definition 15.1.30, we have:

*Example* 16.2.14. The functors of the preceding Example do *not* satisfy the compatibility condition $C^{\bullet}(\mathfrak{h}^l; \operatorname{Ind}_{\mathfrak{h}^k}^{\mathfrak{h}^l}(\cdot)) \cong C^{\bullet}(\mathfrak{h}^k; (\cdot))$, and thus do not canonically extend to the limit category $\mathfrak{h}$-Mod. However, fixing a compact open subalgebra $\mathfrak{h}^0 \subset \mathfrak{h}$ and taking $\mathcal{E}^k = \det(\mathfrak{h}^k/\mathfrak{h}^0)[\dim(\mathfrak{h}^k/\mathfrak{h}^0)]$, we obtain functors

$$C^{\frac{\infty}{2}}(\mathfrak{h}^k, \mathfrak{h}^0; (\cdot)) := C^{\bullet}(\mathfrak{h}^k; (\cdot) \otimes \det(\mathfrak{h}^k/\mathfrak{h}^0)[\dim(\mathfrak{h}^k/\mathfrak{h}^0)]) : \mathfrak{h}^k\text{-Mod} \to \mathrm{Vect}$$

such that we have the desired natural equivalences:

$$\begin{aligned}
C^{\frac{\infty}{2}}(\mathfrak{h}^l, \mathfrak{h}^0; \operatorname{Ind}_{\mathfrak{h}^k}^{\mathfrak{h}^l}(M^k)) &= C^{\bullet}(\mathfrak{h}^l; \operatorname{Ind}_{\mathfrak{h}^k}^{\mathfrak{h}^l}(M^k) \otimes \det(\mathfrak{h}^l/\mathfrak{h}^0)[\dim(\mathfrak{h}^l/\mathfrak{h}^0)]) \\
&\cong C^{\bullet}(\mathfrak{h}^l; \operatorname{Ind}_{\mathfrak{h}^k}^{\mathfrak{h}^l}(M^k) \otimes \det(\mathfrak{h}^l/\mathfrak{h}^k)[\dim(\mathfrak{h}^l/\mathfrak{h}^k)] \otimes \det(\mathfrak{h}^k/\mathfrak{h}^0)[\dim(\mathfrak{h}^k/\mathfrak{h}^0)]) \\
&\cong C^{\bullet}(\mathfrak{h}^k; M^k \otimes \det(\mathfrak{h}^k/\mathfrak{h}^0)[\dim(\mathfrak{h}^k/\mathfrak{h}^0)]) \\
&= C^{\frac{\infty}{2}}(\mathfrak{h}^k, \mathfrak{h}^0; M^k)
\end{aligned}$$

Thus, for $M = (M^k)_{k \in \mathcal{L}} \in \mathfrak{h}$-Mod as in Remark 16.2.10, there are canonical maps

$$C^{\frac{\infty}{2}}(\mathfrak{h}^k, \mathfrak{h}^0; M^k) \cong C^{\frac{\infty}{2}}(\mathfrak{h}^l, \mathfrak{h}^0; \operatorname{Ind}_{\mathfrak{h}^k}^{\mathfrak{h}^l}(\operatorname{oblv}^{kl} M^l)) \to C^{\frac{\infty}{2}}(\mathfrak{h}^l, \mathfrak{h}^0; M^l)$$

defined by the counit of the $(\operatorname{Ind}_{\mathfrak{h}^k}^{\mathfrak{h}^l}, \operatorname{oblv}^{kl})$ adjunction, and we define

$$C^{\frac{\infty}{2}}(\mathfrak{h}, \mathfrak{h}^0; (\cdot)) : \mathfrak{h}\text{-Mod} \to \mathrm{Vect} \qquad \text{by} \qquad M = (M^k)_{k \in \mathcal{L}} \mapsto \operatorname*{colim}_{k \in \mathcal{L}} C^{\frac{\infty}{2}}(\mathfrak{h}^k, \mathfrak{h}^0; M^k) \ ,$$

the desired functor of semi-infinite cohomology with respect to $\mathfrak{h}$ relative to $\mathfrak{h}^0$.

16.2.3. *Equivariant semi-infinite cohomology for Tate Lie algebras filtered by compact open subalgebras.* Let $H$ be a polarizable Tate group indscheme filtered by compact open subgroups as in Example 16.2.7, and moreover fix a compact open subgroup $K \subset H$. Then in analogy with *loc. cit.*, we have:

*Example* 16.2.15. There is a presentation

$$(\mathrm{B}H)^{\wedge}_K = \operatorname{colim}_{k \in \mathcal{L}}(\mathrm{B}H^k)^{\wedge}_K = \operatorname{colim}\left[\ldots \leftarrow (\mathrm{B}H^l)^{\wedge}_K \xleftarrow{\iota^{kl}} (\mathrm{B}H^k)^{\wedge}_K \leftarrow \ldots\right] \ ,$$

noting that $K \subset H^k, H^l$ for $k, l >> 0$. Further, in analogy with Example 16.2.1, for each $k \in \mathcal{L}$, there is a presentation

$$\mathrm{B}(H^k)_K^\wedge = \lim_{i \in \mathcal{I}} \mathrm{B}(H^k/H_i)_{K/H_i}^\wedge = \lim \left[ \dots \to \mathrm{B}(H^k/H_i)_{K/H_i}^\wedge \xrightarrow{\varphi_{ij}} \mathrm{B}(H^k/H_j)_{K/H_j}^\wedge \to \dots \right] .$$

Following examples 16.2.9 and 16.2.2, we have:

*Example* 16.2.16. The category $\mathfrak{h}\text{-Mod}^K = \mathrm{IndCoh}^!((\mathrm{B}H)_K^\wedge) \in \mathrm{DGCat}$ is presented as

$$\mathrm{IndCoh}^!((\mathrm{B}H)_K^\wedge) = \lim_{k \in \mathcal{L}} \mathrm{IndCoh}^!((\mathrm{B}H^k)_K^\wedge) = \lim \left[ \dots \to \mathrm{IndCoh}^!((\mathrm{B}H^l)_K^\wedge) \xrightarrow{\iota^{kl!}} \mathrm{IndCoh}^!((\mathrm{B}H^k)_K^\wedge) \to \dots \right] ,$$

or equivalently

$$(16.2.9) \qquad \mathfrak{h}\text{-Mod}^K = \lim_{k \in \mathcal{L}} \mathfrak{h}^k\text{-Mod}^K = \lim \left[ \dots \to \mathfrak{h}^l\text{-Mod}^K \xrightarrow{\mathrm{oblv}^{kl}} \mathfrak{h}^k\text{-Mod}^K \to \dots \right] .$$

Further, for each $k \in \mathcal{L}$, the category $\mathfrak{h}^k\text{-Mod}^K := \mathrm{IndCoh}^!(\mathrm{B}(H^k)_K^\wedge) \in \mathrm{DGCat}$ is presented as
$$(16.2.10)$$
$$\mathrm{IndCoh}^!(\mathrm{B}(H^k)_K^\wedge) = \mathrm{colim} \left[ \dots \leftarrow \mathrm{IndCoh}^!(\mathrm{B}(H^k/H_i)_{K/H_i}^\wedge) \xleftarrow{\varphi_{ij}^!} \mathrm{IndCoh}^!(\mathrm{B}(H^k/H_j)_{K/H_j}^\wedge) \leftarrow \dots \right]$$

or equivalently, identifying $\mathfrak{h}^k/\mathfrak{h}_i\text{-Mod}^{K/H_i} = \mathrm{IndCoh}^!(\mathrm{B}(H^k/H_i)_{K/H_i}^\wedge)$, as
$$(16.2.11)$$
$$\mathfrak{h}^k\text{-Mod}^K = \mathrm{colim}_{i \in \mathcal{I}} \mathfrak{h}^k/\mathfrak{h}_i\text{-Mod}^{K/H_i} = \mathrm{colim} \left[ \dots \leftarrow \mathfrak{h}^k/\mathfrak{h}_i\text{-Mod}^{K/H_i} \xleftarrow{\mathrm{oblv}_{ij}} \mathfrak{h}^k/\mathfrak{h}_j\text{-Mod}^{K/H_j} \leftarrow \dots \right] .$$

Further, following Example 16.2.4, for each $k \in \mathcal{L}$ we have:

*Example* 16.2.17. The category $\mathfrak{h}^k\text{-Mod}^K \in \mathrm{DGCat}$ admits an equivalent presentation
$$(16.2.12)$$
$$\mathrm{IndCoh}^!((\mathrm{B}H^k)_K^\wedge) = \lim \left[ \dots \to \mathrm{IndCoh}^!(\mathrm{B}(H^k/H_i)_{K/H_i}^\wedge) \xrightarrow{\underline{\mathrm{Hom}}(\omega_{\mathrm{B}(H_j/H_i)_{K/H_i}^\wedge}, \cdot)} \mathrm{IndCoh}^!(\mathrm{B}(H^k/H_j)_{K/H_j}^\wedge) \to \dots \right]$$

or equivalently
$$(16.2.13)$$
$$\mathfrak{h}^k\text{-Mod}^K = \lim_{i \in \mathcal{I}} \mathfrak{h}/\mathfrak{h}_i\text{-Mod}^{K/H_i} = \lim \left[ \dots \to \mathfrak{h}/\mathfrak{h}_i\text{-Mod}^{K/H_i} \xrightarrow{C^\bullet(\mathfrak{h}_j/\mathfrak{h}_i, H_j/H_i; (\cdot))} \mathfrak{h}/\mathfrak{h}_j\text{-Mod}^{K/H_j} \to \dots \right] ,$$

where

$$\underline{\mathrm{Hom}}(\omega_{\mathrm{B}(H_j/H_i)_{K/H_i}^\wedge}, \cdot) = C^\bullet(\mathfrak{h}_j/\mathfrak{h}_i, H_j/H_i; (\cdot)) : \mathfrak{h}/\mathfrak{h}_i\text{-Mod}^{K/H_i} \to \mathfrak{h}/\mathfrak{h}_j\text{-Mod}^{K/H_j}$$

denotes the functor of internal Hom out of the relative dualizing sheaf as defined in Example 15.1.26, which identifies algebraically with the the relative Chevalley-Eilenberg cochains functor for Harish-Chandra pairs.

In particular, following Remark 16.2.13 and Example 16.2.6, we have

*Example* 16.2.18. The presentation of Example 16.2.17 above induces a functor

$$\underline{\mathrm{Hom}}(\omega_{(\mathrm{B}H^k)_K^\wedge}, \cdot) = C^\bullet(\mathfrak{h}^k, K; (\cdot)) : \mathfrak{h}^k\text{-Mod}^K \to \mathrm{Vect} .$$

Finally, following Example 16.2.14, we arrive at our desired definition:

*Definition* 16.2.19. For each $k \in \mathcal{L}$, the $K$-equivariant semi-infinite cohomology functor with respect to $\mathfrak{h}^k$ relative to $\mathfrak{h}^0$ is defined by

$$C^{\frac{\infty}{2}}(\mathfrak{h}^k, \mathfrak{h}^0, K; (\cdot)) := C^\bullet(\mathfrak{h}^k, K; (\cdot) \otimes \det(\mathfrak{h}^k/\mathfrak{h}^0)[\dim(\mathfrak{h}^k/\mathfrak{h}^0)]) : \mathfrak{h}^k\text{-Mod}^K \to \mathrm{Vect} \ .$$

The functor of $K$-equivariant semi-infinite cohomology with respect to $\mathfrak{h}$ relative to $\mathfrak{h}^0$ is defined by

$$C^{\frac{\infty}{2}}(\mathfrak{h}, \mathfrak{h}^0, K; (\cdot)) := \operatorname*{colim}_{k \in \mathcal{L}} C^{\frac{\infty}{2}}(\mathfrak{h}^k, \mathfrak{h}^0, K; (\cdot)) : \mathfrak{h}\text{-Mod}^K \to \mathrm{Vect} \ .$$

## 17. The action of the three dimensional A model on chiral differential operators

In this subsection, we explain the analogues of Theorem 14.2.1 and propositions 14.3.4 and 14.3.9 in the case $Y = N/G$ is a quotient stack, using the results of Section 16. To begin, in analogy with Proposition 13.1.9, we give a compatible presentation of the factorization space $\mathcal{J}^{\mathrm{mer}}(G, N)^\wedge_{\mathcal{J}(G,N)}$:

*Proposition* 17.0.1. There is an equivalence of factorization spaces

$$\mathcal{J}^{\mathrm{mer}}(G, N)^\wedge_{\mathcal{J}(G,N)} = \mathcal{J}(G, N)_{\mathrm{dR}} \times_{\mathcal{J}^{\mathrm{mer}}(G,N)_{\mathrm{dR}}} \mathcal{J}^{\mathrm{mer}}(G, N) \cong \left( \mathcal{J}(G, N)_{\mathrm{dR}} \times_{\mathcal{J}^{\mathrm{mer}}(N)_{\mathrm{dR}}} \mathcal{J}^{\mathrm{mer}}(N) \right) / \mathcal{J}^{\mathrm{mer}}(G) \ .$$

*Remark* 17.0.2. Concretely, over each point $x \in X$ the preceding equivalence is given by

$$\mathcal{J}^{\mathrm{mer}}(G, N)^\wedge_{\mathcal{J}(G,N)_x} = [N/G]_{\mathcal{O},\mathrm{dR}} \times_{[N/G]_{\mathcal{K},\mathrm{dR}}} [N/G]_{\mathcal{K}} \cong \left( [(N_{\mathcal{O}} \times G_{\mathcal{K}})/G_{\mathcal{O}}]_{\mathrm{dR}} \times_{N_{\mathcal{K},\mathrm{dR}}} N_{\mathcal{K}} \right) / G_{\mathcal{K}}$$

*Remark* 17.0.3. Under the preceding equivalence, the fibre sequence of Equation 16.1.1 is given by

$$N_{\mathcal{O},\mathrm{dR}} \times_{N_{\mathcal{K},\mathrm{dR}}} N_{\mathcal{K}} \hookrightarrow \left( [(N_{\mathcal{O}} \times G_{\mathcal{K}})/G_{\mathcal{O}}]_{\mathrm{dR}} \times_{N_{\mathcal{K},\mathrm{dR}}} N_{\mathcal{K}} \right) / G_{\mathcal{K}} \twoheadrightarrow G_{\mathcal{O},\mathrm{dR}} \backslash G_{\mathcal{K},\mathrm{dR}} / G_{\mathcal{K}}$$

and similarly the prequotient by $G_{\mathcal{K}}$ is given by

$$N_{\mathcal{O},\mathrm{dR}} \times_{N_{\mathcal{K},\mathrm{dR}}} N_{\mathcal{K}} \hookrightarrow [(N_{\mathcal{O}} \times G_{\mathcal{K}})/G_{\mathcal{O}}]_{\mathrm{dR}} \times_{N_{\mathcal{K},\mathrm{dR}}} N_{\mathcal{K}} \twoheadrightarrow \mathrm{Gr}_{G,\mathrm{dR}} \ .$$

*Example* 17.0.4. There is a canonical $G_{\mathcal{K}}$ equivariant equivalence

$$\left( (N_{\mathcal{K}})^\wedge_{N_{\mathcal{O}}} \times G_{\mathcal{K},\mathrm{dR}} \right) / G_{\mathcal{O},\mathrm{dR}} = \left( (N_{\mathcal{O},\mathrm{dR}} \times_{N_{\mathcal{K},\mathrm{dR}}} N_{\mathcal{K}}) \times G_{\mathcal{K},\mathrm{dR}}) \right) / G_{\mathcal{O},\mathrm{dR}} \cong [(N_{\mathcal{O}} \times G_{\mathcal{K}})/G_{\mathcal{O}}]_{\mathrm{dR}} \times_{N_{\mathcal{K},\mathrm{dR}}} N_{\mathcal{K}}$$

inducing an equivalence of weak $G_{\mathcal{K}}$ categories

$$\left( \mathrm{IndCoh}^!((N_{\mathcal{K}})^\wedge_{N_{\mathcal{O}}}) \otimes D^!(G_{\mathcal{K}}) \right)^{G_{\mathcal{O}}} \cong D^!(N_{\mathcal{O}} \times G_{\mathcal{K}})^{G_{\mathcal{O}}} \underset{D^!(N_{\mathcal{K}})}{\otimes} \mathrm{IndCoh}^!(N_{\mathcal{K}})$$

where $D^!(N_{\mathcal{K}})$ is equipped with the $\otimes^!$ (symmetric) monoidal structure. In particular, there is an equivalence of the corresponding categories of weak $G_{\mathcal{K}}$ invariants

$$\mathrm{IndCoh}^!(\mathcal{J}^{\mathrm{mer}}(G, N)^\wedge_{\mathcal{J}(G,N)x}) \cong \left( D^!(N_{\mathcal{O}} \times G_{\mathcal{K}})^{G_{\mathcal{O}}} \underset{D^!(N_{\mathcal{K}})}{\otimes} \mathrm{IndCoh}^!(N_{\mathcal{K}}) \right)^{G_{\mathcal{K}},w} \ .$$

Thus, in analogy with the description of Remark 13.4.3, we have the following:

*Remark* 17.0.5. The object $\tilde{\mathcal{D}}(G, N)_x \in \mathrm{IndCoh}^!(\mathcal{J}^{\mathrm{mer}}(G, N)^\wedge_{\mathcal{J}(G,N)x})$ is given, under the equivalence of the preceding Remark, by

$$\left( \mathrm{IndCoh}^!((N_{\mathcal{K}})^\wedge_{N_{\mathcal{O}}}) \otimes D^!(G_{\mathcal{K}}) \right)^{G_{\mathcal{O}} \times (G_{\mathcal{K}}, w)} \overset{\cong}{\longrightarrow} \left( D^!(N_{\mathcal{O}} \times G_{\mathcal{K}})^{G_{\mathcal{O}}} \underset{D^!(N_{\mathcal{K}})}{\otimes} \mathrm{IndCoh}^!(N_{\mathcal{K}}) \right)^{G_{\mathcal{K}},w}$$

$$\tilde{\mathcal{D}}(G, N)_x := \omega_{(N_{\mathcal{K}})^\wedge_{N_{\mathcal{O}}}} \boxtimes \omega_{G_{\mathcal{K}}} \longmapsto \omega_{N_{\mathcal{O}}} \boxtimes \omega_{G_{\mathcal{K}}} \boxtimes_{N_{\mathcal{K}}} \omega_{N_{\mathcal{K}}} \ .$$

Towards establishing Proposition 14.3.4 in the present context, we recall the analogue of Remark 14.3.3:

*Example* 17.0.6. The convolution diagram relating $\mathcal{Z}(G,N)$ and $\mathcal{J}^{\mathrm{mer}}(G,N)^{\wedge}_{\mathcal{J}(G,N)}$ is given by

$$\mathcal{J}^{\mathrm{mer}}(G,N)^{\wedge}_{(3)} = \mathcal{J}(G,N)_{\mathrm{dR}} \times_{\mathcal{J}^{\mathrm{mer}}(G,N)_{\mathrm{dR}}} \mathcal{J}(G,N)_{\mathrm{dR}} \times_{\mathcal{J}^{\mathrm{mer}}(G,N)_{\mathrm{dR}}} \mathcal{J}^{\mathrm{mer}}(G,N)$$

together with the canonical projections

$$\pi_{12} : \mathcal{J}^{\mathrm{mer}}(G,N)^{\wedge}_{(3)} \to \mathcal{Z}(G,N)_{\mathrm{dR}} \qquad \text{and} \qquad \pi_{13}, \pi_{23} : \mathcal{J}^{\mathrm{mer}}(G,N)^{\wedge}_{(3)} \to \mathcal{J}^{\mathrm{mer}}(G,N)^{\wedge}_{\mathcal{J}(G,N)} \ .$$

*Example* 17.0.7. Concretely, the convolution diagram over each point $x \in X$ relating $\mathcal{Z}(G,N)_x$ and $\mathcal{J}^{\mathrm{mer}}(G,N)^{\wedge}_{\mathcal{J}(G,N),x}$ is given by

$$\mathcal{J}^{\mathrm{mer}}(G,N)^{\wedge}_{(3),x} = [N/G]_{\mathcal{O},\mathrm{dR}} \times_{[N/G]_{\mathcal{K},\mathrm{dR}}} [N/G]_{\mathcal{O},\mathrm{dR}} \times_{[N/G]_{\mathcal{K},\mathrm{dR}}} [N/G]_{\mathcal{K}}$$

together with the canonical projections

$$\pi_{12} : \mathcal{J}^{\mathrm{mer}}(G,N)^{\wedge}_{(3)} \to \mathcal{Z}(G,N)_x \qquad \text{and} \qquad \pi_{13}, \pi_{23} : \mathcal{J}^{\mathrm{mer}}(G,N)^{\wedge}_{(3)} \to \mathcal{J}^{\mathrm{mer}}(G,N)^{\wedge}_{\mathcal{J}(G,N),x} \ .$$

The description analogous to that of Proposition 17.0.1 above is given by

*Proposition* 17.0.8. There is an equivalence of factorization spaces

$$\mathcal{J}^{\mathrm{mer}}(G,N)^{\wedge}_{(3)} \cong \left( \mathcal{T}(G,N)_{\mathrm{dR}} \times_{\mathcal{J}^{\mathrm{mer}}(N)_{\mathrm{dR}}} \mathcal{T}(G,N)_{\mathrm{dR}} \times_{\mathcal{J}^{\mathrm{mer}}(N)_{\mathrm{dR}}} \mathcal{J}^{\mathrm{mer}}(N) \right) / \mathcal{J}^{\mathrm{mer}}(G)$$

such that the convolution diagram of Example 17.0.6 above is equivalent to

$$\begin{array}{ccc}
\mathcal{T}(G,N)_{\mathrm{dR}} \times_{\mathcal{J}^{\mathrm{mer}}(N)_{\mathrm{dR}}} \mathcal{T}(G,N)_{\mathrm{dR}} \times_{\mathcal{J}^{\mathrm{mer}}(N)_{\mathrm{dR}}} \mathcal{J}^{\mathrm{mer}}(N) & \xrightarrow{\pi_{13}} & \mathcal{T}(G,N)_{\mathrm{dR}} \times_{\mathcal{J}^{\mathrm{mer}}(N)_{\mathrm{dR}}} \mathcal{J}^{\mathrm{mer}}(N) \\
\downarrow{\scriptstyle \pi_{12}} & & \\
\mathcal{T}(G,N)_{\mathrm{dR}} \times_{\mathcal{J}^{\mathrm{mer}}(N)_{\mathrm{dR}}} \mathcal{J}(N)_{\mathrm{dR}} & & \mathcal{T}(G,N)_{\mathrm{dR}} \times_{\mathcal{J}^{\mathrm{mer}}(N)_{\mathrm{dR}}} \mathcal{J}^{\mathrm{mer}}(N)
\end{array}$$

after passing to the prequotient by $\mathcal{J}^{\mathrm{mer}}(G)$.

*Remark* 17.0.9. Concretely, over each point $x \in X$, the preceding equivalence is given by

$$\mathcal{J}^{\mathrm{mer}}(G,N)^{\wedge}_{(3),x} \cong \left( [(N_{\mathcal{O}} \times G_{\mathcal{K}})/G_{\mathcal{O}}]_{\mathrm{dR}} \times_{N_{\mathcal{K},\mathrm{dR}}} [(N_{\mathcal{O}} \times G_{\mathcal{K}})/G_{\mathcal{O}}]_{\mathrm{dR}} \times_{N_{\mathcal{K},\mathrm{dR}}} N_{\mathcal{K}} \right) / G_{\mathcal{K}}$$

and the convolution diagram is equivalent to

$$\begin{array}{ccc}
[(N_{\mathcal{O}} \times G_{\mathcal{K}})/G_{\mathcal{O}}]_{\mathrm{dR}} \times_{N_{\mathcal{K},\mathrm{dR}}} [(N_{\mathcal{O}} \times G_{\mathcal{K}})/G_{\mathcal{O}}]_{\mathrm{dR}} \times_{N_{\mathcal{K},\mathrm{dR}}} N_{\mathcal{K}} & \xrightarrow{\pi_{13}} & [(N_{\mathcal{O}} \times G_{\mathcal{K}})/G_{\mathcal{O}}]_{\mathrm{dR}} \times_{N_{\mathcal{K},\mathrm{dR}}} N_{\mathcal{K}} \\
\downarrow{\scriptstyle \pi_{12}} & & \\
[(N_{\mathcal{O}} \times G_{\mathcal{K}})/G_{\mathcal{O}}]_{\mathrm{dR}} \times_{N_{\mathcal{K},\mathrm{dR}}} N_{\mathcal{O},\mathrm{dR}} & & [(N_{\mathcal{O}} \times G_{\mathcal{K}})/G_{\mathcal{O}}]_{\mathrm{dR}} \times_{N_{\mathcal{K},\mathrm{dR}}} N_{\mathcal{K}}
\end{array}$$

after passing to the prequotient by $G_{\mathcal{K}}$.

Motivated by the preceding proposition, we make the following definition:

*Definition* 17.0.10. The category $\mathrm{IndCoh}^!(\mathcal{J}^{\mathrm{mer}}(G,N)^{\wedge}_{(3),x}) \in \mathrm{DGCat}$ is defined by

$$\mathrm{IndCoh}^!(\mathcal{J}^{\mathrm{mer}}(G,N)^{\wedge}_{(3),x}) = \left( D^!(N_{\mathcal{O}} \times G_{\mathcal{K}})^{G_{\mathcal{O}}} \underset{D^!(N_{\mathcal{K}})}{\otimes} D^!(N_{\mathcal{O}} \times G_{\mathcal{K}})^{G_{\mathcal{O}}} \underset{D^!(N_{\mathcal{K}})}{\otimes} \mathrm{IndCoh}^!(N_{\mathcal{K}}) \right)^{G_{\mathcal{K}},w}$$

Moreover, we have:

*Example* 17.0.11. There are natural (partially defined, see Warning 17.0.12 below) functors

$$\pi_{12}^\bullet : D^*(\mathcal{Z}(G,N)_x) \to \mathrm{IndCoh}^!(\mathcal{J}^{\mathrm{mer}}(G,N)_{(3),x}^\wedge)$$

$$\pi_{23}^\bullet : \mathrm{IndCoh}^!(\mathcal{J}^{\mathrm{mer}}(G,N)_{\mathcal{J}(G,N)x}^\wedge) \to \mathrm{IndCoh}^!(\mathcal{J}^{\mathrm{mer}}(G,N)_{(3),x}^\wedge)$$

$$\pi_{13}^! : \mathrm{IndCoh}^!(\mathcal{J}^{\mathrm{mer}}(G,N)_{\mathcal{J}(G,N)x}^\wedge) \to \mathrm{IndCoh}^!(\mathcal{J}^{\mathrm{mer}}(G,N)_{(3),x}^\wedge)$$

$$\pi_{13\bullet} : \mathrm{IndCoh}^!(\mathcal{J}^{\mathrm{mer}}(G,N)_{(3),x}^\wedge) \to \mathrm{IndCoh}^!(\mathcal{J}^{\mathrm{mer}}(G,N)_{\mathcal{J}(G,N)x}^\wedge)$$

defined by

$$\pi_{12}^\bullet := (\cdot) \,\tilde{\boxtimes}\, \underline{\mathbb{K}}_{\mathrm{Gr}_G} \boxtimes_{N_{\mathcal{K}}} \mathcal{O}_{N_{\mathcal{K}}} : D_{G_{\mathcal{O}}}^!([(N_{\mathcal{O}} \times G_{\mathcal{K}})/G_{\mathcal{O}}] \times_{N_{\mathcal{K}}} N_{\mathcal{O}}) \to \mathrm{IndCoh}^!(\mathcal{J}^{\mathrm{mer}}(G,N)_{(3),x}^\wedge)$$

$$\pi_{23}^\bullet = \underline{\mathbb{K}}_{N_{\mathcal{O}}} \,\tilde{\boxtimes}\, \underline{\mathbb{K}}_{\mathrm{Gr}_G} \boxtimes_{N_{\mathcal{K}}} (\cdot) \boxtimes_{N_{\mathcal{K}}} (\cdot) : \left( D^!(N_{\mathcal{O}} \times G_{\mathcal{K}})^{G_{\mathcal{O}}} \underset{D^!(N_{\mathcal{K}})}{\otimes} \mathrm{IndCoh}^!(N_{\mathcal{K}}) \right)^{G_{\mathcal{K}},w} \to \mathrm{IndCoh}^!(\mathcal{J}^{\mathrm{mer}}(G,N)_{(3),x}^\wedge)$$

$$\pi_{13}^! = (\cdot) \boxtimes_{N_{\mathcal{K}}} \omega_{N_{\mathcal{O}}} \tilde{\boxtimes}\, \omega_{\mathrm{Gr}_G} \boxtimes_{N_{\mathcal{K}}} (\cdot) : \left( D^!(N_{\mathcal{O}} \times G_{\mathcal{K}})^{G_{\mathcal{O}}} \underset{D^!(N_{\mathcal{K}})}{\otimes} \mathrm{IndCoh}^!(N_{\mathcal{K}}) \right)^{G_{\mathcal{K}},w} \to \mathrm{IndCoh}^!(\mathcal{J}^{\mathrm{mer}}(G,N)_{(3),x}^\wedge)$$

$$\pi_{13\bullet} = (\cdot) \boxtimes_{N_{\mathcal{K}}} \mathrm{p}_{\mathcal{R}(G,N)x*} \boxtimes_{N_{\mathcal{K}}} (\cdot) : \mathrm{IndCoh}^!(\mathcal{J}^{\mathrm{mer}}(G,N)_{(3),x}^\wedge) \to \left( D^!(N_{\mathcal{O}} \times G_{\mathcal{K}})^{G_{\mathcal{O}}} \underset{D^!(N_{\mathcal{K}})}{\otimes} \mathrm{IndCoh}^!(N_{\mathcal{K}}) \right)^{G_{\mathcal{K}},w}$$

*Warning* 17.0.12. As in Warning 13.4.9, the functors $\pi_{12}^\bullet$ and $\pi_{23}^\bullet$ are only partially defined. For example, they are each defined on the full subcategories defined verbatim but in terms of the full subcategory $D_{\mathrm{rh}}^!(N_{\mathcal{O}} \times G_{\mathcal{K}}) \hookrightarrow D^!(N_{\mathcal{O}} \times G_{\mathcal{K}})$ of holonomic $D$ modules, in the sense of Definition B.8.11, which suffices for the application in Theorem 17.0.15 below.

Moreover, the composition

$$(\pi_{12} \times \pi_{23})^\bullet := \Delta^\bullet \circ (\pi_{12}^\bullet \boxtimes \pi_{23}^\bullet) : D^*(\mathcal{Z}(G,N)_x) \otimes \mathrm{IndCoh}^!(\mathcal{J}^{\mathrm{mer}}(G,N)_{\mathcal{J}(G,N)x}^\wedge) \to \mathrm{IndCoh}^!(\mathcal{J}^{\mathrm{mer}}(G,N)_{(3),x}^\wedge)$$

used in Proposition 17.0.13 below is in fact well-defined on the entire categories, so that the proposition holds as stated. This follows from the fact that $\pi_{12} \times \pi_{23}$ is pro-smooth on the prequotient, even though $\pi_{12}$ and $\pi_{23}$ individually are not, so that the (equivariant) $D$ module pullback functor $(\pi_{12} \times \pi_{23})^*$ is well-defined on the entire category.

Following Proposition 14.3.4, analogously to the statement of Proposition 13.4.8 following Proposition 12.2.4, we have:

*Proposition* 17.0.13. The factorization category

$$\mathrm{IndCoh}_{\mathcal{J}^{\mathrm{mer}}(G,N)_{\mathcal{J}(G,N)}^\wedge} \in D_{\mathcal{Z}(G,N)}^\star\text{-}\mathrm{Mod}(\mathrm{Cat}_{\mathrm{un}}^{\mathrm{fact}}(X))$$

is naturally a factorization $\mathbb{E}_1$-module category over $D_{\mathcal{Z}(G,N)}^\star \in \mathrm{Cat}_{\mathbb{E}_1,\mathrm{un}}^{\mathrm{fact}}(X)$, with respect to the convolution module structure $(\cdot) \star (\cdot) : D_{\mathcal{Z}(G,N)} \otimes \mathrm{IndCoh}_{\mathcal{J}^{\mathrm{mer}}(G,N)_{\mathcal{J}(G,N)}^\wedge} \to \mathrm{IndCoh}_{\mathcal{J}^{\mathrm{mer}}(G,N)_{\mathcal{J}(G,N)}^\wedge}$ defined by the composition

$$D_{\mathcal{Z}(G,N)} \otimes^* \mathrm{IndCoh}_{\mathcal{J}^{\mathrm{mer}}(G,N)_{\mathcal{J}(G,N)}^\wedge} \xrightarrow{\pi_{12}^\bullet \boxtimes \pi_{23}^\bullet} \mathrm{IndCoh}_{(\mathcal{J}^{\mathrm{mer}}(G,N)_{(3)}^\wedge)^{\times 2}} \xrightarrow{\Delta^\bullet} \mathrm{IndCoh}_{\mathcal{J}^{\mathrm{mer}}(G,N)_{(3)}^\wedge} \xrightarrow{\pi_{13\bullet}} \mathrm{IndCoh}_{\mathcal{J}^{\mathrm{mer}}(G,N)_{\mathcal{J}(G,N)}^\wedge} .$$

Further, the pushforward functors $\mathrm{p}_{\mathcal{Z}(G,N)*} : D_{\mathcal{Z}(G,N)}^\star \to D_{\mathrm{Ran}_{X,\mathrm{un}}}^{\otimes^!}$ as in Proposition 13.4.8 and $\mathrm{p}_{\mathcal{J}^{\mathrm{mer}}(G,N)_{\mathcal{J}(G,N)}^\wedge \bullet} : \mathrm{IndCoh}_{\mathcal{J}^{\mathrm{mer}}(G,N)_{\mathcal{J}(G,N)}^\wedge} \to D_{\mathrm{Ran}_{X,\mathrm{un}}}$ as in Definition 16.1.24 define unital, lax compatible factorization functors with respect to the above module structure.

Further, following Corollary 14.3.5 and in turn Corollary 8.2.5, we have:

*Corollary* 17.0.14. The pushforward functor $\mathrm{p}_{\mathcal{J}^{\mathrm{mer}}(G,N)^{\wedge}_{\mathcal{J}(G,N)}\check{\bullet}} : \mathrm{IndCoh}_{\mathcal{J}^{\mathrm{mer}}(G,N)^{\wedge}_{\mathcal{J}(G,N)}} \to D_{\mathrm{Ran}_{X,\mathrm{un}}}$ induces a functor

$$(17.0.1) \qquad \tilde{\mathcal{A}}(G,N)\text{-}\mathrm{Mod}(\mathrm{Alg}^{\mathrm{fact}}_{\mathrm{un}}(\mathrm{IndCoh}_{\mathcal{J}^{\mathrm{mer}}(G,N)^{\wedge}_{\mathcal{J}(G,N)}})) \to \mathcal{A}(G,N)\text{-}\mathrm{Mod}(\mathrm{Alg}^{\mathrm{fact}}_{\mathrm{un}}(X)) \ .$$

Finally, following Remark 8.2.7 and in particular Example 8.2.9, we arrive at the statement of Theorem 14.2.1 in the present generality; in turn, the result can be understood as an enhancement of Proposition A.6.12 and in particular Example A.6.17 via Proposition 17.0.13 above:

*Theorem* 17.0.15. The internal variant of chiral differential operators to $N/G$ of Definition 16.1.11 admits a natural module structure

$$\tilde{\mathcal{D}}^{\mathrm{ch}}(G,N) \in \tilde{\mathcal{A}}(G,N)\text{-}\mathrm{Mod}(\mathrm{Alg}^{\mathrm{fact}}_{\mathrm{un}}(\mathrm{IndCoh}^{!}_{\mathcal{J}^{\mathrm{mer}}(G,N)^{\wedge}_{\mathcal{J}(G,N)}}))$$

over the internal variant of the three dimensional A model $\tilde{\mathcal{A}}(G,N) \in \mathrm{Alg}^{\mathrm{fact}}_{\mathbb{E}_1,\mathrm{un}}(D^{\star}_{\mathcal{Z}(G,N)})$ defined in Theorem 13.4.12, with respect to the action of $D^{\star}_{\mathcal{Z}(G,N)} \in \mathrm{Cat}^{\mathrm{fact}}_{\mathbb{E}_1,\mathrm{un}}(X)$ on $\mathrm{IndCoh}_{\mathcal{J}^{\mathrm{mer}}(G,N)^{\wedge}_{\mathcal{J}(G,N)}} \in \mathrm{Cat}^{\mathrm{fact}}_{\mathrm{un}}(X)$ of Proposition 17.0.13 above.

In particular, the chiral differential operators to $N/G$ of Definition 16.1.24 admits a natural module structure

$$\mathcal{D}^{\mathrm{ch}}(G,N) = \mathrm{p}_{\mathcal{J}^{\mathrm{mer}}(G,N)^{\wedge}_{\mathcal{J}(G,N)}\check{\bullet}} \ \tilde{\mathcal{D}}^{\mathrm{ch}}(G,N) \in \mathcal{A}(G,N)\text{-}\mathrm{Mod}(\mathrm{Alg}^{\mathrm{fact}}_{\mathrm{un}}(X))$$

over the three dimensional A model $\mathcal{A}(G,N) \in \mathrm{Alg}^{\mathrm{fact}}_{\mathbb{E}_1,\mathrm{un}}(X)$ defined in Theorem 13.4.12.

*Proof.* Analogously to the proof of Theorem 13.4.12 following Proposition A.6.2, we follow the proof of Proposition A.6.12, in keeping with Remark A.6.18, using the functors of Example 17.0.11 above. Following Example 17.0.11, we use the description of $\omega^{\mathrm{ren}}_{\mathcal{Z}(G,N)}$ in Remark 13.4.3 and that of $\tilde{\mathcal{D}}^{\mathrm{ch}}(G,N)$ in Remark 17.0.5 to compute

$$\mathrm{IndCoh}^{!}(\mathcal{J}^{\mathrm{mer}}(G,N)^{\wedge}_{(3),x}) = \left( D^{!}(N_{\mathcal{O}} \times G_{\mathcal{K}})^{G_{\mathcal{O}}} \underset{D^{!}(N_{\mathcal{K}})}{\otimes} D^{!}(N_{\mathcal{O}} \times G_{\mathcal{K}})^{G_{\mathcal{O}}} \underset{D^{!}(N_{\mathcal{K}})}{\otimes} \mathrm{IndCoh}^{!}(N_{\mathcal{K}}) \right)^{G_{\mathcal{K}},w}$$

$$\pi^{\bullet}_{12}\omega^{\mathrm{ren}}_{\mathcal{Z}(G,N)_x} = \omega_{N_{\mathcal{O}}}\tilde{\boxtimes} \ \omega_{\mathrm{Gr}_G} \boxtimes_{N_{\mathcal{K}}} \mathbb{K}_{N_{\mathcal{O}}} \tilde{\boxtimes} \mathbb{K}_{\mathrm{Gr}_G} \boxtimes_{N_{\mathcal{K}}} \mathcal{O}_{N_{\mathcal{K}}}$$

$$\pi^{\bullet}_{23}\tilde{\mathcal{D}}^{\mathrm{ch}}(G,N)_x = \mathbb{K}_{N_{\mathcal{O}}}\tilde{\boxtimes} \ \mathbb{K}_{\mathrm{Gr}_G} \boxtimes_{N_{\mathcal{K}}} \omega_{N_{\mathcal{O}}}\tilde{\boxtimes} \omega_{\mathrm{Gr}_G} \boxtimes_{N_{\mathcal{K}}} \omega_{N_{\mathcal{K}}}$$

$$\pi^{!}_{13}\tilde{\mathcal{D}}^{\mathrm{ch}}(G,N)_x = \omega_{N_{\mathcal{O}}}\tilde{\boxtimes} \ \omega_{\mathrm{Gr}_G} \boxtimes_{N_{\mathcal{K}}} \omega_{N_{\mathcal{O}}}\tilde{\boxtimes} \omega_{\mathrm{Gr}_G} \boxtimes_{N_{\mathcal{K}}} \omega_{N_{\mathcal{K}}}$$

As in the proof of Theorem 13.4.12, the tensor unit data provides a canonical map

$$\pi^{\bullet}_{12}\omega^{\mathrm{ren}}_{\mathcal{Z}(G,N)_x}\otimes^{\bullet}\pi^{\bullet}_{23}\tilde{\mathcal{D}}^{\mathrm{ch}}(G,N)_x \xrightarrow{\cong} \pi^{!}_{13}\tilde{\mathcal{D}}^{\mathrm{ch}}(G,N)_x \qquad \text{and} \qquad \pi_{13\bullet}\left(\pi^{\bullet}_{12}\omega^{\mathrm{ren}}_{\mathcal{Z}(G,N)_x} \otimes^{\bullet} \pi^{\bullet}_{23}\tilde{\mathcal{D}}^{\mathrm{ch}}(G,N)_x\right) \to \tilde{\mathcal{D}}^{\mathrm{ch}}(G,N)_x \ ,$$

as desired. $\qquad\square$

## 18. The three dimensional holomorphic-B model and its deformations

In this section, we outline a construction of the factorization $\mathbb{E}_1$ algebra $\mathcal{C}(Y) \in \mathrm{Alg}_{\mathbb{E}_1,\mathrm{un}}^{\mathrm{fact}}(X)$ corresponding to the three dimensional holmorphic-B model to $Y$, which arises as the mixed holomorphic-B type twist of a three dimensional $\mathcal{N}=4$ supersymmetric quantum field theory with target space $T^{\vee}Y$. The factorization algebra is defined in keeping with the general format outlined in Section 9.1, by

$$\mathcal{C}(Y)_x := \mathcal{H}\mathrm{om}_{\mathrm{IndCoh}(Y_{\mathcal{K},\mathrm{Dol}})}(\iota_* \omega_{Y_{\mathcal{O}}}, \iota_* \omega_{Y_{\mathcal{O}}}) \cong \Gamma(\mathcal{Z}(Y)_{\mathrm{Dol}}; \omega_{\mathcal{Z}(Y)_{\mathrm{Dol}}/\mathcal{J}(Y)_{\mathrm{Dol}}}),$$

where the line operator category $\mathrm{IndCoh}(Y_{\mathcal{K},\mathrm{Dol}}) \in \mathrm{Cat}_{\mathrm{un}}^{\mathrm{fact}}(X)$ is given by the category of (ind)coherent sheaves on the Dolbeault stack of the meromorphic jet scheme to $Y$; the construction is summarized as in Equation 8.1.4 by the following diagrams in factorization spaces and categories, which for simplicity we denote by their fibre over a fixed point $x \in X$:
(18.0.1)

We also explain that the three dimensional holomorphic-B model to $Y$ admits a deformation to a family of factorization $\mathbb{E}_1$ algebras $\mathcal{C}(Y)^h \in \mathrm{Alg}_{\mathbb{E}_1,\mathrm{un}}^{\mathrm{fact}}(X)_{/(\mathbb{A}^1/\mathbb{G}_m)}$ with generic fibre equivalent to the three dimensional A model to $Y$ constructed above. Moreover, we show that the module structure of chiral differential operators over the three dimensional A model lifts to an action of this deformation on the Rees algebra of chiral differential operators with respect to the PBW filtration.

18.0.1. *Summary.* In Section 18.1 we outline the construction of the factorization $\mathbb{E}_1$ algebra corresponding to the three dimensional holomorphic B model, and in Section 18.2 we construct its deformation described above. In Section 18.3, we give the internal variant of this construction, and in Section 18.4 we explain the induced action on the Rees chiral differential operators.

*Warning* 18.0.1. In keeping with Warning 9.1.1, we do not formulate specific hypotheses on the space $Y$ used in this section, so that the results stated throughout are only an outline of the general expectations. In the present work, we have not established as careful an understanding of the relevant sheaf theory calculations in this example as in the previous sections, but we hope to return to this problem in future work.

### 18.1. The three dimensional holomorphic B-model factorization algebra.

The starting point for the construction of the three dimensional holomorphic-B model is the Dolbeault stack of the meromorphic jet space of $Y$, and its associated factorization category. Analogously to Example 12.1.1, we have:

*Example* 18.1.1. The Doulbeault stack of the meromorphic jet space of $Y$ is the (unital) factorization space

$$\mathcal{J}^{\mathrm{mer}}(Y)_{\mathrm{Dol}} \in \mathrm{PreStk}_{\mathrm{un}}^{\mathrm{fact}}(X)$$

defined by considering the Dolbeault stack, as defined in Example 10.3.9, of the (unital) factorization space $\mathcal{J}^{\mathrm{mer}}(Y)$ of Example 4.2.4.

As in Example 12.1.3, we have:

*Example* 18.1.2. The category of indcoherent sheaves on the Dolbeault stack of the loop space of $Y$ is the (unital) factorization category

$$\mathrm{IndCoh}_{\mathcal{J}^{\mathrm{mer}}(Y)_{\mathrm{Dol}}} \in \mathrm{Cat}^{\mathrm{fact}}_{\mathrm{un}}(X)$$

associated to $\mathcal{J}^{\mathrm{mer}}(Y)_{\mathrm{Dol}} \in \mathrm{PreStk}^{\mathrm{fact}}_{\mathrm{un}}(X)$, following Definition 5.2.7.

Following Example 5.4.3, analogously to Example 12.1.6, we have:

*Example* 18.1.3. The factorization unit object defines a factorization algebra

$$\mathrm{unit}_{\mathcal{J}^{\mathrm{mer}}(Y)_{\mathrm{Dol}}} = \iota_{\mathcal{J}(Y)*}\mathrm{p}^!_{\mathcal{J}(Y)}\omega_{\mathrm{Ran}_{X,\mathrm{un}}} = \iota_{\mathcal{J}(Y)*}\omega_{\mathcal{J}(Y)_{\mathrm{Dol}}} \ \in \mathrm{Alg}^{\mathrm{fact}}_{\mathrm{un}}(\mathrm{IndCoh}_{\mathcal{J}^{\mathrm{mer}}(Y)_{\mathrm{Dol}}})$$

internal to $\mathrm{IndCoh}_{\mathcal{J}^{\mathrm{mer}}(Y)}$, where

$$\mathrm{p}_{\mathcal{J}(Y)} : \mathcal{J}(Y) \to \mathrm{Ran}_{X_{\mathrm{dR}},\mathrm{un}} \qquad \text{and} \qquad \iota_{\mathcal{J}(Y)} : \mathcal{J}(Y) \to \mathcal{J}^{\mathrm{mer}}(Y)$$

are the factorization space structure map for $\mathcal{J}(Y) = \mathrm{unit}_{\mathcal{J}^{\mathrm{mer}}(Y)}$ and the map of unital factorization spaces given by the inclusion of arcs into loops, respectfully, and

$$\mathrm{p}^!_{\mathcal{J}(Y)} : \mathrm{IndCoh}_{\mathrm{Ran}_{X,\mathrm{un}}} \to \mathrm{IndCoh}_{\mathcal{J}(Y)_{\mathrm{Dol}}} \qquad \text{and} \qquad \iota_{\mathcal{J}(Y)*} : \mathrm{IndCoh}_{\mathcal{J}(Y)_{\mathrm{Dol}}} \to \mathrm{IndCoh}_{\mathcal{J}^{\mathrm{mer}}(Y)_{\mathrm{Dol}}}$$

are the induced unital factorization functors.

Now, suppose that $\mathrm{unit}_{\mathcal{J}^{\mathrm{mer}}(Y)_{\mathrm{Dol}}} \in \mathrm{Alg}^{\mathrm{fact}}_{\mathrm{un}}(\mathrm{IndCoh}_{\mathcal{J}^{\mathrm{mer}}(Y)})$ admits internal Hom objects over $X$, in the sense of Definition 8.1.9. Then, following Example 7.1.7 analogously to Definition 12.1.8, we give the following definition of the factorization $\mathbb{E}_1$ algebra describing the three dimensional holomorphic-B model to $Y$:

*Definition* 18.1.4. The three dimensional Holomorphic-B model to $Y$ is the unital factorization $\mathbb{E}_1$-algebra

$$\mathcal{C}(Y) = \mathcal{H}\mathrm{om}_{\mathrm{IndCoh}(\mathcal{J}^{\mathrm{mer}}(Y))}(\mathrm{unit}_{\mathcal{J}^{\mathrm{mer}}(Y)_{\mathrm{Dol}}}, \mathrm{unit}_{\mathcal{J}^{\mathrm{mer}}(Y)_{\mathrm{Dol}}}) \ \in \mathrm{Alg}^{\mathrm{fact}}_{\mathbb{E}_1,\mathrm{un}}(X) \ .$$

## 18.2. Deformation to the three dimensional A model.

There is a natural deformation of the three dimensional holomorphic-B model factorization algebra defined in the preceding subsection to the three dimensional A model factorization of Section 12, as we now explain:

*Example* 18.2.1. The categories of families of factorization spaces and factorization categores over $\mathbb{A}^1/\mathbb{G}_m$ are defined following the general construction of categories of factorization objects in Subsection 2.4 as

$$\mathrm{PreStk}^{\mathrm{fact}}_{\mathrm{un}}(X)_{/(\mathbb{A}^1/\mathbb{G}_m)} = \mathrm{F}^{\mathrm{fact}}_{\mathrm{un}}(X) \qquad \text{and} \qquad \mathrm{Cat}^{\mathrm{fact}}_{\mathrm{un}}(X)_{/(\mathbb{A}^1/\mathbb{G}_m)} = \mathrm{G}^{\mathrm{fact}}_{\mathrm{un}}(X)$$

where the lax symmetric monoidal functors

$$\mathrm{F} := \mathrm{PreStk}_{/(\cdot)\times\mathbb{A}^1/\mathbb{G}_m} : \mathrm{Sch}^{\mathrm{op}}_{\mathrm{aff}} \to \mathrm{Cat}^\times \qquad \text{and} \qquad G := \mathrm{ShvCat}((\cdot)\times\mathbb{A}^1/\mathbb{G}_m) : \mathrm{Sch}^{\mathrm{op}}_{\mathrm{aff}} \to \mathrm{Cat}^\times$$

are the functors which assign to an affine scheme $S \in \mathrm{Sch}_{\mathrm{aff}}$ the category of prestacks over $S\times\mathbb{A}^1/\mathbb{G}_m$ and the category of sheaves of categories over $S\times\mathbb{A}^1/\mathbb{G}_m$, respectively.

*Remark* 18.2.2. Concretely, a family of factorization spaces $\mathcal{Y} \in \mathrm{PreStk}^{\mathrm{fact}}_{\mathrm{un}}(X)_{/(\mathbb{A}^1/\mathbb{G}_m)}$ or factorization categories $\mathrm{C} \in \mathrm{Cat}^{\mathrm{fact}}_{\mathrm{un}}(X)_{/(\mathbb{A}^1/\mathbb{G}_m)}$ over $\mathbb{A}^1/\mathbb{G}_m$ assigns to each $I \in \mathrm{fSet}_{\emptyset}$ a prestack or sheaf of categories

$$\mathcal{Y}_I \in \mathrm{PreStk}_{/(X^I\times\mathbb{A}^1/\mathbb{G}_m)} \qquad \text{or} \qquad \mathrm{C}_I \in \mathrm{ShvCat}(X^I\times\mathbb{A}^1/\mathbb{G}_m) \ ,$$

respectively, together with the usual factorization and unit structure data defined compatibly over $\mathbb{A}^1/\mathbb{G}_m$.

Now, following Example 10.3.9, we have:

*Example* 18.2.3. There is a canonical deformation of the factorization space $\mathcal{J}^{\mathrm{mer}}(Y)_{\mathrm{Dol}} \in \mathrm{PreStk}^{\mathrm{fact}}_{\mathrm{un}}(X)$ of Example 18.1.1 above, given by

$$\mathcal{J}^{\mathrm{mer}}(Y)_{\mathrm{Hdg}} \in \mathrm{PreStk}^{\mathrm{fact}}_{\mathrm{un}}(X)_{/(\mathbb{A}^1/\mathbb{G}_m)}$$

such that the specializations

$$\mathcal{J}^{\mathrm{mer}}(Y)_{\mathrm{Hdg}} \times_{\mathbb{A}^1/\mathbb{G}_m} \{1\} \cong \mathcal{J}^{\mathrm{mer}}(Y)_{\mathrm{dR}} \quad \text{and} \quad \mathcal{J}^{\mathrm{mer}}(Y)_{\mathrm{Hdg}} \times_{\mathbb{A}^1/\mathbb{G}_m} \{0\} = \mathcal{J}^{\mathrm{mer}}(Y)_{\mathrm{Dol}}$$

are equivalent to the factorization spaces of examples 12.1.1 and 18.1.1, respectively.

Correspondingly, following Example 10.3.10 we obtain a deformation of the factorization category of Example 18.1.2 above:

*Example* 18.2.4. The category of indcoherent sheaves on the Hodge stack of the meromorphic jet space of $Y$ is the family of (unital) factorization categories over $\mathbb{A}^1/\mathbb{G}_m$ given by

$$D^{\hbar}_{\mathcal{J}^{\mathrm{mer}}(Y)} := \mathrm{IndCoh}_{\mathcal{J}^{\mathrm{mer}}(Y)_{\mathrm{Hdg}}} \in \mathrm{Cat}^{\mathrm{fact}}_{\mathrm{un}}(X)_{/(\mathbb{A}^1/\mathbb{G}_m)}$$

such that the specializations

$$D^{\hbar}_{\mathcal{J}^{\mathrm{mer}}(Y)}|_{\{1\}} = \mathrm{IndCoh}_{\mathcal{J}^{\mathrm{mer}}(Y)_{\mathrm{dR}}} = D_{\mathcal{J}^{\mathrm{mer}}(Y)} \quad \text{and} \quad D^{\hbar}_{\mathcal{J}^{\mathrm{mer}}(Y)}|_{\{0\}} = \mathrm{IndCoh}_{\mathcal{J}^{\mathrm{mer}}(Y)_{\mathrm{Dol}}},$$

agree with the factorization categories of Examples 12.1.3 and 18.1.2, respectively.

Further, we have the analogous deformation of the factorization unit object of Example 18.1.3 above, with generic fibre given by the factorization unit object of Example 12.1.6:

*Example* 18.2.5. The family of factorization unit objects defines a family of factorization algebras

$$\mathrm{unit}_{\mathcal{J}^{\mathrm{mer}}(Y)_{\mathrm{Hdg}}} = \iota_{\mathcal{J}(Y)*}\mathrm{p}^{!}_{\mathcal{J}(Y)}\omega_{\mathrm{Ran}_{X,\mathrm{un}}} = \iota_{\mathcal{J}(Y)*}\omega_{\mathcal{J}(Y)_{\mathrm{Hdg}}} \in \mathrm{Alg}^{\mathrm{fact}}_{\mathrm{un}}(D^{\hbar}_{\mathcal{J}^{\mathrm{mer}}(Y)})_{\mathbb{A}^1/\mathbb{G}_m}$$

over $\mathbb{A}^1/\mathbb{G}_m$ internal to $D^{\hbar}_{\mathcal{J}^{\mathrm{mer}}(Y)} \in \mathrm{Cat}^{\mathrm{fact}}_{\mathrm{un}}(X)_{/(\mathbb{A}^1/\mathbb{G}_m)}$, where

$$\mathrm{p}^{!}_{\mathcal{J}(Y)} : \mathrm{IndCoh}_{\mathrm{Ran}_{X,\mathrm{un}}} \to D^{\hbar}_{\mathcal{J}(Y)} \quad \text{and} \quad \iota_{\mathcal{J}(Y)*} : D^{\hbar}_{\mathcal{J}(Y)} \to D^{\hbar}_{\mathcal{J}^{\mathrm{mer}}(Y)}$$

are the families of unital factorization functors induced as in Example 18.1.3 above.

Finally, we obtain the desired deformation of the three dimensional holomorphic-B model factorization algebra of Definition 18.1.4 to that of the three dimensional A model as in Definition 12.1.8:

*Definition* 18.2.6. The deformation of the three dimensional holomorphic-B model to the A model is the family of unital factorization $\mathbb{E}_1$ algebras

$$\mathcal{C}(Y)^{\hbar} = \mathcal{H}om_{D^{\hbar}(\mathcal{J}^{\mathrm{mer}}(Y))}(\mathrm{unit}_{\mathcal{J}^{\mathrm{mer}}(Y)_{\mathrm{Hdg}}}, \mathrm{unit}_{\mathcal{J}^{\mathrm{mer}}(Y)_{\mathrm{Hdg}}}) \in \mathrm{Alg}^{\mathrm{fact}}_{\mathbb{E}_1,\mathrm{un}}(X)_{/(\mathbb{A}^1/\mathbb{G}_m)}.$$

*Remark* 18.2.7. The specializations of $\mathcal{C}(Y)^{\hbar} \in \mathrm{Alg}^{\mathrm{fact}}_{\mathbb{E}_1,\mathrm{un}}(X)_{/(\mathbb{A}^1/\mathbb{G}_m)}$ to the generic and central fibres

$$\mathcal{C}(Y)^{\hbar}|_{\{1\}} = \mathcal{A}(Y) \quad \text{and} \quad \mathcal{C}(Y)^{\hbar}|_{\{0\}} = \mathcal{C}(Y)$$

are equivalent to the factorization $\mathbb{E}_1$ algebras of Definitions 12.1.8 and 18.1.4, respectively, as desired.

### 18.3. Internal construction of the three dimensional holomorphic-B model and its deformation.

There is again an internal variant of the above construction of the three dimensional holomorphic-B model and its deformation to the three dimensional A model, following Example 8.1.15 via Remark 8.1.16.

Following Proposition 8.1.1 as in Example 12.2.1, we have:

*Example* 18.3.1. The self fibre product

$$\mathcal{Z}(Y)_{\mathrm{Hdg}} = \mathcal{J}(Y)_{\mathrm{Hdg}} \times_{\mathcal{J}^{\mathrm{mer}}(Y)_{\mathrm{Hdg}}} \mathcal{J}(Y)_{\mathrm{Hdg}} \in \mathrm{PreStk}^{\mathrm{fact}}_{\mathrm{un}}(X)_{/(\mathbb{A}^1/\mathbb{G}_m)}$$

defines a unital factorization space over $X$.

*Remark* 18.3.2. More generally, the iterated fibre products

$$\mathcal{Z}(Y)_{(n),\mathrm{Hdg}} := \mathcal{J}(Y)_{\mathrm{Hdg}} \times_{\mathcal{J}^{\mathrm{mer}}(Y)_{\mathrm{Hdg}}} \times \ldots \times_{\mathcal{J}^{\mathrm{mer}}(Y)_{\mathrm{Hdg}}} \mathcal{J}(Y)_{\mathrm{Hdg}} \in \mathrm{PreStk}^{\mathrm{fact}}_{\mathrm{un}}(X)_{/(\mathbb{A}^1/\mathbb{G}_m)}$$

and projections $\pi_{ij} : \mathcal{Z}(Y)_{(n),\mathrm{Hdg}} \to \mathcal{Z}(Y)_{\mathrm{Hdg}}$, define factorization spaces over $\mathbb{A}^1/\mathbb{G}_m$ and maps of such for each $n \in \mathbb{N}$ and $i, j \in \{1, ..., n\}$.

Following Proposition 8.1.5 as in Proposition 12.2.4, we have:

*Proposition* 18.3.3. The family of factorization categories

$$D^{\hbar,\star}_{\mathcal{Z}(Y)} = \mathrm{IndCoh}^{\star}_{\mathcal{Z}(Y)_{\mathrm{Hdg}}} \in \mathrm{Cat}^{\mathrm{fact}}_{\mathbb{E}_1,\mathrm{un}}(X_{\mathrm{dR}})_{/(\mathbb{A}^1/\mathbb{G}_m)}$$

over $\mathbb{A}^1/\mathbb{G}_m$ is naturally a family of (unital) $\mathbb{E}_1$-factorization categories with respect to the convolution monoidal structure defined by the composition

$$D^{\hbar}_{\mathcal{Z}(Y)} \otimes^* D^{\hbar}_{\mathcal{Z}(Y)} \xrightarrow{\pi^*_{12} \boxtimes \pi^*_{23}} D^{\hbar}_{\mathcal{Z}_{(3)}(Y)^{\times 2}} \xrightarrow{\Delta^*} D^{\hbar}_{\mathcal{Z}_{(3)}(Y)} \xrightarrow{\pi_{13,*}} D^{\hbar}_{\mathcal{Z}(Y)} \ .$$

Further, the pushforward functor $\mathrm{p}_{\mathcal{Z}(Y)*} : D^{\hbar,\star}_{\mathcal{Z}(Y)} \to D^{\otimes^!}_{\mathrm{Ran}_{X,\mathrm{un}}} \otimes \mathrm{D}^b_{\mathrm{fg}}(\mathbb{K}[\hbar])$ defines a family of unital, lax $\mathbb{E}_1$-monoidal factorization functors over $\mathbb{A}^1/\mathbb{G}_m$.

Further, following Corollary 8.1.6 as in Corollary 13.4.10, we have:

*Corollary* 18.3.4. The pushforward functor $\mathrm{p}_{\mathcal{Z}*} : \mathrm{IndCoh}^{\star}_{\mathcal{Z}(Y)_{\mathrm{Hdg}}} \to D_{\mathrm{Ran}_{X,\mathrm{un}}} \otimes \mathrm{D}^b_{\mathrm{fg}}(\mathbb{K}[\hbar])$ induces a functor

$$(18.3.1) \qquad \mathrm{p}_* : \mathrm{Alg}^{\mathrm{fact}}_{\mathbb{E}_1,\mathrm{un}}(D^{\hbar,\star}_{\mathcal{Z}(Y)})_{/(\mathbb{A}^1/\mathbb{G}_m)} \to \mathrm{Alg}^{\mathrm{fact}}_{\mathbb{E}_1,\mathrm{un}}(X)_{/(\mathbb{A}^1/\mathbb{G}_m)} \ .$$

Finally, again following Example 8.1.15, we define the internal variant of the family of factorization algebras over $\mathbb{A}^1/\mathbb{G}_m$ in Definition 18.2.6, which give a deformation from the internal variant of the three dimensional holomorphic-B model to internal variant of three dimensional A model of Definition 12.2.7:

*Definition* 18.3.5. The internal variant of the deformation from the three dimensional holomorphic-B model to the three dimensional A model is the family of factorization $\mathbb{E}_1$ algebras

$$\tilde{\mathfrak{e}}(Y)^h = \mathcal{H}\widetilde{\mathrm{om}}_{D^{\hbar}_{\mathcal{J}^{\mathrm{mer}}(Y)}}(\mathrm{unit}_{\mathcal{J}^{\mathrm{mer}}(Y)_{\mathrm{Hdg}}}, \mathrm{unit}_{\mathcal{J}^{\mathrm{mer}}(Y)_{\mathrm{Hdg}}}) \in \mathrm{Alg}^{\mathrm{fact}}_{\mathbb{E}_1,\mathrm{un}}(D^{\hbar,\star}_{\mathcal{Z}(Y)})_{/(\mathbb{A}^1/\mathbb{G}_m)}$$

internal to $D^{\hbar,\star}_{\mathcal{Z}(Y)} \in \mathrm{Cat}^{\mathrm{fact}}_{\mathbb{E}_1,\mathrm{un}}(X_{\mathrm{dR}})_{/(\mathbb{A}^1/\mathbb{G}_m)}$.

There is also an internal variant of Remark 18.2.7:

*Remark* 18.3.6. The specializations of the family of factorization spaces $\mathcal{Z}(Y)_{\mathrm{Hdg}} \in \mathrm{PreStk}^{\mathrm{fact}}_{\mathrm{un}}(X)_{/(\mathbb{A}^1/\mathbb{G}_m)}$ to the generic and central fibres are given by

$$\mathcal{Z}(Y)_{\mathrm{Hdg}} \times_{\mathbb{A}^1/\mathbb{G}_m} \{1\} \cong \mathcal{Z}(Y)_{\mathrm{dR}} \quad \text{and} \quad \mathcal{Z}(Y)_{\mathrm{Hdg}} \times_{\mathbb{A}^1/\mathbb{G}_m} \{0\} = \mathcal{Z}(Y)_{\mathrm{Dol}}$$

and thus the specializations of the family of factorization $\mathbb{E}_1$ categories $D^{\hbar,\star}_{\mathcal{Z}(Y)} \in \mathrm{Cat}^{\mathrm{fact}}_{\mathbb{E}_1,\mathrm{un}}(X_{\mathrm{dR}})_{/(\mathbb{A}^1/\mathbb{G}_m)}$ are given by

$$D^{\hbar,\star}_{\mathcal{Z}(Y)}|_{\{1\}} = D_{\mathcal{Z}(Y)} \quad \text{and} \quad D^{\hbar,\star}_{\mathcal{Z}(Y)}|_{\{0\}} = \mathrm{IndCoh}^{\star}_{\mathcal{Z}(Y)_{\mathrm{Dol}}} \ .$$

Further, under these identifications, we have

$$\tilde{\mathcal{C}}(Y)^{\hbar}|_{\{1\}} = \tilde{\mathcal{A}}(Y) \ \in D_{\mathcal{Z}(Y)} \quad \text{and} \quad \tilde{\mathcal{C}}(Y)^{\hbar}|_{\{0\}} = \tilde{\mathcal{C}}(Y) \ \in \mathrm{IndCoh}^{\star}_{\mathcal{Z}(Y)_{\mathrm{Dol}}} \ ,$$

justifying the preceding definition.

Following Example 12.2.9, and in turn Example 8.1.15 and Corollary 8.1.6, we have:

*Proposition* 18.3.7. The image of $\tilde{\mathcal{C}}(Y)^{\hbar}$ in $\mathrm{Alg}^{\mathrm{fact}}_{\mathbb{E}_1,\mathrm{un}}(D^{\hbar,\star}_{\mathcal{Z}(Y)})_{/(\mathbb{A}^1/\mathbb{G}_m)}$ under the functor of Equation 18.3.1 above is canonically equivalent to $\mathcal{C}(Y)^{\hbar} \in \mathrm{Alg}^{\mathrm{fact}}_{\mathbb{E}_1,\mathrm{un}}(X)_{/(\mathbb{A}^1/\mathbb{G}_m)}$ of Definition 18.2.6.

18.4. **Action of the deformation on Rees chiral differential operators.** In this subsection, we construct an action of $\mathcal{C}(X)^{\hbar} \in \mathrm{Alg}^{\mathrm{fact}}_{\mathbb{E}_1,\mathrm{un}}(X)_{/(\mathbb{A}^1/\mathbb{G}_m)}$ on the family of factorization algebras $\mathcal{D}^{\mathrm{ch}}(Y)^{\hbar} \in \mathrm{Alg}^{\mathrm{fact}}_{\mathrm{un}}(X)_{/(\mathbb{A}^1/\mathbb{G}_m)}$ underlying the factorization $\mathbb{B}\mathbb{D}^{\hbar}_0$ algebra describing the filtered quantization of chiral differential operators. In particular, the generic fibre is given by the construction outlined in Subsection 14.3. The central fibre is identified with the classical limit of the chiral differential operators factorization algebra, which thus admits an action of the holomorphic-B model factorization algebra.

Following Example 14.3.1 and in turn Example 8.2.1, we have:

*Example* 18.4.1. The fibre product of families of factorization spaces

$$\mathcal{J}^{\mathrm{mer}}_{\hbar}(Y)^{\wedge}_{\mathcal{J}(Y)} := \mathcal{J}(Y)_{\mathrm{Hdg}} \times_{\mathcal{J}^{\mathrm{mer}}(Y)_{\mathrm{Hdg}}} \mathcal{J}^{\mathrm{mer}}(Y) \ \in \mathrm{PreStk}^{\mathrm{fact}}_{\mathrm{un}}(X)_{/(\mathbb{A}^1/\mathbb{G}_m)}$$

defines a family of unital factorization spaces on $X$ over $\mathbb{A}^1/\mathbb{G}_m$, where

$$\mathcal{J}^{\mathrm{mer}}_{\hbar}(Y) = \mathcal{J}^{\mathrm{mer}}(Y) \times \mathbb{A}^1_{\hbar}/\mathbb{G}_m \in \mathrm{PreStk}^{\mathrm{fact}}_{\mathrm{un}}(X)_{/(\mathbb{A}^1/\mathbb{G}_m)} \ .$$

*Remark* 18.4.2. The specializations of $\mathcal{J}^{\mathrm{mer}}_{\hbar}(Y)^{\wedge}_{\mathcal{J}(Y)} \in \mathrm{PreStk}^{\mathrm{fact}}_{\mathrm{un}}(X)_{/(\mathbb{A}^1/\mathbb{G}_m)}$ to the generic and central fibres

$$\mathcal{J}^{\mathrm{mer}}_{\hbar}(Y)^{\wedge}_{\mathcal{J}(Y)} \times_{\mathbb{A}^1/\mathbb{G}_m} \{1\} \cong \mathcal{J}^{\mathrm{mer}}(Y)^{\wedge}_{\mathcal{J}(Y)} \quad \text{and} \quad \mathcal{J}^{\mathrm{mer}}_{\hbar}(Y)^{\wedge}_{\mathcal{J}(Y)} \times_{\mathbb{A}^1/\mathbb{G}_m} \{0\} = \mathcal{J}(Y)_{\mathrm{Dol}} \times_{\mathcal{J}^{\mathrm{mer}}(Y)_{\mathrm{Dol}}} \mathcal{J}^{\mathrm{mer}}(Y)$$

are equivalent to the factorization space of Example 14.3.1 and its Dolbeault analogue, respectively.

Further, following Example 14.3.3, we have:

*Remark* 18.4.3. More generally, the iterated fibre products

$$\mathcal{J}^{\mathrm{mer}}_{\hbar}(Y)^{\wedge}_{(n)} := \mathcal{J}(Y)_{\mathrm{Hdg}} \times_{\mathcal{J}^{\mathrm{mer}}(Y)_{\mathrm{Hdg}}} \mathcal{J}(Y)_{\mathrm{Hdg}} \times_{\mathcal{J}^{\mathrm{mer}}(Y)_{\mathrm{Hdg}}} \ldots \times_{\mathcal{J}^{\mathrm{mer}}(Y)_{\mathrm{Hdg}}} \mathcal{J}^{\mathrm{mer}}_{\hbar}(Y) \in \mathrm{PreStk}^{\mathrm{fact}}_{\mathrm{un}}(X)_{/(\mathbb{A}^1/\mathbb{G}_m)}$$

and projections

$$\pi_{ij} : \mathcal{J}^{\mathrm{mer}}_{\hbar}(Y)^{\wedge}_{(n)} \to \mathcal{J}^{\mathrm{mer}}_{\hbar}(Y)^{\wedge}_{\mathcal{J}(Y)} \qquad \text{for } i \in \{1, ..., n-1\} \text{ and } j = n, \text{ and}$$
$$\pi_{ij} : \mathcal{J}^{\mathrm{mer}}_{\hbar}(Y)^{\wedge}_{(n)} \to \mathcal{Z}(Y)_{\mathrm{Hdg}} \qquad \text{for } i, j \in \{1, ..., n-1\},$$

define families of factorization spaces over $\mathbb{A}^1/\mathbb{G}_m$ and maps of such for each $n \in \mathbb{N}$.

Following Proposition 14.3.4 and in turn Proposition 8.2.4, we have:

*Proposition* 18.4.4. The factorization category

$$\mathrm{IndCoh}_{\mathcal{J}^{\mathrm{mer}}_\hbar(Y)^\wedge_{\mathcal{J}(Y)}} \in D^{h,\bullet}_{\mathcal{Z}(Y)}\text{-Mod}(\mathrm{Cat}^{\mathrm{fact}}_{\mathrm{un}}(X)_{/(\mathbb{A}^1/\mathbb{G}_m)})$$

is naturally a family over $\mathbb{A}^1/\mathbb{G}_m$ of factorization $\mathbb{E}_1$ module categories over $D^{h,\bullet}_{\mathcal{Z}(Y)} \in \mathrm{Cat}^{\mathrm{fact}}_{\mathbb{E}_1,\mathrm{un}}(X)_{/(\mathbb{A}^1/\mathbb{G}_m)}$, with respect to the convolution module structure $(\cdot)\star(\cdot): D^h_{\mathcal{Z}(Y)}\otimes\mathrm{IndCoh}_{\mathcal{J}^\mathrm{mer}_\hbar(Y)^\wedge_{\mathcal{J}(Y)}} \to \mathrm{IndCoh}_{\mathcal{J}^\mathrm{mer}_\hbar(Y)^\wedge_{\mathcal{J}(Y)}}$ defined by the composition

$$D^\hbar_{\mathcal{Z}(Y)}\otimes^*\mathrm{IndCoh}_{\mathcal{J}^\mathrm{mer}_\hbar(Y)^\wedge_{\mathcal{J}(Y)}} \xrightarrow{\pi^\bullet_{12}\boxtimes\pi^\bullet_{23}} \mathrm{IndCoh}_{(\mathcal{J}^\mathrm{mer}_\hbar(Y)_{(3)})^{\times 2}} \xrightarrow{\Delta^\bullet} \mathrm{IndCoh}_{\mathcal{J}^\mathrm{mer}_\hbar(Y)_{(3)}} \xrightarrow{\pi_{13\bullet}} \mathrm{IndCoh}_{\mathcal{J}^\mathrm{mer}_\hbar(Y)^\wedge_{\mathcal{J}(Y)}} .$$

Further, if $Y$ is equipped with a Tate structure, then the pushforward functors $\mathrm{p}_{\mathcal{Z}(Y)*}: D^{h,\bullet}_{\mathcal{Z}(Y)} \to D^\otimes_{\mathrm{Ran}_{X,\mathrm{un}}}\otimes D^b_{\mathrm{fg}}(\mathbb{K}[\hbar])$ as in Proposition 18.3.3 and $\mathrm{p}_{\mathcal{J}^\mathrm{mer}_\hbar(Y)^\wedge_{\mathcal{J}(Y)}*}: \mathrm{IndCoh}_{\mathcal{J}^\mathrm{mer}_\hbar(Y)^\wedge_{\mathcal{J}(Y)}} \to D_{\mathrm{Ran}_{X,\mathrm{un}}}\otimes D^b_{\mathrm{fg}}(\mathbb{K}[\hbar])$ define families over $\mathbb{A}^1/\mathbb{G}_m$ of unital, lax compatible factorization functors with respect to the above module structure.

Thus, following Corollary 14.3.5 and in turn Corollary 8.2.5, we have:

*Corollary* 18.4.5. The pushforward functor $\mathrm{p}_{\mathcal{J}^\mathrm{mer}_\hbar(Y)^\wedge_{\mathcal{J}(Y)}\bullet}: \mathrm{IndCoh}_{\mathcal{J}^\mathrm{mer}_\hbar(Y)^\wedge_{\mathcal{J}(Y)}} \to D_{\mathrm{Ran}_{X,\mathrm{un}}}\otimes D^b_{\mathrm{fg}}(\mathbb{K}[\hbar])$ induces a functor

$$(18.4.1) \qquad \tilde{\mathcal{C}}^\hbar(Y)\text{-Mod}(\mathrm{Alg}^\mathrm{fact}_\mathrm{un}(\mathrm{IndCoh}_{\mathcal{J}^\mathrm{mer}_\hbar(Y)^\wedge_{\mathcal{J}(Y)}})_{/(\mathbb{A}^1/\mathbb{G}_m)}) \to \mathcal{C}(Y)^\hbar\text{-Mod}(\mathrm{Alg}^\mathrm{fact}_\mathrm{un}(X)_{/(\mathbb{A}^1/\mathbb{G}_m)}) .$$

In analogy with Definition 14.3.7, we make the following definition:

*Definition* 18.4.6. The filtered quantization $\tilde{\mathcal{D}}^\mathrm{ch}(Y)^\hbar \in \mathrm{Alg}^\mathrm{fact}_\mathrm{un}(\mathrm{IndCoh}_{\mathcal{J}^\mathrm{mer}_\hbar(Y)^\wedge_{\mathcal{J}(Y)}})_{/(\mathbb{A}^1/\mathbb{G}_m)}$ of the internal variant of chiral differential operators is defined by

$$\tilde{\mathcal{D}}^\mathrm{ch}(Y)^\hbar = \pi^!_{\mathcal{J}(Y)_\mathrm{Hdg}}\omega_{\mathcal{J}(Y)_\mathrm{Hdg}} = \mathrm{p}^!_{\mathcal{J}^\mathrm{mer}_\hbar(Y)^\wedge_{\mathcal{J}(Y)}}\omega_{\mathrm{Ran}_{X,\mathrm{un}}} .$$

The preceding definition is justified by the following:

*Proposition* 18.4.7. If $Y$ admits a Tate structure, than the image of $\tilde{\mathcal{D}}^\mathrm{ch}(Y)^\hbar \in \mathrm{Alg}^\mathrm{fact}_\mathrm{un}(\mathrm{IndCoh}_{\mathcal{J}^\mathrm{mer}_\hbar(Y)^\wedge_{\mathcal{J}(Y)}})_{/(\mathbb{A}^1/\mathbb{G}_m)}$ under the pushforward factorization functor

$$\mathrm{p}_{\mathcal{J}^\mathrm{mer}_\hbar(Y)^\wedge_{\mathcal{J}(Y)}\bullet}\tilde{\mathcal{D}}^\mathrm{ch}(Y)^\hbar \cong \mathcal{D}^\mathrm{ch}(Y)^\hbar \in \mathrm{Alg}^\mathrm{fact}_{\mathbb{B}\mathbb{D}_0,\mathrm{un}}(X)$$

is equivalent to the filtered quantization of chiral differential operators.

Further, following Theorem 17.0.15, we have:

*Theorem* 18.4.8. The filtered quantization of the internal variant of chiral differential operators admits a natural module structure

$$\tilde{\mathcal{D}}^\mathrm{ch}(Y)^\hbar \in \tilde{\mathcal{C}}(Y)^\hbar\text{-Mod}(\mathrm{Alg}^\mathrm{fact}_\mathrm{un}(\mathrm{IndCoh}_{\mathcal{J}^\mathrm{mer}_\hbar(Y)^\wedge_{\mathcal{J}(Y)}})_{/(\mathbb{A}^1/\mathbb{G}_m)})$$

over the internal variant of the deformation from the three dimensional holomorphic-B model to A model $\tilde{\mathcal{C}}(Y)^\hbar \in \mathrm{Alg}^\mathrm{fact}_{\mathbb{E}_1,\mathrm{un}}(D^{h,\bullet}_{\mathcal{Z}(Y)})_{/(\mathbb{A}^1/\mathbb{G}_m)}$ of Definition 18.3.5.

In particular, the filtered quantization of chiral differential operators admits a natural module structure

$$\mathcal{D}^\mathrm{ch}(Y)^\hbar \in \mathcal{C}(Y)^\hbar\text{-Mod}(\mathrm{Alg}^\mathrm{fact}_\mathrm{un}(X)_{/(\mathbb{A}^1/\mathbb{G}_m)})$$

over the deformation from the three dimensional holomorphic-B model to A model $\mathcal{C}(Y)^\hbar \in \mathrm{Alg}^\mathrm{fact}_{\mathbb{E}_1,\mathrm{un}}(X)_{/(\mathbb{A}^1/\mathbb{G}_m)}$ of Definition 18.2.6.

## 19. The four dimensional holomorphic-B model and chiral quantization

In this section, we outline a construction of the factorization $\mathbb{E}_2$ algebra $\mathcal{F}(Y) \in \mathrm{Alg}^{\mathrm{fact}}_{\mathbb{E}_2,\mathrm{un}}(X)$ corresponding to the four dimensional holomorphic-B model to $Y$, which arises as the mixed holomorphic-B type twist of a four dimensional $\mathcal{N} = 2$ supersymmetric quantum field theory with target space $T^\vee Y$, introduced in [Kap06]. Moreover, we outline an identification between Tate structures on $Y$ and $S^1$ equivariant structures on $\mathcal{F}(Y)$, and correspondingly an identification of the (two-periodic) filtered quantization of a factorization algebra on $X$, corresponding to $\mathcal{F}(Y) \in \mathrm{Alg}^{\mathrm{fact}}_{\mathbb{E}_2^{S^1},\mathrm{un}}(X)$ under the equivalence of Proposition I-25.3.1, with the (two-periodic) Rees chiral differential operators $\mathcal{D}^{\mathrm{ch}}(Y)_u$ to $Y$. This gives a mathematical variant of the main construction of [BLL+15], which constructs chiral algebras from the observables of four dimensional $\mathcal{N} = 2$ theories. Finally, we outline an identification of the negative cyclic chains $\mathrm{CC}^-_\bullet(\mathcal{F}(Y)) \in \mathrm{Alg}^{\mathrm{fact}}_\bullet(X)_{\mathbb{K}[u]}$ of $\mathcal{F}(Y)$ with a two-periodic variant of the deformation $\mathcal{C}(Y)^u$ of the three dimensional holomorphic-B model to the A model, such that the module structure induced by the equivariant cigar reduction principle of Section I-25.4 identifies with the two-periodic variant of that constructed in the preceding section. This explains the relationship between the predictions of [BLL+15] and those of [CG18] explained in the preceding sections.

The factorization algebra $\mathcal{F}(Y) \in \mathrm{Alg}^{\mathrm{fact}}_{\mathbb{E}_2,\mathrm{un}}(X)$ is defined in keeping with the general format outlined in section 9.1, by

$$\mathcal{F}(Y)_x = \mathrm{Hom}_{\mathrm{IndCoh}(\mathcal{Z}(Y)_x)}(\mathrm{u}_{\mathcal{Z}(Y)_x\bullet}\mathcal{O}_{\mathcal{J}(Y)_x}, \mathrm{u}_{\mathcal{Z}(Y)_x\bullet}\mathcal{O}_{\mathcal{J}(Y)_x}) ,$$

where the line operator category $\mathrm{IndCoh}(\mathcal{Z}(Y)) \in \mathrm{Cat}^{\mathrm{fact}}_{\mathbb{E}_1,\mathrm{un}}(X)$ is given by the category of (ind)coherent sheaves on the space $\mathcal{Z}(Y) = Y_\mathcal{K} \times_{Y_\mathcal{O}} Y_\mathcal{K}$; for example, in the case $Y = \mathrm{B}G$, this is precisely the derived coherent Satake category studied in [CW19]. The construction of the factorization algebra is summarized as in Equation 8.1.4 by the following diagrams in factorization spaces and categories, which for simplicity we denote by their fibre over a fixed point $x \in X$:

(19.0.1)

19.0.1. *Summary.* In Section 19.1 we outline the construction of the factorization algebra describing the four dimensional holomorphic-B model, in section 19.2 we describe the internal variant of the construction, and in section 19.3 we discuss the results related to the equivariant cigar reduction principle of Section I-25.4 in this example.

*Warning* 19.0.1. In keeping with Warning 9.1.1, we do not formulate specific hypotheses on the space $Y$ used in this section, so that the results stated throughout are only an outline of the general expectations. In the present work, we have not established as careful an understanding of the relevant sheaf theory calculations in this example as in the previous sections, but we hope to return to this problem in future work; we do believe that the theory of indcoherent sheaves on renormalizable prestacks developed in Section 6 of [Ras20b] suffices for the case $Y = N/G$, following the constructions of sections 16 and 17.

**19.1. The four dimensional holomorphic-B model factorization algebra.** To begin, we have the following analogue of Example 12.2.1 in the coherent (as opposed to de Rham) setting:

*Example* 19.1.1. The self fibre product

$$\mathcal{Z}(Y) = \mathcal{J}(Y) \times_{\mathcal{J}^{\mathrm{mer}}(Y)} \mathcal{J}(Y) \ \in \mathrm{PreStk}^{\mathrm{fact}}_{\mathrm{un}}(X)$$

defines a unital factorization space over $X$.

*Remark* 19.1.2. Concretely, the prestacks $\mathcal{Z}(Y)_I \in \mathrm{PreStk}_{/X^I}$ over $X^I$ assigned to each $I \in \mathrm{fSet}$ and the fibre $\mathcal{Z}(Y)_x \in \mathrm{PreStk}$ over $x \in X$ are given by

$$\mathcal{Z}(Y)_I = \mathcal{J}(Y)_I \times_{\mathcal{J}^{\mathrm{mer}}(Y)_I} \mathcal{J}(Y)_I \qquad \text{and} \qquad \mathcal{Z}(Y)_x = Y_{\mathcal{O}_x} \times_{Y_{\mathcal{K}_x}} Y_{\mathcal{O}_x} \ .$$

Similarly, as in Remark 18.3.2, we have:

*Remark* 19.1.3. More generally, the iterated fibre products

$$\mathcal{Z}(Y)_{(n)} := \mathcal{J}(Y) \times_{\mathcal{J}^{\mathrm{mer}}(Y)} \times \ldots \times_{\mathcal{J}^{\mathrm{mer}}(Y)} \mathcal{J}(Y) \in \mathrm{PreStk}^{\mathrm{fact}}_{\mathrm{un}}(X)$$

and projections $\pi_{ij} : \mathcal{Z}(Y)_{(n)} \to \mathcal{Z}(Y)$, define factorization spaces and maps of such for each $n \in \mathbb{N}$ and $i, j \in \{1, ..., n\}$.

*Remark* 19.1.4. Following Remark 4.2.5, $\mathcal{Z}(Y) \in \mathrm{PreStk}^{\mathrm{fact}}_{\mathrm{un}}(X_{\mathrm{dR}})$ is canonically defined over $X_{\mathrm{dR}}$.

Following Proposition 8.1.5, analogously to Proposition 12.2.4 but in the coherent setting, we have:

*Proposition* 19.1.5. The factorization category

$$\mathrm{IndCoh}^{\star}_{\mathcal{Z}(Y)} \in \mathrm{Cat}^{\mathrm{fact}}_{\mathbb{E}_1, \mathrm{un}}(X_{\mathrm{dR}})$$

is naturally a (unital) $\mathbb{E}_1$-factorization category with respect to the convolution monoidal structure defined by the composition

$$\mathrm{IndCoh}_{\mathcal{Z}(Y)} \otimes^* \mathrm{IndCoh}_{\mathcal{Z}(Y)} \xrightarrow{\pi^{\star}_{12} \boxtimes \pi^{\star}_{23}} \mathrm{IndCoh}_{\mathcal{Z}(Y)^{\times 2}_{(3)}} \xrightarrow{\Delta^{\bullet}} \mathrm{IndCoh}_{\mathcal{Z}(Y)_{(3)}} \xrightarrow{\pi_{13, \bullet}} \mathrm{IndCoh}_{\mathcal{Z}(Y)} \ .$$

Further, the pushforward functor $\mathrm{p}_{\mathcal{Z}(Y)_{\bullet}} : \mathrm{IndCoh}^{\star}_{\mathcal{Z}(Y)} \to \mathcal{D}^{\otimes^!}_{\mathrm{Ran}_{X, \mathrm{un}}}$ is a unital, lax $\mathbb{E}_1$-monoidal factorization functor.

Analogously to Example 12.2.6, we have:

*Example* 19.1.6. Let $G$ be a reductive algebraic group and $Y = \mathrm{B}G$ be its classifying stack, so that the factorization space above is given by

$$\mathcal{Z}(Y) \cong \mathcal{J}(G) \backslash \mathrm{Gr}_{G, \mathrm{Ran}_{X, \mathrm{un}}} \ .$$

The resulting $\mathbb{E}_1$ factorization category

$$\mathrm{IndCoh}_{\mathcal{Z}(Y)} \cong \mathrm{IndCoh}_{\mathcal{J}(G) \backslash \mathrm{Gr}_{G, \mathrm{Ran}_{X, \mathrm{un}}}} \in \mathrm{Cat}^{\mathrm{fact}}_{\mathbb{E}_1, \mathrm{un}}(X)$$

is the derived (ind)coherent Satake category, closely related to that considered in [CW19] and references therein; concretely, its fibre category over $x \in X$ is given by

$$\mathrm{IndCoh}_{\mathcal{Z}(Y), x} = \mathrm{IndCoh}(\mathrm{Gr}_{G, x})^{G_{\mathcal{O}, w}} \ .$$

*Example* 19.1.7. The unit correspondence for $\mathcal{Z}(Y) \in \mathrm{PreStk}_{\mathrm{un}}^{\mathrm{fact}}(X)$ is given by

$$
\begin{array}{ccc}
& \mathcal{J}(Y) & \\
{}^{\mathrm{p}_{\mathcal{J}(Y)}}\swarrow & & \searrow^{\mathrm{u}_{\mathcal{Z}(Y)}} \\
\mathrm{Ran}_{X,\mathrm{un}} & & \mathcal{Z}(Y)
\end{array}
\qquad \text{or concretely} \qquad
\begin{array}{ccc}
& Y_{\mathcal{O}} & \\
{}^{\mathrm{p}_{\mathcal{J}(Y)}}\swarrow & & \searrow^{\mathrm{u}_{\mathcal{Z}(Y)_x}} \\
\mathrm{pt} & & Y_{\mathcal{O}} \times_{Y_{\mathcal{K}}} Y_{\mathcal{O}}
\end{array}
$$

over each point $x \in X$.

*Example* 19.1.8. The factorization $\mathbb{E}_1$ unit object $\mathrm{unit}_{\mathrm{IndCoh}_{\mathcal{Z}(Y)}} \in \mathrm{Alg}_{\mathbb{E}_1,\mathrm{un}}^{\mathrm{fact}}(\mathrm{IndCoh}_{\mathcal{Z}(Y)}^{\star})$ is given by

$$
\mathrm{unit}_{\mathrm{IndCoh}_{\mathcal{Z}(Y)}} = \mathrm{u}_{\mathcal{Z}(Y)\bullet}\mathrm{p}_{\mathcal{J}(Y)}^{\bullet}\mathcal{O}_{\mathrm{Ran}_{X,\mathrm{un}}} = \mathrm{u}_{\mathcal{Z}(Y)\bullet}\mathcal{O}_{\mathcal{J}(Y)} \in \mathrm{Alg}_{\mathbb{E}_1,\mathrm{un}}^{\mathrm{fact}}(\mathrm{IndCoh}_{\mathcal{Z}(Y)}^{\star}) \ .
$$

*Remark* 19.1.9. Concretely, the $\mathbb{E}_1$ algebra object internal to $\mathrm{IndCoh}(\mathcal{Z}(Y)_I)^{\star}$ assigned to each $I \in \mathrm{fSet}_{\oslash}$, and the $\mathbb{E}_1$ algebra object of the fibre category over each point $x \in X$, are given by

$$
\mathrm{unit}_{\mathrm{IndCoh}_{\mathcal{Z}(Y)},I} = \mathrm{u}_{\mathcal{Z}(Y)_I\bullet}\mathrm{p}_{\mathcal{J}(Y)_I}^{\bullet}\mathcal{O}_{X^I} = \mathrm{u}_{\mathcal{Z}(Y)_I\bullet}\mathcal{O}_{\mathcal{J}(Y)_I} \qquad \in \mathrm{Alg}_{\mathbb{E}_1}(\mathrm{IndCoh}(\mathcal{Z}(Y)_I)^{\star}), \text{ and}
$$
$$
\mathrm{unit}_{\mathrm{IndCoh}_{\mathcal{Z}(Y)},x} = \mathrm{u}_{\mathcal{Z}(Y)_x\bullet}\mathrm{p}_{\mathcal{J}(Y)_x}^{\bullet}\mathcal{O}_{\mathrm{pt}} = \mathrm{u}_{\mathcal{Z}(Y)_x\bullet}\mathcal{O}_{\mathcal{J}(Y)_x} \qquad \in \mathrm{Alg}_{\mathbb{E}_1}(\mathrm{IndCoh}(\mathcal{Z}(Y)_x)^{\star}) \ .
$$

*Definition* 19.1.10. The four dimensional holomorphic B model factorization $\mathbb{E}_2$ algebra is defined by

$$
\mathcal{F}(Y) = \mathcal{H}\mathrm{om}_{\mathrm{IndCoh}(\mathcal{Z}(Y))}(\mathrm{unit}_{\mathrm{IndCoh}_{\mathcal{Z}(Y)}}, \mathrm{unit}_{\mathrm{IndCoh}_{\mathcal{Z}(Y)}}) \ \in \mathrm{Alg}_{\mathbb{E}_2,\mathrm{un}}^{\mathrm{fact}}(X) \ .
$$

*Remark* 19.1.11. There is a natural (factorization) $\mathbb{E}_2$ structure on $\mathcal{F}(Y) \in \mathrm{Alg}_{\mathbb{E}_2,\mathrm{un}}^{\mathrm{fact}}(X)$, defined as in Remark 10.1.6. Concretely, $\mathcal{F}(Y)$ has an $\mathbb{E}_1$ algebra structure given by composition of endomorphisms, and a second compatible (factorization) $\mathbb{E}_1$ algebra structure as a space of maps between (factorization) $\mathbb{E}_1$ algebra objects internal to the category $\mathrm{IndCoh}_{\mathcal{Z}(Y)}^{\star} \in \mathrm{Cat}_{\mathbb{E}_1,\mathrm{un}}^{\mathrm{fact}}(X_{\mathrm{dR}})$.

*Remark* 19.1.12. Concretely, following Remark 8.1.10, the $\mathbb{E}_2$ algebra internal to $\mathcal{D}(X^I)$ assigned to each $I \in \mathrm{fSet}_{\oslash}$, and the $\mathbb{E}_2$ algebra in Vect over each $x \in X$, are given by

$$
\mathcal{F}(Y)_I = \mathrm{p}_{\mathcal{Z}(Y)_I\bullet}\underline{\mathrm{Hom}}_{\mathrm{IndCoh}(\mathcal{Z}(Y)_I)}(\mathrm{u}_{\mathcal{Z}(Y)_I\bullet}\mathcal{O}_{\mathcal{J}(Y)_I}, \mathrm{u}_{\mathcal{Z}(Y)_I\bullet}\mathcal{O}_{\mathcal{J}(Y)_I}) \qquad \in \mathrm{Alg}_{\mathbb{E}_2}(\mathcal{D}(X^I)) \ , \text{ and}
$$
$$
\mathcal{F}(Y)_x = \mathrm{Hom}_{\mathrm{IndCoh}(\mathcal{Z}(Y)_x)}(\mathrm{u}_{\mathcal{Z}(Y)_x\bullet}\mathcal{O}_{\mathcal{J}(Y)_x}, \mathrm{u}_{\mathcal{Z}(Y)_x\bullet}\mathcal{O}_{\mathcal{J}(Y)_x}) \qquad \in \mathrm{Alg}_{\mathbb{E}_2}(\mathrm{Vect}) \ .
$$

## 19.2. Internal construction of the four dimensional holomorphic-B model. Following Subsection 8.1, there is again an internal variant of the preceding construction:

*Example* 19.2.1. The self fibre product

$$
\zeta(Y) = \mathcal{J}(Y) \times_{\mathcal{Z}(Y)} \mathcal{J}(Y) \ \in \mathrm{PreStk}_{\mathrm{un}}^{\mathrm{fact}}(X)
$$

defines a unital factorization space over $X$.

*Remark* 19.2.2. Concretely, the prestacks $\zeta(Y)_I \in \mathrm{PreStk}_{/X^I}$ over $X^I$ assigned to each $I \in \mathrm{fSet}$ and the fibre $\zeta(Y)_x \in \mathrm{PreStk}$ over $x \in X$ are given by

$$
\zeta(Y)_I = \mathcal{J}(Y)_I \times_{\mathcal{Z}(Y)_I} \mathcal{J}(Y)_I \qquad \text{and} \qquad \zeta(Y)_x = Y_{\mathcal{O}_x} \times_{\mathcal{Z}(Y)_x} Y_{\mathcal{O}_x} \ .
$$

*Remark* 19.2.3. More generally, the iterated fibre products

$$
\zeta(Y)_{(n)} := \mathcal{J}(Y) \times_{\mathcal{Z}(Y)} \times \ldots \times_{\mathcal{Z}(Y)} \mathcal{J}(Y) \in \mathrm{PreStk}_{\mathrm{un}}^{\mathrm{fact}}(X)
$$

and projections $\pi_{ij} : \zeta(Y)_{(n)} \to \zeta(Y)$, define factorization spaces and maps of such for each $n \in \mathbb{N}$ and $i, j \in \{1, ..., n\}$.

*Remark* 19.2.4. Following Remark 4.2.5, $\zeta(Y) \in \mathrm{PreStk}_{\mathrm{un}}^{\mathrm{fact}}(X_{\mathrm{dR}})$ is canonically defined over $X_{\mathrm{dR}}$.

Again following Proposition 8.1.5, as well generalizing Proposition 10.2.3 in the factorization setting, we have:

*Proposition* 19.2.5. The factorization category

$$\mathrm{IndCoh}_{\zeta(Y)}^{\star} \in \mathrm{Cat}_{\mathbb{E}_2,\mathrm{un}}^{\mathrm{fact}}(X_{\mathrm{dR}})$$

is naturally a (unital) $\mathbb{E}_2$-factorization category with respect to the convolution monoidal structure defined by the composition

$$\mathrm{IndCoh}_{\zeta(Y)} \otimes^* \mathrm{IndCoh}_{\zeta(Y)} \xrightarrow{\pi_{12}^{\bullet}\boxtimes\pi_{23}^{\bullet}} \mathrm{IndCoh}_{\zeta_{(3)}(Y)^{\times 2}} \xrightarrow{\Delta^{\bullet}} \mathrm{IndCoh}_{\zeta(Y)_{(3)}} \xrightarrow{\pi_{13\bullet}} \mathrm{IndCoh}_{\zeta(Y)} .$$

Further, the pushforward functor $\mathrm{p}_{\zeta(Y)\bullet} : \mathrm{IndCoh}_{\zeta(Y)}^{\star} \to D_{\mathrm{Ran}_{X,\mathrm{un}}}^{\otimes^!}$ is a unital, lax $\mathbb{E}_2$-monoidal factorization functor.

Following Corollary 8.1.6, we have:

*Corollary* 19.2.6. The pushforward functor $\mathrm{p}_{\zeta(Y)\bullet} : \mathrm{IndCoh}_{\zeta(Y)}^{\star} \to D_{\mathrm{Ran}_{X,\mathrm{un}}}^{\otimes^!}$ induces a functor

$$(19.2.1) \qquad\qquad \mathrm{Alg}_{\mathbb{E}_2,\mathrm{un}}^{\mathrm{fact}}(\mathrm{IndCoh}_{\zeta(Y)}^{\star}) \to \mathrm{Alg}_{\mathbb{E}_2,\mathrm{un}}^{\mathrm{fact}}(X) .$$

Following Example 8.1.15, we define the internal variant of the four dimensional holomorphic-B model factorization $\mathbb{E}_2$ algebra:

*Definition* 19.2.7. The internal variant of the four dimensional holomorphic-B model to $Y$ is the factorization $\mathbb{E}_2$ algebra

$$\tilde{\mathcal{F}}(Y) = \tilde{\mathcal{H}om}_{\mathrm{IndCoh}(\mathcal{Z}(Y))}(\mathrm{unit}_{\mathrm{IndCoh}_{\mathcal{Z}(Y)}}, \mathrm{unit}_{\mathrm{IndCoh}_{\mathcal{Z}(Y)}}) = \mathcal{O}_{\mathcal{J}(Y)} \underset{\mathcal{Z}(Y)}{\boxtimes} \omega_{\mathcal{J}(Y)/\mathcal{Z}(Y)} \cong \omega_{\zeta(Y)/\mathcal{J}(Y)} \in \mathrm{Alg}_{\mathbb{E}_2,\mathrm{un}}^{\mathrm{fact}}(\mathrm{IndCoh}_{\zeta(Y)}^{\star})$$

internal to $\mathrm{IndCoh}_{\zeta(Y)}^{\star} \in \mathrm{Cat}_{\mathbb{E}_2,\mathrm{un}}^{\mathrm{fact}}(X_{\mathrm{dR}})$.

The preceding definition is justified by the following:

*Proposition* 19.2.8. The image of $\tilde{\mathcal{F}}(Y) \in \mathrm{Alg}_{\mathbb{E}_2,\mathrm{un}}^{\mathrm{fact}}(\mathrm{IndCoh}_{\zeta(Y)}^{\star})$ under the functor of Equation 19.2.1 above is canonically equivalent to $\mathcal{F}(Y) \in \mathrm{Alg}_{\mathbb{E}_2,\mathrm{un}}^{\mathrm{fact}}(X)$ of Definition 19.1.10.

### 19.3. The equivariant cigar compactification principle for the four dimensional holmorphic-B model.
In this subsection, we apply the results of Subsection I-25.4 to the example of the four dimensional holomorphic-B model. In particular, we deduce the equivalence of the existence of an $S^1$ equivariant structure on $\mathcal{F}(Y) \in \mathrm{Alg}_{\mathbb{E}_2,\mathrm{un}}^{\mathrm{fact}}(X)$ with the existence of a Tate structure on $Y$, and identify the corresponding factorization $\mathbb{BD}_0^u$ algebra with the two-periodic filtered quantization of chiral differential operators on $Y$. We recommend reviewing the analogous application of the results of Subsection I-25.2 to the three dimensional B model, explained in Subsection 11.4, before reading the present subsection.

*Definition* 19.3.1. The topological loop space of the meromorphic jet space $\mathcal{L}\mathcal{J}^{\mathrm{mer}}(Y) \in \mathrm{PreStk}_{\mathrm{un}}^{\mathrm{fact}}(X)$ is the factorization space defined by

$$\mathcal{L}\mathcal{J}^{\mathrm{mer}}(Y) = \mathcal{J}^{\mathrm{mer}}(Y) \times_{\mathcal{J}^{\mathrm{mer}}(Y)^{\times 2}} \mathcal{J}^{\mathrm{mer}}(Y) ,$$

following the construction of fibre products of factorization spaces in Example 8.2.1. The topological loop space of the jet space $\mathcal{L}\mathcal{J}(Y) \in \mathrm{PreStk}_{\mathrm{un}}^{\mathrm{fact}}(X)$ is defined analogously.

*Remark* 19.3.2. The factorization space $\mathcal{L}\mathcal{J}^{\mathrm{mer}}(Y)$ of the preceding definition is a shifted variant of $\mathcal{J}^{\mathrm{mer}}(Y)_{\mathrm{Dol}} \in \mathrm{PreStk}_{\mathrm{un}}^{\mathrm{fact}}(X)$ the Dolbeault stack of the meromorphic jet space of $Y$, as defined in Example 18.1.1.

*Proposition* 19.3.3. There is an equivalence of factorization prestacks

$$\zeta(Y) = \mathcal{J}(Y) \times_{\mathcal{Z}(Y)} \mathcal{J}(Y) \cong \mathcal{L}\mathcal{J}(Y) \times_{\mathcal{L}\mathcal{J}^{\mathrm{mer}}(Y)} \mathcal{J}^{\mathrm{mer}}(Y) \in \mathrm{PreStk}_{\mathrm{un}}^{\mathrm{fact}}(X)$$

where $\zeta(Y) \in \mathrm{PreStk}_{\mathrm{un}}^{\mathrm{fact}}(X)$ is as in Example 19.2.1 and $\mathcal{L}\mathcal{J}(Y)$ and $\mathcal{L}\mathcal{J}^{\mathrm{mer}}(Y)$ are as in Definition 19.3.1 above.

*Proof.* For simplicity, consider first the statement over a single point $x \in X$, which is written explicitly in Equation 19.3.1 below. In this case, each of the spaces admits the required maps such that, by the universal property of fibre products, it maps to the other; the induced maps are evidently inverse equivalences. This identification canonically extends in a factorization compatible way. $\square$

*Remark* 19.3.4. Concretely, the equivalence of the preceding Proposition is given over each point $x \in X$ by

$$(19.3.1) \qquad \zeta(Y)_x = Y_{\mathcal{O}} \underset{Y_{\mathcal{O}} \times Y_{\mathcal{K}}}{\times} Y_{\mathcal{O}} \cong (Y_{\mathcal{O}} \times_{Y_{\mathcal{O}}^{\times 2}} Y_{\mathcal{O}}) \underset{Y_{\mathcal{K}} \times_{Y_{\mathcal{K}}^{\times 2}} Y_{\mathcal{K}}}{\times} Y_{\mathcal{K}} = \mathcal{L}Y_{\mathcal{O}} \times_{\mathcal{L}Y_{\mathcal{K}}} Y_{\mathcal{K}} .$$

*Corollary* 19.3.5. There is an equivalence of factorization categories

$$\mathrm{IndCoh}_{\zeta(Y)} = \mathrm{IndCoh}_{\mathcal{J}(Y) \times_{\mathcal{Z}(Y)} \mathcal{J}(Y)} \cong \mathrm{IndCoh}_{\mathcal{L}\mathcal{J}(Y) \times_{\mathcal{L}\mathcal{J}^{\mathrm{mer}}(Y)} \mathcal{J}^{\mathrm{mer}}(Y)} \in \mathrm{Cat}_{\mathrm{un}}^{\mathrm{fact}}(X) .$$

*Example* 19.3.6. Under the equivalence of the preceding corollary, the internal variant of the four dimensional holomorphic-B model is given by

$$\tilde{\mathcal{F}}(Y) = \mathcal{O}_{\mathcal{J}(Y)} \underset{\mathcal{Z}(Y)}{\boxtimes} \omega_{\mathcal{J}(Y)/\mathcal{Z}(Y)} \cong \omega_{\mathcal{L}\mathcal{J}(Y)/\mathcal{J}(Y)} \underset{\mathcal{L}\mathcal{J}^{\mathrm{mer}}(Y)}{\boxtimes} \omega_{\mathcal{J}^{\mathrm{mer}}(Y)}^{-1} \in \mathrm{Alg}_{\mathrm{un}}^{\mathrm{fact}}(\mathrm{IndCoh}_{\zeta(Y)}) .$$

*Remark* 19.3.7. Concretely, the equivalence of the calculation of the preceding example is given over each point $x \in X$ by

$$\tilde{\mathcal{F}}(Y)_x = \mathcal{O}_{Y_{\mathcal{O}}} \underset{\mathcal{Z}(Y)_x}{\boxtimes} \omega_{Y_{\mathcal{O}}/\mathcal{Z}(Y)_x} \cong \omega_{\mathcal{L}Y_{\mathcal{O}}/Y_{\mathcal{O}}} \underset{\mathcal{L}Y_{\mathcal{K}}}{\boxtimes} \omega_{Y_{\mathcal{K}}}^{-1} .$$

From the latter description of the factorization space $\zeta(Y) \in \mathrm{PreStk}_{\mathrm{un}}^{\mathrm{fact}}(X)$, we have the following:

*Corollary* 19.3.8. There is a canonical deformation of $\zeta(Y)$ to a family of factorization spaces

$$\zeta^u(Y) := \mathcal{L}^u\mathcal{J}(Y) \times_{\mathcal{L}^u\mathcal{J}^{\mathrm{mer}}(Y)} \mathcal{J}^{\mathrm{mer}}(Y) \in \mathrm{PreStk}_{\mathrm{un}}^{\mathrm{fact}}(X)_{/(\mathbb{A}^1/\mathbb{G}_m)}$$

over $\mathbb{A}^1/\mathbb{G}_m$ with generic and central fibres given by the factorization spaces

$$\zeta^u(Y)|_{\{1\}} = \mathcal{J}^{\mathrm{mer}}(Y)_{\mathcal{J}(Y)}^{\wedge} \qquad \text{and} \qquad \zeta^u(Y)|_{\{0\}} = \zeta(Y) \in \mathrm{PreStk}_{\mathrm{un}}^{\mathrm{fact}}(X) .$$

*Corollary* 19.3.9. There is a canonical deformation of $\mathrm{IndCoh}_{\zeta(Y)}$ to a family of factorization categories

$$\mathrm{IndCoh}_{\zeta^u(Y)} \in \mathrm{Cat}_{\mathrm{un}}^{\mathrm{fact}}(X)_{/(\mathbb{A}^1/\mathbb{G}_m)}$$

with generic and central fibre given by the factorization categories

$$\mathrm{IndCoh}_{\zeta^u(Y)}|_{\{1\}} = \mathrm{IndCoh}_{\mathcal{J}^{\mathrm{mer}}(Y)_{\mathcal{J}(Y)}^{\wedge}} \qquad \text{and} \qquad \mathrm{IndCoh}_{\zeta^u(Y)}|_{\{0\}} = \mathrm{IndCoh}_{\zeta(Y)} \in \mathrm{Cat}_{\mathrm{un}}^{\mathrm{fact}}(X) .$$

Next, towards stating the main result of this subsection, we give the following description of the negative cyclic chains on $\mathcal{F}(Y) \in \mathrm{Alg}_{\mathbb{E}_2,\mathrm{un}}^{\mathrm{fact}}(X)$:

*Theorem* 19.3.10. There is a natural equivalence

$$\mathrm{CC}_\bullet^-(\mathcal{F}(Y)) \cong \mathcal{H}om_{\mathrm{IndCoh}(\mathcal{L}^u\mathcal{J}^{\mathrm{mer}}(Y))}(\mathrm{u}_\bullet\omega_{\mathcal{L}^u\mathcal{J}(Y)}, \mathrm{u}_\bullet\omega_{\mathcal{L}^u\mathcal{J}(Y)}) \in \mathrm{Alg}_{\mathbb{E}_1,\mathrm{un}}^{\mathrm{fact}}(X)_{/(\mathbb{A}^1/\mathbb{G}_m)} \ .$$

*Proof.* Analogously to the proof of 11.4.4, there is a natural commutative diagram

$$
\begin{array}{ccc}
\mathrm{Cat}_{\mathbb{E}_1,\mathrm{un}}^{\mathrm{fact}}(X) \xrightarrow{\mathcal{E}nd_{(\cdot)}(\mathrm{unit})} \mathrm{Alg}_{\mathbb{E}_2,\mathrm{un}}^{\mathrm{fact}}(X) & & \mathrm{IndCoh}(\mathcal{Z}(Y)) \longmapsto \mathcal{F}(Y) \\
\Big\downarrow \mathrm{CC}_\bullet^- \qquad\qquad \Big\downarrow \mathrm{CC}_\bullet^- & \text{under which} & \Big\downarrow \qquad\qquad \Big\downarrow \\
\mathrm{Cat}_{\mathrm{un}}^{\mathrm{fact}}(X)_{/\mathbb{K}[u]} \xrightarrow{\mathcal{E}nd_{(\cdot)}(\mathrm{unit})} \mathrm{Alg}_{\mathbb{E}_1}^{\mathrm{fact}}(X)_{/\mathbb{K}[u]} & & \mathrm{IndCoh}(\mathcal{L}^u\mathcal{J}^{\mathrm{mer}}(Y)) \longmapsto \mathcal{E}nd_{\mathrm{IndCoh}(\mathcal{L}^u\mathcal{J}^{\mathrm{mer}}(Y))}(\mathrm{u}_\bullet\omega_{\mathcal{L}^u\mathcal{J}(Y)})
\end{array}
$$

where we have again assumed an extension of the homological variant of Theorem 5.3 of [BZFN10] to compute $\mathrm{CC}_\bullet^-(\mathrm{IndCoh}(\mathcal{Z}(Y))) \cong \mathrm{IndCoh}(\mathcal{L}^u\mathcal{J}^{\mathrm{mer}}(Y))$. □

*Remark* 19.3.11. The family of factorization $\mathbb{E}_1$ algebras

$$\mathcal{C}(Y)^u := \mathcal{H}om_{\mathrm{IndCoh}(\mathcal{L}^u\mathcal{J}^{\mathrm{mer}}(Y))}(\mathrm{u}_\bullet\omega_{\mathcal{L}^u\mathcal{J}(Y)}, \mathrm{u}_\bullet\omega_{\mathcal{L}^u\mathcal{J}(Y)}) \in \mathrm{Alg}_{\mathbb{E}_1,\mathrm{un}}^{\mathrm{fact}}(X)_{/(\mathbb{A}^1/\mathbb{G}_m)}$$

in the preceding Theorem is a two-periodic analogue of the family of factorization $\mathbb{E}_1$ algebras giving the deformation from the three dimensional holomorphic-B model to A model as in Definition 18.2.6. In particular, the generic fibre is given by the three dimensional A model factorization $\mathbb{E}_1$ algebra

$$\mathcal{C}(Y)^u|_{\{1\}} = \mathcal{A}(Y) \in \mathrm{Alg}_{\mathbb{E}_1,\mathrm{un}}^{\mathrm{fact}}(X) \ .$$

We now state the main result of this subsection:

*Theorem* 19.3.12. The existence of an $S^1$ equivariant structure on $\mathcal{F}(Y) \in \mathrm{Alg}_{\mathbb{E}_2,\mathrm{un}}^{\mathrm{fact}}(X)$ is equivalent to the existence of a Tate structure on $Y$. Further, under the equivalence of Proposition I-25.4.3, the $S^1$ equivariant structure corresponds to

$$\mathcal{F}(Y)^u = \mathrm{p}_{\mathcal{J}^{\mathrm{mer}}(Y)\bullet}\mathrm{q}_{\mathcal{J}^{\mathrm{mer}}(Y)}^!\mathrm{u}_\bullet\omega_{\mathcal{L}^u\mathcal{J}(Y)} \in \mathrm{CC}_\bullet^-(\mathcal{F}(Y))\text{-Mod}(\mathrm{Alg}_{\mathrm{un}}^{\mathrm{fact}}(X)_{/(\mathbb{A}^1/\mathbb{G}_m)})$$

with module structure given by composition with endomorphisms, as in Example 7.1.10, under the equivalence of Proposition 19.3.10 above.

We give only a sketch of the proof, in keeping with Warning 19.0.1: By Proposition I-25.4.3, an $S^1$ equivariant structure on $\mathcal{F}(Y) \in \mathrm{Alg}_{\mathbb{E}_2,\mathrm{un}}^{\mathrm{fact}}(X)$ is equivalent to a deformation of $\mathcal{F}(Y)_0$ to a module factorization algebra $\mathcal{F}(Y)_u \in \mathrm{CC}_\bullet^-(\mathcal{F}(Y))\text{-Mod}(\mathrm{Alg}_{\mathrm{un}}^{\mathrm{fact}}(X)_{/(\mathbb{A}^1/\mathbb{G}_m)})$, which corresponds to a deformation of the internal object $\widetilde{\mathcal{F}}(Y) \in \mathrm{Alg}_{\mathrm{un}}^{\mathrm{fact}}(\zeta(Y))$ to an object of $\mathrm{IndCoh}(\zeta^u(Y))$.

From the description of Example 19.3.6, this appears to require a deformation of the object $\omega_{\mathcal{L}\mathcal{J}(Y)/\mathcal{J}(Y)} \in \mathrm{IndCoh}(\mathcal{L}\mathcal{J}(Y))$ to an object of $\mathrm{IndCoh}(\mathcal{L}^u\mathcal{J}(Y))$, which is a priori obstructed in analogy with Remark 10.3.14. A Tate structure on $Y$ is exactly what is required to identify

$$\omega_{\mathcal{L}\mathcal{J}(Y)/\mathcal{J}(Y)} \underset{\mathcal{L}\mathcal{J}^{\mathrm{mer}}(Y)}{\overset{\boxtimes}{}} \omega_{\mathcal{J}^{\mathrm{mer}}(Y)}^{-1} \cong \omega_{\mathcal{L}\mathcal{J}(Y)} \underset{\mathcal{L}\mathcal{J}^{\mathrm{mer}}(Y)}{\overset{\boxtimes}{}} \text{``}\omega_{\mathcal{J}(Y)/\mathcal{J}^{\mathrm{mer}}(Y)}\text{''}$$

so that the resulting object $\omega_{\mathcal{L}\mathcal{J}(Y)} \in \mathrm{IndCoh}(\mathcal{L}\mathcal{J}(Y))$ can be deformed to an object of $\mathrm{IndCoh}(\mathcal{L}^u\mathcal{J}(Y))$, again in analogy with Remark 10.3.14.

*Remark* 19.3.13. The underlying family of factorization algebras

$$\mathcal{F}(Y)^u = \mathrm{p}_{\mathcal{J}^{\mathrm{mer}}(Y)\bullet}\mathrm{q}_{\mathcal{J}^{\mathrm{mer}}(Y)}^!\mathrm{u}_\bullet\omega_{\mathcal{L}^u\mathcal{J}(Y)} \in \mathrm{Alg}_{\mathrm{un}}^{\mathrm{fact}}(X)_{/(\mathbb{A}^1/\mathbb{G}_m)}$$

in the preceding Theorem is a two-periodic analogue of the graded quantization of chiral differential operators to $Y$. Under this identification, following Remark 19.3.11 above, the module structure in the preceding Theorem is the two-periodic analogue of that in Theorem 18.4.8.

# Chapter 3
# Appendices

<div align="center">APPENDIX A. SHEAF THEORY</div>

Let $X$ be a smooth variety of dimension $d_X$ over $\mathbb{K} = \mathbb{C}$ or a field of characteristic 0. We write $\mathcal{O}_X$ for the sheaf of regular functions, $\mathcal{D}_X$ for the sheaf of differential operators, $\Theta_X$ for the tangent sheaf, $\Omega^1_X$ for the sheaf of Kahler differentials, $\Omega^{d_X}_X$ for the sheaf of sections of the canonical bundle, and $\omega_X = \Omega^{d_X}_X[d_X]$ for the dualizing sheaf on $X$. Let $\mathrm{Sh}_z(X)$ denote the category of complexes of sheaves of $\mathbb{K}$-modules on $X$ in the Zariski topology.

*Warning* A.0.1. Note that we assume $X$ is a smooth variety throughout, and only define the category of $D$ modules on more general spaces in Supappendix A.5.

A.1. $\mathcal{O}$-**module conventions.** Let $\mathrm{D}(\mathcal{O}_X)$ be the DG category of complexes of $\mathcal{O}_X$-modules, $\mathrm{QCoh}(X)$ and $\mathrm{Coh}(X)$ be the full subcategories of complexes with quasi-coherent and coherent cohomology sheaves, and $\mathrm{Perf}(X)$ the subcategory of bounded complexes with finitely generated cohomology sheaves. The category $\mathrm{D}(\mathcal{O}_X)$ is symmetric monoidal with respect to the tensor product $\otimes_{\mathcal{O}_X}$, with unit object $\mathcal{O}_X$, and $\mathrm{QCoh}(X)$, $\mathrm{Coh}(X)$, and $\mathrm{Perf}(X)$ are monoidal subcategories.

*Definition* A.1.1. Let $f : X \to Y$ a map of schemes. The inverse and direct image functors are

$$f^\bullet : \mathrm{D}(\mathcal{O}_Y) \to \mathrm{D}(\mathcal{O}_X) \qquad f^\bullet \mathcal{F} = f^{-1}F \otimes_{f^{-1}\mathcal{O}_Y} \mathcal{O}_X \qquad \text{and} \qquad f_\bullet : \mathrm{D}(\mathcal{O}_X) \to \mathrm{D}(\mathcal{O}_Y) \qquad f_\bullet \mathcal{F} = f_\bullet \mathcal{F},$$

where $f_\bullet : \mathrm{Sh}_z(X) \to \mathrm{Sh}_z(Y)$ and $f^{-1} : \mathrm{Sh}_z(Y) \to \mathrm{Sh}_z(X)$ are the usual direct and inverse image functors on sheaves of $\mathbb{K}$-modules.

*Remark* A.1.2. Note that $f^\bullet$ preserves quasicoherence, as does $f_\bullet$ for quasicompact, quasiseperated maps. We define the global sections functor by $\Gamma = \pi_\bullet : \mathrm{D}(\mathcal{O}_X) \to \mathrm{Vect}$ where $\pi : X \to \mathrm{pt}$.

*Definition* A.1.3. Let $\mathcal{F}, \mathcal{G} \in \mathrm{D}(\mathcal{O}_X)$. The internal hom object in $\mathrm{D}(\mathcal{O}_X)$ is

$$\underline{\mathrm{Hom}}_{\mathcal{O}_X}(\mathcal{F}, \mathcal{G}) \in \mathrm{D}(\mathcal{O}_X) \qquad \text{by} \qquad \underline{\mathrm{Hom}}_{\mathcal{O}_X}(\mathcal{F}, \mathcal{G})(U) := \mathrm{Hom}_{\mathcal{O}_X|_U}(\mathcal{F}|_U, \mathcal{G}|_U) .$$

*Remark* A.1.4. For $\mathcal{H} \in \mathrm{D}(\mathcal{O}_X)$, we have

$$\mathrm{Hom}(\mathcal{H}, \underline{\mathrm{Hom}}_{\mathcal{O}_X}(\mathcal{F}, \mathcal{G})) \cong \mathrm{Hom}(\mathcal{H} \otimes_{\mathcal{O}_X} \mathcal{F}, \mathcal{G}) .$$

In particular, the space of homomorphisms is given by the space of sections of the internal hom object

$$\mathrm{Hom}(\mathcal{F}, \mathcal{G}) = \mathrm{Hom}(\mathcal{O}_X, \underline{\mathrm{Hom}}_{\mathcal{O}_X}(\mathcal{F}, \mathcal{G})) = \Gamma(X, \underline{\mathrm{Hom}}_{\mathcal{O}_X}(\mathcal{F}, \mathcal{G})).$$

*Remark* A.1.5. For $\mathcal{F} \in \mathrm{Coh}(X)$ coherent and $\mathcal{G} \in \mathrm{QCoh}(X)$ quasi-coherent, the object $\underline{\mathrm{Hom}}_{\mathcal{O}_X}(\mathcal{F}, \mathcal{G}) \in \mathrm{QCoh}(X)$ is quasi-coherent. If $\mathcal{F}, \mathcal{G} \in \mathrm{Coh}(X)$ are both coherent, then $\underline{\mathrm{Hom}}_{\mathcal{O}_X}(\mathcal{F}, \mathcal{G}) \in \mathrm{Coh}(X)$ is also coherent.

*Definition* A.1.6. The duality functor on coherent $\mathcal{O}_X$-modules is defined by

$$(-)^\vee : \mathrm{Coh}(X) \to \mathrm{Coh}(X) \qquad \text{by} \qquad \mathcal{F} \mapsto \mathcal{F}^\vee := \underline{\mathrm{Hom}}_{\mathcal{O}_X}(\mathcal{F}, \mathcal{O}_X).$$

*Remark* A.1.7. There are canonical isomorphisms $\underline{\mathrm{Hom}}_{\mathcal{O}_X}(\mathcal{F}, \mathcal{G}) \cong \mathcal{G} \otimes_{\mathcal{O}_X} \mathcal{F}^\vee$ and $(\mathcal{F}^\vee)^\vee \cong \mathcal{F}$.

A.2. **$D$-module conventions.** Let $D^l(X)$ and $D^r(X)$ be the concrete DG categories of complexes of left and right $\mathcal{D}_X$-modules which are quasicoherent as $\mathcal{O}_X$-modules, and let $D_c^l(X)$ and $D_c^r(X)$ denote the full sub DG categories of complexes with cohomology that is coherent as a module over $\mathcal{D}_X$.

*Example* A.2.1. The sheaf of regular functions $\mathcal{O}_X \in D^l(X)^\heartsuit$ has the structure of a left $D$ module, given by the defining action of the sheaf of differential operators $\mathcal{D}_X$ on $\mathcal{O}_X$.

More generally, a left $D$ module (or a complex of such) $M \in D^l(X)$ on $X$ is given by a quasicoherent sheaf (or a complex of such) $M \in \mathrm{QCoh}(X)$, together with a flat connection, that is, $\nabla \in \mathrm{Hom}_{\mathrm{Sh}_\mathbb{K}(X)}(M, \Omega_X^1 \otimes_{\mathcal{O}_X} M)$ such that

- $\nabla_\theta(fs) = \theta(f)s + f\nabla_\theta(s)$ , and
- $\nabla_{[\theta_1,\theta_2]}s = [\nabla_{\theta_1}, \nabla_{\theta_2}]s$ ,

where $\theta, \theta_1, \theta_2 \in \Theta_X$, $f \in \mathcal{O}_X$, and $s \in M$. The first condition is that $\nabla$ defines a connection, and the second that $\nabla$ is flat.

*Example* A.2.2. The sheaf of sections of the canonical bundle $\Omega_X^{d_X} \in D^r(X)^\heartsuit$ is the protypical example of a right $D_X$ module, with action of vector fields given by $\theta(\eta) = -\mathrm{Lie}_\theta(\eta)$ for $\theta \in \Theta_X$ and $\eta \in \Omega_X^{d_X}$.

*Remark* A.2.3. There is a canonical equivalence of the categories $D^l(X)$ and $D^r(X)$

$$D^l(X) \underset{(-)^l}{\overset{(-)^r}{\rightleftarrows}} D^r(X) \qquad \text{defined by} \qquad \begin{cases} M \mapsto M^l := M \otimes_{\mathcal{O}_X} \omega_X^\vee & \text{for } M \in D^r(X) \quad \text{and} \\ L \mapsto L^r := \omega_X \otimes_{\mathcal{O}_X} L & \text{for } L \in D^l(X). \end{cases}$$

We write $D(X)$ for the abstract DG category given by the common value of $D^r(X)$ and $D^l(X)$ under this identification, and $D_c(X)$ for the full sub DG category corresponding to $D_c^r(X)$ and $D_c^l(X)$, which are also identified under this equivalence. $D^r(X)$ and $D^l(X)$ both have natural forgetful functors to $\mathrm{QCoh}(X)$, which are intertwined by tensoring with $\omega_X$. This perspective is summarized in the following diagram:

$$\begin{array}{ccc} D^l(X) & \xrightarrow[\simeq]{\omega_X} & D^r(X) \\ \downarrow{o^l} & & \downarrow{o^r} \\ \mathrm{QCoh}^l(X) & \xrightarrow[\simeq]{\omega_X} & \mathrm{QCoh}^r(X) \end{array} \qquad \text{so that} \qquad \begin{array}{c} D(X) \\ {}^{o^l}\downarrow \quad \searrow^{o^r} \\ \mathrm{QCoh}(X) \xrightarrow[\simeq]{\omega_X} \mathrm{QCoh}(X) \end{array} \quad ,$$

where $\mathrm{QCoh}^l(X)$ and $\mathrm{QCoh}^r(X)$ are just the category $\mathrm{QCoh}(X)$

*Remark* A.2.4. Throughout, when defining a functor involving (potentially several copies of) the category $D(X)$, we will prescribe the values of the functor in terms of a particular choice of realization $D^r(X)$ or $D^l(X)$ for each copy of $D(X)$, with the extension to all other choices of concrete realizations of $D(X)$ implicitly specified via the above equivalence.

*Remark* A.2.5. Note that the above equivalence is exact up to a cohomological degree shift of $d_X = \dim_\mathbb{K} X$, so that the category $D(X)$ inherits two different t-structures, which differ only by this shift. We choose to preference the right t structure, and all statements about exactness of functors involving $D(X)$ will be given in these terms. This t-structure will be the one which corresponds to the perverse t-structure on constructible sheaves under the Riemann-Hilbert correspondence. In particular, under this identification $\omega_X \in D(X)$ is the dualizing sheaf, $\omega_X[-d_X] \in D(X)^\heartsuit$ is the IC sheaf, and $\underline{\mathbb{K}}_X := \omega_X[-2d_X] \in D(X)$ is the constant sheaf.

*Definition* A.2.6. The $\otimes^!$ monoidal structure on $D(X)$ is $\otimes^! : D(X)^{\otimes 2} \to D(X)$ defined by

$$\otimes^! : D^l(X) \times D^l(X) \to D^l(X) \qquad M \otimes^! N = M \otimes_{\mathcal{O}_X} N \quad \text{with} \quad P(m \otimes n) = Pm \otimes n + m \otimes Pn \,,$$

for $P \in D_X$.

*Remark* A.2.7. This formula agrees with the usual definition of the tensor product of connections, and tensor products of flat connections are flat. The corresponding functor $\otimes^! : D^r(X) \otimes D^r(X) \to D^r(X)$ is given by $M \otimes^! N = M \otimes_{\mathcal{O}_X} N \otimes_{\mathcal{O}_X} \omega_X^\vee$.

Let $\mathbb{1} \in D(X)$ denote the tensor unit, and note $o^l(\mathbb{1}) = \mathcal{O}_X$ and $o^r(\mathbb{1}) = \omega_X$. We will often use just $\otimes$ to denote this symmetric monoidal structure on $D(X)$.

*Definition* A.2.8. Let $f : X \to Y$ be a map of smooth varieties. The inverse image functor $f^! : D(Y) \to D(X)$ is defined by

$$f^! : D^l(Y) \to D^l(X) \qquad f^!(M) = f^\bullet(M) \quad \text{equipped with the pullback flat connection.}$$

*Remark* A.2.9. This functor is symmetric monoidal with respect to $\otimes^!$, and in particular maps the tensor unit $\mathbb{1}_Y$ to $\mathbb{1}_X$.

*Remark* A.2.10. The corresponding functor $f^! : D^r(Y) \to D^r(X)$ is given by

$$f^!(M) = f^\bullet(M \otimes_{\mathcal{O}_Y} \omega_Y^\vee) \otimes_{\mathcal{O}_X} \omega_X \cong f^\bullet M \otimes_{\mathcal{O}_X} \omega_{X/Y} \,.$$

*Definition* A.2.11. The exterior product is defined by

$$\boxtimes : D(X) \times D(Y) \to D(X \times Y) \qquad \text{by} \qquad M \boxtimes N = \pi_X^! M \otimes \pi_Y^! N \,,$$

for $\pi_X : X \times Y \to X, \pi_Y : X \times Y \to Y$.

*Remark* A.2.12. Note that

$$M \otimes N = \Delta^!(M \boxtimes N)$$

for $M, N \in D(X)$ and $\Delta : X \to X \times X$ the diagonal embedding.

*Definition* A.2.13. Let $f : X \to Y$ again be a map of smooth varieties. The direct image functor is

$$f_* : D^r(X) \to D^r(Y) \qquad f_*(M) = f_\bullet(M \otimes_{\mathcal{D}_X} \mathcal{D}_{X \to Y}) \qquad \text{for} \qquad \mathcal{D}_{X \to Y} := f^! \mathcal{D}_Y \in (\mathcal{D}_X, f^{-1}\mathcal{D}_Y)\text{-Mod}$$

where $\mathcal{D}_{X \to Y} = f^! \mathcal{D}_Y \in D^l(X)$ is defined in terms of $\mathcal{D}_Y \in D^l(Y)$ as a left module, so that the additional $(\mathcal{D}_Y, \mathcal{D}_Y)$-bimodule structure on $\mathcal{D}_Y$ equips $\mathcal{D}_{X \to Y}$ with the structure of a $(\mathcal{D}_X, f^{-1}\mathcal{D}_Y)$-bimodule.

*Definition* A.2.14. The de Rham cochains functor is $C_{\mathrm{dR}}^\bullet := \pi_* : D(X) \to \mathrm{Vect}_{\mathbb{K}}$, where $\pi : X \to \mathrm{pt}$. The de Rham chains functor is $C_\bullet^{\mathrm{dR}} := \pi_! : D(X) \to \mathrm{Vect}_{\mathbb{K}}$.

*Remark* A.2.15. Note that the de Rham cochain and chains functors are calculated as

$$C_{\mathrm{dR}}^\bullet : D^r(X) \to \mathrm{Vect} \qquad\qquad C_{\mathrm{dR}}^\bullet(X; M) = \pi_\bullet(M \otimes_{D_X} \mathcal{O}_X)$$

$$C_{\mathrm{dR}}^\bullet : D^l(X) \to \mathrm{Vect} \qquad\qquad C_{\mathrm{dR}}^\bullet(X; M) = \pi_\bullet(\omega_X \otimes_{D_X} M) \,.$$

*Definition* A.2.16. The sheaf internal hom functor

$$\underline{\mathrm{Hom}}_{D(X)}(\cdot, \cdot) : D(X)^{\mathrm{op}} \times D(X) \to \mathrm{Sh}_z(X) \qquad \text{by} \qquad \underline{\mathrm{Hom}}_{D(X)}(M, N)(U) = \mathrm{Hom}_{D(U)}(j^! M, j^! N) \,,$$

for $M, N \in D(X)$, where $j : U \to X$ is the open embedding.

*Remark* A.2.17. Note that

$$\Gamma \circ \underline{\mathrm{Hom}} = \mathrm{Hom} : D(X)^{\mathrm{op}} \times D(X) \to \mathrm{Vect} \qquad\qquad \Gamma(X, \underline{\mathrm{Hom}}_{D(X)}(M, N)) = \mathrm{Hom}_{D(X)}(M, N) \ .$$

*Definition* A.2.18. The duality functor $\mathbb{D} : D_{\mathrm{c}}(X)^{\mathrm{op}} \to D_{\mathrm{c}}(X)$ is defined by

$$\mathbb{D} : D_{\mathrm{c}}^r(X)^{\mathrm{op}} \to D_{\mathrm{c}}^l(X) \qquad\qquad \mathbb{D}(M) = \underline{\mathrm{Hom}}_{D^r(X)}(M, \mathcal{D}_X) \ ,$$

where $\mathcal{D}_X \in D^r(X)$ is considered as a $(\mathcal{D}_X, \mathcal{D}_X)$-bimodule so that $\mathbb{D}(M)$, which is a priori an object in $\mathrm{Sh}_z(X)$, defines an object of $D^l(X)$ as desired.

*Remark* A.2.19. Note that $\mathbb{D}$ preserves coherence, but if $M$ is not coherent, then the resulting object of $\mathcal{D}_X$-Mod is not in general quasicoherent as an object of $\mathrm{D}(\mathcal{O}_X)$.

*Definition* A.2.20. The genuine internal hom functor $\mathcal{H}om(\cdot, \cdot) : D_{\mathrm{c}}(X)^{\mathrm{op}} \otimes D(X) \to D(X)$ is defined by

$$\mathcal{H}om(\cdot, \cdot) : D_{\mathrm{c}}^r(X)^{\mathrm{op}} \times D^l(X) \to D^l(X) \qquad\qquad \mathcal{H}om(M, N) = \underline{\mathrm{Hom}}_{D^r(X)}(M, N \otimes_{\mathcal{O}_X} \mathcal{D}_X) \ ,$$

where $N \otimes_{\mathcal{O}_X} \mathcal{D}_X \in D^r(X)$ is considered as a $(\mathcal{D}_X, \mathcal{D}_X)$-bimodule so that $\mathcal{H}om(M, N) \in D^l(X)$ as above.

*Remark* A.2.21. Note that

$$C_{\mathrm{dR}}^\bullet \circ \mathcal{H}om = \mathrm{Hom} : D_{\mathrm{c}}(X)^{\mathrm{op}} \times D(X) \to \mathrm{Vect} \qquad\qquad C_{\mathrm{dR}}^\bullet(X, \mathcal{H}om_{D(X)}(M, N)) = \mathrm{Hom}_{D(X)}(M, N) \ .$$

Further, we have

$$\mathcal{H}om(\cdot, \cdot) = \mathbb{D}(\cdot) \otimes^! (\cdot) : D_{\mathrm{c}}(X)^{\mathrm{op}} \times D(X) \to D(X) \ ,$$

and in particular $\mathcal{H}om(\cdot, \mathbb{1}) = \mathbb{D} : D(X)^{\mathrm{op}} \to D(X)$; we could equivalently take this as the definition of $\mathcal{H}om$.

*Remark* A.2.22. The pushforward and pullback functors $f_*$ and $f^!$ above were defined on the entire category $D(X)$, but their putative adjoints can not always be defined. In general, the best we can do is the following: Let $f : X \to Y$ again be a map of smooth varieties, and let $D_{\mathrm{c}}^{f^!}(Y)$ be the full subcategory of objects $M \in D_{\mathrm{c}}(Y)$ such that $f^! \mathbb{D} M \in D_{\mathrm{c}}(X)$ is coherent, and similarly $D_{\mathrm{c}}^{f_*}(X)$ be the full subcategory of objects $M \in D_{\mathrm{c}}(X)$ such that $f_* \mathbb{D} M \in D_{\mathrm{c}}(Y)$ is coherent. Then we define

$$f^* := \mathbb{D} f^! \mathbb{D} : D_{\mathrm{c}}^{f^!}(Y) \to D_{\mathrm{c}}(X) \qquad\qquad f_! := \mathbb{D} f_* \mathbb{D} : D_{\mathrm{c}}^{f_*}(X) \to D_{\mathrm{c}}(Y) \ .$$

In various situations, these definitions simplify to more useful ones, as in the following propositions.

*Proposition* A.2.23. Let $f : X \to Y$ be a smooth map of relative dimension $d = d_X - d_Y$ of smooth varieties. Then $f^! : D_{\mathrm{c}}(Y) \to D_{\mathrm{c}}(X)$ preserves coherence, so that $f^* : D_{\mathrm{c}}(Y) \to D_{\mathrm{c}}(X)$ is defined. Moreover, in this case $f^* = f^![-2d]$, and we have a natural isomorphism

$$\mathrm{Hom}_{D(X)}(f^* M, N) \cong \mathrm{Hom}_{D(Y)}(M, f_* N)$$

of functors $D_{\mathrm{c}}(X) \times D(Y) \to \mathrm{Vect}$.

*Proposition* A.2.24. Let $f : X \to Y$ be a proper map of smooth varieties. Then $f_* : D_{\mathrm{c}}(X) \to D_{\mathrm{c}}(Y)$ preserves coherence, so that $f_! : D_{\mathrm{c}}(X) \to D_{\mathrm{c}}(Y)$ is defined. Moreover, in this case $f_! = f_*$ and we have a natural isomorphism

$$\mathrm{Hom}_{D(Y)}(f_! M, N) \cong \mathrm{Hom}_{D(X)}(M, f^! N)$$

of functors $D_{\mathrm{c}}(X) \times D(Y) \to \mathrm{Vect}$.

A.3. **The six functors formalism.** Let $D_{\mathrm{rh}}(X)$ the full subcategory of $D(X)$ on bounded complexes with regular holonomic cohomology modules.

*Theorem* A.3.1. There are functors

$$\otimes^! : D_{\mathrm{rh}}(X) \times D_{\mathrm{rh}}(X) \to D_{\mathrm{rh}}(X) \qquad \mathbb{D} : D_{\mathrm{rh}}(X)^{\mathrm{op}} \xrightarrow{\cong} D_{\mathrm{rh}}(X) \qquad \mathcal{H}\mathrm{om}_X : D_{\mathrm{rh}}(X)^{\mathrm{op}} \times D_{\mathrm{rh}}(X) \to D_{\mathrm{rh}}(X) ,$$

and for $f : X \to Y$ natural adjunctions,

$$f^* : D_{\mathrm{rh}}(Y) \rightleftarrows D_{\mathrm{rh}}(X) : f_* \qquad f_! : D_{\mathrm{rh}}(X) \rightleftarrows D_{\mathrm{rh}}(Y) : f^! .$$

Moreover, these satisfy:

- for $f : X \to Y$ smooth of relative dimension $d$, $f^* = f^![-2d]$ as in A.2.23 above;
- for $f : X \to Y$ proper, $f_* = f_!$ as in A.2.24 above;
- for $f : X \to Y$, there are natural equivalences $\mathbb{D}_Y f_* \cong f_! \mathbb{D}_Y$, $\mathbb{D}_X f^* \cong f^! \mathbb{D}_Y$;
- $\otimes^!$ defines a symmetric monoidal structure on $D_{\mathrm{rh}}(X)$; and
- for $f : X \to Y$, there are natural equivalences

$$f_!(M \otimes f^* N) \cong f_!(M) \otimes N \qquad \mathcal{H}\mathrm{om}_Y(f_! M, N) \cong f_* \mathcal{H}\mathrm{om}_X(M, f^! N) \qquad f^! \mathcal{H}\mathrm{om}_Y(M, N) \cong \mathcal{H}\mathrm{om}_X(f^* M, f^! N) .$$

- For a Cartesian square

(A.3.1)
$$\begin{array}{ccc} Z & \xrightarrow{\tilde{g}} & X \\ \downarrow{\tilde{f}} & & \downarrow{f} \\ Y & \xrightarrow{g} & W \end{array} \qquad \text{there is a natural isomorphism} \qquad \tilde{f}_* \circ \tilde{g}^! \cong g^! \circ f_* .$$

*Definition* A.3.2. The $\otimes^*$ tensor structure $\otimes^* : D_{\mathrm{rh}}(X)^{\times 2} \to D_{\mathrm{rh}}(X)$ is defined by $(M, N) \mapsto \Delta^*(M \boxtimes N)$.

*Definition* A.3.3. The constant $D$ module on $X$ is $\underline{\mathbb{K}}_X = \pi^* \underline{\mathbb{K}}_{\mathrm{pt}} \in D_{\mathrm{rh}}(X)$, where $\pi : X \to \mathrm{pt}$ and $\underline{\mathbb{K}}_{\mathrm{pt}} \in D(\mathrm{pt})$ is the object corresponding to $\mathbb{K} \in \mathrm{Vect}_{\mathbb{K}} \cong D(\mathrm{pt})$.

*Example* A.3.4. If $X$ is smooth of dimension $d_X$, then $\underline{\mathbb{K}}_X \cong \omega_X[-2d_X]$ is a shift of the dualizing sheaf, as in Remark A.2.5.

*Definition* A.3.5. The de-Rham (Borel-Moore) (co)chains on $X$ are

$$C_\bullet(X) = \pi_! \pi^! \underline{\mathbb{K}}_{\mathrm{pt}} \qquad C_\bullet^{\mathrm{BM}}(X) = \pi_* \pi^! \underline{\mathbb{K}}_{\mathrm{pt}} \qquad C^\bullet(X) = \pi_* \pi^* \underline{\mathbb{K}}_{\mathrm{pt}} \qquad C_{\mathrm{c}}^\bullet(X) = \pi_! \pi^* \underline{\mathbb{K}}_{\mathrm{pt}} \quad \in \quad D(\mathrm{pt})$$

where $\pi : X \to \mathrm{pt}$ and $\underline{\mathbb{K}}_{\mathrm{pt}} \in D(\mathrm{pt})$ are as in the preceding definition.

*Remark* A.3.6. Note that we use cochain complexes throughout, and do not use the convention of reversing the grading on homology. Thus, classes in homology which are correspond to higher dimensional cycles geometrically contribute to the homology groups of lower (cohomological) degree.

*Example* A.3.7. For $X$ a smooth variety of dimension $d_X$, we have an isomorphism $C_\bullet^{\mathrm{BM}}(X) \cong C_{\mathrm{dR}}^\bullet(X)[2d_X]$.

*Remark* A.3.8. The unit and counit of the above adjunctions for $f : X \to Y$ give canonical maps

$$A \to f_* f^* A \qquad \text{and} \qquad f_! f^! A \to A .$$

Applying these to $A = \underline{\mathbb{K}}_Y$ and $A = \omega_Y$, respectively, and composing with $\pi_*$ for $\pi : Y \to \mathrm{pt}$, we obtain maps

$$f^* : C^\bullet(Y) \to C^\bullet(X) \qquad \text{and} \qquad f_* : C_\bullet(X) \to C_\bullet(Y)$$

of objects in $D(\mathrm{pt}) = \mathrm{Vect}$, as expected for usual chains and cochains.

If $f$ is proper, then we similarly have maps

$$f^* : C^\bullet_c(Y) \to C^\bullet_c(X) \qquad \text{and} \qquad f_* : C^{\mathrm{BM}}_\bullet(X) \to C^{\mathrm{BM}}_\bullet(Y) \ ,$$

while if $f$ is smooth of relative dimension $d = d_X - d_Y$, then we have maps

$$f^* : C^{\mathrm{BM}}_\bullet(Y) \to C^{\mathrm{BM}}_\bullet(X)[-2d] \qquad \text{and} \qquad f_* : C^\bullet_c(X) \to C^\bullet_c(Y)[-2d] \ .$$

Finally, if $f$ is proper and smooth of relative dimension $d$, then we have maps

$$f^* : C_\bullet(Y) \to C_\bullet(X)[-2d] \qquad \text{and} \qquad f_* : C^\bullet(X) \to C^\bullet(Y)[-2d] \ .$$

A.4. **The de Rham functor and Riemann-Hilbert correspondence.** Throughout this section, let $X$ be a smooth, finite dimensional variety over $\mathbb{K} = \mathbb{C}$, and let $\Omega^\bullet_X \in \mathrm{Sh}_z(X)$ denote the algebraic de Rham complex, viewed as a complex of sheaves with the usual de Rham differential.

*Remark* A.4.1. Each $\Omega^i_X \in \mathcal{O}_X\text{-Mod}$ is a coherent $\mathcal{O}_X$ module, but the differential on $\Omega^\bullet_X$ is not $\mathcal{O}_X$ linear. Rather, the de Rham differential $d_{\mathrm{dR}} \in \mathrm{Diff}(\Omega^i, \Omega^{i+1})$ is a differential operator, so the de Rham complex can equivalently be described in terms of the induced complex of $D$ modules $\Omega^\bullet_{X,\mathcal{D}} = \Omega^\bullet_X \otimes_{\mathcal{O}_X} \mathcal{D}_X \in D^r(X)$, recalling $\mathrm{Diff}(\mathcal{F}, \mathcal{G}) = \mathrm{Hom}_{D(X)}(\mathcal{F}_D, \mathcal{G}_D)$.

*Proposition* A.4.2. There is a natural quasiisomorphism $\Omega^\bullet_{X,\mathcal{D}} \xrightarrow{\cong} \omega_X[-2d_X] = \underline{\mathbb{K}}_X \in D^r(X)$.

*Example* A.4.3. For $M \in D(X)$, applying this resolution to the calculation of de Rham cochains of $M^l \in D^l(X)$ following A.2.15 yields

$$C^\bullet_{\mathrm{dR}}(X; A) = \pi_\bullet(\omega_X \otimes_{D_X} M) \cong \Gamma(X, \Omega^\bullet_X \otimes_{\mathcal{O}_X} M^l)[2d_X]$$

where $\Omega^\bullet_X \otimes_{\mathcal{O}_X} M \in \mathrm{Sh}_z(X)$ denotes the usual de Rham complex with coefficients in a complex of $\mathcal{O}_X$ modules with flat connection. For $X$ smooth and projective, this is calculated by $\Gamma(X^{\mathrm{an}}, \Omega^\bullet_{X^{\mathrm{an}}} \otimes_{\mathcal{O}_{X^{\mathrm{an}}}} M^{\mathrm{an}})$ where $\Omega^\bullet_{X^{\mathrm{an}}} \otimes_{\mathcal{O}_{X^{\mathrm{an}}}} M^{\mathrm{an}} \in \mathrm{Sh}(X^{\mathrm{an}}; \mathbb{C})$ denotes the analytic variant of the above de Rham complex.

*Definition* A.4.4. The analytic de Rham functor is

$$\mathrm{dR} : D(X) \to \mathrm{D}^b(X) \qquad \text{defined by} \qquad \mathrm{dR}(M) = \Omega^\bullet_{X^{\mathrm{an}}} \otimes_{\mathcal{O}_{X^{\mathrm{an}}}} M^{l,\mathrm{an}}[2d_x] \ ,$$

for each $M \in D(X)$, where $\mathrm{D}^b(X) = \mathrm{D}^b(X(\mathbb{C}))$ denotes the derived category of sheaves on $X$ in the analytic topology, as in subappendix I-A.5.

Let $\mathrm{D}^b_c(X) = \mathrm{D}^b_c(X(\mathbb{C}); \mathbb{C})$ denote the bounded derived category constructible sheaves on $X$ in the analytic topology, as defined in *loc. cit.*.

*Theorem* A.4.5. The de Rham functor restricts to a derived equivalence $\mathrm{dR} : D_{\mathrm{rh}}(X) \to \mathrm{D}^b_c(X)$.

Moreover, it naturally intertwines the six functors operations stated in Theorems A.3.1 and I-A.5.3, with the caveat that it intertwines the $\otimes^*$ tensor structure on $D_{\mathrm{rh}}(X)$ of Definition A.3.2 with that of Theorem I-A.5.3.

*Remark* A.4.6. The object $\mathcal{O}_X \in D^l(X)$ or equivalently $\omega_X \in D^r(X)$ corresponds to

$$\mathrm{dR}(\mathcal{O}_X) = \Omega^\bullet_X[2n] \cong \underline{\mathbb{K}}_X[2n]$$

which is the dualizing sheaf in $\mathrm{D}^b_c(X)$. Equivalently, the objects $\underline{\mathbb{K}}_X = \mathcal{O}_X[-2n]$ and $\underline{\mathrm{IC}}_X = \mathcal{O}_X[-n]$ correspond to the constant sheaf $\mathbb{K}_X$ and the intersection cohomology sheaf $\mathrm{IC}_X = \mathbb{K}_X[n]$, in keeping with Remark A.2.5.

A.5. $\Omega_X^\bullet$-**modules and the de Rham stack.** In this subappendix, we explain an alternate description of the category of $D$ modules on a smooth variety defined above, and give the definition of the category of $D$ modules on more general spaces via this approach, following [] and references therein.

*Remark* A.5.1. The de Rham functor $\mathrm{dR} : D(X) \to \mathrm{D}^b(X)$ of Definition A.4.4 can equivalently be expressed in terms of the sheaf internal Hom functor

$$\mathrm{dR} = \underline{\mathrm{Hom}}_{D(X)}(\underline{\mathbb{K}}_X, \cdot) : D(X) \to \mathrm{D}^b(X) \ .$$

Further, the image of the constant sheaf $\underline{\mathbb{K}}_X = \mathcal{O}_X[-2d_X] \in D^l(X)$ is given by the de Rham complex

$$\mathrm{dR}(\underline{\mathbb{K}}_X) = \underline{\mathrm{Hom}}_{D(X)}(\underline{\mathbb{K}}_X, \underline{\mathbb{K}}_X) \cong \Omega_X^\bullet \quad \in \mathrm{Alg}_{\mathbb{E}_1}(\mathrm{D}^b(X)) \ ,$$

such that the associative algebra structure corresponding to composition of endomorphisms is given by the usual wedge product of differential forms. Thus, the de Rham functor naturally lifts to a functor

$$\tilde{\mathrm{dR}} : D(X) \to \Omega_X^\bullet\text{-Mod} := \Omega_X^\bullet\text{-Mod}(\mathrm{D}^b(X)) \ ,$$

where the latter is the category of (DG) modules over $\Omega_X^\bullet$ internal to $\mathrm{D}^b(X)$.

*Remark* A.5.2. This functor defines an equivalence between the category of $D$ modules and that of sheaves of modules over $\Omega_X^\bullet$, given appropriate definitions of the two categories; we recall the variant of this setup which will be convenient for our purposes in definitions A.5.9 and A.5.10 below.

*Remark* A.5.3. In addition to identifying the algebra objects $\underline{\mathbb{K}}_X$ and $\Omega_X^\bullet$, the de Rham functor sends the object $D_X \in D^r(X)^*$ given by differential operators under right multiplication to the structure sheaf $D_X \mapsto \Omega_X^\bullet \otimes_{\mathcal{O}_X} D_X^l[2d_X] \cong \mathcal{O}_X \in \Omega_X^\bullet\text{-Mod}$ by Proposition A.4.2, equipped with the natural augmentation module structure.

*Remark* A.5.4. For $f : X \to Y$ a map of smooth varieties, the functors $f^!$ and $f_*$ of $D$ module pullback and pushforward as in definitions A.2.8 and A.2.13 are intertwined with the functors

$$f^! : \Omega_Y^\bullet\text{-Mod} \to \Omega_X^\bullet\text{-Mod} \qquad\qquad f^! M = \Omega_X^\bullet \otimes_{f^{-1}\Omega_Y^\bullet} f^{-1}(M)$$

$$f_* : \Omega_X^\bullet\text{-Mod} \to \Omega_Y^\bullet\text{-Mod} \qquad\qquad f_* M = f_\bullet M \ ,$$

where $f^{-1} : \mathrm{D}^b(Y) \to \mathrm{D}^b(X)$ and $f_\bullet : \mathrm{D}^b(X) \to \mathrm{D}^b(Y)$ denote the sheaf-theoretic inverse and direct image.

*Remark* A.5.5. The functors $f^!$ and $f_*$ of the preceding remark agree with the usual first principles definition of the quasicoherent inverse and direct image functors $f^\bullet$ and $f_\bullet$, as recalled in Definition A.1.1, for a locally ringed space $(X, \Omega_X^\bullet)$ with 'structure sheaf' given by the de Rham complex, viewed as a sheaf of commutative DG algebras.

*Remark* A.5.6. The locally ringed space $(X, \Omega_X^\bullet)$ of the preceding remark is a model for the quotient stack of $X$ by the formal groupoid integrating the tangent sheaf $\Theta_X$, viewed as a Lie algebroid on $X$ under usual Lie bracket of vector fields and anchor map given by the identity; the sheaf of commutative DG algebras $\Omega_X^\bullet = C_{\mathrm{CE}}^\bullet(\Theta_X)$ is precisely the sheaf of Chevalley-Eilenberg cochains on the Lie algebroid $\Theta_X$ (with coefficients in the trivial module). Further, in these terms, the equivalence of categories between $D$ modules and $\Omega_X^\bullet$ modules can be understood as an instance of filtered Koszul duality between the Chevalley-Eilenberg cochains and the universal enveloping algebra of the Lie algebroid $\Theta_X$.

The stack modelled by the locally ringed space of the preceding remark is called the de Rham stack of $X$, and denoted $X_{\mathrm{dR}}$. We now give the formal definition of the de Rham stack of a general prestack, and give a careful statement of our preferred variant of the equivalence of Remark A.5.2, following [GR14a].

Throughout, let $\mathcal{Y} \in \mathrm{PreStk}$ be a prestack, as in Definition B.1.1.

*Definition* A.5.7. The de Rham prestack $\mathcal{Y}_{\mathrm{dR}} \in \mathrm{PreStk}$ of $\mathcal{Y}$ is defined by

$$\mathcal{Y}_{\mathrm{dR}} : \mathrm{DGSch}_{\mathrm{aff}} \to \mathrm{Grpd} \qquad S \mapsto \mathcal{Y}^{(\mathrm{cl,red}\,S)} \ .$$

*Remark* A.5.8. This construction is natural in $\mathcal{Y}$, and thus upgrades to a functor $(\cdot)_{\mathrm{dR}} : \mathrm{PreStk} \to \mathrm{PreStk}$.

*Definition* A.5.9. The category of left $D$ modules on $\mathcal{Y}$ is defined by

$$D^l(\mathcal{Y}) = \mathrm{QCoh}^\bullet(\mathcal{Y}_{\mathrm{dR}})$$

where the category of quasicoherent sheaves on prestacks is as in Definition B.4.2.

*Definition* A.5.10. For $\mathcal{Y} \in \mathrm{PreStk}_{\mathrm{laft}}$ a prestack locally almost of finite type in the sense of 1.3.1 of [GR14a], the category of right $D$ modules on $\mathcal{Y}$ is defined by

$$D^r(\mathcal{Y}) = \mathrm{IndCoh}^!(\mathcal{Y}_{\mathrm{dR}})$$

where the category of indcoherent sheaves is as defined in 10.2.3 of [Gai13].

*Warning* A.5.11. We will make limited use of the definition of $D$ modules on prestacks given here; we will by default use the definitions given in Appendix B.

*Warning* A.5.12. Throughout the remainder of this subsection, we will assume $\mathcal{Y} \in \mathrm{PreStk}_{\mathrm{laft}}$ is a prestack locally almost of finite type whenever discussing $D^r(\mathcal{Y})$, though see also the following remark.

*Remark* A.5.13. For $X \in \mathrm{Sch}_{\mathrm{ft}} \hookrightarrow \mathrm{PreStk}_{(\mathrm{laft})}$ a smooth, finite type variety, these definitions recover the usual categories of left and right $D$ modules of subappendix A.2; see e.g. Subsection 5.5 of [GR14a]. This equivalence, which justifies the preceding general definitions, is the also the analogue of the equivalence of Remark A.5.2 above.

*Remark* A.5.14. We recall constructions of the category of indcoherent sheaves on certain classes of prestacks of infinite type in Appendix B.5, following [Ras15b]; we also use the preceding definition in this context, when applicable.

*Remark* A.5.15. For each $\mathcal{Y} \in \mathrm{PreStk}$, there is a natural quotient map $\mathrm{q}_{\mathcal{Y}} : \mathcal{Y} \to \mathcal{Y}_{\mathrm{dR}}$, defining a natural transformation $\mathbb{1}_{\mathrm{PreStk}} \to (\cdot)_{\mathrm{dR}}$ of endofunctors of $\mathrm{PreStk}$.

*Definition* A.5.16. The left and right forgetful functors are defined by the quasicoherent and indcoherent pullback along $\mathrm{q}_{\mathcal{Y}} : \mathcal{Y} \to \mathcal{Y}_{\mathrm{dR}}$, respectfully:

$$\mathrm{oblv}^l := \mathrm{q}_{\mathcal{Y}}^\bullet : D^l(\mathcal{Y}) \to \mathrm{QCoh}(\mathcal{Y}) \qquad \mathrm{oblv}^r := \mathrm{q}_{\mathcal{Y}}^! : D^r(\mathcal{Y}) \to \mathrm{IndCoh}(\mathcal{Y}) \ .$$

*Remark* A.5.17. The right forgetful functor $\mathrm{oblv}^r$ of the preceding definition is t-exact, while the left analogue $o^l$ is not, in keeping with our choice to preference the right t structure on $D$ modules, as in Remark A.2.5. Further in keeping with this, in the case when $\mathcal{Y}$ is not a scheme, we will by default consider the sheaf of $\mathcal{O}$ modules underlying a $D$ module to be indcoherent rather than quasicoherent.

*Remark* A.5.18. In keeping with Remark A.5.5, we define various operations on $D$ modules in terms of those on indcoherent sheaves, following [GR14a] and chapter 4 of [GR17b].

*Proposition* A.5.19. The functor $\mathrm{oblv}^r = \mathrm{q}^!_{\mathcal{Y}} : D^r(\mathcal{Y}) \to \mathrm{IndCoh}(\mathcal{Y})$ admits a left adjoint

$$(\cdot)_D := \mathrm{q}_{\mathcal{Y}_\bullet} : \mathrm{IndCoh}(\mathcal{Y}) \to D^r(\mathcal{Y}) \ ,$$

following 3.3 of [GR14a].

*Remark* A.5.20. Concretely, for finite type scheme $X \in \mathrm{Sch}_{\mathrm{ft}}$, the functor of the preceding proposition is given by induction...

*Remark* A.5.21. For $f : \mathcal{X} \to \mathcal{Y}$ a map of (locally almost of finite type) prestacks, the functoriality of the de Rham stack construction together with that of $\mathrm{QCoh}^\bullet$ ($\mathrm{IndCoh}^!$) defines a pullback functor $D(\mathcal{Y}) \to D(\mathcal{X})$.

*Definition* A.5.22. The inverse image functor on $D$ modules is defined by

$$f^! := f^\bullet_{\mathrm{dR}} : D^l(\mathcal{Y}) \to D^l(\mathcal{X}) \qquad f^! := f^!_{\mathrm{dR}} : D^r(\mathcal{X}) \to D^r(\mathcal{Y}) \ ,$$

as in the preceding remark.

*Proposition* A.5.23. There is a natural commutative diagram

$$(A.5.1) \qquad \begin{array}{ccc} D^r(\mathcal{Y}) & \xrightarrow{f^!} & D^r(\mathcal{X}) \\ {\scriptstyle \mathrm{oblv}^r_{\mathcal{Y}}} \downarrow & & \downarrow {\scriptstyle \mathrm{oblv}^r_{\mathcal{X}}} \\ \mathrm{IndCoh}(\mathcal{Y}) & \xrightarrow{f^!} & \mathrm{IndCoh}(\mathcal{X}) \end{array} \qquad .$$

*Remark* A.5.24. For $f : \mathcal{X} \to \mathcal{Y}$ an (ind) nil-schematic map of prestacks in the sense of 4.1.3.5 of [GR17b], the corresponding map on de Rham stacks is (ind) inf-schematic, and defines a natural pushforward functor $D^r(\mathcal{X}) \to D^r(\mathcal{Y})$, following *loc. cit.*.

*Definition* A.5.25. For $f : \mathcal{X} \to \mathcal{Y}$ (ind)nil-schematic, the direct image functor on $D$ modules is defined by the ind-coherent direct image functor

$$f_* = f_{\mathrm{dR}, \bullet} : D^r(\mathcal{X}) \to D^r(\mathcal{Y}) \ ,$$

following 4.2.1.3 of [GR17b].

*Proposition* A.5.26. There is a natural commutative diagram

$$(A.5.2) \qquad \begin{array}{ccc} \mathrm{IndCoh}(\mathcal{X}) & \xrightarrow{f_\bullet} & \mathrm{IndCoh}(\mathcal{Y}) \\ {\scriptstyle (\cdot)_{D_{\mathcal{X}}}} \downarrow & & \downarrow {\scriptstyle (\cdot)_{D_{\mathcal{Y}}}} \\ D^r(\mathcal{X}) & \xrightarrow{f_*} & D^r(\mathcal{Y}) \end{array} \qquad .$$

A.6. **Convolution.** In this subsection, we review the formalism of convolution algebras in geometric representation theory, following the approach to Springer theory of the well-known textbook of Chriss and Ginzburg [CG11], for example. In our applications of interest, we will use analogous constructions in the setting of (quasi and ind)coherent sheaves which will rely more heavily on derived algebraic geometry. We follow the results of Ben-Zvi, Francis, and Nadler [BZFN10], as well as Gaitsgory and Rozenblyum [GR17a, GR17b], on this topic.

*Warning* A.6.1. Throughout this subsection, we assume for simplicity all $D$ modules are regular holonomic, and abbreviate $D_{\mathrm{rh}}(X)$ by just $D(X)$, in contrast with our general conventions.

Let $f : X \to Y$ be a map of finite type schemes, $Z = X \times_Y X$ be the self fibre product, $\Delta : X \to Z$ the diagonal inclusion, and consider the convolution diagram

$$Z_{(3)} := X \times_Y X \times_Y X \xrightarrow{\pi_{13}} X \times_Y X \quad .$$

$$\downarrow \pi_{12} \qquad \nearrow \pi_{23}$$

$$X \times_Y X \qquad\qquad X \times_Y X$$

*Proposition* A.6.2. Let $A \in D(Z)$ be a (regular holonomic) $D$ module on $Z$, equipped with a map

$$m : \pi_{13*}(\pi_{12}^* A \otimes^* \pi_{23}^* A) \to A \quad \text{in} \quad D(Z) \ ,$$

together with a coherent choice of its higher arity analogues as in Remark A.6.3, and a map $u : \Delta_* \mathbb{K}_X \to A$ such that the map induced as in Equation A.6.1 below is the identity. Then $C^\bullet(Z, A) \in \mathrm{Alg}_{\mathbb{E}_1}(\mathrm{Vect}_\mathbb{K})$ defines an associative algebra.

*Proof.* The unit of the $(\pi^*, \pi_*)$ adjunction for holonomic $D$ modules give a map $A \to \pi_* \pi^* A$, so that under the above hypotheses we have maps

$$A^{\boxtimes 2} \to (\pi_{12,*} \pi_{12}^* A) \boxtimes (\pi_{23,*} \pi_{23}^* A) \cong (\pi_{12} \times \pi_{23})_* \pi_{12}^* A \boxtimes \pi_{23}^* A \qquad \text{in } D(Z^{\times 2}) \ ,$$

$$\pi_{12}^* A \boxtimes \pi_{23}^* A \to \Delta_* \Delta^* (\pi_{12}^* A \boxtimes \pi_{23}^* A) \cong \Delta_* (\pi_{12}^* A \otimes^* \pi_{23}^* A) \qquad \text{in } D(Z_{(3)}^{\times 2}) \ .$$

Together with the map $m : \pi_{13*}(\pi_{12}^* A \otimes^* \pi_{23}^* A) \to A$, the induced maps on de Rham coch ains give

$$C^\bullet(Z; A)^{\otimes \mathbb{K} 2} \to C^\bullet(Z_{(3)}^{\times 2}, \pi_{12}^* A \boxtimes \pi_{23}^* A) \to C^\bullet(Z_{(3)}, \pi_{12}^* A \otimes^* \pi_{23}^* A) \xrightarrow{\pi_{*}(m)} C^\bullet(Z, A) \ .$$

If the map $m$ is specified together with a compatible, coherent choice of higher multiplication maps as required, then this induces a homotopy coherent associative product structure.

If $u : \Delta_* \mathbb{K}_X \to A$ is such that map $C^\bullet(A) \to C^\bullet(A)$ induced as above by

$$(\mathrm{A.6.1}) \qquad A \to \pi_{12,*} \pi_{12}^* A \qquad \text{and} \qquad \pi_{13*}(\pi_{12}^* A \otimes^* \pi_{23}^* \mathbb{K}_Z) \xrightarrow{m \circ (\mathbb{1} \otimes u)} A$$

is the identity, then the induced map $u : \mathbb{K} = C^\bullet(\mathrm{pt}) \to C^\bullet(Z, \mathbb{K}) \to C^\bullet(Z, A)$ defines a unit for the above associative product.                                                                                   $\square$

*Remark* A.6.3. For each $n \in \mathbb{N}$ there is a higher arity analogue of the convolution diagram

$$Z_{(n)} = X \times_Y X \times_Y \ldots \times_Y X \qquad \text{with maps} \qquad \pi_{ij} : Z_{(n)} \to Z$$

defined for each $i \neq j \in \{1, ..., n\}$. To specify a dg associative algebra in a homotopy coherent way, we require maps

$$\pi_{1n*}\left( \pi_{12}^* A \otimes^* \pi_{23}^* A \otimes^* \ldots \otimes^* \pi_{(n-1)n}^* A \right) \to A$$

for each $n \in \mathbb{N}$ together with natural compatibility data with iterated compositions of lower arity multiplication maps. We will mostly ignore this subtlety throughout.

*Remark* A.6.4. Suppose the projection maps $\pi_{ij}$ in the convolution diagram are proper. Then we can equivalently formulate the required data as maps

$$\pi_{12}^* A \otimes^* \pi_{23}^* A \otimes^* \ldots \otimes^* \pi_{(n-1)n}^* A \to \pi_{1n}^! A \ .$$

We will restrict to this setting for the following two examples.

*Remark* A.6.5. Let $M \in D(X)$, $N \in D(W)$, $f : X \to Y$, $g : W \to Y$ proper, $Z = X \times_Y W$, and consider the internal Hom object

$$\mathcal{H}om^*(g_*N, f_*M) = f_*M \otimes^!_Y \mathbb{D}_Y(g_*N) \quad \in \quad D(Y)$$

for the $\otimes^*$ tensor structure, where $\otimes^!_Y : D(Y)^{\otimes 2} \to D(Y)$ denotes the usual ! tensor structure on $Y$. Then by base change along

$$(A.6.2) \quad \begin{array}{ccc} X \times_Y W & \xrightarrow{\pi_X^{\times 2}} & X \times W \\ \downarrow{\scriptstyle \tilde{f}} & & \downarrow{\scriptstyle f \times g} \\ Y & \xrightarrow{\Delta_Y} & Y^{\times 2} \end{array} \quad \text{we have} \quad \mathcal{H}om^*(g_*N, f_*M) \cong \tilde{f}_*((\pi_X^! M) \otimes^!_Z (\pi_W^! \mathbb{D}N)) \ .$$

Thus, we obtain a lift of the internal hom object to $D(X \times_Y W)$, which we denote

$$\widetilde{\mathcal{H}om}(g_*N, f_*M) = M \boxtimes_Y \mathbb{D}N := (\pi_X^! N) \otimes^!_{X \times_Y W} (\pi_W^! \mathbb{D}N) \quad \in \quad D(X \times_Y W) \ ,$$

where we have introduced the shorthand $\boxtimes_Y : D(X) \times D(W) \to D(X \times_Y W)$, which agrees with the usual exterior product in the case $Y = \mathrm{pt}$.

*Example* A.6.6. Consider the special case $X = W$ and $M = N$ in the previous example and let $A = \mathcal{H}om(f_*M, f_*M) \in D(Z)$. Then there is a canonical multiplication map defined by

$$\begin{aligned} \pi_{12}^* A \otimes^* \pi_{23}^* A &= (M \boxtimes_Y \mathbb{D}M \boxtimes_Y \mathbb{K}_X) \otimes^*_{Z_{(3)}} (\mathbb{K}_X \boxtimes_Y M \boxtimes_Y \mathbb{D}M) \\ &= M \boxtimes (\mathbb{D}M \otimes^*_X M) \boxtimes_Y \mathbb{D}M \\ &\to M \boxtimes_Y \omega_X \boxtimes_Y \mathbb{D}M = \pi_{13}^! A \ , \end{aligned}$$

where the map is given by the $\otimes^*$ tensor structure duality pairing $\mathbb{D}M \otimes^* M \to \omega_X$. The associative algebra structure on $C^\bullet(Z, A) = \mathrm{Hom}_{D(Y)}(f_*M, f_*M)$ induced by Proposition A.6.2 agrees with that given by composition of maps.

*Example* A.6.7. Consider the special case of the previous example where $X$ is smooth of dimension $d_X$, and $M = \omega_X \in D(X)$. Then

$$\mathrm{Hom}_{D(Y)}(f_*\omega_X, f_*\omega_X) \cong \mathrm{Hom}_{D(X)}(\omega_X, f^! f_*\omega_X) \cong \mathrm{Hom}_{D(X)}(\omega_X, \pi_{X*}\pi_X^!\omega_X) \cong H^{\mathrm{BM}}_\bullet(Z)[-2d_X]$$

so that the latter inherits the structure of an algebra under composition of maps. Indeed, in this case we have

$$\widetilde{\mathcal{H}om}(f_*\omega_X, f_*\omega_X) = \pi_X^! \omega_X \otimes^! \pi_X^! \mathbb{D}\omega_X \cong \omega_Z[-2d_X] \ ;$$

this is precisely the usual format for convolution algebras, as in e.g. [CG11].

This calculation can be summarized by the diagram in spaces, and induced diagram of categories, given by

(A.6.3)

*Remark* A.6.8. The proposition A.6.2 and the remarks and examples which follow it above evidently apply in much greater generality then that in which they are stated; all that is required is a sheaf theory admitting a (partially defined) six functors formalism, in the sense of subsection A.3, satisfying the adjunctions used to construct each of the relevant maps and check their compatibilities.

In particular, we use the above paradigm to deduce the existence of associative algebra structures in situations where the spaces $Z$ and $Z_{(n)}$ in question are infinite type schemes or stacks, and/or where we use a sheaf theory other than that of regular holonomic $D$ modules; the proof of the preceding proposition reduces the proofs of these more general statements to checking the existence, and adjointness and base-change properties, of the analogous pushforward and pullback functors between the relevant categories of sheaves in these more general settings.

We now recall a slight conceptual enhancement of the construction of Proposition A.6.2 above:

*Proposition* A.6.9. The category of $D$ modules $D(Z)^\star \in \mathrm{DGCat}_{\mathbb{E}_1}$ on $Z$ is naturally a monoidal category with respect to the convolution monoidal structure $(\cdot) \star (\cdot) : D(Z) \otimes D(Z) \to D(Z)$ defined by the composition

$$D(Z) \otimes D(Z) \xrightarrow{\pi_{12}^* \boxtimes \pi_{23}^*} D(Z_{(3)}^{\times 2}) \xrightarrow{\Delta^*} D(Z_{(3)}) \xrightarrow{\pi_{13,*}} D(Z) \ .$$

Further, the de Rham cochains functor $C^\bullet(Z, \cdot) : D(Z)^\star \to \mathrm{Vect}_{\mathbb{K}}$ is lax monoidal with respect to the convolution monoidal structure on $D(Z)$. In particular, there exist natural maps

$$C_{\mathrm{dR}}^\bullet(Z, M) \otimes_{\mathbb{K}} C_{\mathrm{dR}}^\bullet(Z, N) \to C^\bullet(Z, M \otimes^\bullet N)$$

for each $M, N \in D(Z)$.

*Proposition* A.6.10. An object $A \in D(Z)$ together with data required in the hypotheses of Proposition A.6.2 defines an associative algebra object $A \in \mathrm{Alg}_{\mathbb{E}_1}(D(Z)^\star)$.

Moreover, the associative algebra structure on $C^\bullet(Z, A) \in \mathrm{Alg}_{\mathbb{E}_1}(\mathrm{Vect}_{\mathbb{K}})$ constructed in *loc. cit.* is the natural one induced on the image of this algebra object under the lax monoidal functor $C^\bullet(Z, \cdot)$.

*Remark* A.6.11. The propositions A.6.9 and A.6.10 also apply in much greater generality than that of the hypotheses stated at the begining of this subsection, as in Remark A.6.8 above.

A.6.1. *Convolution constructions for modules, module categories, and module objects.* We now give the analogues of propositions A.6.9 and A.6.10 required for constructing module categories and module objects over the monoidal categories and algebras above. In particular, we deduce the analogue of A.6.2 which gives a construction of modules over the algebra $C^\bullet(Z, A) \in \mathrm{Alg}_{\mathbb{E}_1}(\mathrm{Vect}_{\mathbb{K}})$.

Let $f : X \to Y$ a map of finite type schemes and $Z = X \times_Y X$ as above such that $\pi_{ij} : Z_{(3)} \to Z$ is proper for each $i \neq j \in \{1, 2, 3\}$, exactly as in Proposition A.6.2.

Further, let $g : W \to Y$, $Z^W = X \times_Y W$, and consider the convolution diagram

$$Z_{(3)}^W := X \times_Y X \times_Y W \xrightarrow{\pi_{13}} X \times_Y W \ .$$

$$\downarrow{\scriptstyle \pi_{12}} \qquad \nearrow{\scriptstyle \pi_{23}}$$

$$X \times_Y X \qquad\qquad X \times_Y W$$

*Proposition* A.6.12. Let $A \in D(Z)$ be as in Proposition A.6.2, and $R \in D(Z^W)$ a $D$ module on $Z^W$, equipped with a map

$$\rho : \pi_{13*}(\pi_{12}^* A \otimes^\bullet \pi_{23}^* R) \to R \quad \text{in} \quad D(Z^W)$$

and a coherent choice of its higher arity analogues as in Remark A.6.3. Then $C^\bullet(Z^W, R) \in C^\bullet(Z, A)$-Mod$(\mathrm{Vect}_\mathbb{K})$ defines a module for the associative algebra $C^\bullet(Z, A) \in \mathrm{Alg}_{\mathbb{E}_1}(\mathrm{Vect}_\mathbb{K})$ defined in Proposition A.6.2.

*Proposition* A.6.13. The category of $D$ modules on $Z^W$ is naturally a module category

$$D(Z^W) \in D(Z)^\star\text{-Mod}(\mathrm{DGCat}_{\mathrm{cont}})$$

for the monoidal category $D(Z)^\star$, with structure functor $(\cdot) \star (\cdot) : D(Z) \otimes D(Z^W) \to D(Z^W)$ defined by the composition

$$D(Z) \otimes D(Z^W) \xrightarrow{\pi_{12}^* \boxtimes \pi_{23}^*} D((Z_{(3)}^W)^{\times 2}) \xrightarrow{\Delta^*} D(Z_{(3)}^W) \xrightarrow{\pi_{13,*}} D(Z^W) \ .$$

Further, the de Rham cochains functors $C^\bullet(Z, \cdot) : D(Z) \to \mathrm{Vect}$ and $C^\bullet(Z^W, \cdot) : D(Z^W) \to \mathrm{Vect}$ are lax compatible with the above module structure, in the sense that there are natural maps

$$C^\bullet(Z, A) \otimes_\mathbb{K} C^\bullet(Z^W, R) \to C^\bullet(Z^W, A \star R)$$

for each $A \in D(Z)$ and $R \in D(Z^W)$.

*Proposition* A.6.14. Let $A \in D(Z)$ and $M \in D(Z^W)$ be as in Proposition A.6.12, so that by Proposition A.6.10 we have $A \in \mathrm{Alg}_{\mathbb{E}_1}(D(Z)^\star)$. Then $M \in A$-Mod$(D(Z^W))$ is a module object internal to the $D(Z)^\star$ module category $D(Z^W)$ given in Proposition A.6.13.

Further, the induced $C^\bullet(Z, A)$ module structure on $C^\bullet(Z^W, M)$ constructed in Proposition A.6.12 is the natural one induced on the images of the above algebra and module objects under the lax compatible functors $C^\bullet$.

*Remark* A.6.15. In analogy with Remark A.6.4, in the case where the convolution diagram maps $\pi_{ij}$ are proper, we can reformulate the required data in terms of maps $\pi_{12}^* A \otimes^* \pi_{23}^* R \to \pi_{13}^! R$ and their higher arity analogues, and we restrict to this setting for the following examples:

*Example* A.6.16. Let $M \in D(X)$, $N \in D(W)$, and following examples A.6.5 and A.6.6 above, let

$$A = \tilde{\mathcal{H}\mathrm{om}}(f_* M, f_* M) = M \boxtimes_Y \mathbb{D}M \ \in D(Z) \ , \text{ and analogously}$$

$$R = \tilde{\mathcal{H}\mathrm{om}}(g_* N, f_* M) = M \boxtimes_Y \mathbb{D}N \ \in D(Z^W) \ .$$

Then $R \in A$-Mod$(D(Z^W))$ with module structure map $\rho : \pi_{12}^* A \otimes^* \pi_{23}^* R \to \pi_{13}^! R$ defined by

$$\pi_{12}^* A \otimes^* \pi_{23}^* R = (M \boxtimes_Y \mathbb{D}M \boxtimes_Y \mathbb{K}_X) \otimes_{Z_{(3)}}^* (\mathbb{K}_X \boxtimes_Y M \boxtimes_Y \mathbb{D}N)$$

$$= M \boxtimes (\mathbb{D}M \otimes_X^* M) \boxtimes_Y \mathbb{D}N$$

$$\to M \boxtimes_Y \omega_X \boxtimes_Y \mathbb{D}N = \pi_{13}^! R \ ,$$

where the map is given by the $\otimes^*$ tensor structure duality pairing $\mathbb{D}M \otimes^* M \to \omega_X$.

The $C^\bullet(Z, A) = \mathrm{Hom}_{D(Y)}(f_* M, f_* M)$ module structure on $C^\bullet(Z^W, M) = \mathrm{Hom}(N, M)$ induced by Proposition A.6.14 agrees with that given by composition of maps.

*Example* A.6.17. Following Example A.6.7, consider the special case of the previous example where $W$ is smooth of dimension $d_W$, $M = \omega_X \in D(X)$ and $N = \omega_W \in D(W)$. Then

$$\mathrm{Hom}_{D(Y)}(g_* \omega_W, f_* \omega_X) \cong \mathrm{Hom}_{D(W)}(\omega_W, g^! f_* \omega_X) \cong \mathrm{Hom}_{D(W)}(\omega_W, \pi_{W*} \pi_X^! \omega_X) \cong H_\bullet^{\mathrm{BM}}(Z^W)[-2d_W]$$

so that the latter inherits the structure of an module over $\mathrm{Hom}_{D(Y)}(f_* \omega_X, f_* \omega_X) \cong H_\bullet^{\mathrm{BM}}(Z)[-2d_X]$ under composition of maps. Indeed, in this case we have

$$\tilde{\mathcal{H}\mathrm{om}}(g_* \omega_W, f_* \omega_X) = \pi_X^! \omega_X \otimes^! \pi_W^! \mathbb{D}\omega_W \cong \omega_{Z^W}[-2d_W] \ ;$$

this is preciesly the usual format for modules over convolution algebras, as in e.g. [CG11].

This calculation can be summarized by the diagram in spaces, and induced diagram of categories, given by

(A.6.4)

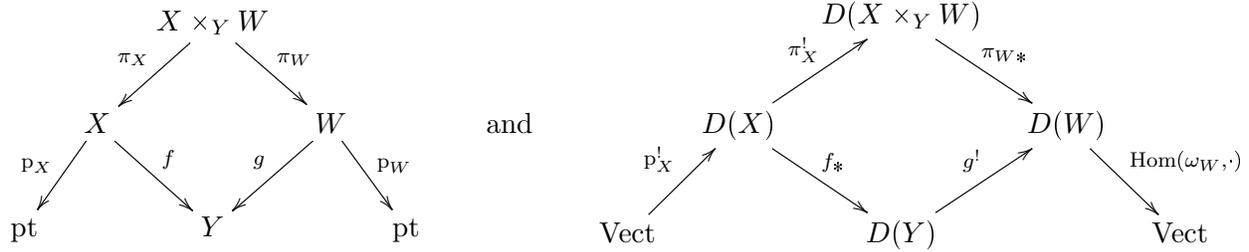

*Remark* A.6.18. The propositions and examples of this subsection also apply in much greater generality than that in which they are stated here, as explained in remakrs A.6.8 and A.6.11.

## Appendix B. Sheaf theory on infinite dimensional varieties

In this section, we recall some aspects of the foundations of sheaf theory on infinite dimensional varieties and stacks given in [GR14a, GR14b, Gai15, GR17a, GR17b, Ras15a, Ras15b, Ras20b] and references therein. None of the material presented here is original.

### B.1. **Prestacks.**

In this subsection, we recall the notion of prestack, which will feature as the ambient category of spaces in which we phrase the formal categorical aspects of constructions involving spaces.

*Definition* B.1.1. A prestack is a functor $\mathcal{Y} : \mathrm{DGSch}_{\mathrm{aff}}^{\mathrm{op}} \to \mathrm{Grpd}$ from the category of affine DG schemes to the category of groupoids.

Let PreStk denote the category of prestacks.

*Remark* B.1.2. An object $\mathcal{Y} \in \mathrm{PreStk}$ is thought of as the functor of points of a putative stack, without any representability criteria; for $S \in \mathrm{Sch}_{\mathrm{aff}}$, $\mathcal{Y}(S) =: \mathrm{Maps}(S, \mathcal{Y})$ is thought of as the groupoid of $S$ points of the space $\mathcal{Y}$.

*Remark* B.1.3. The Yoneda lemma defines a fully faithful embedding $\mathrm{DGSch} \hookrightarrow \mathrm{PreStk}$, identifying sets with the full subcategory of discrete groupoids.

*Definition* B.1.4. For $n \in \mathbb{N}$, an affine DG scheme $S \in \mathrm{DGSch}_{\mathrm{aff}}$ is called $n$-coconnective if $S = \mathrm{Spec}\,(A)$ for $A$ in cohomological degree $\leqslant -n$. Similarly, $S \in \mathrm{DGSch}_{\mathrm{aff}}$ is called eventually coconnective if it is $n$-coconnective for some $n$, and classical if it is 0-coconnective.

Let $^{\leqslant n}\mathrm{DGSch}_{\mathrm{aff}}$ denote the full subcategory of $\mathrm{DGSch}_{\mathrm{aff}}$ on $n$-coconnective objects, $^{<\infty}\mathrm{DGSch}_{\mathrm{aff}}$ the full subcategory on eventually coconnective objects, and note the category of classical affine DG schemes $^{\leqslant 0}\mathrm{DGSch}_{\mathrm{aff}} = \mathrm{Sch}_{\mathrm{aff}}$ is tautologically equivalent to that of usual affine schemes.

*Definition* B.1.5. A prestack $\mathcal{Y} \in \mathrm{PreStk}$ is called $n$-coconnective if it is defined by left Kan extension from a functor $\mathcal{Y} : ^{\leqslant n}\mathrm{DGSch}_{\mathrm{aff}} \to \mathrm{Grpd}$. Similarly, a prestack $\mathcal{Y} \in \mathrm{PreStk}$ is called eventually coconnective if it is $n$-coconnective for some $n$, and classical if it is 0-coconnective.

Let $^{\leqslant n}\mathrm{PreStk}$ denote the full subcategory on $n$-coconnective prestacks, $^{<\infty}\mathrm{PreStk}$ the full subcategory on eventually coconnective prestacks, and $^{\mathrm{cl}}\mathrm{PreStk}$ the full subcategory on classical prestacks.

### B.2. **Pro-finite type schemes.**

*Warning* B.2.1. For simplicity of notation, abstract limit and colimit diagram categories in this Appendix will be assumed to be filtered, unless stated otherwise.

*Remark* B.2.2. The category of $\mathbb{K}$ algebras is closed under colimits, so that the category of affine schemes is closed under limits. Thus, an object which is presented as a limit of finite type affine schemes is again simply an (affine) scheme. In fact, every $\mathbb{K}$ algebra is canonically equivalent to the union of its finite type subalgebras, so that $\mathrm{Sch}_{\mathrm{aff}} \cong \mathrm{Pro}(\mathrm{Sch}_{\mathrm{aff,ft}})$.

*Remark* B.2.3. More generally, the category of schemes is closed under colimits of diagrams with all structure maps affine. Further, Noetherian approximation (see [TT07], and regarding applications to the present context, Section 3.2 of [Ras15b] and references therein) implies that every quasi-compact, quasi-separated scheme can be presented as a limit of a diagram in finite type schemes with affine structure maps.

The canonical functor $\mathrm{Pro}^{\mathrm{aff}}(\mathrm{Sch}_{\mathrm{ft}}) \to \mathrm{Sch}$ from the category of such limits to schemes is fully faithful with essential image the quasi-compact, quasi-separated schemes.

For notational simplicity, we make the following (non-standard) definition:

*Definition* B.2.4. A pro-finite type (DG) scheme is a quasi-compact, quasi-seperated (and eventually coconnective, DG) scheme.

Let $\mathrm{Sch}_{\mathrm{pft}}$ and $\mathrm{DGSch}_{\mathrm{pft}}$ denote the categories of pro-finite type (DG) schemes.

*Remark* B.2.5. In summary, each object $X \in \mathrm{Sch}_{\mathrm{pft}}$ admits a presentation

$$X = \lim_{i \in \mathcal{I}} X_i = \lim \left[ \ldots \to X_i \xrightarrow{\varphi_{ij}} X_j \to \ldots \right] \qquad \text{so that}$$

with $X_i \in \mathrm{Sch}_{\mathrm{ft}}$ finite type and $\varphi_{ij}$ affine for each $i, j \in \mathcal{I}$.

*Example* B.2.6. The infinite dimensional affine space $\mathbb{A}^\infty \in \mathrm{Sch}_{\mathrm{pft}}$ over $\mathbb{K}$ is defined by

$$\mathbb{A}^\infty = \lim_n \mathbb{A}^n = \lim \left[ \ldots \twoheadrightarrow \mathbb{A}^{n+1} \twoheadrightarrow \mathbb{A}^n \twoheadrightarrow \ldots \twoheadrightarrow \mathbb{A}^1 \right]$$

with structure maps given by projecting out the last coordinate. It follows

$$\mathbb{A}^\infty = \mathrm{Spec} \ \mathbb{K}[t_1, t_2, ...] \qquad \text{where} \qquad \mathbb{K}[t_1, t_2, ...] = \mathrm{colim}_n \mathbb{K}[t_1, ..., t_n] \ .$$

*Example* B.2.7. Let $\mathbb{D} = \mathrm{Spf} \ \mathcal{O}$ denote the formal punctured disk, where $\mathcal{O} = \mathbb{K}[[z]]$. The jet scheme, or arc scheme

$$\mathcal{J}(\mathbb{A}^d)_0 := \mathbb{A}^d_{\mathcal{O}} = \mathrm{Maps}(\mathbb{D}, \mathbb{A}^d) \ \in \mathrm{Sch}_{\mathrm{pft}}$$

of $\mathbb{A}^d$ is a pro-finite type scheme as follows: The description $\mathbb{D} = \mathrm{colim}_i \mathbb{D}_i$ corresponding to $\mathbb{K}[[z]] = \lim_i \mathbb{K}[z]/z^i$ induces the presentation

$$\mathrm{Maps}(\mathbb{D}, \mathbb{A}^d) = \lim_i \mathrm{Maps}(\mathbb{D}_i, \mathbb{A}^d) = \lim \left[ \ldots \twoheadrightarrow \mathrm{Maps}(\mathbb{D}_i, \mathbb{A}^d) \xrightarrow{\mathrm{rest}_{ij}} \mathrm{Maps}(\mathbb{D}_j, \mathbb{A}^d) \twoheadrightarrow \ldots \right]$$

where each term is given by

$$\mathrm{Maps}(\mathbb{D}_i, \mathbb{A}^d) = \mathrm{Hom}(\mathbb{K}[t_1, ..., t_d], \mathbb{K}[z]/z^i) = \mathbb{A}^d_{\mathcal{O}}/z^i \mathbb{A}^d_{\mathcal{O}} \cong \mathbb{A}^{di}$$

so that $\mathbb{A}^d_{\mathcal{O}}$ indeed defines a pro-finite type (affine) scheme. In the case $d = 1$, $\mathbb{A}^1_{\mathcal{O}}$ is isomorphic to $\mathbb{A}^\infty \in \mathrm{Sch}_{\mathrm{pft}}$, as defined in the preceding example.

*Example* B.2.8. Let $Y$ be a finite type affine scheme. The jet scheme, or arc scheme

$$\mathcal{J}(Y)_0 = Y_{\mathcal{O}} = \mathrm{Maps}(\mathbb{D}, Y) \ \in \mathrm{Sch}_{\mathrm{pft}}$$

of $Y$ is a pro-finite type scheme as above, presented by

$$\mathrm{Maps}(\mathbb{D}, Y) = \lim_i \mathrm{Maps}(\mathbb{D}_i, Y) = \lim \left[ \ldots \twoheadrightarrow \mathrm{Maps}(\mathbb{D}_i, Y) \xrightarrow{\mathrm{rest}_{ij}} \mathrm{Maps}(\mathbb{D}_j, Y) \twoheadrightarrow \ldots \right] \ .$$

Note that a closed embedding $Y \hookrightarrow \mathbb{A}^n$ induces an embedding $Y_{\mathcal{O}} \hookrightarrow \mathbb{A}^n_{\mathcal{O}}$ giving an equivalent pro-finite type presentation.

### B.3. Indschemes.

*Remark* B.3.1. In contrast to Remark B.2.2 above, the category of rings is not closed under limits, and the category of (affine) schemes is not closed under colimits.

*Definition* B.3.2. A (DG) indscheme is a prestack $X = \operatorname{colim}_k X^k \in \mathrm{PreStk}$ presented as a colimit of a diagram in pro-finite type (DG) schemes $X^k \in \mathrm{Sch}_{\mathrm{pft}}$ ($X^k \in \mathrm{DGSch}_{\mathrm{pft}}$) with structure maps given by closed embeddings.

Let IndSch and DGIndSch denote the categories of (DG) indschemes.

*Remark* B.3.3. Concretely, in analogy with Remark B.2.5, we have

$$X = \operatorname*{colim}_{k \in \mathcal{L}} X^k = \operatorname{colim}\left[\ldots \leftarrow X^l \xleftarrow{\iota^{kl}} X^k \leftarrow \ldots\right] \qquad \text{so that}$$

$$\ldots \longrightarrow X^k \xrightarrow{\iota^{kl}} X^l \longrightarrow \ldots$$

where $X^k \in \mathrm{Sch}_{\mathrm{pft}}$ is a pro-finite type scheme and $\iota^{kl} : X^k \hookrightarrow X^l$ is a closed embedding for each $k, l \in \mathcal{L}$. In particular, applying the description of Remark B.2.5 to each $X^k \in \mathrm{Sch}_{\mathrm{pft}}$, we have a presentation

$$(\mathrm{B.3.1}) \qquad X = \operatorname*{colim}_{k \in \mathcal{L}} \lim_{i \in \mathcal{I}} X_i^k \qquad \text{where} \qquad 
\begin{array}{ccc}
X_i^k & \xrightarrow{\varphi_{ij}^k} & X_j^k \\
\iota_i^{kl} \downarrow & & \downarrow \iota_j^{kl} \\
X_i^l & \xrightarrow{\varphi_{ij}^l} & X_j^l
\end{array}$$

is the general component of the bi-diagram, with $X_i^k \in \mathrm{Sch}_{\mathrm{ft}}$ finite type, $\varphi_{ij}^k$ affine, and $\iota_i^{kl}$ a closed embedding, for each combination of $i, j \in \mathcal{I}$ and $k, l \in \mathcal{L}$.

The existence of such a presentation is called local compactness in [KV06], where they do not require the schemes $X^k \in \mathrm{Sch}$ to be pro-finite type in the definition of indscheme.

*Example* B.3.4. The affine Grassmannian $\mathrm{Gr}_G$ is an indscheme, via the presentation $\mathrm{Gr}_G = \operatorname{colim}_{\lambda \in \Lambda^+}(\mathrm{Gr}_G^{\leq \lambda})$. Note $\mathrm{Gr}_G$ is ind-finite type and ind-proper.

*Example* B.3.5. Let $\mathbb{D}^\circ = \mathbb{D}\backslash\{0\}$ be the formal punctured disk and $\mathcal{K} = \mathbb{K}((x))$. The meromorphic jet scheme, or algebraic loop scheme

$$\mathcal{J}^{\mathrm{mer}}(\mathbb{A}^d)_0 = \mathbb{A}_{\mathcal{K}}^d = \mathrm{Maps}(\mathbb{D}^\circ, \mathbb{A}^d) \in \mathrm{IndSch}$$

of $\mathbb{A}^d$ is an indscheme as follows: The description $\mathbb{K}((x)) = \operatorname{colim}_k \mathbb{K}((x))^{<k}$, where $\mathbb{K}((x))^{<k}$ is space of laurent polynomials with poles bounded by $k$, induces the presentation

$$\mathrm{Maps}(\mathbb{D}^\circ, \mathbb{A}^d) = \operatorname{colim}_k \mathrm{Maps}(\mathbb{D}^\circ, \mathbb{A}^d)^{<k} = \operatorname{colim}\left[\ldots \hookrightarrow \mathrm{Maps}(\mathbb{D}^\circ, \mathbb{A}^d)^{<k} \xrightarrow{\iota^{kl}} \mathrm{Maps}(\mathbb{D}^\circ, \mathbb{A}^d)^{<l} \hookrightarrow \ldots\right]$$

where each subscheme is given by the space of maps with poles bounded by $k$

$$\mathrm{Maps}(\mathbb{D}^\circ, \mathbb{A}^d)^{<k} := (\mathbb{K}((x))^{<k})^{\times d} = \lim_i (\mathbb{K}[x^{\pm 1}]^{<k}/x^i)^{\times d} \quad \in \mathrm{Sch}_{\mathrm{pft}}$$

which is presented as a pro-finite type scheme as in Example B.2.7 above.

In summary, we obtain a presentation

$$\mathbb{A}^d_{\mathcal{K}} = \operatorname*{colim}_{k \in \mathcal{L}} \lim_{i \in \mathcal{I}} \ \operatorname{Maps}(\mathbb{D}^\circ_i, \mathbb{A}^d)^{<k} \qquad \text{where}$$

$$
\begin{array}{ccc}
\operatorname{Maps}(\mathbb{D}^\circ_i, \mathbb{A}^d)^{<k} & \xrightarrow{\operatorname{rest}^k_{ij}} & \operatorname{Maps}(\mathbb{D}^\circ_j, \mathbb{A}^d)^{<k} \\
{\scriptstyle \iota^{kl}_i} \downarrow & & \downarrow {\scriptstyle \iota^{kl}_j} \\
\operatorname{Maps}(\mathbb{D}^\circ_i, \mathbb{A}^d)^{<l} & \xrightarrow{\operatorname{rest}^l_{ij}} & \operatorname{Maps}(\mathbb{D}^\circ_j, \mathbb{A}^d)^{<l}
\end{array}
$$

is the general component of the bidiagram, and $\operatorname{Maps}(\mathbb{D}^\circ_i, \mathbb{A}^d)^{<k} := (\mathbb{K}[x^{\pm 1}]^{<k}/x^i)^{\times d}$. Note that $\mathbb{A}^d_{\mathcal{K}}$ is a genuinely ind-pro finite type object; it is not ind-finite type.

*Example* B.3.6. Let $Y$ be a finite type affine scheme. The meromorphic jet scheme, or algebraic loop scheme

$$\mathcal{J}^{\operatorname{mer}}(Y)_0 = Y_{\mathcal{K}} = \operatorname{Maps}(\mathbb{D}^\circ, Y) \ \in \operatorname{IndSch}$$

of $Y$ is an indscheme, as follows: as in Example B.2.8, a close dembedding $Y \in \mathbb{A}^d$ induces an embedding $Y_{\mathcal{K}} \hookrightarrow \mathbb{A}^d_{\mathcal{K}}$, the the restriction of the indscheme presentation of the latter in the preceding example presents $Y_{\mathcal{K}}$ as an indscheme.

B.3.1. *Reasonable indschemes.*

*Definition* B.3.7. Let $X$ be an indscheme. A subscheme $Y \hookrightarrow X$ is called a reasonable subscheme if $Y$ is a closed, pro-finite type subscheme, and for any closed subscheme $Y' \hookrightarrow X$ containing $Y$, the closed embedding $Y \hookrightarrow Y'$ is finitely presented.

An indscheme $X$ is called reasonable if it is the colimit of its reasonable subschemes.

Let $\operatorname{IndSch}_{\operatorname{reas}}$ denote the full subcategory of reasonable indschemes.

*Remark* B.3.8. Heuristically, a subscheme is reasonable if it is maximally infinite dimensional (up to finite codimension) among all closed subschemes, and the indscheme is reasonable if it can be approximated by such subschemes.

*Example* B.3.9. Let $X$ be an ind-finite type indscheme. Then any (necessarily finite type) closed subscheme is reasonable, and thus $X$ is reasonable.

*Example* B.3.10. Let $\mathcal{T}$ be the total space of a pro-finite type vector bundle over an ind-finite type indscheme $X$. Then the restriction of $\mathcal{T}$ to a closed subscheme of $X$ is a reasonable subscheme of $\mathcal{T}$. An infinite codimension subbundle of $\mathcal{T}$ restricted to a closed subscheme of $X$ is not reasonable.

More generally, in the case of a DG indscheme, we make the following definition:

*Definition* B.3.11. A DG indscheme is callled reasonable if it can be presented as a colimit $X = \operatorname{colim}_k X^k$ of pro-finite type DG schemes $X^k \in \operatorname{DGSch}_{\operatorname{pft}}$ under almost finitely presented closed embeddings.

Let $\operatorname{DGIndSch}_{\operatorname{reas}}$ denote the full subcategory of reasonable DG indschemes, and define the category of reasonable DG schemes by $\operatorname{DGSch}_{\operatorname{reas}} = \operatorname{DGIndSch}_{\operatorname{reas}} \cap \operatorname{DGSch}_{\operatorname{pft}}$ be the category of reasonable DG schemes.

B.4. **QCoh on prestacks.**

*Remark* B.4.1. The functoriality properties of the theory of quasicoherent sheaves on affine schemes, recalled in Appendix A, imply that there is a functor
(B.4.1)
$$\operatorname{QCoh}^\bullet : \operatorname{DGSch}^{\operatorname{op}}_{\operatorname{aff}} \to \operatorname{DGCat}_{\operatorname{cont}} \qquad X \mapsto \operatorname{QCoh}(X) \qquad [f : X \to Y] \mapsto [f^\bullet : \operatorname{QCoh}(Y) \to \operatorname{QCoh}(X)] \,.$$

*Definition* B.4.2. The functor $\mathrm{QCoh}^\bullet : \mathrm{PreStk}^{\mathrm{op}} \to \mathrm{DGCat}_{\mathrm{cont}}$ is the right Kan extension of $\mathrm{QCoh}^\bullet : \mathrm{DGSch}_{\mathrm{aff}}^{\mathrm{op}} \to \mathrm{DGCat}_{\mathrm{cont}}$ along $\mathrm{DGSch}_{\mathrm{aff}}^{\mathrm{op}} \to \mathrm{PreStk}^{\mathrm{op}}$.

*Remark* B.4.3. Concretely, by pointwise evaluation of the opposite left Kan extension as a weighted colimit, an object of $\mathcal{F} \in \mathrm{QCoh}^\bullet(\mathcal{Y})$ is given by an assignment

$$S \mapsto [(\cdot)^\bullet(\mathcal{F}) : \ \mathrm{Maps}(S, \mathcal{Y}) \longrightarrow \mathrm{QCoh}(S)\ ] \qquad [\varphi : S \to T] \mapsto \ \mathrm{Maps}(T, \mathcal{Y}) \xrightarrow{(\cdot)^\bullet(\mathcal{F})} \mathrm{QCoh}(T) \ ,$$

$$\text{``} [f : S \to \mathcal{Y}] \longmapsto f^\bullet(\mathcal{F}) \text{``} \qquad\qquad \downarrow {\scriptstyle (\cdot)\circ\varphi} \qquad {\scriptstyle \varphi^\bullet}\downarrow$$

$$\mathrm{Maps}(S, \mathcal{Y}) \xrightarrow{(\cdot)^\bullet(\mathcal{F})} \mathrm{QCoh}(S)$$

defined for each $S \in \mathrm{DGSch}_{\mathrm{aff}}$ and each morphism $\varphi : S \to T$ of affine schemes.

*Proposition* B.4.4. Let $f : \mathcal{X} \to \mathcal{Y}$ be a schematic, quasicompact map of prestacks. Then there exists a natural continuous right adjoint $f_\bullet : \mathrm{QCoh}(\mathcal{X}) \to \mathrm{QCoh}(\mathcal{Y})$ to $f^\bullet$, defining a functor

$$\mathrm{QCoh} : \mathrm{PreStk}_{\mathrm{sch\text{-}qc}} \to \mathrm{DGCat}_{\mathrm{cont}} \qquad \mathcal{Y} \mapsto \mathrm{QCoh}(\mathcal{Y}) \qquad [f : \mathcal{X} \to \mathcal{Y}] \mapsto [f_\bullet : \mathrm{QCoh}(\mathcal{X}) \to \mathrm{QCoh}(\mathcal{Y})] \,.$$

*Remark* B.4.5. More generally, the above constructions satisfy base-change, and can be extended to a functor from the correspondence category $\mathrm{PreStk}_{(\mathrm{corr};\mathrm{all},\mathrm{sch\text{-}qc})} \to \mathrm{DGCat}_{\mathrm{cont}}$, where the former is as defined in Example 5.2.4. This can also be lifted to a 2-categorical variant as in section 5.3.2 of [GR17a], though we do not work directly with this definition.

## B.5. IndCoh **on pro-finite type DG schemes and reasonable DG indschemes.** Throughout, let $X, Y \in \mathrm{DGSch}_{\mathrm{pft}}$ be pro-finite type (which we have defined to mean quasicompact, quasiseperated, and eventually coconnective) DG schemes.

*Definition* B.5.1. The category $\mathrm{IndCoh}(X) \in \mathrm{DGCat}$ of indcoherent sheaves on $X \in \mathrm{DGSch}_{\mathrm{pft}}$ is defined as the ind completion of the category $\mathrm{Coh}(X)$ of coherent sheaves on the DG scheme $X$.

*Proposition* B.5.2. For $f : X \to Y$ an arbitrary map, there is a natural functor $f_\bullet : \mathrm{IndCoh}(X) \to \mathrm{IndCoh}(Y)$, defining a functor

$$\mathrm{IndCoh} : \mathrm{DGSch}_{\mathrm{pft}} \to \mathrm{DGCat}_{\mathrm{cont}} \qquad X \mapsto \mathrm{IndCoh}(X) \qquad [f : X \to Y] \mapsto [f_\bullet : \mathrm{IndCoh}(X) \to \mathrm{IndCoh}(Y)] \,.$$

*Proposition* B.5.3. For $f : X \to Y$ flat, there is a natural, continuous left adjoint $f^\bullet : \mathrm{IndCoh}(Y) \to \mathrm{IndCoh}(X)$ to $f_\bullet$, defining a functor

$$\mathrm{IndCoh} : \mathrm{DGSch}_{\mathrm{pft},\mathrm{flat}}^{\mathrm{op}} \to \mathrm{DGCat}_{\mathrm{cont}} \qquad X \mapsto \mathrm{IndCoh}(X) \qquad [f : X \to Y] \mapsto [f^\bullet : \mathrm{IndCoh}(Y) \to \mathrm{IndCoh}(X)] \,,$$

where $\mathrm{DGSch}_{\mathrm{pft},\mathrm{flat}}$ denotes the wide subcategory of $\mathrm{DGSch}_{\mathrm{pft}}$ with flat maps.

*Remark* B.5.4. More generally, following [GR17a], the above constructions satisfy base-change, and can be extended to a functor from the correspondence category $\mathrm{Sch}_{\mathrm{pft},(\mathrm{corr};\mathrm{flat},\mathrm{all})} \to \mathrm{DGCat}_{\mathrm{cont}}$, in the notation of Example 5.2.4. This can also be lifted to a 2-categorical variant, though we do not work directly with this definition.

More generally, we make the following definition:

*Remark* B.5.5. Concretely, this means that for $X = \lim_i X_i \in \mathrm{DGSch}_{\mathrm{pft}}$ presented as a limit of finite type schemes $X_i \in \mathrm{DGSch}_{\mathrm{ft}}$ under affine maps, we have

$$\mathrm{IndCoh}^!(X) = \operatorname*{colim}_i \mathrm{IndCoh}(X_i) = \operatorname*{colim} \left[ \dots \leftarrow \mathrm{IndCoh}(X_i) \xleftarrow{\varphi_{ij}^!} \mathrm{IndCoh}(X_j) \leftarrow \dots \right] \ , \text{ and}$$

$$\mathrm{IndCoh}^\bullet(X) = \lim_i \mathrm{IndCoh}(X_i) = \lim \left[ \dots \to \mathrm{IndCoh}(X_i) \xrightarrow{\varphi_{ij*}} \mathrm{IndCoh}(X_j) \to \dots \right] \ ,$$

where we have used the canonical identifications $\mathrm{IndCoh}^!(X_i) = \mathrm{IndCoh}^\bullet(X_i) = \mathrm{IndCoh}(X_i)$ for $X_i \in \mathrm{DGSch}_{\mathrm{ft}}$. Further, the remaining content of Remark B.6.3 applies similarly here, *mutatis mutandis*.

Recall the notion of reasonable DG indscheme from Subappendix B.3.1. We define categories of indcoherent sheaves on reasonable DG indschemes, following [Ras20b] and references therein, as follows:

*Definition* B.5.6. The functor $\mathrm{IndCoh}^\bullet : \mathrm{DGIndSch}_{\mathrm{reas}} \to \mathrm{DGCat}_{\mathrm{cont}}$ is defined as the left Kan extension of $\mathrm{IndCoh}^\bullet : \mathrm{DGSch}_{\mathrm{reas}} \to \mathrm{DGCat}_{\mathrm{cont}}$ along $\mathrm{DGSch}_{\mathrm{reas}} \to \mathrm{DGIndSch}_{\mathrm{reas}}$.

Similarly, the functor $\mathrm{IndCoh}^!_{\mathrm{reas}} \to \mathrm{DGCat}_{\mathrm{cont}}$ is defined as the right Kan extension of $\mathrm{IndCoh}^! : \mathrm{DGSch}^{\mathrm{op}}_{\mathrm{reas}} \to \mathrm{DGCat}_{\mathrm{cont}}$ along $\mathrm{DGSch}^{\mathrm{op}}_{\mathrm{reas}} \to \mathrm{DGIndSch}^{\mathrm{op}}_{\mathrm{reas}}$.

*Remark* B.5.7. Concretely, for $X = \mathrm{colim}_k X^k$ presented as a colimit of reasonable, pro-finite type DG schemes under almost finitely presented closed embeddings, there is a presentation of $\mathrm{IndCoh}^\bullet(X) \in \mathrm{DGCat}$ as a colimit

$$\mathrm{IndCoh}^\bullet(X) = \underset{k}{\mathrm{colim}}\, \mathrm{IndCoh}^\bullet(X^k) = \mathrm{colim}\left[\ldots \leftarrow \mathrm{IndCoh}^\bullet(X^l) \overset{\iota^{kl}_\bullet}{\leftarrow} \mathrm{IndCoh}^\bullet(X^k) \leftarrow \ldots\right].$$

Similarly, there is a presentation of $\mathrm{IndCoh}^!(X) \in \mathrm{DGCat}$ as a limit

$$\mathrm{IndCoh}^!(X) = \underset{k}{\mathrm{lim}}\, \mathrm{IndCoh}^!(X^k) = \mathrm{lim}\left[\ldots \to \mathrm{IndCoh}^!(X^k) \overset{\iota^{kl,!}}{\longrightarrow} \mathrm{IndCoh}^!(X^l) \to \ldots\right].$$

*Proposition* B.5.8. For $X, Y \in \mathrm{DGSch}_{\mathrm{reas}}$ and $f : X \to Y$ proper and almost finitely presented, there is a natural, continuous right adjoint $f^! : \mathrm{IndCoh}^\bullet(Y) \to \mathrm{IndCoh}^\bullet(X)$ to $f_\bullet$, defining a functor

$$\mathrm{IndCoh} : \mathrm{DGSch}^{\mathrm{op}}_{\mathrm{pft,prop}} \to \mathrm{DGCat}_{\mathrm{cont}} \qquad X \mapsto \mathrm{IndCoh}^\bullet(X) \qquad [f : X \to Y] \mapsto [f^! : \mathrm{IndCoh}^\bullet(Y) \to \mathrm{IndCoh}^\bullet(X)]$$

where $\mathrm{DGSch}_{\mathrm{pft,prop}}$ denotes the wide subcategory of $\mathrm{DGSch}_{\mathrm{pft}}$ with proper maps.

*Remark* B.5.9. There is a presentation of $\mathrm{IndCoh}^\bullet(X)$ as a limit

$$\mathrm{IndCoh}^\bullet(X) = \underset{k}{\mathrm{lim}}\, \mathrm{IndCoh}^\bullet(X^k) = \mathrm{lim}\left[\ldots \to \mathrm{IndCoh}^\bullet(X^k) \overset{\iota^{kl,!}}{\longrightarrow} \mathrm{IndCoh}^\bullet(X^l) \to \ldots\right],$$

as follows from the preceding proposition, by passing to right adjoints in the sense of 1.8.4.2 in [GR17a] in the presentation of Remark B.5.7.

B.6. *D* **modules on prestacks.** We would like to extend the $D$ module formalism of Subappendix A.2 to pro-finite type schemes and indschemes. We begin by discussing those aspects of the theory which follow from formal categorical constructions.

*Remark* B.6.1. The functoriality properties of the theory of $D$ modules on finite type schemes recalled in Appendix A.2 imply that we have functors:

(B.6.1) $\qquad D^! : \mathrm{Sch}^{\mathrm{op}}_{\mathrm{aff,ft}} \to \mathrm{DGCat}_{\mathrm{cont}} \qquad X \mapsto D(X) \qquad [f : X \to Y] \mapsto [f^! : D(Y) \to D(X)]$

(B.6.2) $\qquad D^* : \mathrm{Sch}_{\mathrm{aff,ft}} \to \mathrm{DGCat}_{\mathrm{cont}} \qquad X \mapsto D(X) \qquad [f : X \to Y] \mapsto [f_* : D(X) \to D(Y)]$

To begin, we extend these functors to the category of affine schemes:

*Definition* B.6.2. The functor $D^! : \mathrm{Sch}^{\mathrm{op}}_{\mathrm{aff}} \to \mathrm{DGCat}_{\mathrm{cont}}$ is the left Kan extension of $D^! : \mathrm{Sch}^{\mathrm{op}}_{\mathrm{aff,ft}} \to \mathrm{DGCat}_{\mathrm{cont}}$ along $\mathrm{Sch}_{\mathrm{aff,ft}} \to \mathrm{Sch}_{\mathrm{aff}}$.

The functor $D^* : \mathrm{Sch}_{\mathrm{aff}} \to \mathrm{DGCat}_{\mathrm{cont}}$ is the right Kan extension of $D^* : \mathrm{Sch}_{\mathrm{aff,ft}} \to \mathrm{DGCat}_{\mathrm{cont}}$ along $\mathrm{Sch}_{\mathrm{aff,ft}} \to \mathrm{Sch}_{\mathrm{aff}}$.

*Remark* B.6.3. Concretely, this means that for $X = \lim X_i \in \mathrm{Sch}_{\mathrm{aff}}$ the presentation of $X$ as a limit of finite type affine schemes $X_i$, we have

$$D^!(X) = \operatorname*{colim}_i D(X_i) = \operatorname{colim}\left[\ldots \leftarrow D(X_i) \xleftarrow{\varphi^!_{ij}} D(X_j) \leftarrow \ldots\right] \ .$$

In particular, there are canonical functors $\varphi^!_i : D(X_i) \to D(X)$ which are by definition the pull back along the canonical map $\varphi_i : X \to X_i$. Heuristically, a typical object in $D^!(X)$ is of the form $\varphi^! M$ for $M \in D(Y)$ with $Y$ a finite type affine scheme, for some $\varphi : X \to Y$.

Similarly, for $X = \lim_i X_i \in \mathrm{Sch}_{\mathrm{aff}}$ as above, we have

$$D^*(X) = \lim_i D(X_i) = \lim\left[\ldots \to D(X_i) \xrightarrow{\varphi_{ij*}} D(X_j) \to \ldots\right] \ .$$

In particular, there are canonical functors $\varphi_{i,*} : D(X) \to D(X_i)$ which are by definition the push-forward along $\varphi_i : X \to X_i$. An object in $D^*(X)$ is specified by a system of $D$ modules $M_i \in D(X_i)$ together with identifications $\varphi_{i,j,*} M_j \xrightarrow{\cong} M_i$.

There are analogous interpretations in terms of finite dimensional approximations for the induced functors $f_* : D(X) \to D(Y)$ and $f^! : D(Y) \to D(X)$ for each map $f : X \to Y$ of affine schemes.

*Definition* B.6.4. The functor $D^! : {}^{\mathrm{cl}}\mathrm{PreStk}^{\mathrm{op}} \to \mathrm{DGCat}_{\mathrm{cont}}$ is the right Kan extension of $D^! : \mathrm{Sch}^{\mathrm{op}}_{\mathrm{aff}} \to \mathrm{DGCat}_{\mathrm{cont}}$ along $\mathrm{Sch}^{\mathrm{op}}_{\mathrm{aff}} \to {}^{\mathrm{cl}}\mathrm{PreStk}^{\mathrm{op}}$.

The functor $D^* : {}^{\mathrm{cl}}\mathrm{PreStk} \to \mathrm{DGCat}_{\mathrm{cont}}$ is the left Kan extension of $D^* : \mathrm{Sch}_{\mathrm{aff}} \to \mathrm{DGCat}_{\mathrm{cont}}$ along $\mathrm{Sch}_{\mathrm{aff}} \to {}^{\mathrm{cl}}\mathrm{PreStk}$.

*Remark* B.6.5. Concretely, by pointwise evaluation of the opposite left Kan extension as a weighted colimit, an object of $M \in D^!(\mathcal{Y})$ is given by an assignment

$$S \mapsto \left[(\cdot)^!(M) : \ \mathrm{Maps}(S, \mathcal{Y}) \longrightarrow D(S)\right] \qquad [\varphi : S \to T] \mapsto \begin{array}{c} \mathrm{Maps}(T, \mathcal{Y}) \xrightarrow{(\cdot)^!(M)} D(T) \\ {\scriptstyle (\cdot)\circ\varphi}\downarrow \qquad \downarrow {\scriptstyle \varphi^!} \\ \mathrm{Maps}(S, \mathcal{Y}) \xrightarrow{(\cdot)^!(M)} D(S) \end{array} \ ,$$

$$\text{``} [f : S \to \mathcal{Y}] \longmapsto f^!(M) \text{``}$$

defined for each $S \in \mathrm{Sch}_{\mathrm{aff}}$ and each morphism $\varphi : S \to T$ of affine schemes.

The analogous interpretation of $D^*(\mathcal{Y})$ identifies it with the dual category to $D^!(\mathcal{Y})$, whenever the latter is dualizeable. Heuristically, a typical object of $D^*(\mathcal{Y})$ is of the form $\varphi_* M$ for $M \in D(Y)$ with $Y$ an affine scheme, for some $\varphi : Y \to \mathcal{Y}$.

There are analogous interpretations in terms of probe affine schemes for the functors $f_* : D^*(\mathcal{X}) \to D^*(\mathcal{Y})$ and $f^! : D^!(\mathcal{Y}) \to D^!(\mathcal{X})$ for a morphism $f : \mathcal{X} \to \mathcal{Y}$ of prestacks.

*Remark* B.6.6. For a general prestack $\mathcal{Y} \in {}^{\mathrm{cl}}\mathrm{PreStk}$, the map $\pi : \mathcal{Y} \to \mathrm{pt}$ yields a dualizing sheaf object $\omega_{\mathcal{Y}} = \pi^! \underline{\mathbb{K}}_{\mathrm{pt}} \in D^!(\mathcal{Y})$ and a de Rham cohomology functor $H^\bullet_{\mathrm{dR}} = \pi_* : D^*(X) \to \mathrm{Vect}_{\mathbb{K}}$. However, without any identification between $D^!(X)$ and $D^*(X)$, this data can not be used to define a cohomology theory; some genuine geometric structures are required.

### B.7. $D$ **modules on pro-finite type schemes.**

Now, we restrict our attention to pro-finite type schemes. Recall that Noetherian approximation identifies the category $\mathrm{Sch}_{\mathrm{pft}}$ of pro-finite type schemes with a full subcategory $\mathrm{Pro}^{\mathrm{aff}}(\mathrm{Sch}_{\mathrm{ft}})$ of $\mathrm{Pro}(\mathrm{Sch}_{\mathrm{ft}})$ on objects presented by diagrams with affine structure maps.

*Definition* B.7.1. The functor $\tilde{\mathcal{D}}^! : \mathrm{Sch}_{\mathrm{pft}}^{\mathrm{op}} \to \mathrm{DGCat}_{\mathrm{cont}}$ is defined as the left Kan extension of $D^! : \mathrm{Sch}_{\mathrm{ft}}^{\mathrm{op}} \to \mathrm{DGCat}_{\mathrm{cont}}$ along $\mathrm{Sch}_{\mathrm{ft}} \hookrightarrow \mathrm{Sch}_{\mathrm{pft}}$.

Similarly, the functor $\tilde{\mathcal{D}}^* : \mathrm{Sch}_{\mathrm{pft}} \to \mathrm{DGCat}_{\mathrm{cont}}$ is defined as the right Kan extension of $D^* : \mathrm{Sch}_{\mathrm{ft}} \to \mathrm{DGCat}_{\mathrm{cont}}$ along $\mathrm{Sch}_{\mathrm{ft}} \hookrightarrow \mathrm{Sch}_{\mathrm{pft}}$.

*Proposition* B.7.2. The canonical natural transformations $\tilde{\mathcal{D}}^! \to D^!$ and $\mathcal{D}^* \to \tilde{\mathcal{D}}^*$ are isomorphisms.

*Remark* B.7.3. Concretely, this means that for $X = \lim_i X_i \in \mathrm{Sch}_{\mathrm{pft}}$ presented as a limit of finite type schemes $X_i \in \mathrm{Sch}_{\mathrm{ft}}$ under affine maps, we have

$$\mathcal{D}^!(X) = \underset{i}{\mathrm{colim}}\, D(X_i) = \mathrm{colim}\left[ \ldots \leftarrow D(X_i) \xleftarrow{\varphi_{ij}^!} D(X_j) \leftarrow \ldots \right] \text{, and}$$

$$\mathcal{D}^*(X) = \underset{i}{\lim}\, D(X_i) = \lim\left[ \ldots \to D(X_i) \xrightarrow{\varphi_{ij*}} D(X_j) \to \ldots \right] ,$$

where we have used the canonical identifications $D^!(X_i) = D^*(X_i) = D(X_i)$ for $X_i \in \mathrm{Sch}_{\mathrm{ft}}$. Further, the remaining content of Remark B.6.3 applies similarly here, *mutatis mutandis*.

*Example* B.7.4. Let $X \in \mathrm{Sch}_{\mathrm{pft}}$ be a pro-finite type scheme. The constant sheaf $\underline{\mathbb{K}}_X \in D^*(X)$ is defined by

$$``\varphi_{i*}(\underline{\mathbb{K}}_X)" := \underset{j \in I_{/i}}{\mathrm{colim}}\, \varphi_{ij*}\underline{\mathbb{K}}_{X_j} \qquad \text{noting} \qquad \varphi_{ij*}``\varphi_{i*}(\underline{\mathbb{K}}_X)" \xrightarrow{\cong} ``\varphi_{j*}(\underline{\mathbb{K}}_X)" .$$

In particular, we have

$$C_{\mathrm{dR}}^\bullet(X; \underline{\mathbb{K}}) := \pi_* \underline{\mathbb{K}}_X \cong \underset{i}{\mathrm{colim}}\, C_{\mathrm{dR}}^\bullet(X_i; \mathbb{K}) .$$

*Definition* B.7.5. Let $f : X \to Y$ be a finitely presented map of pro-finite type schemes. The functors

$$f_{*,!} : \mathcal{D}^!(X) \to D^!(Y) \qquad \text{and} \qquad f^{!,*} : D^*(Y) \to D^*(X)$$

are defined as the colimit and limit of the functors $f_*$ and $f^!$ defined on finite type, affine approximations in the $\tilde{\mathcal{D}}^!$ and $\tilde{\mathcal{D}}^*$ presentations.

*Proposition* B.7.6. For $f : X \to Y$ a proper (in particular, finitely presented) map of pro-finite type schemes, the functor $f_{*,!}$ is canonically left adjoint to $f^!$ on $D^!$, and similarly $f^{*,!}$ is canonically right adjoint to $f_*$ on $D^*$.

*Proposition* B.7.7. For $f : X \to Y$ a smooth, finitely presented map of pro-finite type schemes, the functor $f^{!,*}[-2d_{X/Y}]$ is canonically left adjoint to $f_*$ on $D^*$, and similarly $f_{*,!}$ is canonically right adjoint to $f^![-2d_{X/Y}]$ on $D^!$.

Here $d_{X/Y} : X \to \mathbb{Z}$ is the rank of $\Omega_{X/Y}^1$, a locally constant function on $X$; see definition B.7.8 below.

*Definition* B.7.8. A locally constant function $d : X \to \mathbb{Z}$ defined as a map of indschemes, where $\mathbb{Z} = \sqcup_{n \in \mathbb{Z}} \mathrm{pt}$.

For $f : X \to Y$ a map of finite type schemes, the relative dimension of $f$ is the locally constant function $d_{X/Y} = \dim_X - \dim_Y \circ f$ on $X$.

For $f : X \to Y$ a finitely presented map of placid schemes, the relative dimension of $f$ is the locally constant function $d_{X/Y} : X \to \mathbb{Z}$ defined by Noetherian approximation.

*Definition* B.7.9. A presentation $X = \lim_i X_i$ of a scheme $X \in \mathrm{Sch}$ is called placid if the index category is filtered, each $X_i \in \mathrm{Sch}_{\mathrm{ft}}$ is finite type, and the structure maps are all affine and smooth. The scheme $X$ is called placid if it admits such a presentation.

*Remark* B.7.10. A placid scheme is in particular profinite type.

*Proposition* B.7.11. Let $X = \lim_i X_i$ be a placid presentation. There are canonical equivalences

$$D^!(X) = \lim_i D^!(X_i) \qquad \text{under} \qquad \varphi_{ij(*,\mathrm{ren})} := \varphi_{ij(*,!)}[-2d_{X_i/X_j}] : D^!(X_i) \to D^!(X_j) \quad , \text{and}$$

$$D^*(X) = \operatorname*{colim}_i D^*(X_i) \qquad \text{under} \qquad \varphi^*_{ij} := \varphi^{!,*}_{ij}[-2d_{X_i/X_j}] : D^*(X_j) \to D^*(X_i) \ .$$

In particular, there are canonical functors $\varphi_{i(*,\mathrm{ren})} : D^!(X) \to D^!(X_i)$ and $\varphi^*_i : D^*(X_i) \to D^*(X)$.

*Definition* B.7.12. Let $X$ be a placid, pro-finite type scheme with placid presentation $X = \lim_i X_i$. The renormalized dualizing sheaf $\omega^{\mathrm{ren}}_X \in D^*(X)$ is defined by

$$\omega^{\mathrm{ren}}_X = \varphi^*_i \omega_{X_i}[-2d_{X_i}] \qquad \text{noting} \qquad \varphi^*_i \omega_{X_i}[-2d_{X_i}] = \varphi^*_i \varphi^{!,*}_{ij} \omega_{X_j}[-2d_{X_i}] = \varphi^*_j \omega_{X_j}[-2d_{X_j}] \ .$$

*Remark* B.7.13. The object $\omega^{\mathrm{ren}}_X \in D^*(X)$ is canonically independent of the placid presentation.

*Example* B.7.14. For $X = \lim_i X_i$ a placid presentation with each $X_i$ a smooth scheme, $\omega^{\mathrm{ren}}_X \cong \underline{\mathbb{K}}_X$.

*Definition* B.7.15. Let $X$ be a placid, pro-finite type scheme. The renormalized Borel-Moore homology of $X$ is $H^{\mathrm{BM,ren}}_\bullet(X; \mathbb{K}) = H^\bullet_{\mathrm{dR}}(X; \omega^{\mathrm{ren}}_X)$.

## B.8. $D$ **modules on ind-schemes.** The definition B.7.1, Proposition B.7.2 and Remark B.7.3, and propositions B.7.5, B.7.6, and B.7.7 extend to the setting of indschemes:

*Proposition* B.8.1. Let $X = \operatorname{colim}_k X^k$ be a presentation of an indscheme. There are canonical equivalences

$$\mathcal{D}^!(X) = \lim_{k \in \mathcal{L}} D^!(X^k) = \lim \left[ \ldots \to D(X^l) \xrightarrow{\iota^{kl!}} D(X^k) \to \ldots \right] \text{, and}$$

$$\mathcal{D}^*(X) = \operatorname*{colim}_{k \in \mathcal{L}} D^*(X^k) = \operatorname{colim} \left[ \ldots \leftarrow D(X^l) \xleftarrow{\iota^{kl}_*} D(X^j) \leftarrow \ldots \right] \ .$$

Moreover, for $f : X \to Y$ an ind-finitely presented map of indschemes, there are natural functors

$$f_{*,!} : \mathcal{D}^!(X) \to D^!(Y) \qquad \text{and} \qquad f^{!,*} : D^*(Y) \to D^*(X)$$

satisfying the usual adjunctions for $f$ ind-proper or smooth.

Recall the notion of reasonable indscheme from Subappendix B.3.1. Following again [Ras15b] and references therein, we have the following descriptions of $D$ modules on reasonable indschemes:

*Proposition* B.8.2. Let $X = \operatorname{colim}_k X^k$ be a reasonable indscheme presented as the colimit of its reasonable subschemes. There are canonical equivalences

$$D^!(X) = \operatorname*{colim}_k D^!(X^k) = \operatorname{colim} \left[ \ldots \leftarrow D(X^l) \xleftarrow{\iota^{kl}_{*,!}} D(X^k) \leftarrow \ldots \right] \text{, and}$$

$$D^*(X) = \lim_k D^*(X^k) = \lim \left[ \ldots \to D(X^l) \xrightarrow{\iota^{kl(!,*)}} D(X^k) \to \ldots \right] \ .$$

In particular, there are canonical functors $\iota^k_{(*,!)} : D^!(X^k) \to D^!(X)$ and $\iota^{k(!,*)} : D^*(X) \to D^*(X^k)$.

*Definition* B.8.3. An indscheme $X$ is placid if $X$ is reasonable, and every reasonable subscheme of $X$ is placid.

*Remark* B.8.4. An indscheme $X$ is placid if and only if it is presented as a colimit $X = \operatorname{colim}_k X^k$ under closed embeddings where each $X_k$ is placid and a reasonable subscheme of $X$.

*Remark* B.8.5. A placid indscheme $X$ does not admit a canonical choice of renormalized dualizing sheaf $\omega_X^{\mathrm{ren}} \in D^*(X)$ in general. Following B.8.2, the natural candidate is given by assigning $\iota^{k(!,*)}\omega_X^{\mathrm{ren}} = \omega_{X^k}^{\mathrm{ren}}$, where the latter is the renormalized dualizing sheaf of the placid, reasonable subscheme $X^k \in \mathrm{Sch}_{\mathrm{pft}}$. The resulting putative object fails to be well-defined because the inclusions of distinct reasonable subschemes are potentially positive codimension, so that the shifts in the definitions of the renormalized dualizing sheaves on distinct reasonable subschemes will disagree.

However, since inclusions of reasonable subschemes are always of finite codimension, we can fix a choice of reasonable subscheme $X^0$ to renormalize 'relative to', in the sense that we define the restriction of the renormalized dualizing sheaf to other reasonable subschemes as their renormalized dualizing sheaf shifted by their (finite) dimension relative to $X^0$.

Following the preceding remark, we introduce the following data:

*Definition* B.8.6. A dimension theory on a placid indscheme $X$ is an assignment $\tau$ of a locally constant function $\tau_k : X^k \to \mathbb{Z}$ to each reasonable subscheme $X^k \hookrightarrow X$, such that for any inclusion of reasonable subschemes $\iota^{kl} : X^k \hookrightarrow X^l$, we have $\tau_k = \tau_l \circ \iota^{kl} + d_{X^k/X^l}$.

*Example* B.8.7. There is a canonical dimension theory on any placid scheme $X$, defined by the condition $\tau_X = 0$.

*Example* B.8.8. Let $X$ be an indscheme of ind-finite type. There is a canonical dimension theory on $X$ defined by $\tau_k = \dim_{X_k}$.

*Example* B.8.9. Let $X, X'$ be reasonable indschemes equipped with dimension theories $\tau, \tilde{\tau}$. Then $X \times \tilde{X}$ inherits a canonical dimension theory defined by the fact that to products $X^k \times \tilde{X}^k$ of reasonable subschemes it assigns $\tau_k \circ \mathrm{p} + \tilde{\tau}_k \circ \tilde{\mathrm{p}}$.

*Definition* B.8.10. Let $X = \mathrm{colim}_k X^k$ be a placid indscheme equipped with a dimension theory $\tau$. The $\tau$-renormalized dualizing sheaf $\omega_X^{\tau\text{-}\mathrm{ren}} \in D^*(X)$ is defined by

$$\text{``}\iota^{k(!,*)}(\omega_X^{\tau\text{-}\mathrm{ren}})\text{''} := \omega_{X^k}^{\mathrm{ren}}[2\tau_k] \qquad \text{noting} \qquad \iota_{kl}^{!,*}(\omega_{X^l}^{\mathrm{ren}}[2\tau_l]) = \omega_{X^k}^{\mathrm{ren}}[2d_{X^k/X^l} + 2\tau_l \circ \iota^{kl}] = \omega_{X^k}^{\mathrm{ren}}[2\tau_k] \,,$$

for any inclusion of reasonable subschemes $\iota^{kl} : X^k \hookrightarrow X^l$.

B.8.1. *Holonomic D modules on indschemes.*

*Definition* B.8.11. Let $X$ be an indscheme. The categories $D_{\mathrm{rh}}^*(X)$ and $D_{\mathrm{rh}}^!(X)$ of holonomic $D$ modules are defined by iterated Kan extension of the functors from $\mathrm{Sch}_{\mathrm{ft}}$, as in Definitions B.6.2 and B.6.4.

*Remark* B.8.12. There are canonical functors

$$D_{\mathrm{rh}}^!(X) \to D^!(X) \qquad \text{and} \qquad D_{\mathrm{rh}}^*(X) \to D^*(X) \,,$$

which intertwine the functors $f^!$ and $f_*$, respectively, for $f : X \to Y$ an arbitrary map of indschemes.

*Definition* B.8.13. A map $f : X \to Y$ of reasonable indschemes is called reasonable if there exists a cofinal subsystem $Y = \mathrm{colim}_k Y^k$ of reasonable subschemes $Y^k \hookrightarrow Y$ such that each $X \times_Y Y^k \hookrightarrow X$ is reasonable for each $k$.

*Example* B.8.14. Any finitely presented map of reasonable indschemes is reasonable.

*Proposition* B.8.15. Let $f : X \to Y$ be a reasonable map of reasonable indschemes. Then the partially defined left adjoint adjoint $f^*$ to $f_* : D^*(X) \to D^*(Y)$ is defined on holonomic objects, inducing $f^* : D_{\mathrm{rh}}^*(Y) \to D_{\mathrm{rh}}^*(X)$.

Let $f : X \to Y$ be a finitely presented map of placid indschemes. Then the partially defined left adjoint $f_!$ to $f^! : D^!(Y) \to D^!(X)$ is defined on holonomic objects, inducing $f_! : D_{\mathrm{rh}}^!(X) \to D_{\mathrm{rh}}^!(Y)$.

*Email address*: dbutson@perimeterinstitute.ca